\documentclass{amsart}

\usepackage{times,graphicx,amsthm}
\usepackage[percent]{overpic}
\usepackage{varioref}
\usepackage{placeins}
\usepackage[matrix,arrow]{xy}
\usepackage{array}
\usepackage{longtable}
\usepackage{booktabs} 
\usepackage{url}
\usepackage{fullpage}
\usepackage{dcolumn}
\usepackage{multicol}
\usepackage{natbib}
\usepackage{fp}
\usepackage{units}

\renewcommand{\subsection}[1]{\medskip\noindent\textbf{#1.} }

\def\marginpar#1{}   

\let\lbl=\label
\def\label#1{\lbl{#1}\ifinner\else\marginpar{\ref{#1} #1}\ignorespaces\fi}

\newcommand\<{\langle}
\newcommand\cross{\times}
\newcommand\dcsd{\operatorname{dcsd}}
\newcommand\Thi{\operatorname{Thi}}
\newcommand\PThi{\operatorname{Thi}_p}
\newcommand\PRop{\operatorname{Rop}_p}
\newcommand\Rop{\operatorname{Rop}}

\newcommand\MinRad{\operatorname{MinRad}}

\providecommand{\abs}[1]{\lvert#1\rvert} 
\providecommand{\norm}[1]{\lVert#1\rVert} 

\renewcommand{\phi}{\varphi}

 \providecommand{\abs}[1]{\lvert#1\rvert} 
 \providecommand{\norm}[1]{\lVert#1\rVert}

\newtheorem{theorem}{Theorem}[section]

\theoremstyle{definition}

\newtheorem{definition}{Definition}

\graphicspath{{./figs/}}

{\makeatletter
 \gdef\xxxmark{%
   \expandafter\ifx\csname @mpargs\endcsname\relax 
     \expandafter\ifx\csname @captype\endcsname\relax 
       \marginpar{xxx}
     \else
       xxx 
     \fi
   \else
     xxx 
   \fi}
 \gdef\xxx{\@ifnextchar[\xxx@lab\xxx@nolab}
 \long\gdef\xxx@lab[#1]#2{{\bf [\xxxmark #2 ---{\sc #1}]}}
 \long\gdef\xxx@nolab#1{{\bf [\xxxmark #1]}}
}

\newcommand{\SimpleChainRRUB}{41.7086588}

\begin{document}

\bibliographystyle{plain}

\title{Self-contact sets for 50 tightly knotted and linked tubes}
\date{June 12, 2004; Revised: \today}

\author{Ted Ashton}

\author{Jason Cantarella}

\author{Michael Piatek}

\author{Eric Rawdon}

\begin{abstract}
We report on new numerical computations of the set of self-contacts in
tightly knotted tubes of uniform circular cross-section. Such contact
sets have been obtained before for the trefoil and figure eight knots
by simulated annealing --- we use constrained gradient-descent to
provide new self-contact sets for those and 48 other knot and link types. 
The minimum length of all unit diameter tubes in a given knot or link type
is called the ropelength of that class of curves. Our computations
yield improved upper bounds for the ropelength of all knots and links
with 9 or fewer crossings except the trefoil.
\end{abstract}


\maketitle

\section{Introduction}
The study of knots as abstract topological objects has inspired a
great deal of fascinating mathematics. But knots are no less
interesting as physical structures tied in flexible ropes 
and pulled tight. While it is intuitively clear that tight knots
organize tension and contact forces to bind tightly and resist
unravelling, the details of their structure remain mysterious. Even
today, there is no explicit mathematical description of any tight
knot.

We define the \emph{thickness} of a space curve to be the supremal
(largest) radius of any embedded tubular neighborhood of the curve,
and the \emph{ropelength} of a curve to be the quotient of its length
and thickness. Ropelength provides a scale-invariant way to measure
the total flexibility of a given length of rope. It has been
established that there is a curve of minimum ropelength in each knot
and link type $L$ (\cite{MR2002m:74035,MR2003j:57010,MR2003h:58014}) and the
ropelength of that curve is called the ropelength $\Rop(L)$ of the
knot or link type. These minimum ropelength curves are called
\emph{tight} or \emph{ideal} knots. In this paper, we give some
results from our numerical computations of the shapes of tight knots
and links.

Over the past ten years, many authors have found approximate shapes
for tight knots and links using numerical
methods~\cite{MR1702021,MR1702022,mglob,baranska}. We follow previous
authors in defining a version of ropelength, $\PRop$, for space
polygons and optimizing this polygonal ropelength functional over the
space of polygons in different knot and link types. But while most other others
rely on simulated annealing, we introduce the use of constrained
gradient descent for ropelength optimization.

This new method has allowed us to significantly expand the scope and accuracy
of existing computations of tight knots. In particular, it seems that our
method has succeeded at the challenging task of resolving the set of
self-contacts in all knot and link types with nine and fewer crossings (212
types in all).  In the process, we have produced a new table of upper bounds
for the ropelengths of these knot and link types which improves upon previous
results. These improvements range from $0.05\%$ for the figure eight
knot (compared to the bound of~\cite{mglob}) to more than
$8.11\%$ for the $9_{20}$ knot (compared to the bound
of~\cite{MR2034393}).  For links, these seem to be the first upper
bounds reported for almost all of the link types we consider.

This dataset is likely to be useful in the study of tight knots, so we
provide here an early view of our results.  This research announcement will be
followed by an expanded paper, ``Tightening Knots with Constrained Gradient
Descent", which describes our methods and results in
detail. The filesize limitations of the arXiv forced us to truncate
the data section of this posting, and to use comparatively low-quality
image files for the three-dimensional views of tight knots and links
in Appendix B. A higher-quality view of these images is provided at
\url{http://www.cs.washington.edu/homes/piatek/contact_table/}.

\section{Background material}
We start by defining ropelength more precisely. Suppose $\gamma$ is a
$C^1$ space curve, parameterized by arclength. The maximum diameter of
an embedded tube around $\gamma$ is controlled by two phenomena:
``self-contacts'' of the tube formed when sections of the curve far
away in arclength approach each other in space, and ``kinks'' formed
by points of high curvature on $\gamma$.  These effects are shown
below in Figure~\ref{fig:thick}.
\begin{figure}[h]
  \centering
  \includegraphics{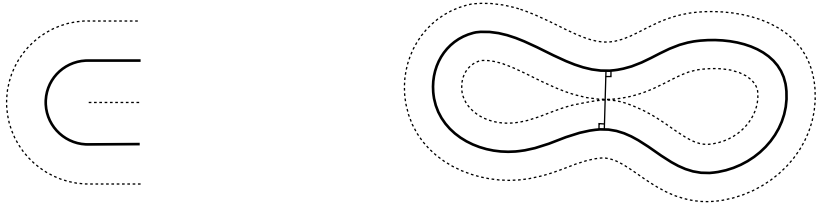}
  \caption{The thickness of a smooth curve $\gamma$ is controlled by
    curvature (as in the left picture), and self-contacts of the tube
    around $\gamma$ (as in the right picture).} 
  \label{fig:thick}
\end{figure}

To make this understanding precise, we need a few definitions. We can define
the self-distance map of $\gamma$ by 
\begin{equation*}
  d(s,t) = \norm{ \gamma(s) - \gamma(t) }.
\end{equation*}
We then have
\begin{definition}The set $\dcsd(\gamma)$ of \emph{doubly-critical
    self-distances} of $\gamma$ is the set of critical points for
  $d(s,t)$ with $s \neq t$. Taking the partial derivatives of $d$, we
  see that $(s,t) \in \dcsd(\gamma)$ if and only if
  \begin{equation*}
    \< \gamma(s) - \gamma(t), \gamma\,'(s) \rangle = 0 \text{ and } 
    \< \gamma(s) - \gamma(t), \gamma\,'(t) \rangle = 0.
  \end{equation*}
\end{definition}
Denoting the curvature of $\gamma$ at $s$ by $\kappa(s)$, Litherland et
al.~proved
\begin{theorem}
  \label{thm:lsdr}
  \cite{MR99k:57025} The thickness of $\gamma$ is the minimum of 
 \begin{equation*}
   \min_s \frac{1}{\kappa(s)} \text{ and } \min_{(s,t) \in
   \dcsd(\gamma)} \frac{d(s,t)}{2}.
 \end{equation*}
\end{theorem}
We can now define the primary object of interest in our computations:
\begin{definition}
  The \emph{self-contact set} or \emph{strut set} of a space curve
  $\gamma$ is the set of $(s,t) \in \dcsd(\gamma)$ with
  $\norm{\gamma(s) - \gamma(t)} = 2\Thi(\gamma)$.
\end{definition}
The term ``strut'', borrowed from tensegrity theory, comes from the
fact that minimum-length chords in $\dcsd(\gamma)$ ``hold the curve
apart from itself'' as the knot tightens.

In our polygonal knot-tightening problem, we will replace the curve
$\gamma$ with a space polygon $\mathcal{V}$ with vertices $v_1, \dots,
v_{V}$ and edges $e_1, \dots, e_{V}$. To define the polygonal
thickness $\PThi(\mathcal{V})$ of $\mathcal{V}$, we will need an idea
of curvature. 
\begin{definition}
  \label{def:minrad}
  The minimum radius of curvature (or $\MinRad$) of $\mathcal{V}$ at
  $v_i$ is given by the radius of the unique circle tangent to both of
  the edges which meet at $v_i$ and passing through the midpoint of
  the shorter one.
  
  If $\theta_i$ is the turning angle of $\mathcal{V}$ at $v_i$, then
  $\MinRad(v_i)$ is given by
  \begin{equation*}
    \MinRad(v_i) =
    \frac{\min\{\abs{e_{i-1}},\abs{e_i}\}}{2\,\tan(\nicefrac{\theta_i}{2})}.
  \end{equation*}
\end{definition}
We will also need a definition of doubly-critical self-distances:
\begin{definition}
  Let $\dcsd(\mathcal{V})$ be the set of $(p,q)$ on $\mathcal{V}$ with
  $p \neq q$ which are local minima of the self-distance function.
\end{definition}
We now define a thickness measure for polygons:
\begin{definition}
  \label{def:polythi}
  The thickness $\PThi(\mathcal{V})$ of a space polygon $\mathcal{V}$
  without self-intersections is given by the minimum of
  \begin{equation*}
    \min_i \MinRad(v_i)  \quad \text{ and } \quad \min_{(p,q) \in
      \dcsd(\mathcal{V})} \frac{d(p,q)}{2}.
  \end{equation*}
\end{definition}
The polygonal ropelength $\PRop$ of $\mathcal{V}$ is then the quotient of the
length of $\mathcal{V}$ and $\PThi(\mathcal{V})$. As expected, when polygons
$\mathcal{V}_n$ with increasing numbers of edges are inscribed in a space curve
$\gamma$, it is known that $\PRop(\mathcal{V}_n) \rightarrow \Rop(\gamma)$
under some mild geometric hypotheses~\cite{rawdonthesis,MR1702029,MR2034393}. 
We define self-contact sets for polygons like we do for smooth curves.

\section{Quantities computed and how the computation was validated}

Our algorithm minimizes the length of $\mathcal{V}$ subject to a
family of constraints derived from the distance and $\MinRad$
functions in Definition~\ref{def:polythi}. We report the polygonal
ropelengths of the minimized configurations in summary form in
Appendix A. We also report preliminary computations of the ropelengths
of piecewise-smooth curves obtained from our polygons by replacing the
corners of the polygon with small circular arcs.  These provide
tentative upper bounds on the ropelengths of the knot and link
types listed. We also provide self-contact sets for $50$ of these
knots and links.

Though our ropelength bounds for knots improve on those of previous
authors, this is an uncertain measure of their quality. After all,
there is no way to check the accuracy of an upper bound for the
ropelength of any knot, since no exact value for the minimum ropelength of  any knot type is known.
On the other hand, the minimum ropelength is known or explicitly conjectured for some
link types~\cite{MR2003h:58014,cfksw}.  A comparison of our results to some of these
known examples appears below.
\FPeval{twooneresult}{100*((25.1388 - 25.132741228718)/(25.132741228718))}
\FPround{\twooneresult}{\twooneresult}{2} 

\FPeval{sixthreetworesult}{100*((58.0145 - 58.006)/(58.006))}
\FPround{\sixthreetworesult}{\sixthreetworesult}{2} 

\FPeval{simplechainresult}{100*((\SimpleChainRRUB - 41.6991)/(41.6991))}
\FPround{\simplechainresult}{\simplechainresult}{2} 

\begin{center}
  \begin{tabular}{p{1in}p{1.3in}p{1.3in}p{1.3in}}
    & \includegraphics[width=1in]{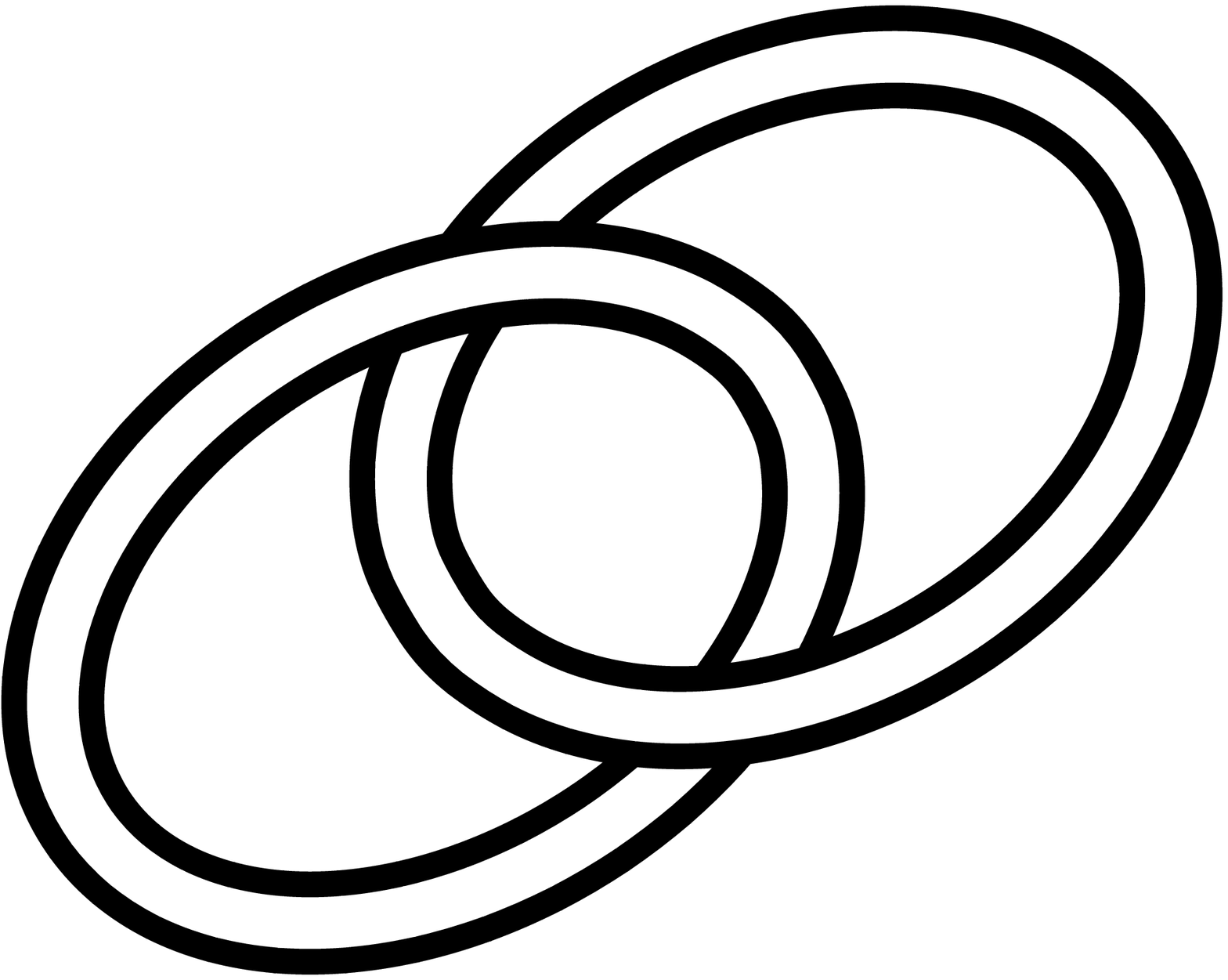} &
    \includegraphics[width=1in]{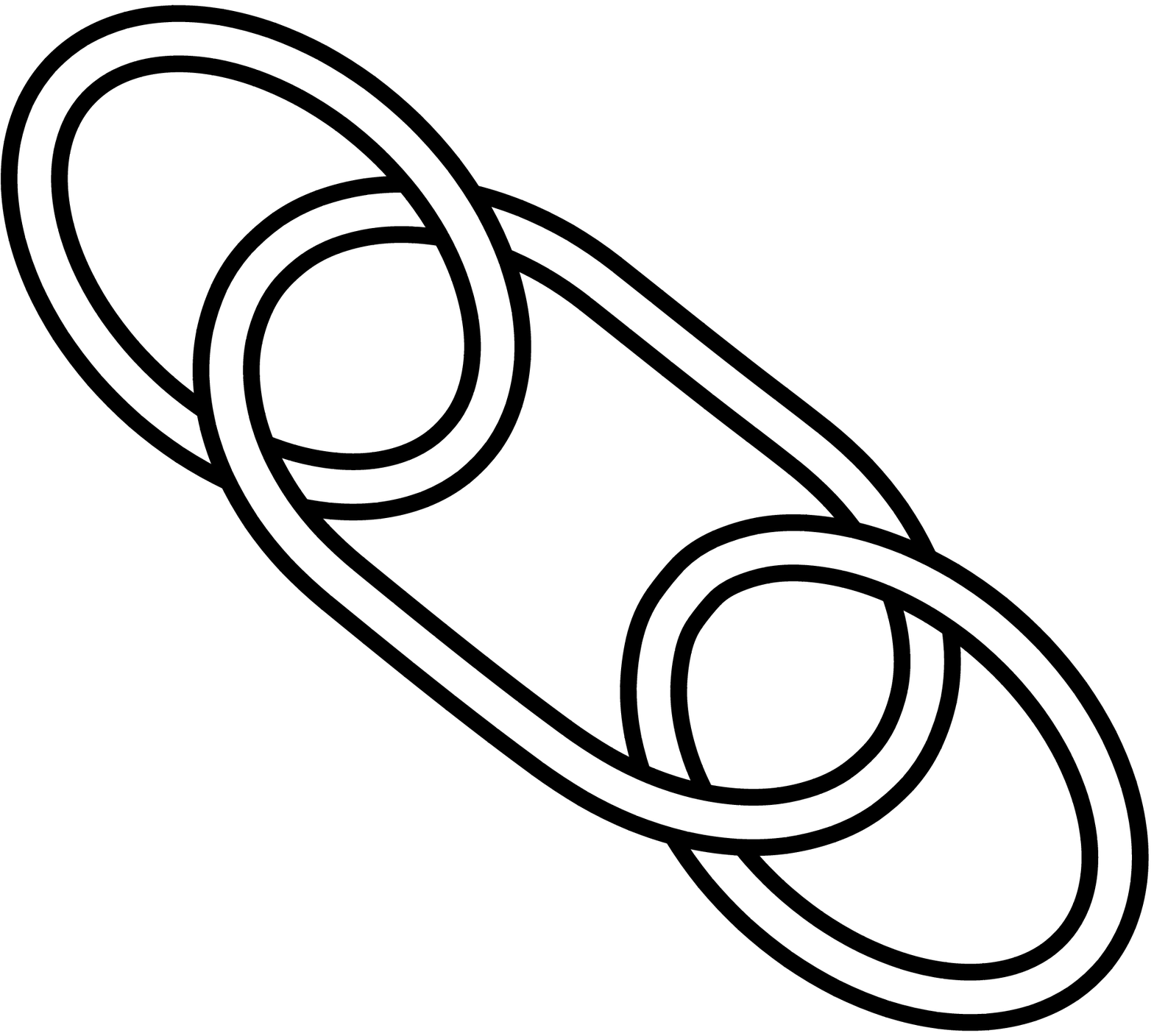} &
    \includegraphics[width=1in]{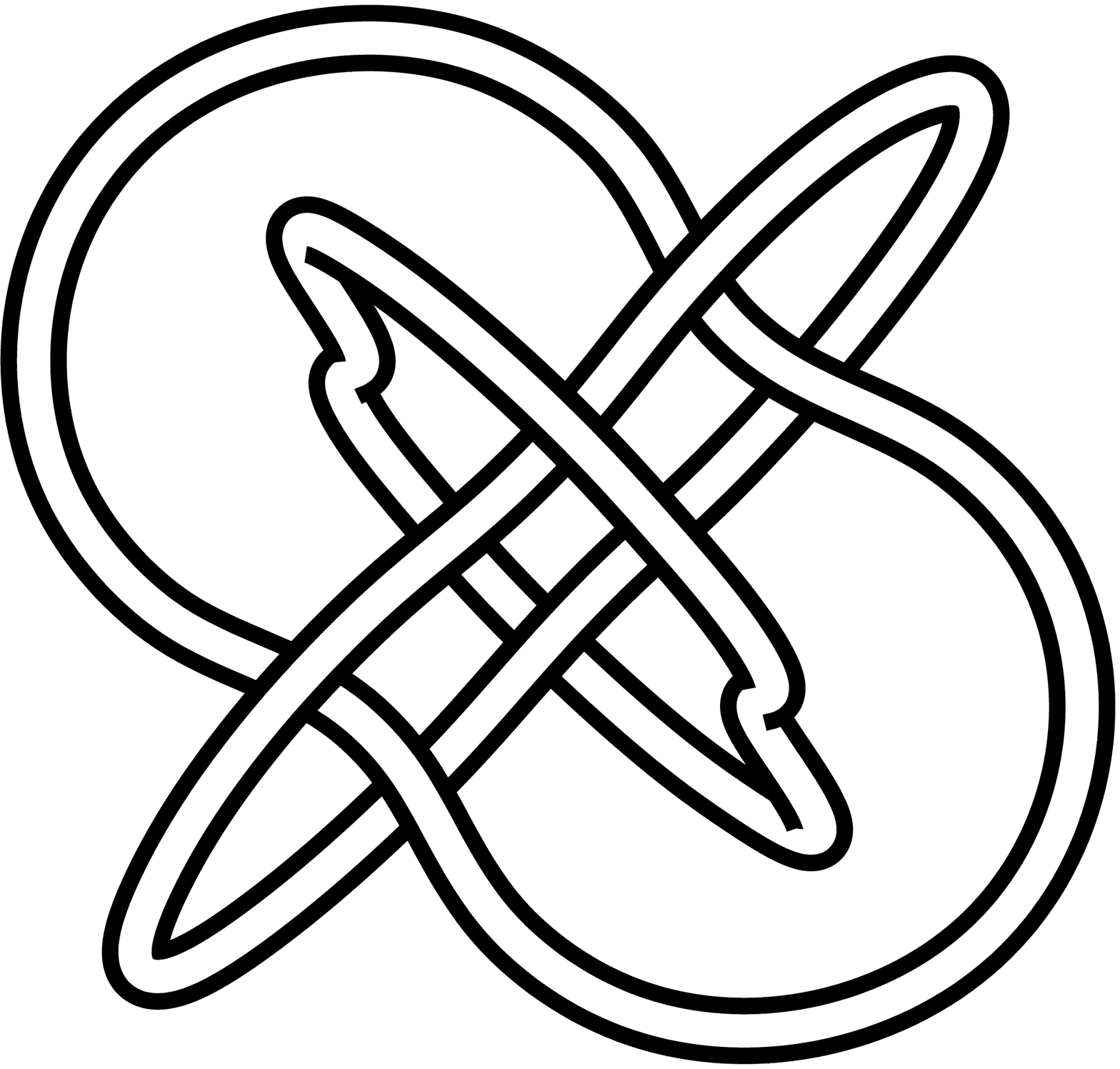} \\ \midrule
    Link name & Hopf link ($2_1$) & $3$-link chain ($2_1 \# 2_1$) & Borromean rings ($6^3_{2}$) \\
    Edges 		& $216$		& $384$		& $630$	\\
    Polygon length 	& $25.1439$	& $41.7131$	& $58.0300$\\
    Upper bound         & $25.1388$ & $\SimpleChainRRUB$   &    $58.0145$ \\
    Smooth length 	& $8\pi$		& $12\pi+4$	& $58.006$\\
    Relative error 	& $\twooneresult\%$	& $\simplechainresult\%$	& $\sixthreetworesult \%$ \\
    \bottomrule
  \end{tabular}
\end{center}
The results above lead us to suspect that most of our smooth ropelength bounds are within 1 to 2 hundreths of a percent of the corresponding minimum ropelength values, but we must be cautious.  The ropelength ``landscape'' for knots is quite complicated, and
we have already discovered a number of local ropelength minima for $8$ and
$9$-crossing knots which are very different from the (apparant) global
ropelength minima for these knot types. If our gradient descent algorithm has
been trapped by one of these local minima, the figures we report might be an
accurate computation of the ropelength of the local minimizer, but far off from
the true minimum ropelength value for the knot type.

We checked our computation of the self-contact sets by comparing our
computed self contact set for the Borromean rings to the contact set
for the ropelength-critical configuration provided by~\cite{cfksw}.
The results appear in the Figure below.

\begin{center}
  \begin{overpic}[width=4.5in]{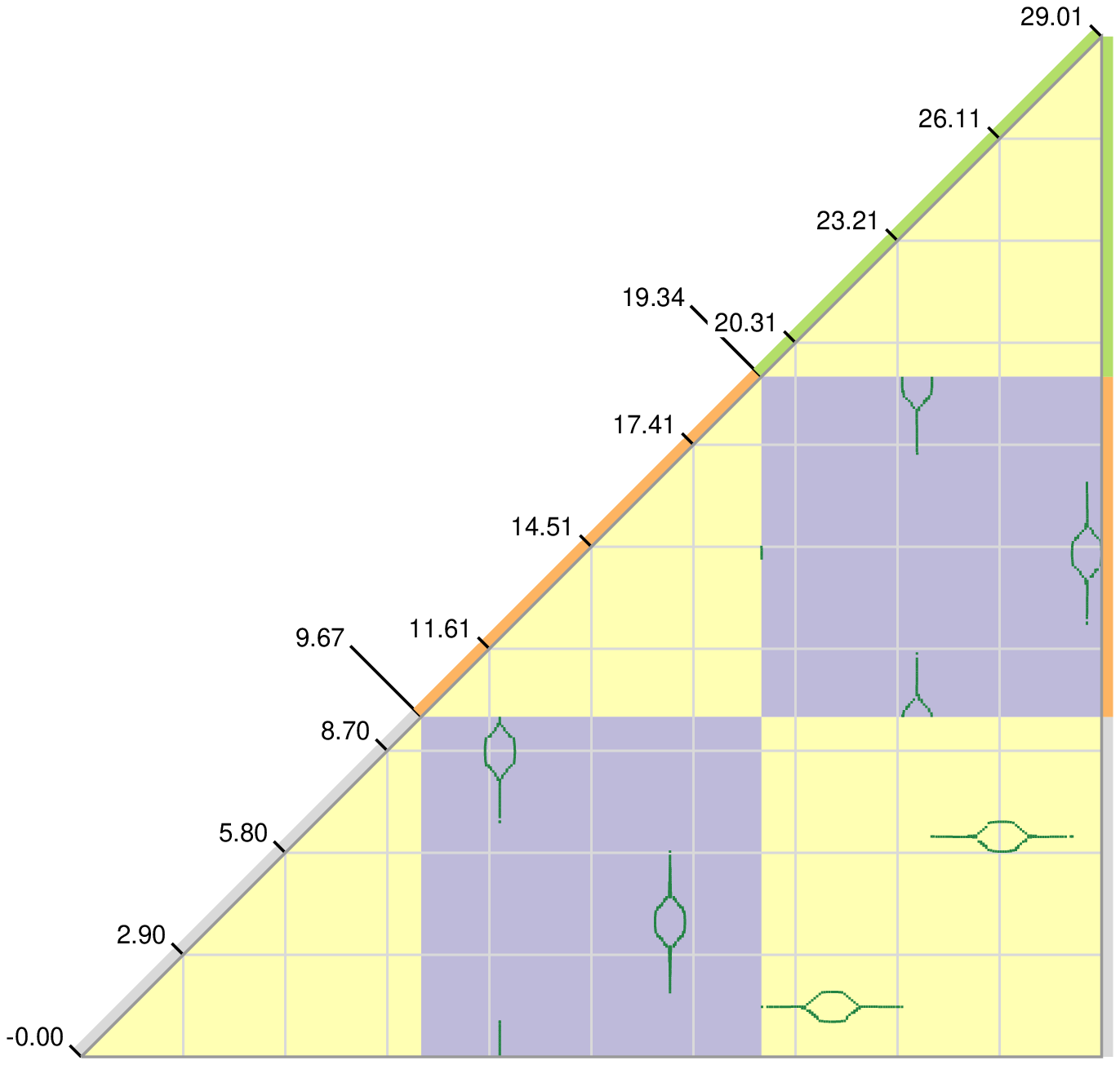}
  \end{overpic}\\
  \vspace{0.25in}
  \begin{overpic}[width=4.5in]{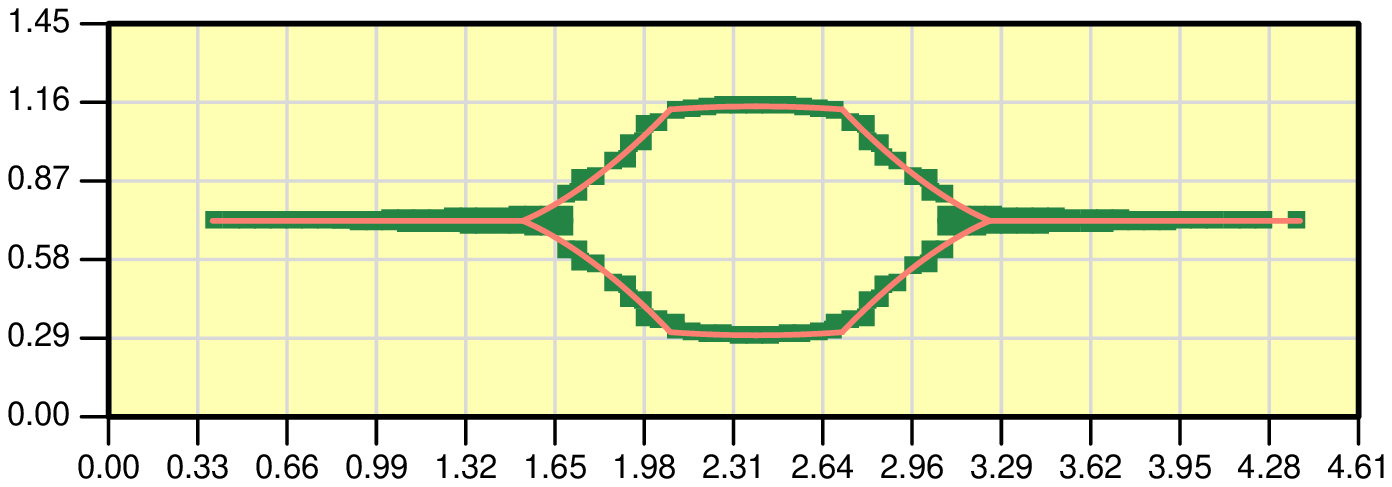}
  \end{overpic}
\end{center}

 The polygonal configuration of the Borromean rings discovered by our
 software has total length about $29.01$. The top plot shows the
 lower triangle of the square $[0,29.01] \cross [0,29.01]$
 representing pairs of arclength values $(s,t)$ on these polygonal
 rings. To describe the position of a point on this $3$-component link
 with a single arclength value, we use the convention that arclength
 values in $[0,9.67)$ refer to points on the first component, values
   in $[9.67,19.34)$ refer to points on the second component, and
     values in $[19.34,29.01)$ refer to points on the third component
       of the link. These breaks are reinforced by the colored bands
       running up the diagonal of the plot, which correspond to the
       colors of the different components of our tight links in the 3d
       renderings of Appendix B.

       The $s$ and $t$ position of points on the plot is indicated
       by the labels running up the diagonal, where the three
       special values representing breaks between components are
       lifted away from the other labels. The breaks between
       components are also shown on the plot by the alternating
       checkerboard pattern of the background colors. 
       
       We then locate every pair $(s,t) \in \dcsd(\mathcal{V})$ with
       $\norm{\mathcal{V}(s) - \mathcal{V}(t)} \leq \PThi(\mathcal{V})
       + 10^{-5}$. In the context of our gradient descent algorithm,
       these points represent active distance constraints. Our code
       computes Lagrange multipliers for all these active constraints,
       which represent self-contact forces borne by these tube
       contacts in the tight knot. We center a dark green box at every
       such $(s,t)$ value with a nonzero Lagrange multiplier.  The
       size of the box represents the average\footnote{Our final
         polygons are almost-equilateral, so this is a good
         approximation of the lengths of the edges incident to the pair
         $(s,t)$.}  edgelength of the polygon $\mathcal{V}$. We chose
       this size to represent the expected error in computed
       self-contact positions introduced by approximating a $C^1$
       ropelength-minimizer by the polygon $\mathcal{V}$.

       As we see from the plot, no tube around a component of the link
       is in contact with itself (so the three cream-colored triangles
       near the diagonal are empty). But each of the components makes
       contact with the other two, as shown by the boxes plotted in
       the purple and cream-colored rectangles forming the remainder
       of the plot. We can see that the contacts break up naturally
       into ``lantern-shaped'' structures.
       
       This link has been studied by Cantarella, Fu, Kusner, Sullivan,
       and Wrinkle, who provide a ropelength-critical configuration
       in~\cite{cfksw}. In the bottom plot, we compare one ``lantern''
       formed by $608$ of these boxes to the self-contact set
       predicted by these authors, which is represented by a red line.
       In this plot, the arclength distances labelled on the $s$ and
       $t$ axes do not correspond to a region of the plot above, but
       merely indicate the scale of the plot. We can see that the
       agreement between theory and computation is generally within
       one edgelength.
       
       Appendix B contains similar plots of our computed self-contact
       sets for $50$ knots and links from our collection of minimized
       examples. To the left of each self-contact plot, we provide a
       3d rendering of the corresponding tight shape. For links, the
       diagonal and right-hand side of the triangular self-contact
       plot are colored gray, orange and green in correspondance with
       the colors of the components in the rendering. The start of
       each component in the rendering is denoted by the tube coming
       to a point. The black bands on the tubes correspond to the
       arclength tick marks on the plot at right. The table also
       contains the polygonal ropelength of the knot (top number,
       slightly higher) and the corresponding smooth ropelength upper
       bound (bottom number, slightly lower), as well as the number of
       edges in the configuration plotted.

\section{Conclusions}

The major contribution of Appendix A is the provision of ropelength
figures for links. This allows us to check the accuracy of numerical
ropelength minimizations against theoretical results for the first
time. We are happy to report that our method passes this test for the
cases we examined.

The pictures in Appendix B are considerably more evocative.  It is
evident from first inspection that the contact sets of tight knots and
links seem to contain a fairly small number of commonly repeated
patterns. Some of these, such as the ``steps'' pattern first seen in
the trefoil knot, change shape from knot to knot. But others (such as
the ``winged'' pattern in the figure eight knot or the ``lantern''
shape seem in the Borromean rings) seem to remain remarkably
consistent throughout our computations. The reader may notice many
other examples as well. Isolating and understanding some of these
structures could provide us with a ``construction kit'' for tight
knots, much like the analysis of the simple chain
in~\cite{MR2003h:58014} led to the construction of an infinite family
of tight links.

\section{Acknowledgements} We would like to acknowledge many helpful
conversations with our colleagues during the years we spent working on
this project. In particular, Mark Peletier and Bob Planque, John
Sullivan, Erik Demaine and Bob Connelly, Herbert Edelsbrunner, and Joe
Fu all provided helpful insights and advice at various points in the
process. We were supported by NSF grants DMS-02-04826 (to Cantarella
and Fu), DMS-00-89927 (the University of Georgia VIGRE grant), and
DMS-03-11010 (to Rawdon).

\bibliography{drl,cantarella,strutpost}


\noindent
Ted Ashton and Jason Cantarella \\
Department of Mathematics, University of Georgia, Athens GA \\

\noindent
Michael Piatek \\
Department of Computer Science, University of Washington, Seattle WA \\

\noindent
Eric Rawdon \\
Department of Mathematics, Duquesne University, Pittsburgh, PA

\newpage

\appendix
\section{Table of polygonal ropelengths and ropelength upper bounds}

The table shows new polygonal ropelengths and ropelength upper bounds for $212$
knots and links with $9$ or fewer crossings. From left to right, there are four
columns: the name of the knot or link, the polygonal ropelength $\PRop$, the
corresponding ropelength upper bound $\Rop$, and the previous best ropelength
upper bound we could find in the literature together with the percentage improvement in ropelength. For these last figures, we used the papers~\cite{baranska,mglob,rawdonpc}. 

\newcommand{\rowmark}
{\addlinespace[1ex]
 \addlinespace[1ex]}

\begin{multicols}{2}
\begin{tabular}{l!{\hspace{-2ex}}c!{\hspace{1ex}}c!{\hspace{1ex}}l}
\toprule
  \hspace{-1ex} Link \hspace{1ex} & $\PRop$ & $\Rop$ & Previous  \\  
\midrule
   $ 2_{1}^{2} $ & $ 25.1439 $ & $ 25.1388 $ &    \\ 
   \rowmark
  $ 3_{1} $ & $ 32.7490 $ & $ 32.7448 $ & $ 32.7433864 \,\, ( - \%)$ \\ 
   \rowmark
  $ 4_{1} $ & $ 42.0997 $ & $ 42.0928 $ & $ 42.1158845 \,\, ( 0.05 \%)$ \\   \rowmark
  $ 4_{1}^{2} $ & $ 40.0247 $ & $ 40.0169 $ &    \\ 
  \rowmark
  $ 5_{1} $ & $ 47.2156 $ & $ 47.2016 $ & $ 47.51 \,\, ( 0.64 \%)$ \\ 
  $ 5_{2} $ & $ 49.4840 $ & $ 49.4704 $ & $ 49.73 \,\, ( 0.52 \%)$ \\ 
\rowmark
  $ 5_{1}^{2} $ & $ 49.7874 $ & $ 49.7723 $ &    \\ 
    \rowmark
  $ 6_{1} $ & $ 56.7316 $ & $ 56.7150 $ & $ 57.11 \,\, ( 0.69 \%)$ \\ 
  $ 6_{2} $ & $ 57.0451 $ & $ 57.0271 $ & $ 57.44 \,\, ( 0.71 \%)$ \\ 
  $ 6_{3} $ & $ 57.8602 $ & $ 57.8435 $ & $ 58.48 \,\, ( 1.08 \%)$ \\ 
\rowmark
  $ 6_{1}^{2} $ & $ 54.4068 $ & $ 54.3893 $ &    \\ 
  $ 6_{2}^{2} $ & $ 56.7132 $ & $ 56.7028 $ &    \\ 
  $ 6_{3}^{2} $ & $ 58.1161 $ & $ 58.1044 $ &    \\ 
\rowmark
  $ 6_{1}^{3} $ & $ 57.8334 $ & $ 57.8170 $ &    \\ 
  $ 6_{2}^{3} $ & $ 58.0300 $ & $ 58.0145 $ &    \\ 
  $ 6_{3}^{3} $ & $ 50.5865 $ & $ 50.5745 $ &    \\ 
\rowmark
  $ 7_{1} $ & $ 61.4319 $ & $ 61.4109 $ & $ 61.89 \,\, ( 0.77 \%)$ \\ 
  $ 7_{2} $ & $ 63.9165 $ & $ 63.8956 $ & $ 65.36 \,\, ( 2.24 \%)$ \\ 
  $ 7_{3} $ & $ 63.9539 $ & $ 63.9327 $ & $ 64.35 \,\, ( 0.64 \%)$ \\ 
  $ 7_{4} $ & $ 64.2960 $ & $ 64.2724 $ & $ 65.63 \,\, ( 2.06 \%)$ \\ 
  $ 7_{5} $ & $ 65.2802 $ & $ 65.2609 $ & $ 65.70 \,\, ( 0.66 \%)$ \\ 
  $ 7_{6} $ & $ 65.7183 $ & $ 65.7012 $ & $ 66.17 \,\, ( 0.7 \%)$ \\ 
  $ 7_{7} $ & $ 65.6316 $ & $ 65.6108 $ & $ 66.09 \,\, ( 0.72 \%)$ \\ 
\rowmark
  $ 7_{1}^{2} $ & $ 64.2585 $ & $ 64.2353 $ &    \\ 
  $ 7_{2}^{2} $ & $ 65.0467 $ & $ 65.0274 $ &    \\ 
  $ 7_{3}^{2} $ & $ 65.3743 $ & $ 65.3561 $ &    \\ 
  $ 7_{4}^{2} $ & $ 65.0971 $ & $ 65.0759 $ &    \\ 
  $ 7_{5}^{2} $ & $ 66.2400 $ & $ 66.2186 $ &    \\ 
  $ 7_{6}^{2} $ & $ 66.3494 $ & $ 66.3372 $ &    \\ 
  $ 7_{7}^{2} $ & $ 55.5451 $ & $ 55.5311 $ &    \\ 
  $ 7_{8}^{2} $ & $ 57.8043 $ & $ 57.7948 $ &    \\ 
\rowmark
  $ 7_{1}^{3} $ & $ 65.8275 $ & $ 65.8090 $ &    \\ 
\rowmark
  $ 8_{1} $ & $ 71.0484 $ & $ 71.0241 $ & $ 71.44 \,\, ( 0.58 \%)$ \\ 
  $ 8_{2} $ & $ 71.4327 $ & $ 71.4107 $ & $ 71.91 \,\, ( 0.69 \%)$ 
\end{tabular}

\begin{tabular}{l!{\hspace{-2ex}}c!{\hspace{1ex}}c!{\hspace{1ex}}l}
\toprule
  \hspace{-1ex} Link \hspace{1ex} & $\PRop$ & $\Rop$ & Previous  \\  
\midrule
$ 8_{3} $ & $ 71.1880 $ & $ 71.1655 $ & $ 71.56 \,\, ( 0.55 \%)$ \\ 
  $ 8_{4} $ & $ 72.0301 $ & $ 72.0049 $ & $ 72.41 \,\, ( 0.55 \%)$ \\ 
  $ 8_{5} $ & $ 72.2100 $ & $ 72.1878 $ & $ 72.70 \,\, ( 0.7 \%)$ \\ 
  $ 8_{6} $ & $ 72.5005 $ & $ 72.4791 $ & $ 72.93 \,\, ( 0.61 \%)$ \\ 
  $ 8_{7} $ & $ 72.2447 $ & $ 72.2204 $ & $ 72.63 \,\, ( 0.56 \%)$ \\ 
  $ 8_{8} $ & $ 73.3533 $ & $ 73.3334 $ & $ 73.88 \,\, ( 0.73 \%)$ \\ 
  $ 8_{9} $ & $ 72.4717 $ & $ 72.4461 $ & $ 72.96 \,\, ( 0.7 \%)$ \\ 
  $ 8_{10} $ & $ 73.4279 $ & $ 73.4095 $ & $ 73.86 \,\, ( 0.6 \%)$ \\ 
  $ 8_{11} $ & $ 73.5029 $ & $ 73.4802 $ & $ 76.70 \,\, ( 4.19 \%)$ \\ 
   $ 8_{12} $ & $ 74.0291 $ & $ 74.0098 $ & $ 74.61 \,\, ( 0.8 \%)$ \\ 
  $ 8_{13} $ & $ 72.8291 $ & $ 72.8045 $ & $ 73.29 \,\, ( 0.66 \%)$ \\ 
  $ 8_{14} $ & $ 73.9226 $ & $ 73.8991 $ & $ 74.93 \,\, ( 1.37 \%)$ \\ 
  $ 8_{15} $ & $ 74.3344 $ & $ 74.3134 $ & $ 74.82 \,\, ( 0.67 \%)$ \\ 
  $ 8_{16} $ & $ 74.9213 $ & $ 74.8962 $ & $ 75.47 \,\, ( 0.76 \%)$ \\ 
  $ 8_{17} $ & $ 74.5276 $ & $ 74.5071 $ & $ 75.08 \,\, ( 0.76 \%)$ \\ 
  $ 8_{18} $ & $ 74.9420 $ & $ 74.9252 $ & $ 75.44 \,\, ( 0.68 \%)$ \\ 
  $ 8_{19} $ & $ 61.0734 $ & $ 61.0430 $ & $ 61.35 \,\, ( 0.5 \%)$ \\ 
  $ 8_{20} $ & $ 63.1530 $ & $ 63.1146 $ & $ 64.11 \,\, ( 1.55 \%)$ \\ 
  $ 8_{21} $ & $ 65.5504 $ & $ 65.5298 $ & $ 65.91 \,\, ( 0.57 \%)$ \\ 
  \rowmark
  $ 8_{1}^{2} $ & $ 68.5198 $ & $ 68.4884 $ &    \\ 
  $ 8_{2}^{2} $ & $ 71.1823 $ & $ 71.1587 $ &    \\ 
  $ 8_{3}^{2} $ & $ 72.7498 $ & $ 72.7291 $ &    \\ 
  $ 8_{4}^{2} $ & $ 72.6102 $ & $ 72.5908 $ &    \\ 
  $ 8_{5}^{2} $ & $ 74.0039 $ & $ 73.9826 $ &    \\ 
  $ 8_{6}^{2} $ & $ 73.2932 $ & $ 73.2502 $ &    \\ 
  $ 8_{7}^{2} $ & $ 74.4165 $ & $ 74.3885 $ &    \\ 
  $ 8_{8}^{2} $ & $ 73.7849 $ & $ 73.7702 $ &    \\ 
  $ 8_{9}^{2} $ & $ 74.0620 $ & $ 74.0386 $ &    \\ 
  $ 8_{10}^{2} $ & $ 73.6890 $ & $ 73.6684 $ &    \\ 
  $ 8_{11}^{2} $ & $ 73.0115 $ & $ 72.9899 $ &    \\ 
  $ 8_{12}^{2} $ & $ 74.0194 $ & $ 73.9140 $ &    \\ 
  $ 8_{13}^{2} $ & $ 74.1685 $ & $ 74.1501 $ &    \\ 
  $ 8_{14}^{2} $ & $ 73.7000 $ & $ 73.6775 $ &    \\ 
  $ 8_{15}^{2} $ & $ 64.3305 $ & $ 64.3086 $ &    \\ 
  $ 8_{16}^{2} $ & $ 66.8434 $ & $ 66.8315 $ &    \\
\rowmark
 $ 8_{1}^{3} $ & $ 72.2883 $ & $ 72.2649 $ &    \\ 
  $ 8_{2}^{3} $ & $ 72.9544 $ & $ 72.9360 $ &    \\ 
 $ 8_{3}^{3} $ & $ 74.9366 $ & $ 74.9139 $ &    \\ 
  $ 8_{4}^{3} $ & $ 77.8544 $ & $ 77.8314 $ &    \\ 
  $ 8_{5}^{3} $ & $ 73.4286 $ & $ 73.4061 $ &    \\ 
  $ 8_{6}^{3} $ & $ 74.7680 $ & $ 74.7468 $ &    \\ 
  $ 8_{7}^{3} $ & $ 60.6065 $ & $ 60.5888 $ &    \\ 
  \end{tabular}

\begin{tabular}{l!{\hspace{-2ex}}c!{\hspace{1ex}}c!{\hspace{1ex}}l}
\toprule
  \hspace{-1ex} Link \hspace{1ex} & $\PRop$ & $\Rop$ & Previous  \\  
\midrule
  
  $ 8_{8}^{3} $ & $ 65.0637 $ & $ 65.0444 $ &    \\ 
  $ 8_{9}^{3} $ & $ 70.1904 $ & $ 70.1810 $ &    \\ 
  $ 8_{10}^{3} $ & $ 68.9823 $ & $ 68.9694 $ &    \\ 
  \rowmark
  $ 8_{1}^{4} $ & $ 75.2901 $ & $ 75.2677 $ &    \\ 
  $ 8_{2}^{4} $ & $ 67.4772 $ & $ 67.4571 $ &    \\ 
  $ 8_{3}^{4} $ & $ 66.4140 $ & $ 66.4046 $ &    \\ 
  \rowmark
  $ 9_{1} $ & $ 75.7507 $ & $ 75.7252 $ & $ 76.43 \,\, ( 0.92 \%)$ \\ 
  $ 9_{2} $ & $ 79.3059 $ & $ 79.2794 $ & $ 79.92 \,\, ( 0.8 \%)$ \\ 
  $ 9_{3} $ & $ 78.5867 $ & $ 78.5591 $ & $ 79.05 \,\, ( 0.62 \%)$ \\ 
  $ 9_{4} $ & $ 78.4117 $ & $ 78.3861 $ & $ 78.84 \,\, ( 0.57 \%)$ \\ 
  $ 9_{5} $ & $ 79.7877 $ & $ 79.7586 $ & $ 80.32 \,\, ( 0.69 \%)$ \\ 
  $ 9_{6} $ & $ 80.1326 $ & $ 80.0822 $ & $ 80.65 \,\, ( 0.7 \%)$ \\ 
  $ 9_{7} $ & $ 80.6597 $ & $ 80.6357 $ & $ 82.65 \,\, ( 2.43 \%)$ \\ 
   $ 9_{8} $ & $ 80.5651 $ & $ 80.5384 $ & $ 81.14 \,\, ( 0.74 \%)$ \\ 
  $ 9_{9} $ & $ 79.9574 $ & $ 79.9323 $ & $ 80.85 \,\, ( 1.13 \%)$ \\ 
  $ 9_{10} $ & $ 79.8161 $ & $ 79.7946 $ & $ 80.33 \,\, ( 0.66 \%)$ \\ 
  $ 9_{11} $ & $ 80.3650 $ & $ 80.3418 $ & $ 81.98 \,\, ( 1.99 \%)$ \\ 
  $ 9_{12} $ & $ 80.1275 $ & $ 80.1047 $ & $ 80.71 \,\, ( 0.74 \%)$ \\ 
  $ 9_{13} $ & $ 80.6525 $ & $ 80.6246 $ & $ 81.33 \,\, ( 0.86 \%)$ \\ 
  $ 9_{14} $ & $ 80.1512 $ & $ 80.1259 $ & $ 80.73 \,\, ( 0.74 \%)$ \\ 
  $ 9_{15} $ & $ 82.1665 $ & $ 82.1396 $ & $ 82.70 \,\, ( 0.67 \%)$ \\ 
  $ 9_{16} $ & $ 80.1446 $ & $ 80.1211 $ & $ 80.67 \,\, ( 0.68 \%)$ \\ 
  $ 9_{17} $ & $ 80.5792 $ & $ 80.5537 $ & $ 81.90 \,\, ( 1.64 \%)$ \\ 
  $ 9_{18} $ & $ 81.6230 $ & $ 81.5960 $ & $ 82.68 \,\, ( 1.31 \%)$ \\ 
  $ 9_{19} $ & $ 82.1828 $ & $ 82.1593 $ & $ 82.72 \,\, ( 0.67 \%)$ \\ 
  $ 9_{20} $ & $ 80.2543 $ & $ 80.2288 $ & $ 87.31 \,\, ( 8.11 \%)$ \\ 
  $ 9_{21} $ & $ 81.1336 $ & $ 81.1098 $ & $ 81.64 \,\, ( 0.64 \%)$ \\ 
  $ 9_{22} $ & $ 81.0714 $ & $ 81.0464 $ & $ 81.60 \,\, ( 0.67 \%)$ \\ 
  $ 9_{23} $ & $ 81.3164 $ & $ 81.2898 $ & $ 81.84 \,\, ( 0.67 \%)$ \\ 
  $ 9_{24} $ & $ 80.9933 $ & $ 80.9701 $ & $ 81.54 \,\, ( 0.69 \%)$ \\ 
  $ 9_{25} $ & $ 81.1873 $ & $ 81.1612 $ & $ 81.85 \,\, ( 0.84 \%)$ \\ 
  $ 9_{26} $ & $ 80.9352 $ & $ 80.9137 $ & $ 81.94 \,\, ( 1.25 \%)$ \\ 
  $ 9_{27} $ & $ 81.9473 $ & $ 81.9092 $ & $ 83.21 \,\, ( 1.56 \%)$ \\ 
  $ 9_{28} $ & $ 81.5694 $ & $ 81.5467 $ & $ 82.25 \,\, ( 0.85 \%)$ \\ 
  $ 9_{29} $ & $ 81.8708 $ & $ 81.8470 $ & $ 83.45 \,\, ( 1.92 \%)$ \\ 
  $ 9_{30} $ & $ 81.8567 $ & $ 81.8344 $ & $ 82.46 \,\, ( 0.75 \%)$ \\ 
 $ 9_{31} $ & $ 81.7436 $ & $ 81.6915 $ & $ 82.22 \,\, ( 0.64 \%)$ \\ 
  $ 9_{32} $ & $ 81.5775 $ & $ 81.5555 $ & $ 82.34 \,\, ( 0.95 \%)$ \\ 
  $ 9_{33} $ & $ 82.8451 $ & $ 82.7989 $ & $ 83.37 \,\, ( 0.68 \%)$ \\ 
  $ 9_{34} $ & $ 82.3213 $ & $ 82.2744 $ & $ 82.99 \,\, ( 0.86 \%)$ \\ 
  $ 9_{35} $ & $ 79.2495 $ & $ 79.2216 $ & $ 80.85 \,\, ( 2.01 \%)$ \\ 
  $ 9_{36} $ & $ 81.0579 $ & $ 81.0297 $ & $ 81.57 \,\, ( 0.66 \%)$ \\ 
  $ 9_{37} $ & $ 81.5845 $ & $ 81.5562 $ & $ 82.10 \,\, ( 0.66 \%)$ \\ 
  $ 9_{38} $ & $ 81.8119 $ & $ 81.7909 $ & $ 82.43 \,\, ( 0.77 \%)$ \\ 
  $ 9_{39} $ & $ 81.9490 $ & $ 81.9266 $ & $ 85.55 \,\, ( 4.23 \%)$ \\ 
  $ 9_{40} $ & $ 81.7008 $ & $ 81.6806 $ & $ 82.67 \,\, ( 1.19 \%)$ \\ 
  $ 9_{41} $ & $ 81.4929 $ & $ 81.4399 $ & $ 82.11 \,\, ( 0.81 \%)$ \\ 
  $ 9_{42} $ & $ 69.6133 $ & $ 69.5939 $ & $ 70.02 \,\, ( 0.6 \%)$ \\
 $ 9_{43} $ & $ 71.7062 $ & $ 71.6863 $ & $ 72.20 \,\, ( 0.71 \%)$ \\  
\end{tabular}

\begin{tabular}{l!{\hspace{-2ex}}c!{\hspace{1ex}}c!{\hspace{1ex}}l}
\toprule
  \hspace{-1ex} Link \hspace{1ex} & $\PRop$ & $\Rop$ & Previous  \\  
\midrule
 
  $ 9_{44} $ & $ 71.6516 $ & $ 71.6305 $ & $ 72.23 \,\, ( 0.82 \%)$ \\ 
  $ 9_{45} $ & $ 74.9154 $ & $ 74.8959 $ & $ 75.51 \,\, ( 0.81 \%)$ \\ 
  $ 9_{46} $ & $ 68.6579 $ & $ 68.6369 $ & $ 69.35 \,\, ( 1.02 \%)$ \\ 
  $ 9_{47} $ & $ 75.1289 $ & $ 75.0875 $ & $ 75.61 \,\, ( 0.69 \%)$ \\ 
  $ 9_{48} $ & $ 74.2918 $ & $ 74.2477 $ & $ 74.94 \,\, ( 0.92 \%)$ \\ 
   $ 9_{49} $ & $ 74.0530 $ & $ 74.0127 $ & $ 74.50 \,\, ( 0.65 \%)$ \\ 
   \rowmark
  $ 9_{1}^{2} $ & $ 78.7014 $ & $ 78.6731 $ &    \\ 
  $ 9_{2}^{2} $ & $ 79.5525 $ & $ 79.5259 $ &    \\ 
  $ 9_{3}^{2} $ & $ 79.9725 $ & $ 79.9476 $ &    \\ 
  $ 9_{4}^{2} $ & $ 78.7248 $ & $ 78.6967 $ &    \\ 
  $ 9_{5}^{2} $ & $ 79.6807 $ & $ 79.6549 $ &    \\ 
  $ 9_{6}^{2} $ & $ 81.1029 $ & $ 81.0791 $ &    \\ 
  $ 9_{7}^{2} $ & $ 81.1705 $ & $ 81.1460 $ &    \\ 
  $ 9_{8}^{2} $ & $ 81.0999 $ & $ 81.0723 $ &    \\ 
  $ 9_{9}^{2} $ & $ 80.3391 $ & $ 80.3181 $ &    \\ 
  $ 9_{10}^{2} $ & $ 80.3964 $ & $ 80.3693 $ &    \\ 
  $ 9_{11}^{2} $ & $ 82.0513 $ & $ 82.0271 $ &    \\ 
  $ 9_{12}^{2} $ & $ 81.9983 $ & $ 81.9738 $ &    \\ 
  $ 9_{13}^{2} $ & $ 79.3586 $ & $ 79.3319 $ &    \\ 
  $ 9_{14}^{2} $ & $ 80.7420 $ & $ 80.7171 $ &    \\ 
  $ 9_{15}^{2} $ & $ 80.5805 $ & $ 80.5550 $ &    \\ 
  $ 9_{16}^{2} $ & $ 81.4190 $ & $ 81.3966 $ &    \\ 
 $ 9_{17}^{2} $ & $ 81.2589 $ & $ 81.2307 $ &    \\ 
  $ 9_{18}^{2} $ & $ 82.2739 $ & $ 82.2429 $ &    \\ 
  $ 9_{19}^{2} $ & $ 79.4899 $ & $ 79.4628 $ &    \\ 
  $ 9_{20}^{2} $ & $ 80.4176 $ & $ 80.3956 $ &    \\ 
  $ 9_{21}^{2} $ & $ 81.2545 $ & $ 81.2356 $ &    \\ 
  $ 9_{22}^{2} $ & $ 81.1318 $ & $ 81.1058 $ &    \\ 
  $ 9_{23}^{2} $ & $ 80.4552 $ & $ 80.4290 $ &    \\ 
  $ 9_{24}^{2} $ & $ 82.5780 $ & $ 82.5207 $ &    \\ 
  $ 9_{25}^{2} $ & $ 81.8205 $ & $ 81.7924 $ &    \\ 
  $ 9_{26}^{2} $ & $ 82.1200 $ & $ 82.1015 $ &    \\ 
  $ 9_{27}^{2} $ & $ 81.3442 $ & $ 81.2046 $ &    \\ 
  $ 9_{28}^{2} $ & $ 81.3869 $ & $ 81.3740 $ &    \\ 
  $ 9_{29}^{2} $ & $ 82.1780 $ & $ 82.1567 $ &    \\ 
  $ 9_{30}^{2} $ & $ 82.2789 $ & $ 82.2517 $ &    \\ 
  $ 9_{31}^{2} $ & $ 80.6125 $ & $ 80.5936 $ &    \\ 
  $ 9_{32}^{2} $ & $ 81.4352 $ & $ 81.4122 $ &    \\ 
  $ 9_{33}^{2} $ & $ 82.2167 $ & $ 82.1951 $ &    \\ 
 $ 9_{34}^{2} $ & $ 81.8743 $ & $ 81.8552 $ &    \\ 
  $ 9_{35}^{2} $ & $ 81.3504 $ & $ 81.3284 $ &    \\ 
  $ 9_{36}^{2} $ & $ 80.7280 $ & $ 80.7103 $ &    \\ 
  $ 9_{37}^{2} $ & $ 81.9353 $ & $ 81.9101 $ &    \\ 
  $ 9_{38}^{2} $ & $ 82.7480 $ & $ 82.7206 $ &    \\ 
  $ 9_{39}^{2} $ & $ 81.9251 $ & $ 81.8973 $ &    \\ 
  $ 9_{40}^{2} $ & $ 82.0090 $ & $ 81.9822 $ &    \\ 
  $ 9_{41}^{2} $ & $ 83.6167 $ & $ 83.5979 $ &    \\ 
  $ 9_{42}^{2} $ & $ 83.6493 $ & $ 83.6292 $ &    \\ 
  $ 9_{43}^{2} $ & $ 66.3264 $ & $ 66.3049 $ &    \\ 
   $ 9_{44}^{2} $ & $ 72.2518 $ & $ 72.2255 $ &    \\ 
\end{tabular}

\begin{tabular}{l!{\hspace{-2ex}}c!{\hspace{1ex}}c!{\hspace{1ex}}l}
\toprule
  \hspace{-1ex} Link \hspace{1ex} & $\PRop$ & $\Rop$ & Previous  \\  
\midrule
 $ 9_{45}^{2} $ & $ 71.3593 $ & $ 71.3490 $ &    \\ 
  $ 9_{46}^{2} $ & $ 73.9646 $ & $ 73.9377 $ &    \\ 
  $ 9_{47}^{2} $ & $ 69.9602 $ & $ 69.9432 $ &    \\ 
  $ 9_{48}^{2} $ & $ 73.6781 $ & $ 73.6545 $ &    \\ 
  $ 9_{49}^{2} $ & $ 66.0806 $ & $ 66.0658 $ &    \\ 
  $ 9_{50}^{2} $ & $ 69.3690 $ & $ 69.3483 $ &    \\ 
  $ 9_{51}^{2} $ & $ 70.5942 $ & $ 70.5699 $ &    \\ 
  $ 9_{52}^{2} $ & $ 72.9916 $ & $ 72.9685 $ &    \\ 
  $ 9_{53}^{2} $ & $ 68.0369 $ & $ 68.0305 $ &    \\ 
  $ 9_{54}^{2} $ & $ 71.0385 $ & $ 71.0185 $ &    \\ 
  $ 9_{55}^{2} $ & $ 73.8423 $ & $ 73.8217 $ &    \\ 
  $ 9_{56}^{2} $ & $ 75.2439 $ & $ 75.2245 $ &    \\ 
  $ 9_{57}^{2} $ & $ 73.8408 $ & $ 73.8194 $ &    \\ 
  $ 9_{58}^{2} $ & $ 74.1911 $ & $ 74.1697 $ &    \\ 
  $ 9_{59}^{2} $ & $ 73.0481 $ & $ 73.0305 $ &    \\ 
  $ 9_{60}^{2} $ & $ 73.5734 $ & $ 73.5553 $ &    \\ 
  $ 9_{61}^{2} $ & $ 69.3978 $ & $ 69.3840 $ &    \\ 
  \rowmark
  $ 9_{1}^{3} $ & $ 81.2897 $ & $ 81.2323 $ &    \\ 
  $ 9_{2}^{3} $ & $ 82.4507 $ & $ 82.4004 $ &    \\ 
  $ 9_{3}^{3} $ & $ 82.3127 $ & $ 82.2861 $ &    \\ 
  $ 9_{4}^{3} $ & $ 82.5449 $ & $ 82.5208 $ &    \\ 
  $ 9_{5}^{3} $ & $ 80.7974 $ & $ 80.7714 $ &    \\ 
  $ 9_{6}^{3} $ & $ 81.0235 $ & $ 81.0054 $ &    \\ 
  $ 9_{7}^{3} $ & $ 82.1986 $ & $ 82.1535 $ &    \\ 
  $ 9_{8}^{3} $ & $ 81.1631 $ & $ 81.1408 $ &    \\ 
  $ 9_{9}^{3} $ & $ 81.6200 $ & $ 81.5735 $ &    \\ 
  $ 9_{10}^{3} $ & $ 82.3446 $ & $ 82.3259 $ &    \\ 
  $ 9_{11}^{3} $ & $ 82.0323 $ & $ 82.0137 $ &    \\ 
  $ 9_{12}^{3} $ & $ 82.6345 $ & $ 82.5740 $ &    \\ 
  $ 9_{13}^{3} $ & $ 72.2098 $ & $ 72.2008 $ &    \\ 
  $ 9_{14}^{3} $ & $ 74.5697 $ & $ 74.5492 $ &    \\ 
  $ 9_{15}^{3} $ & $ 74.3877 $ & $ 74.3655 $ &    \\ 
  $ 9_{16}^{3} $ & $ 75.0664 $ & $ 75.0430 $ &    \\ 
  $ 9_{17}^{3} $ & $ 74.2972 $ & $ 74.2779 $ &    \\ 
  $ 9_{18}^{3} $ & $ 72.5059 $ & $ 72.4741 $ &    \\ 
  $ 9_{19}^{3} $ & $ 72.7143 $ & $ 72.6859 $ &    \\ 
  $ 9_{20}^{3} $ & $ 76.3557 $ & $ 76.1829 $ &    \\ 
  $ 9_{21}^{3} $ & $ 74.9369 $ & $ 74.9212 $ &    \\ 
  \rowmark
  $ 9_{1}^{4} $ & $ 85.5620 $ & $ 85.5115 $ &    \\ 
\end{tabular}

\end{multicols}

\newpage

\graphicspath{{data/texfiles/}}
\section{Table of Self-Contact Sets}
\begin{center}
\noindent
\begin{minipage}[t]{6in}
  \vspace{2mm}
    \begin{overpic}[width=2.65in]{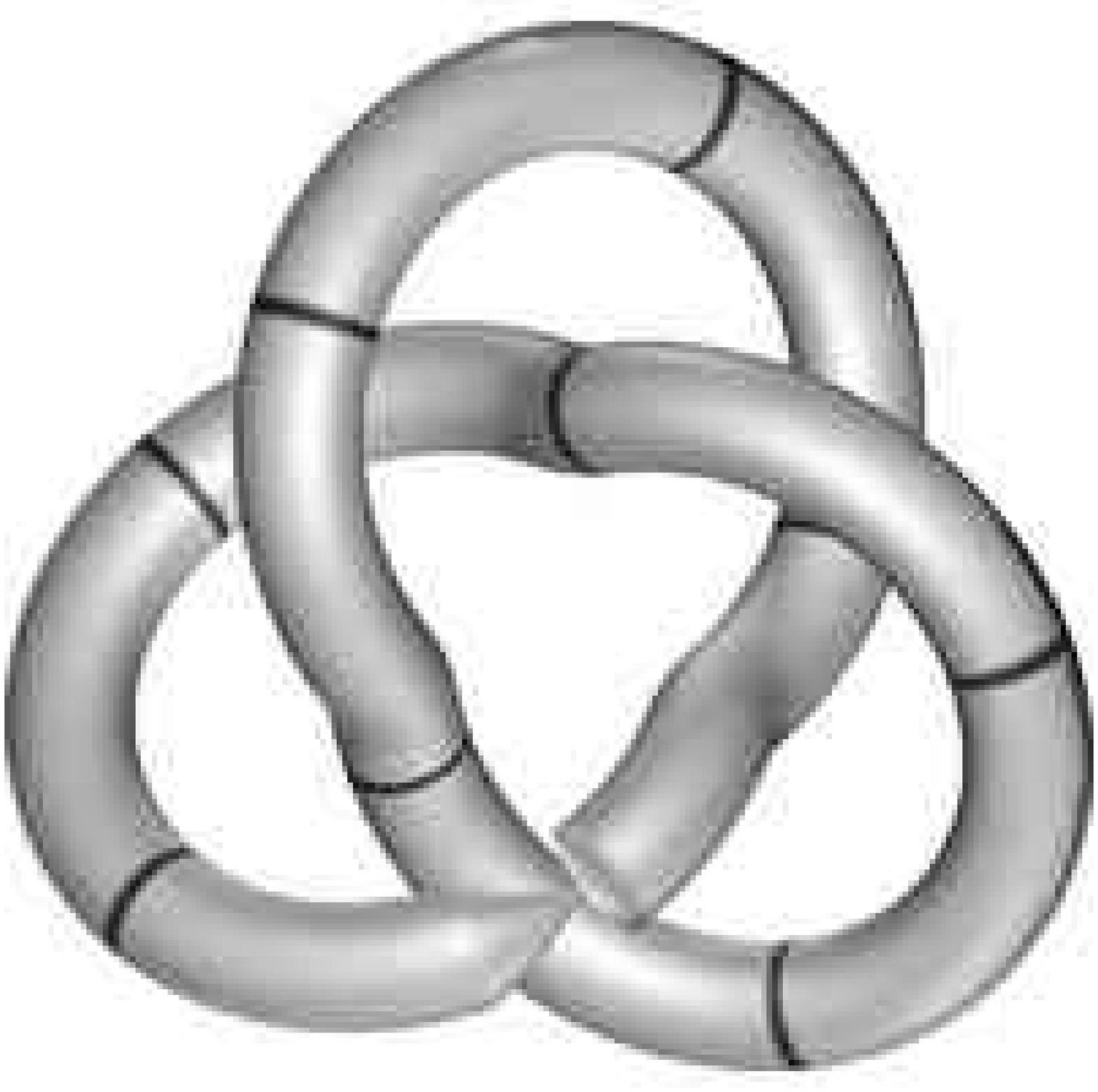}
        \put(-10,90){\large{$3_{1}$}}
    \end{overpic}
      \hspace{7mm}
    \begin{overpic}[width=2.65in]{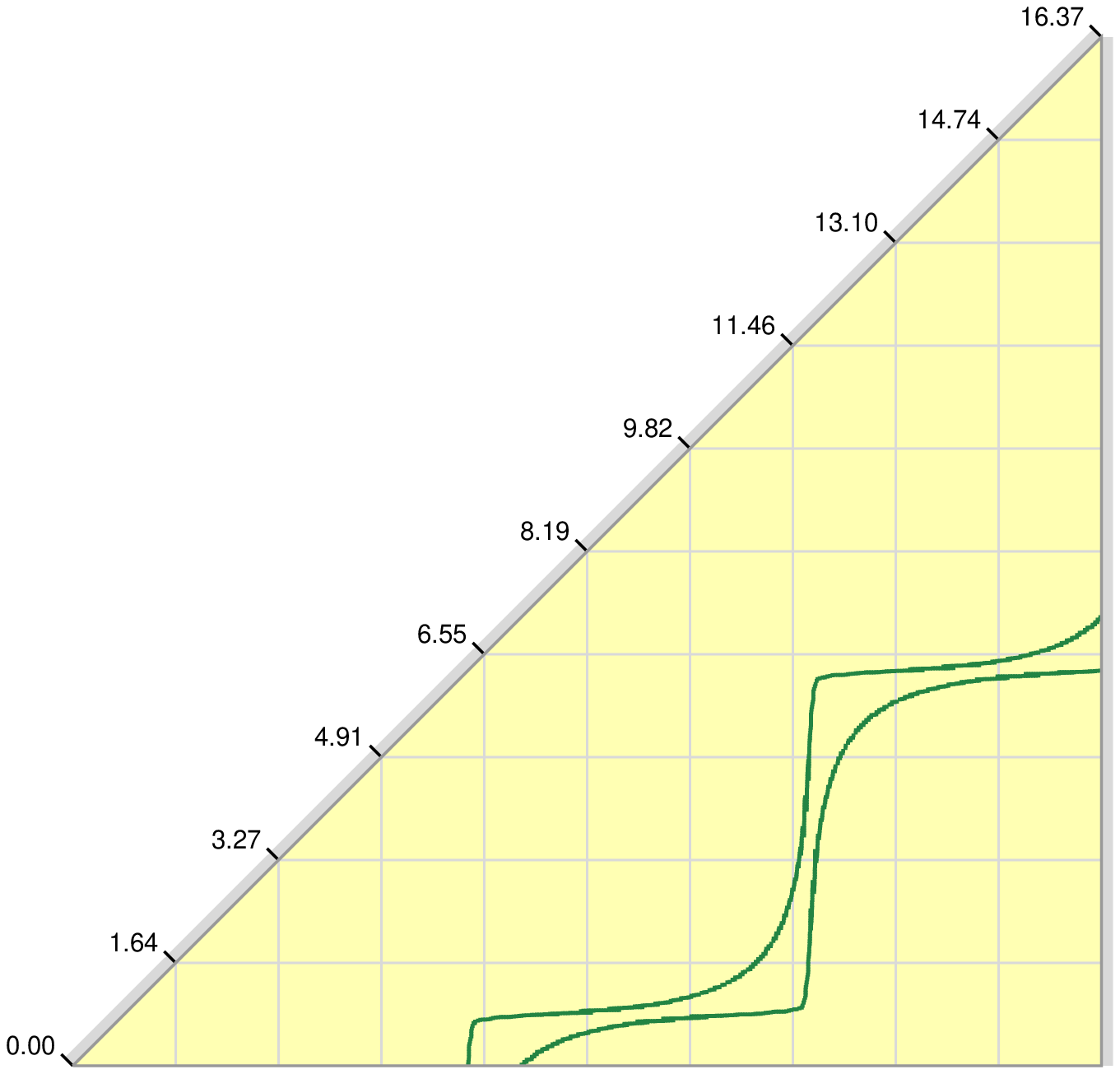}
        \put(8,94){\scriptsize{Poly$(K)$: $32.75$}}
        \put(8,89){\scriptsize{Bound: $32.75$}}
        \put(8,84){\scriptsize{Vertices: $400$}}
    \end{overpic}
\end{minipage} 
\hfill
\begin{minipage}[t]{6in}
  \vspace{2mm}
    \begin{overpic}[height=2.65in]{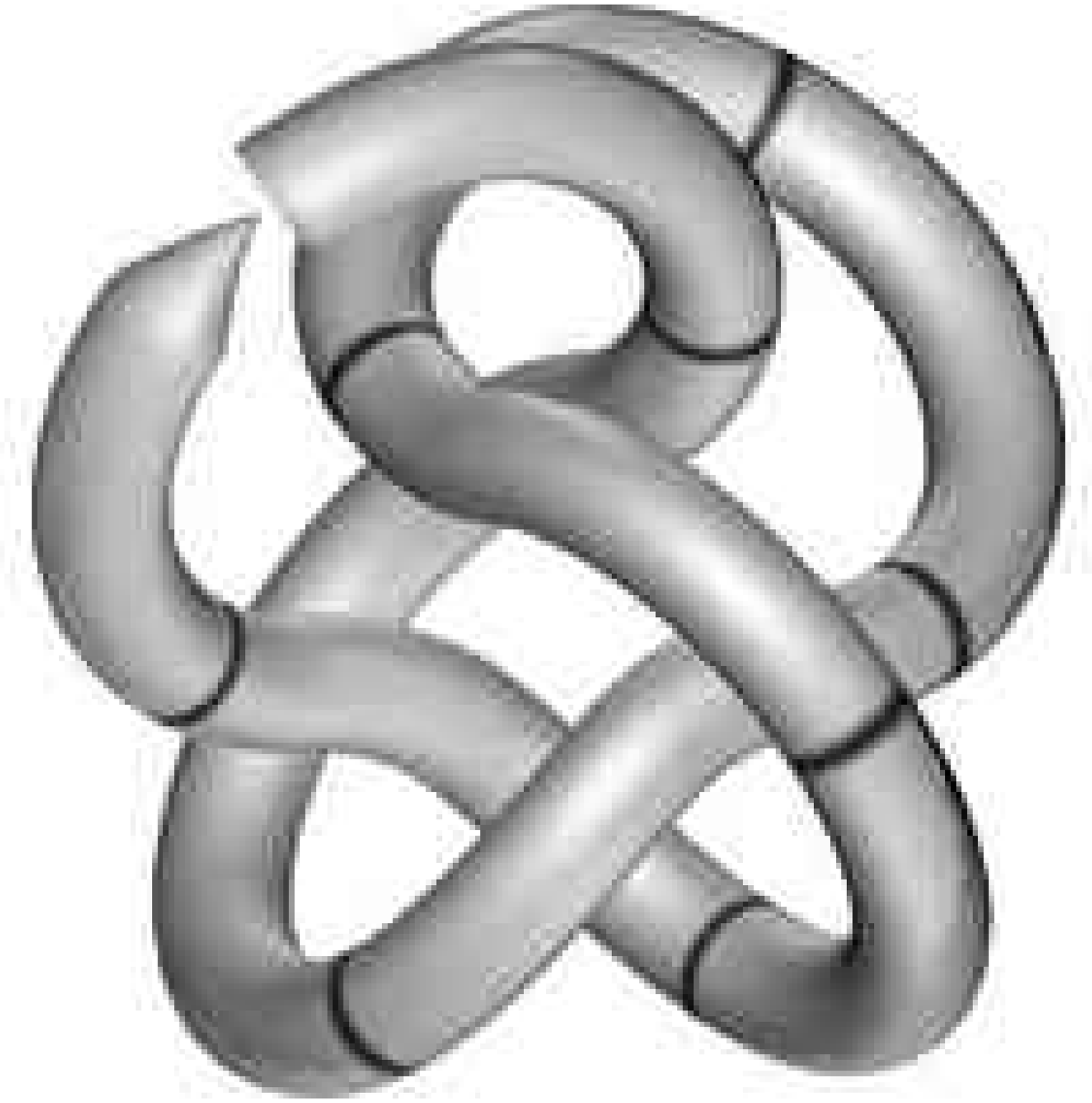}
        \put(-10,90){\large{$4_{1}$}}
    \end{overpic}
      \hspace{7mm}
    \begin{overpic}[width=2.65in]{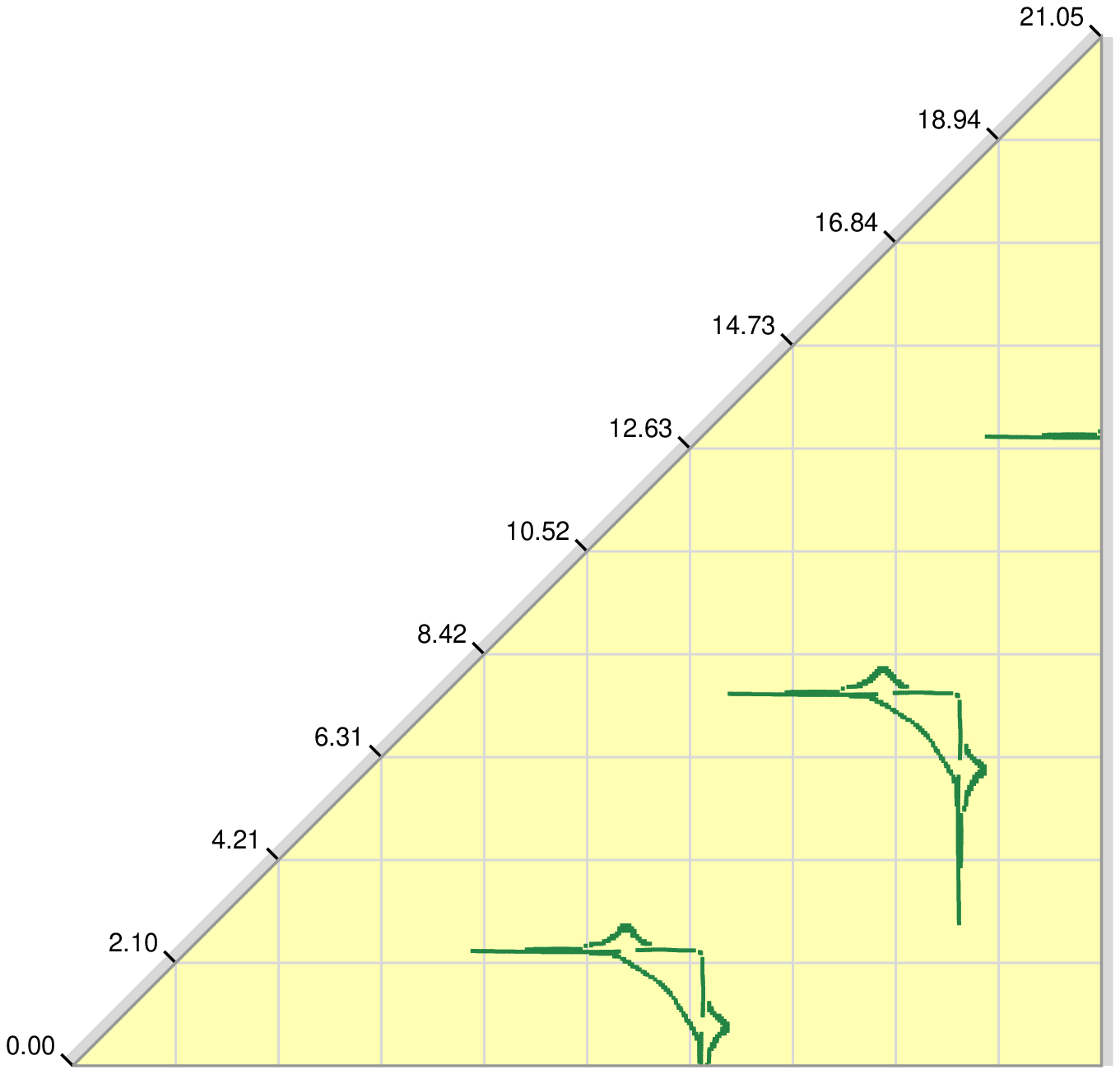}
        \put(8,94){\scriptsize{$42.10$}}
        \put(8,89){\scriptsize{$42.10$}}
        \put(8,84){\scriptsize{$400$}}
    \end{overpic}
\end{minipage} 
\hfill
\begin{minipage}[t]{6in}
  \vspace{2mm}
    \begin{overpic}[width=2.65in]{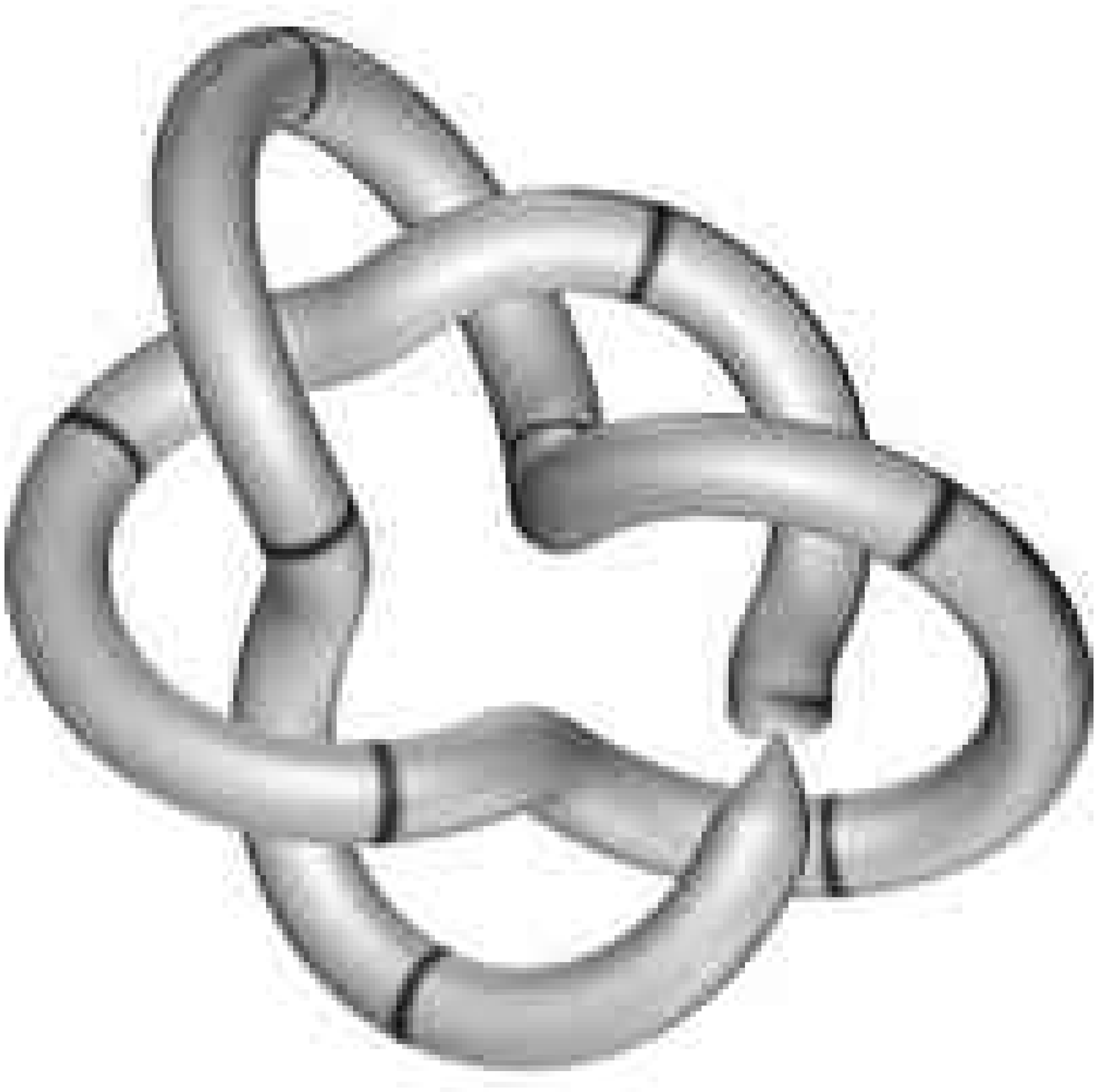}
        \put(-10,90){\large{$5_{1}$}}
    \end{overpic}
      \hspace{7mm}
    \begin{overpic}[width=2.65in]{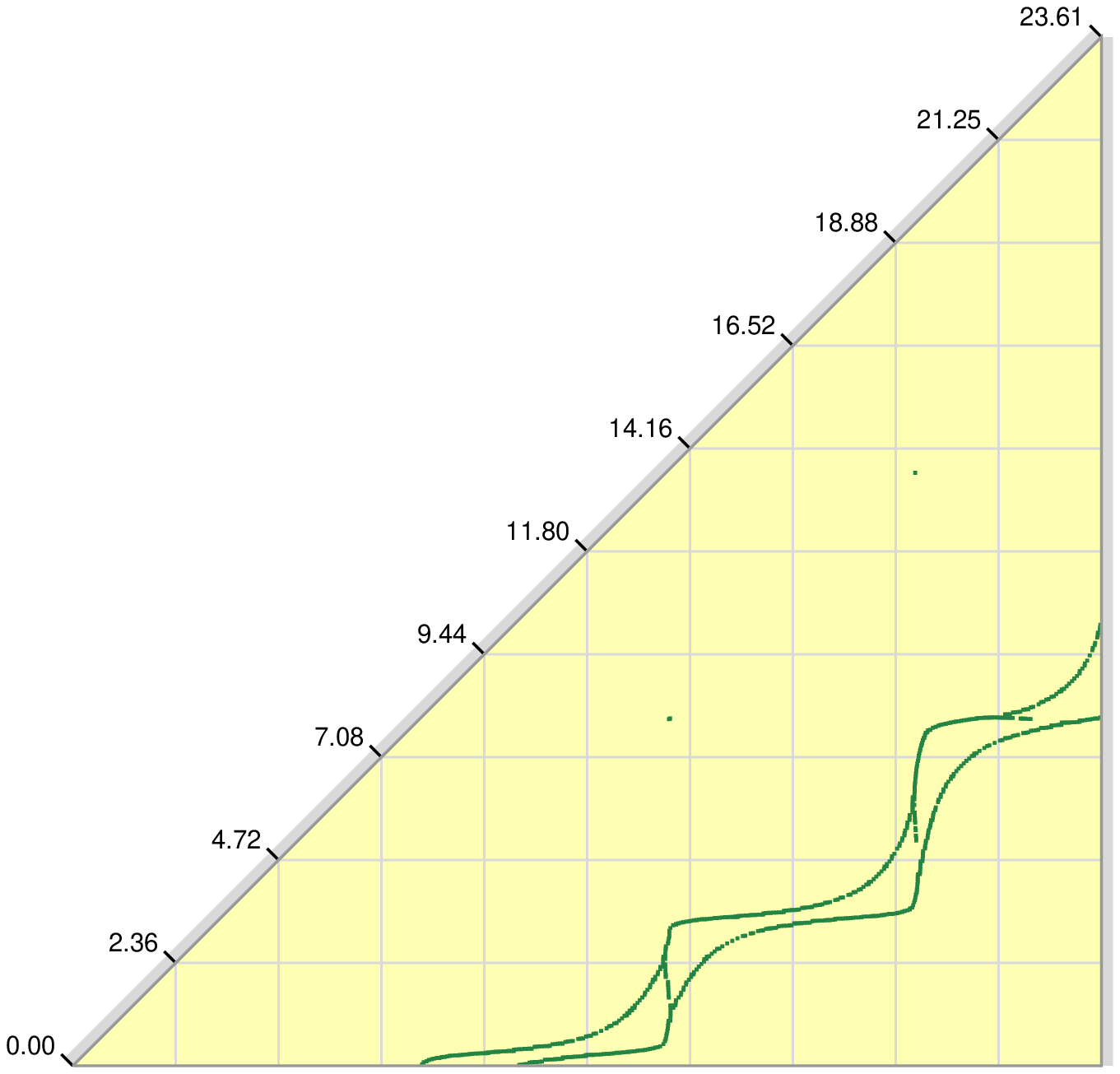}
        \put(8,94){\scriptsize{$47.22$}}
        \put(8,89){\scriptsize{$47.21$}}
        \put(8,84){\scriptsize{$400$}}
    \end{overpic}
\end{minipage} 
\hfill
\clearpage
\pagebreak
\begin{figure}
\begin{minipage}[t]{6in}
  \vspace{2mm}
    \begin{overpic}[width=2.8in]{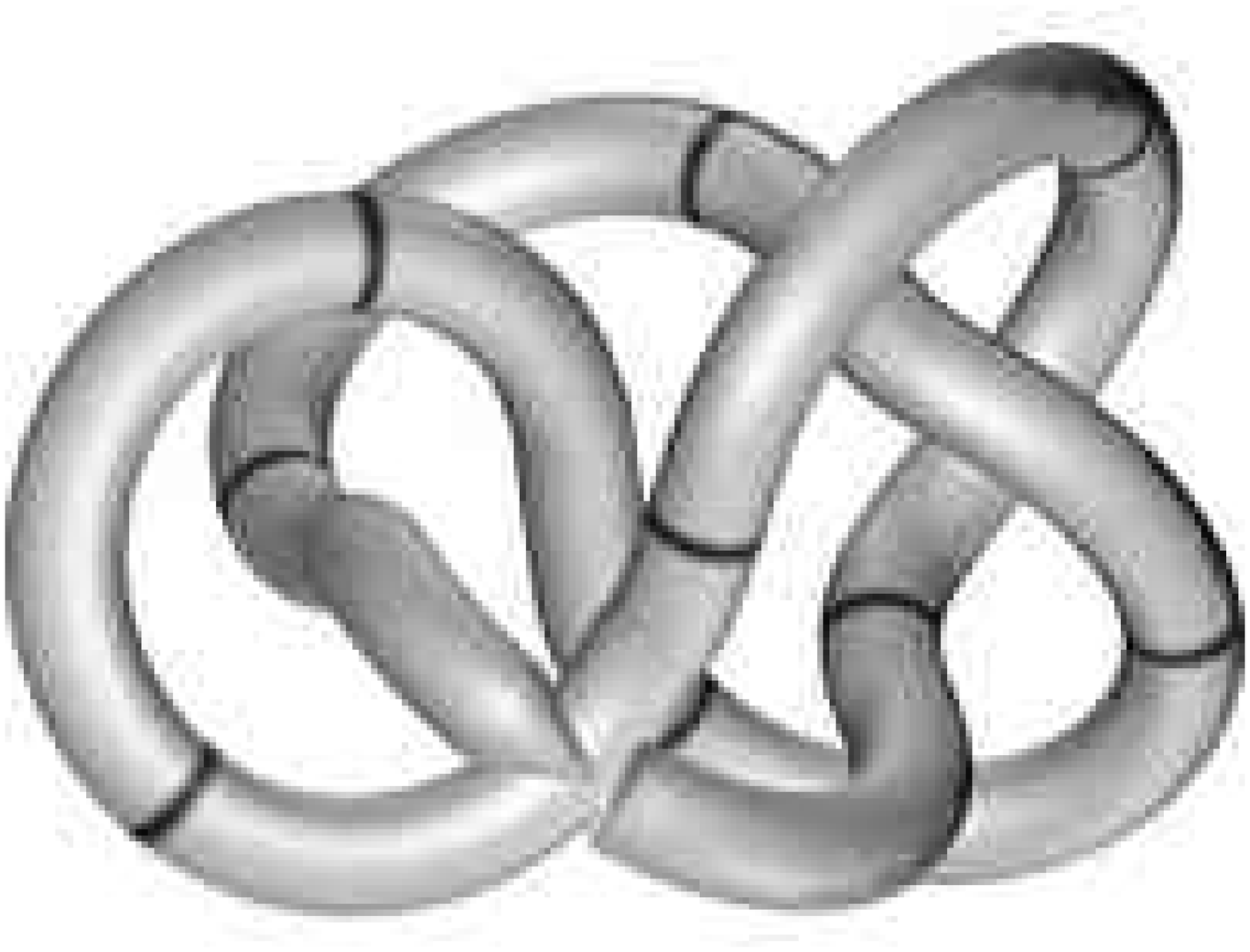}
        \put(-10,90){\large{$5_{2}$}}
    \end{overpic}
      \hspace{7mm}
    \begin{overpic}[width=2.8in]{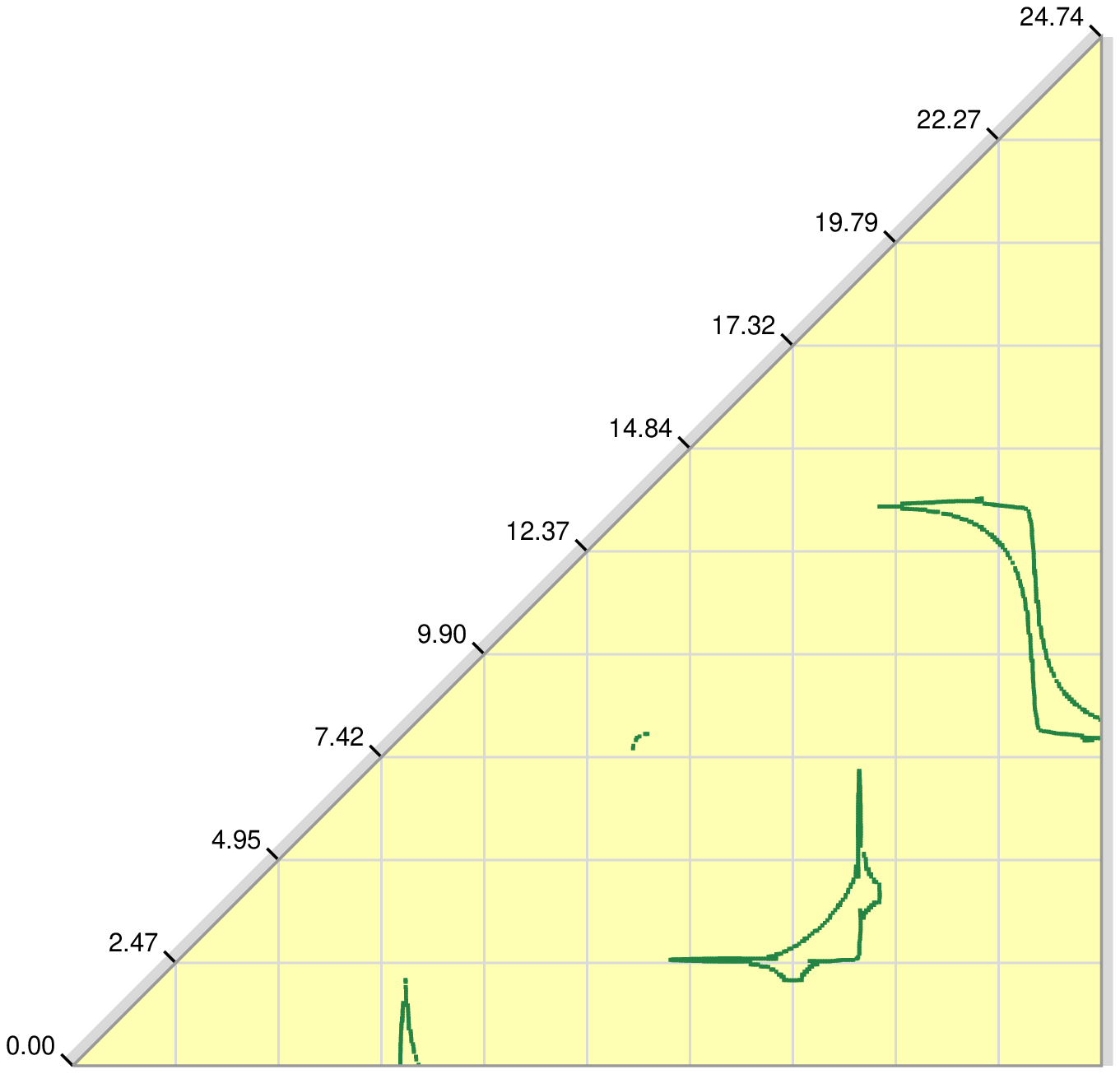}
        \put(8,94){\scriptsize{$49.49$}}
        \put(8,89){\scriptsize{$49.48$}}
        \put(8,84){\scriptsize{$400$}}
    \end{overpic}
\end{minipage} 
\hfill
\begin{minipage}[t]{6in}
  \vspace{2mm}
    \begin{overpic}[height=2.8in]{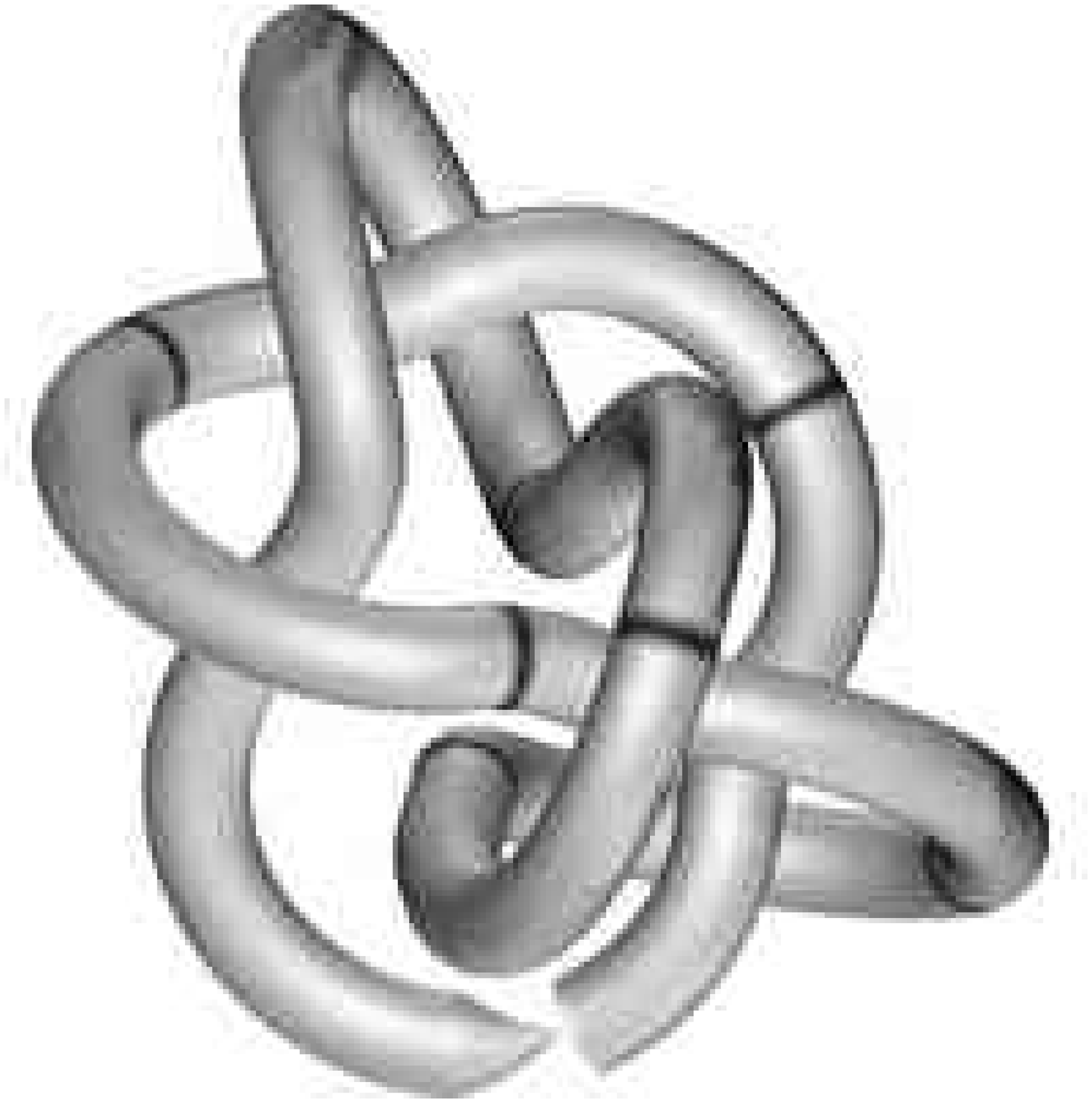}
        \put(-10,90){\large{$6_{2}$}}
    \end{overpic}
      \hspace{7mm}
    \begin{overpic}[width=2.8in]{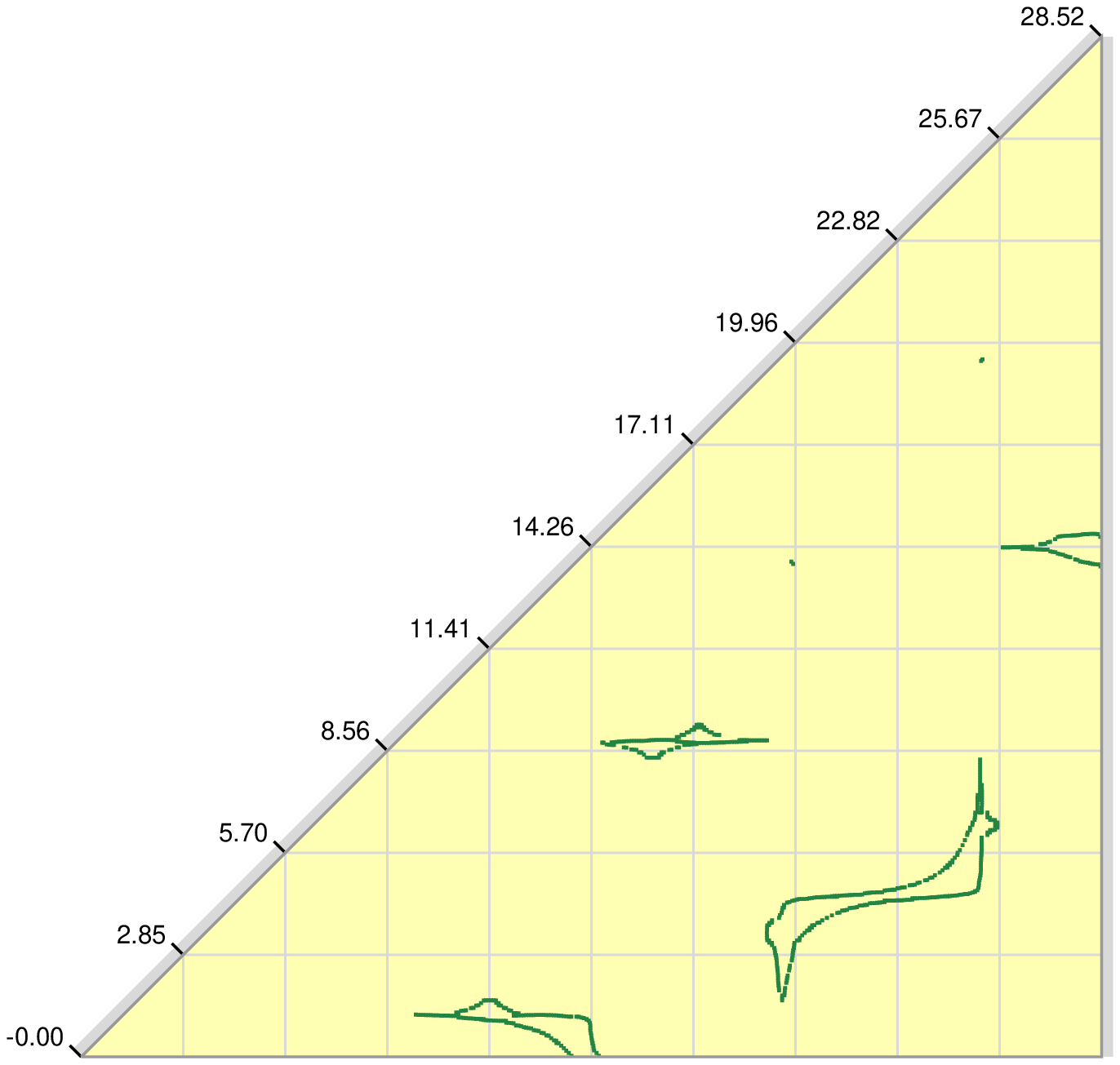}
        \put(8,94){\scriptsize{$57.05$}}
        \put(8,89){\scriptsize{$57.03$}}
        \put(8,84){\scriptsize{$407$}}
    \end{overpic}
\end{minipage} 
\hfill
\begin{minipage}[t]{6in}
  \vspace{2mm}
    \begin{overpic}[height=2.8in]{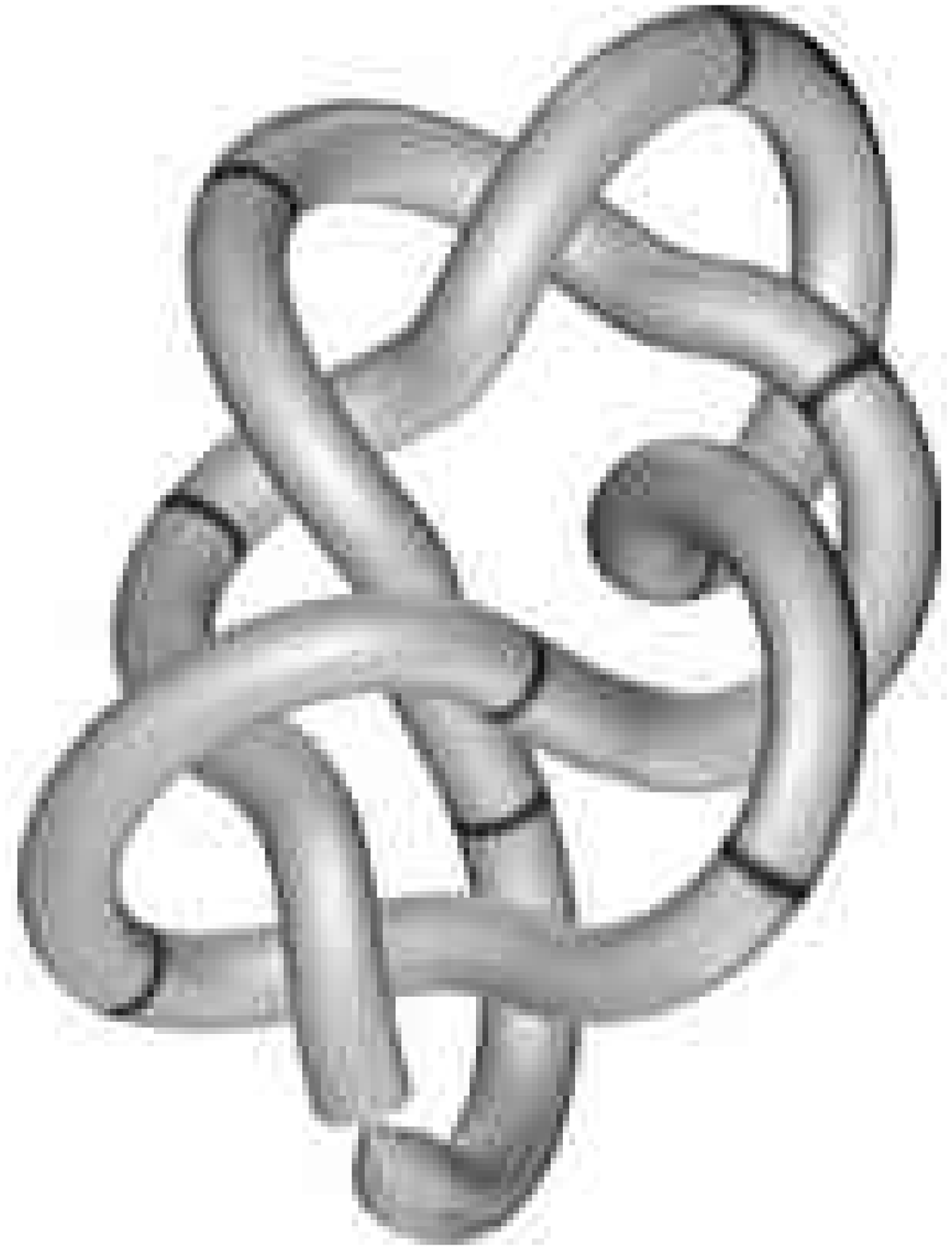}
        \put(-10,90){\large{$7_{2}$}}
    \end{overpic}
      \hspace{7mm}
    \begin{overpic}[width=2.8in]{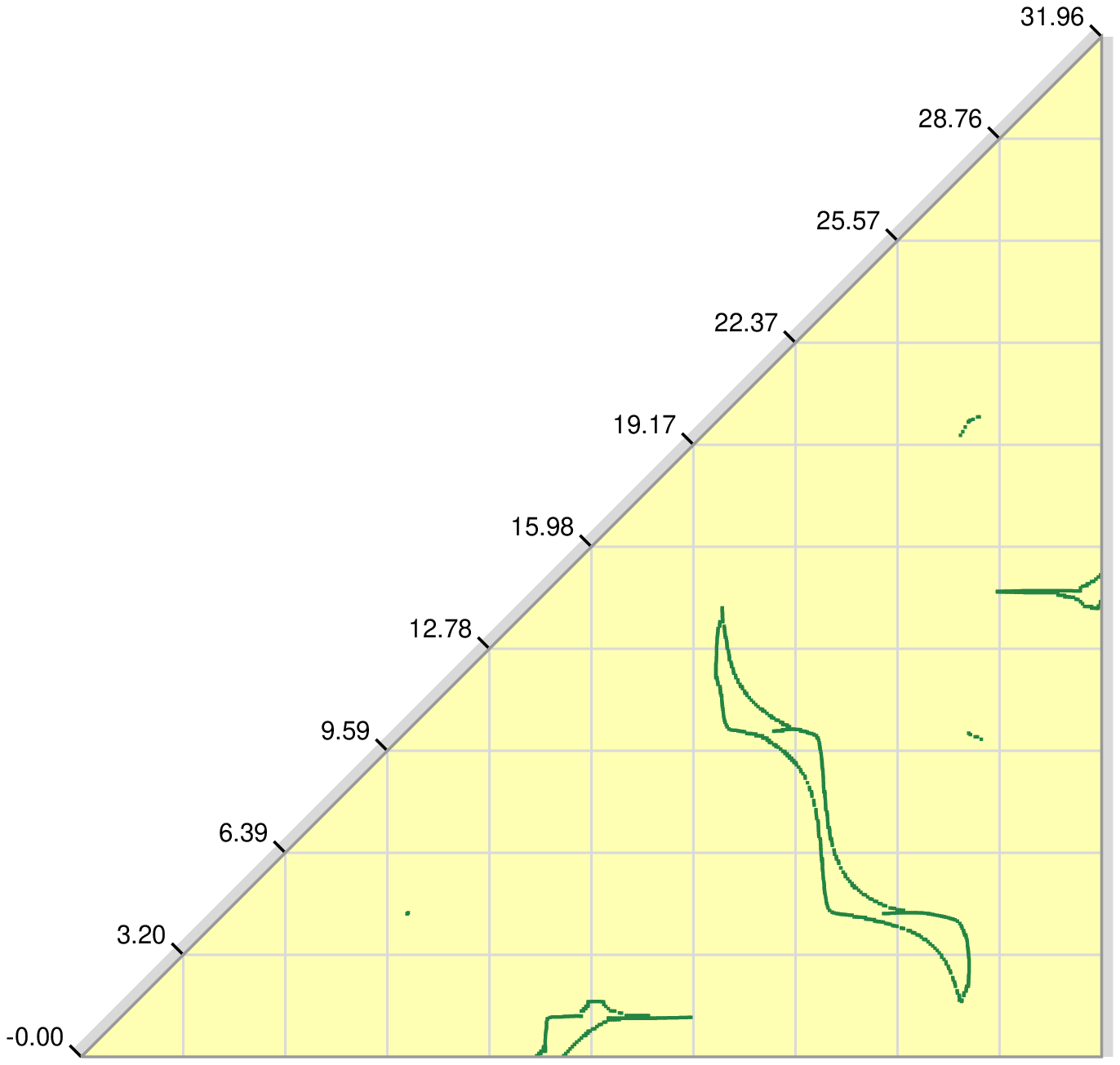}
        \put(8,94){\scriptsize{$63.92$}}
        \put(8,89){\scriptsize{$63.90$}}
        \put(8,84){\scriptsize{$456$}}
    \end{overpic}
\end{minipage} 
\hfill
\end{figure}
\clearpage
\pagebreak
\begin{figure}
\begin{minipage}[t]{6in}
  \vspace{2mm}
    \begin{overpic}[width=2.8in]{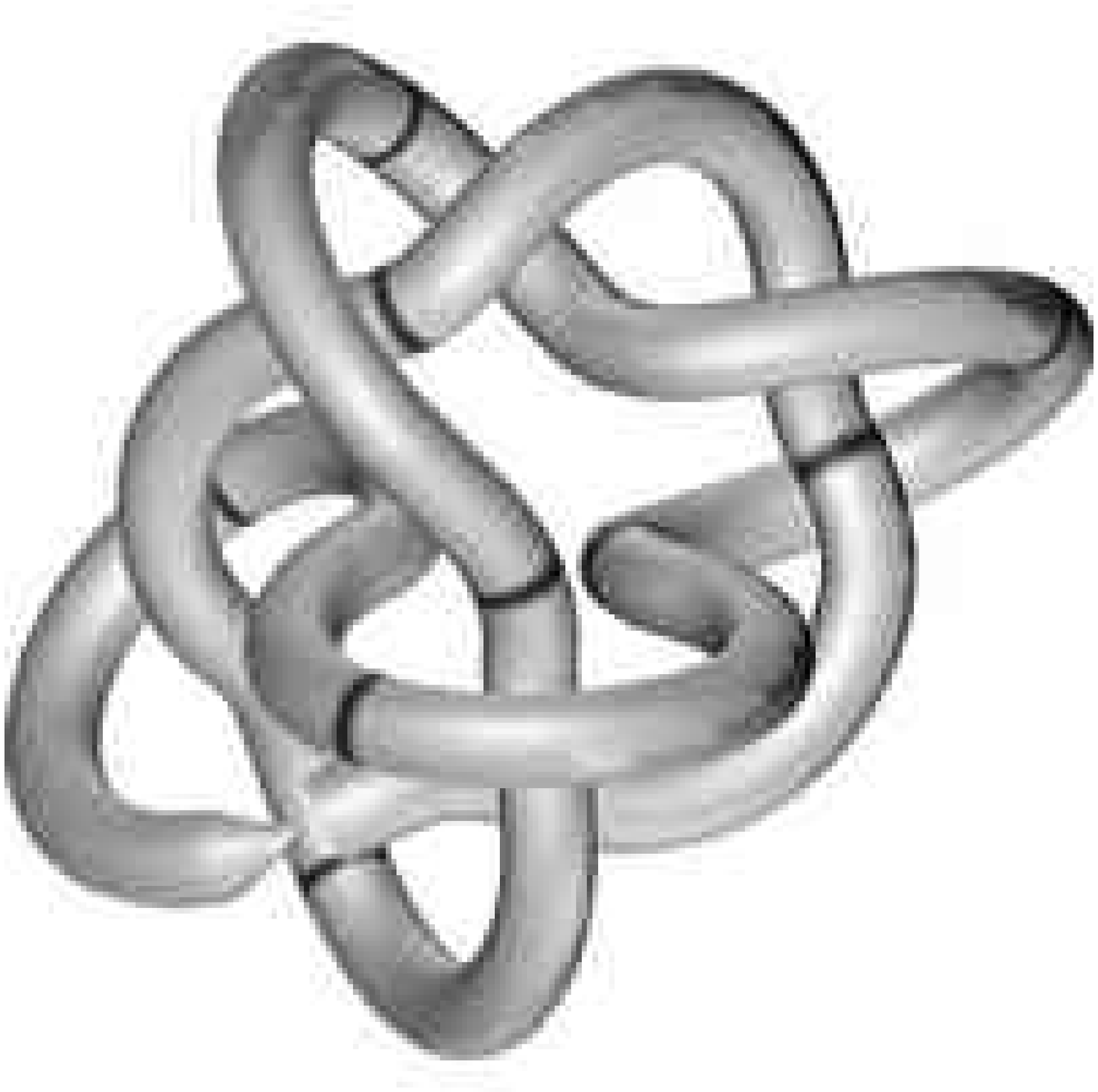}
        \put(-10,90){\large{$8_{4}$}}
    \end{overpic}
      \hspace{7mm}
    \begin{overpic}[width=2.8in]{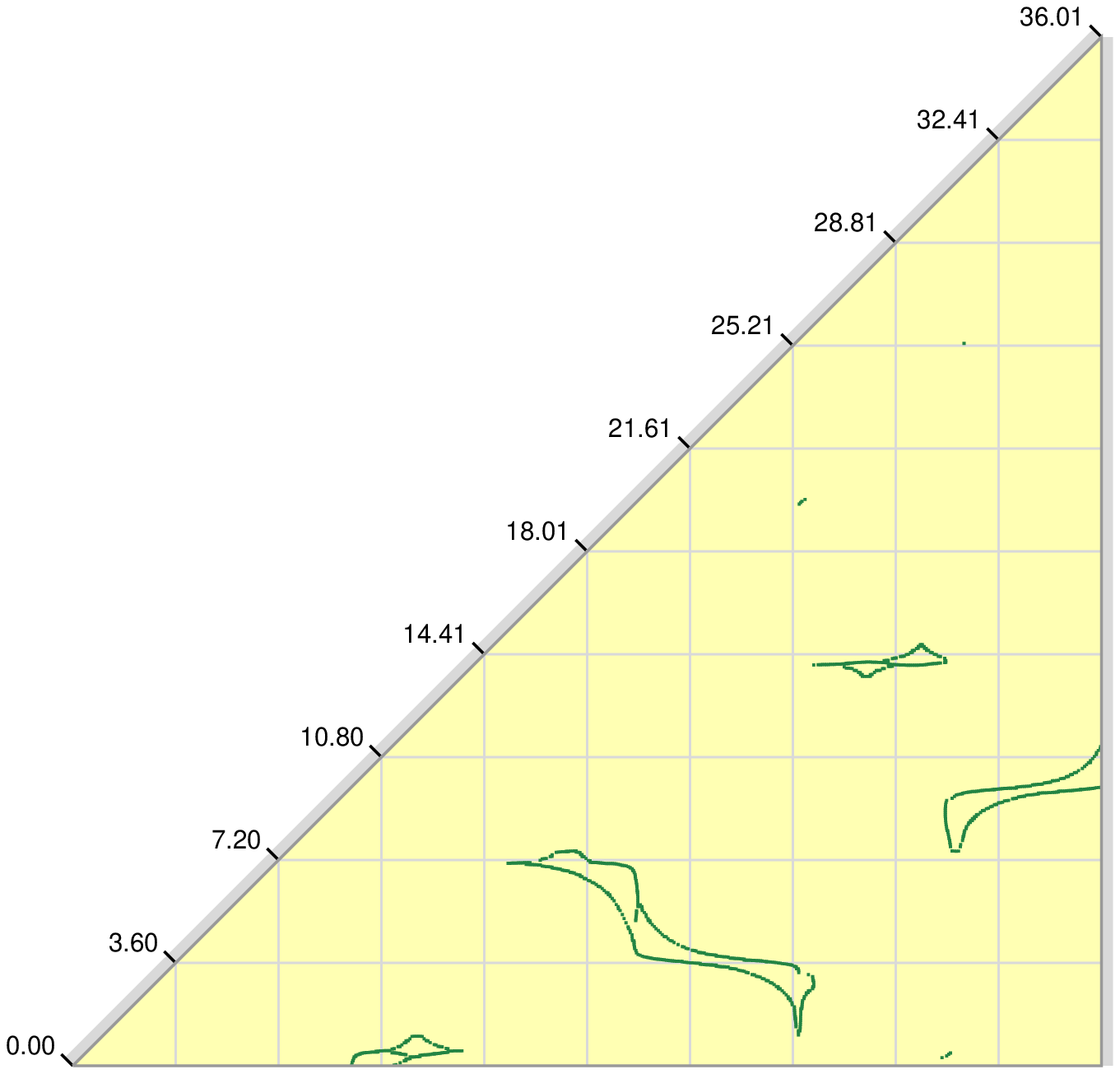}
        \put(8,94){\scriptsize{$72.04$}}
        \put(8,89){\scriptsize{$72.01$}}
        \put(8,84){\scriptsize{$514$}}
    \end{overpic}
\end{minipage} 
\hfill
\begin{minipage}[t]{6in}
  \vspace{2mm}
    \begin{overpic}[height=2.8in]{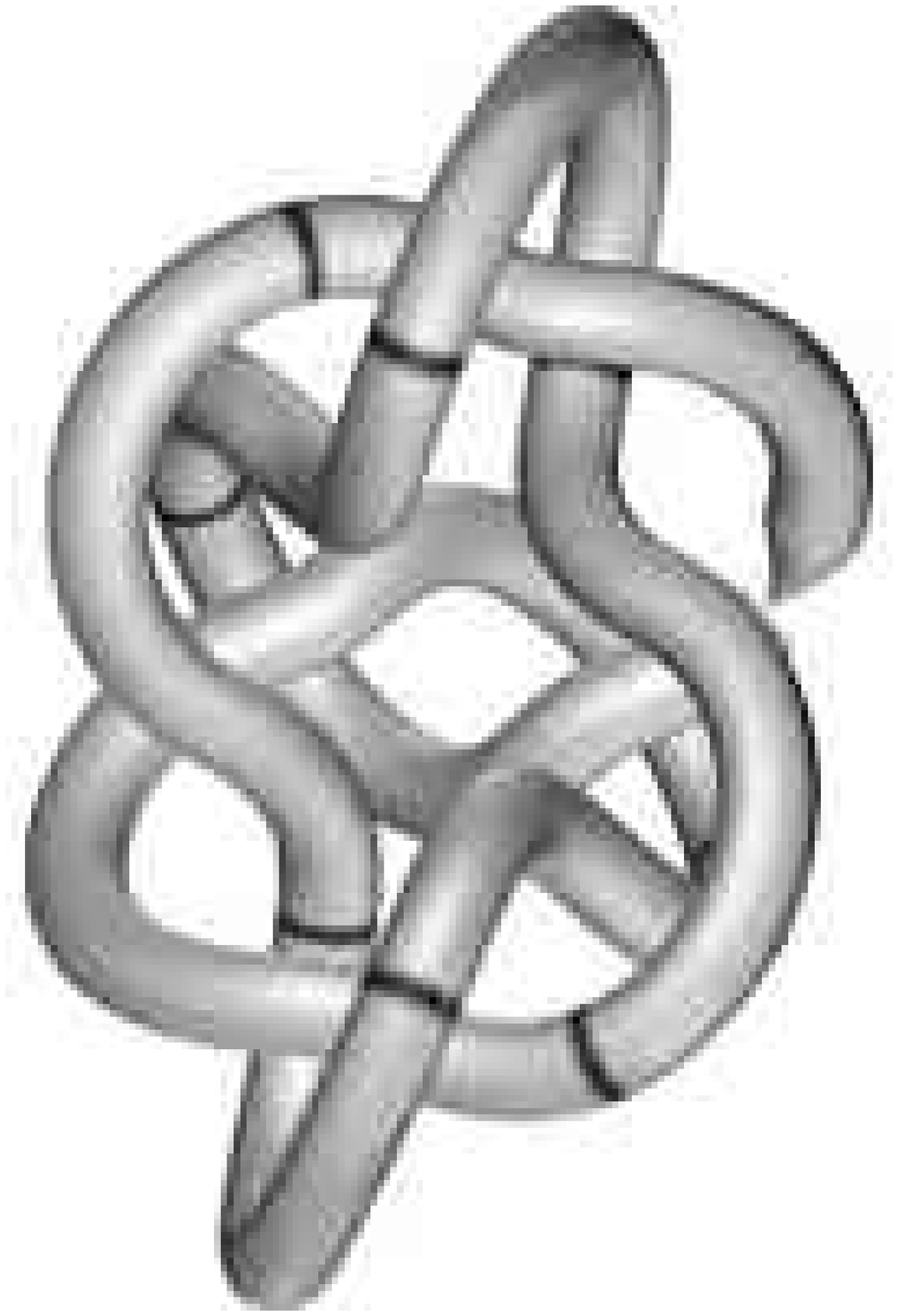}
        \put(-10,90){\large{$8_{5}$}}
    \end{overpic}
      \hspace{7mm}
    \begin{overpic}[width=2.8in]{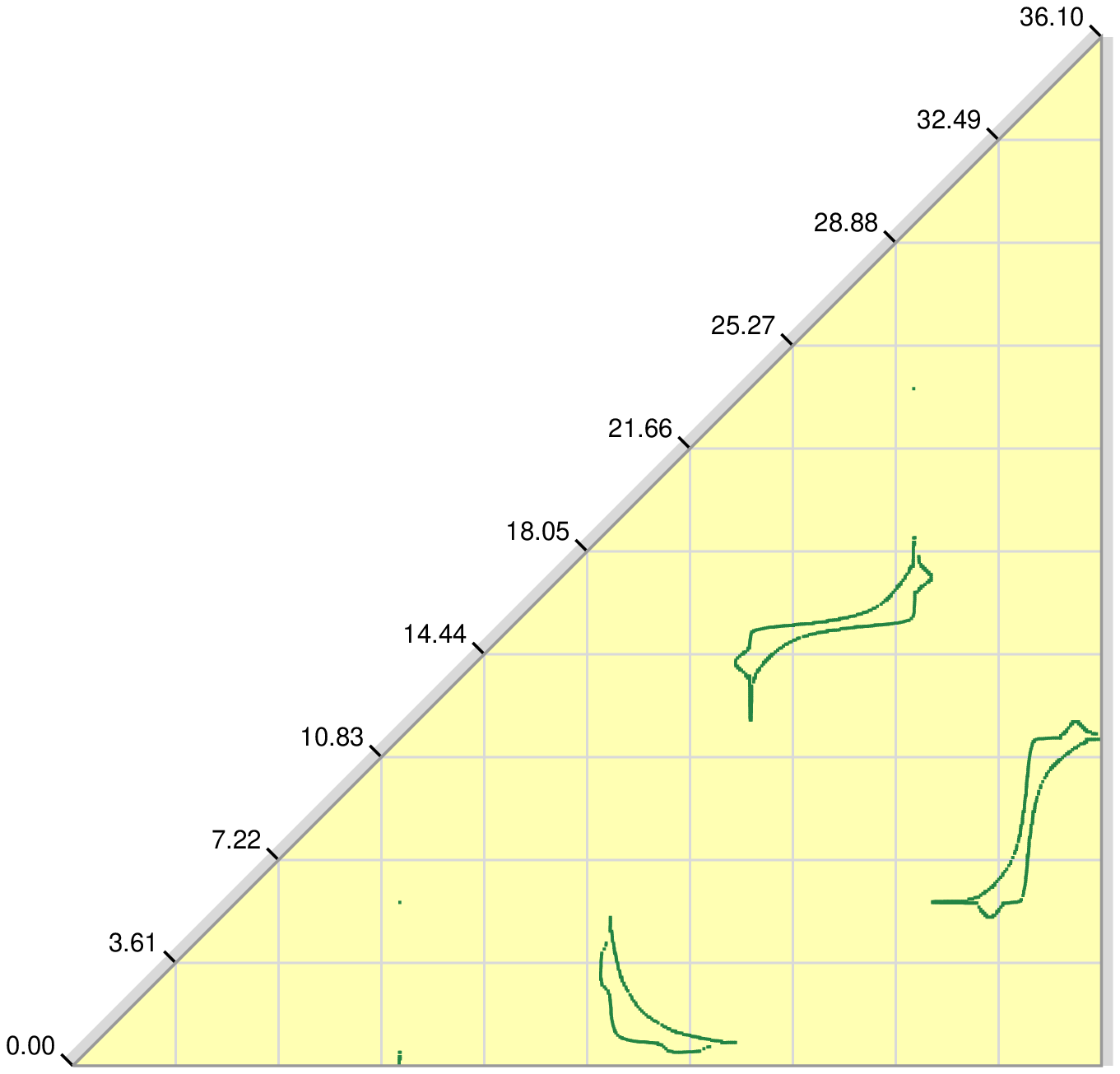}
        \put(8,94){\scriptsize{$72.21$}}
        \put(8,89){\scriptsize{$72.19$}}
        \put(8,84){\scriptsize{$516$}}
    \end{overpic}
\end{minipage} 
\hfill
\begin{minipage}[t]{6in}
  \vspace{2mm}
    \begin{overpic}[height=2.8in]{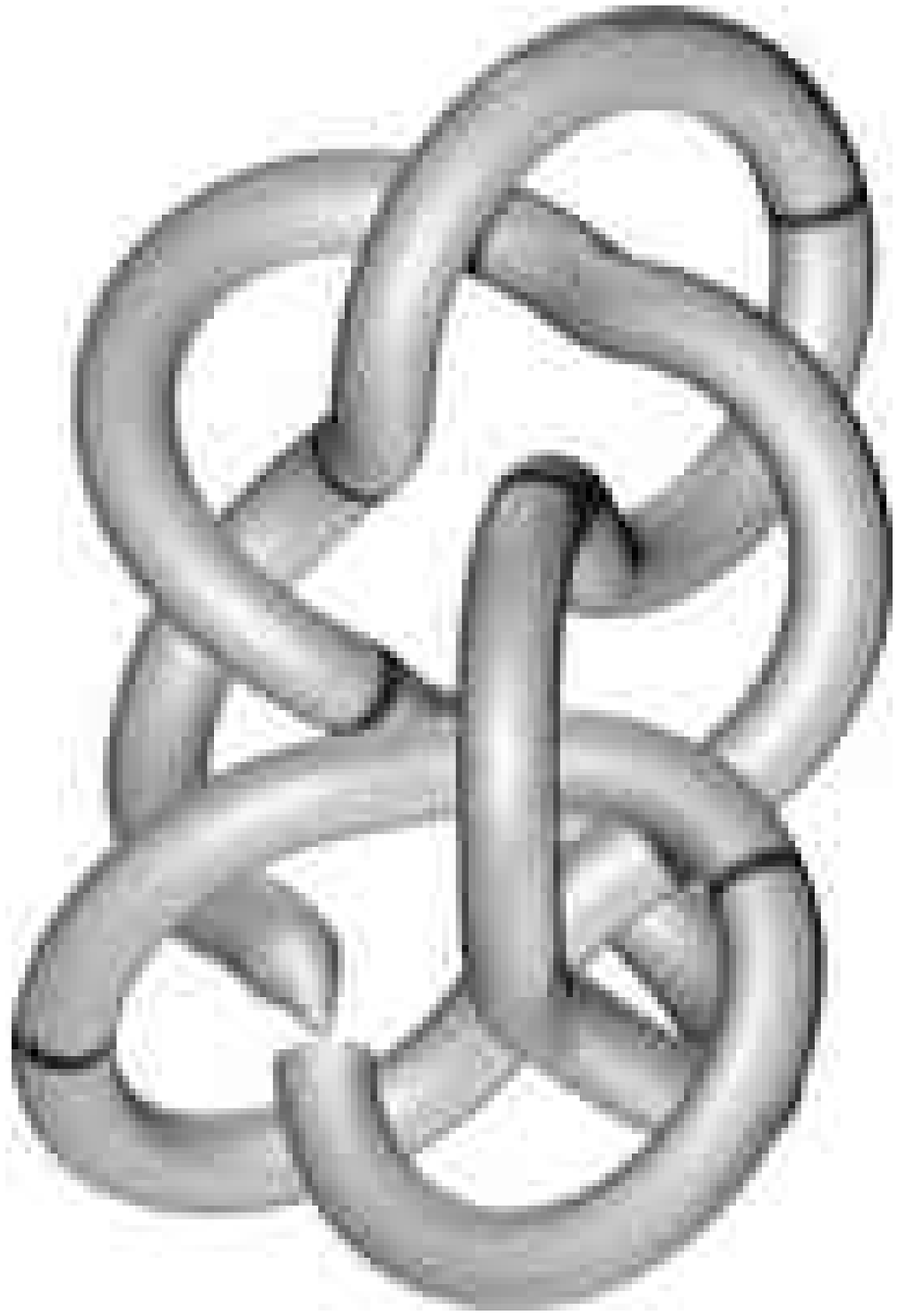}
        \put(-10,90){\large{$8_{7}$}}
    \end{overpic}
      \hspace{7mm}
    \begin{overpic}[width=2.8in]{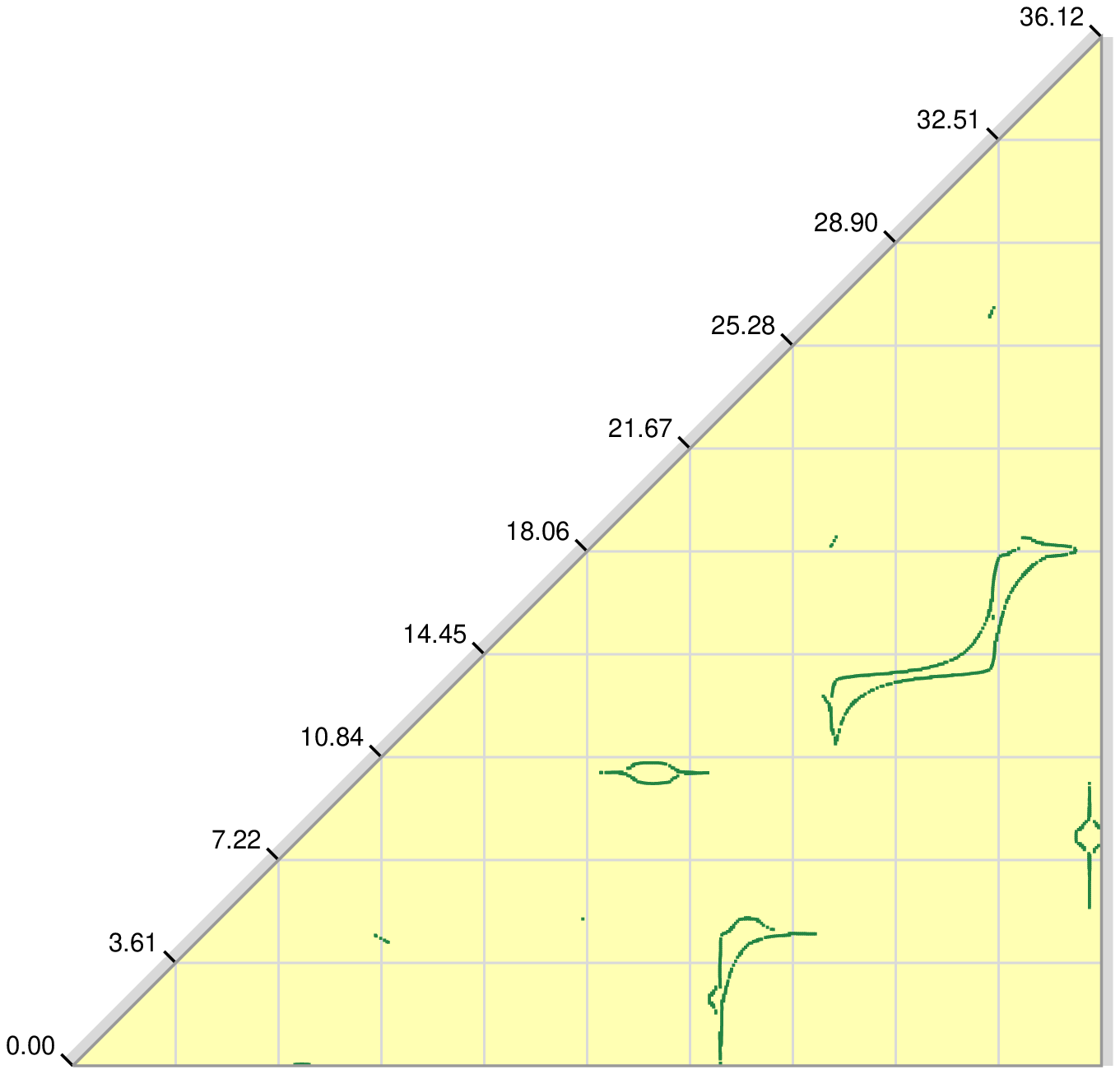}
        \put(8,94){\scriptsize{$72.25$}}
        \put(8,89){\scriptsize{$72.23$}}
        \put(8,84){\scriptsize{$516$}}
    \end{overpic}
\end{minipage} 
\hfill
\end{figure}
\clearpage
\pagebreak
\begin{figure}
\begin{minipage}[t]{6in}
  \vspace{2mm}
    \begin{overpic}[height=2.8in]{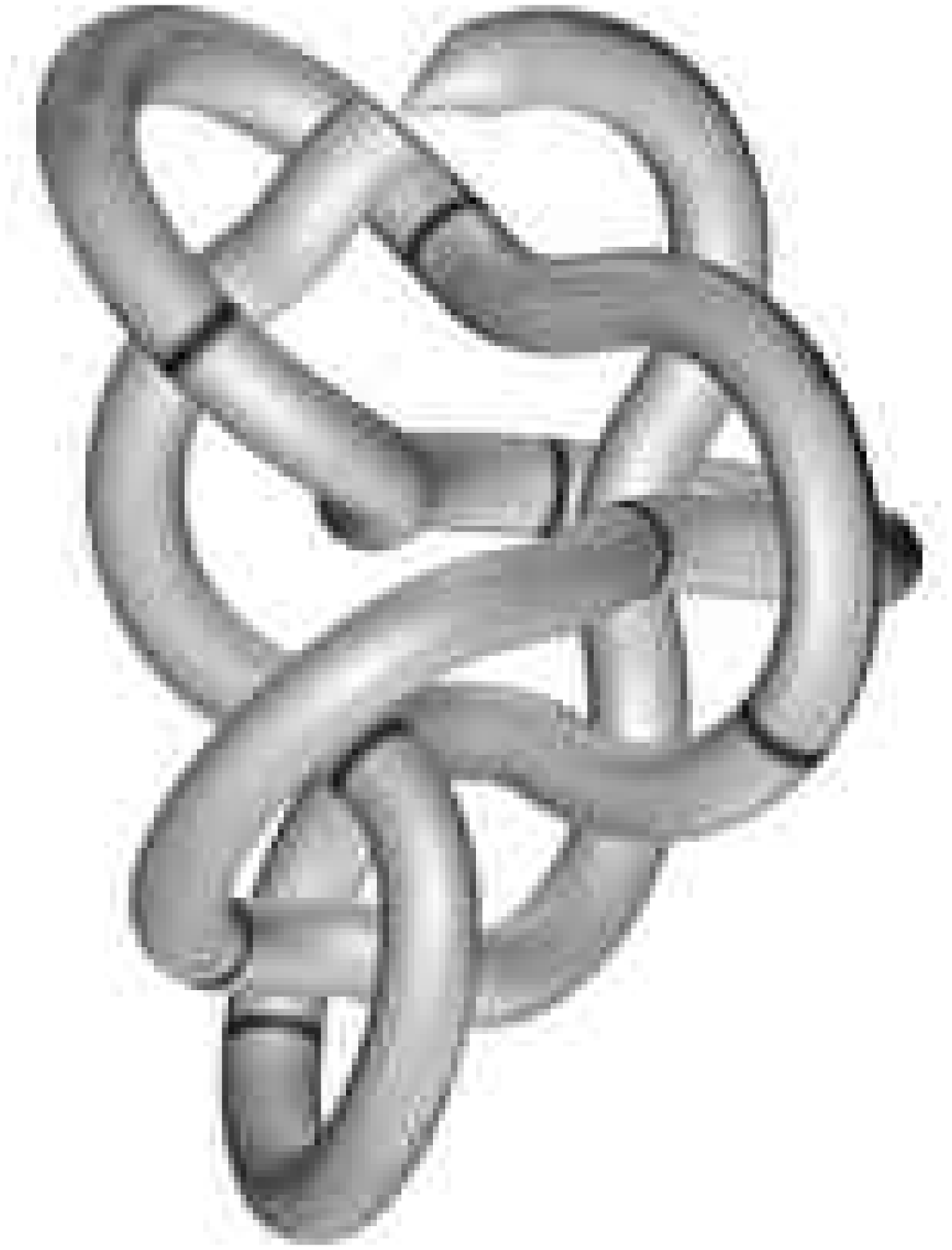}
        \put(-10,90){\large{$8_{9}$}}
    \end{overpic}
      \hspace{7mm}
    \begin{overpic}[width=2.8in]{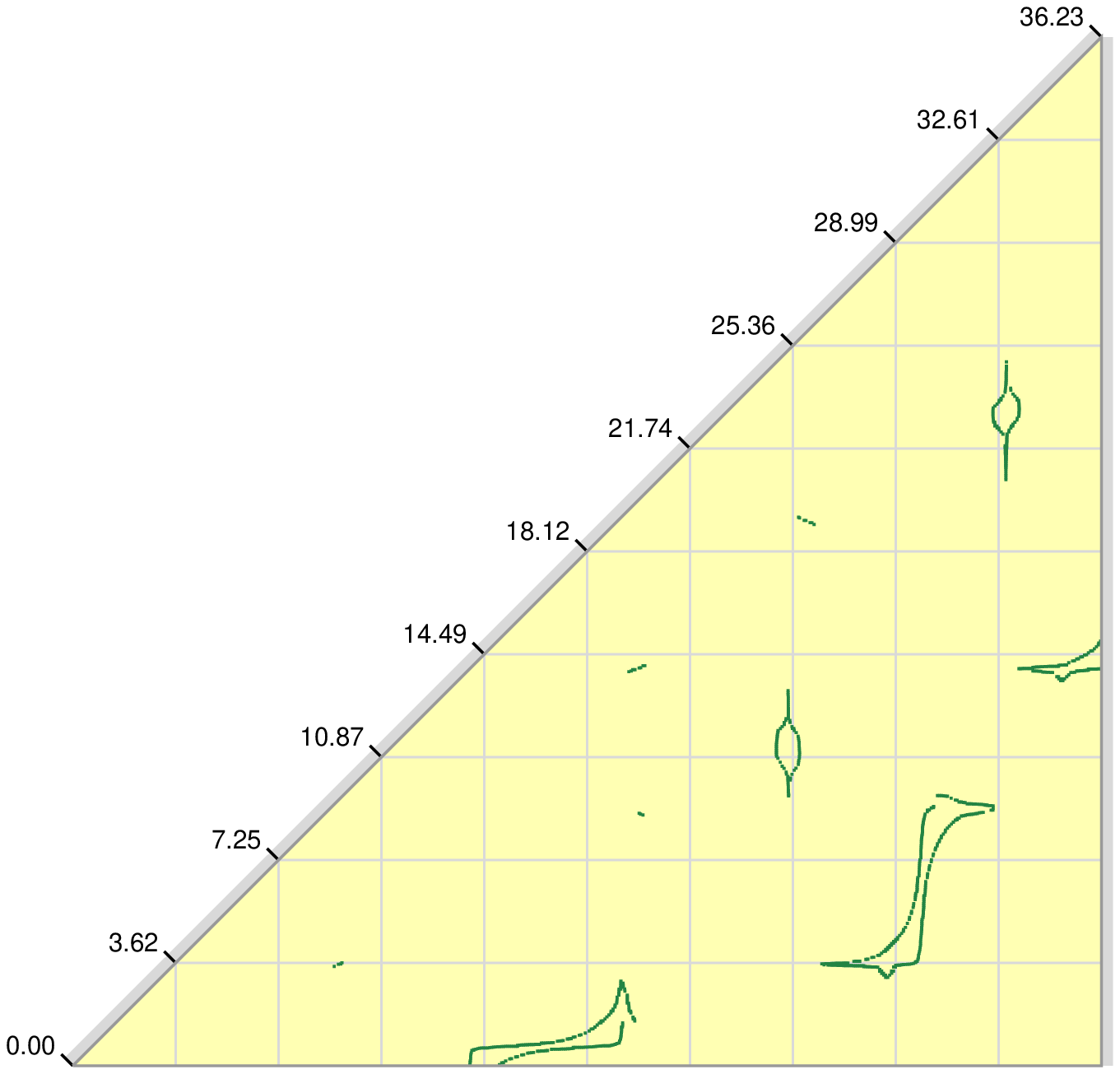}
        \put(8,94){\scriptsize{$72.48$}}
        \put(8,89){\scriptsize{$72.45$}}
        \put(8,84){\scriptsize{$517$}}
    \end{overpic}
\end{minipage} 
\hfill
\begin{minipage}[t]{6in}
  \vspace{2mm}
    \begin{overpic}[width=2.8in]{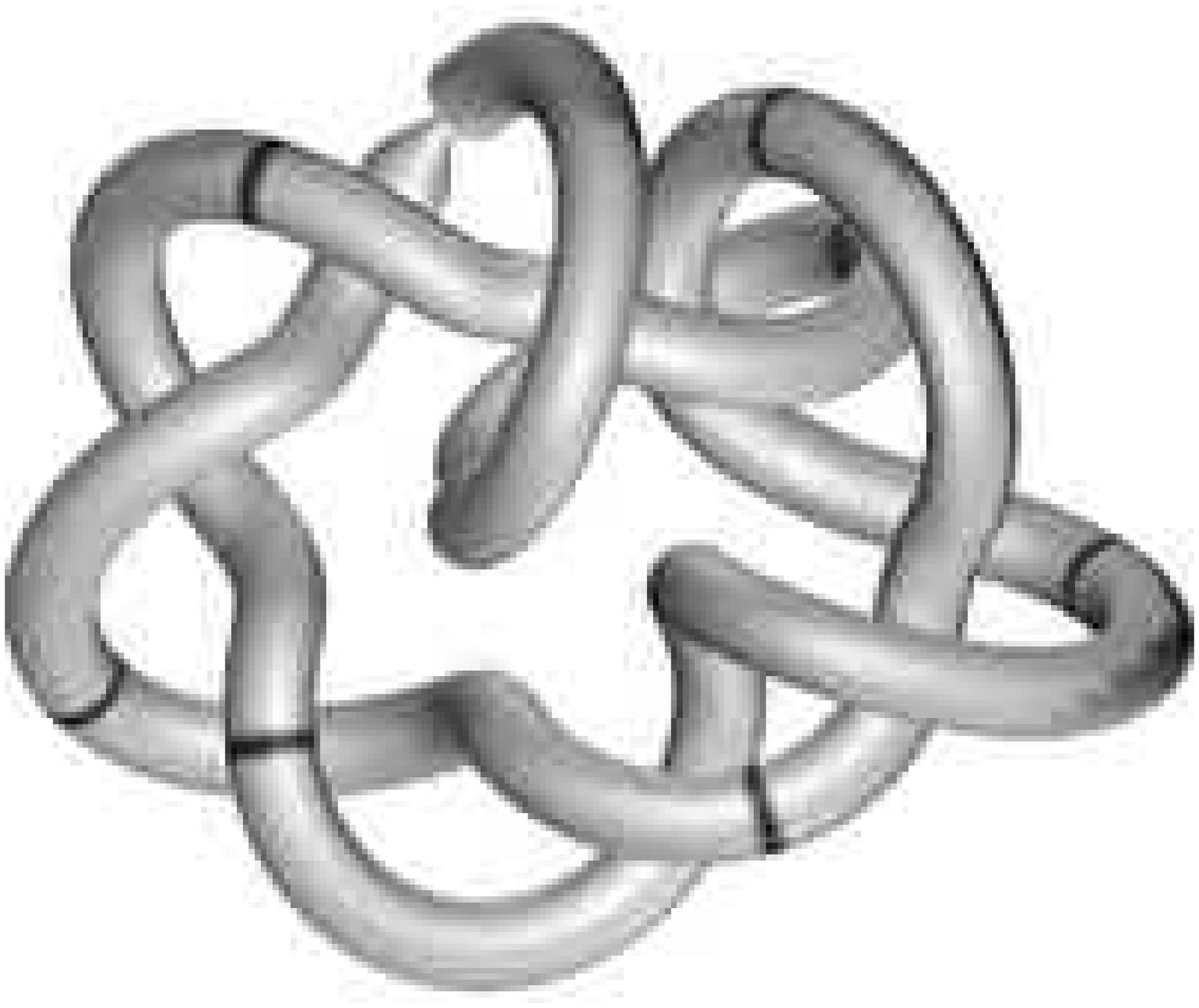}
        \put(-10,90){\large{$9_{2}$}}
    \end{overpic}
      \hspace{7mm}
    \begin{overpic}[width=2.8in]{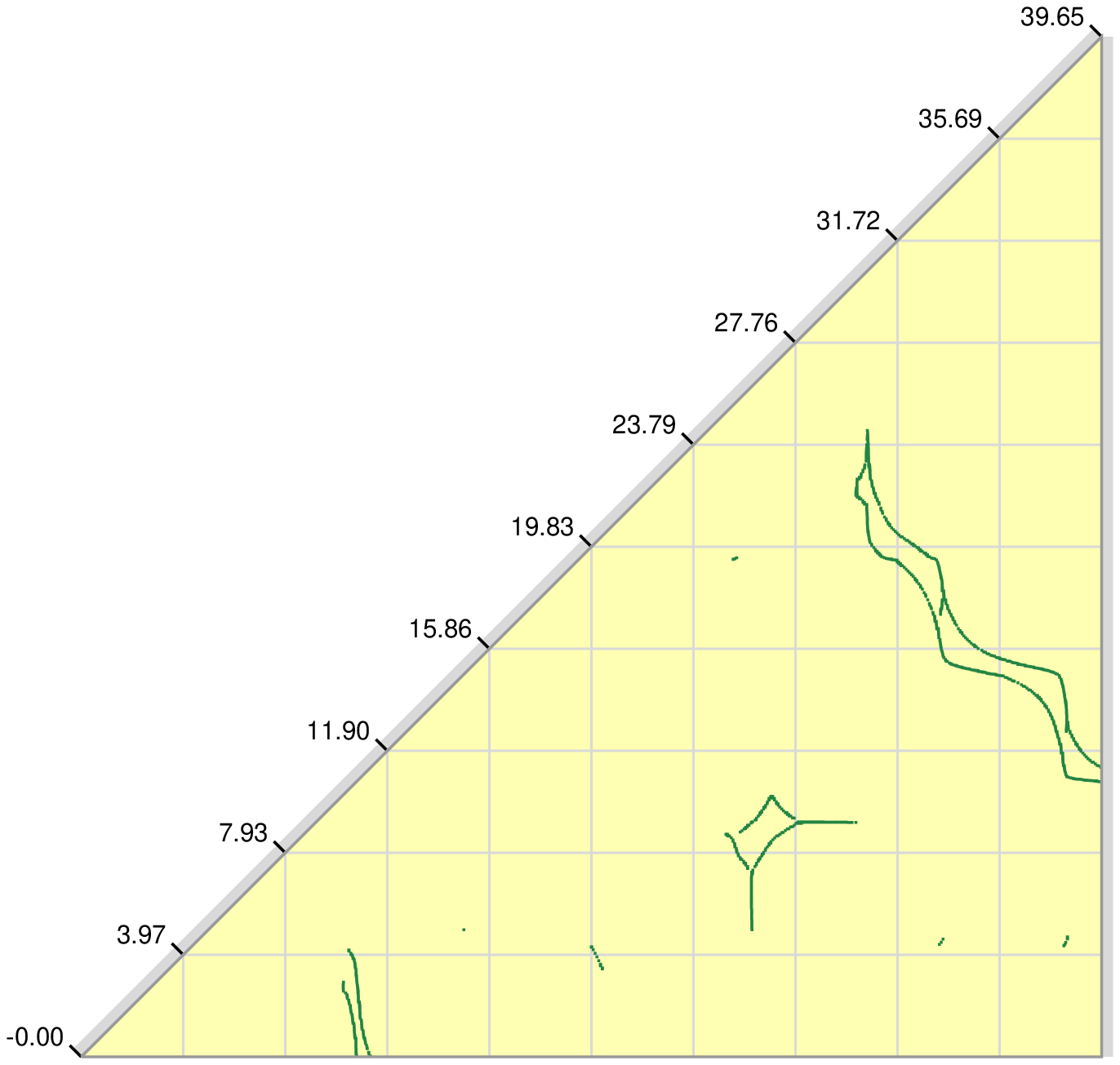}
        \put(8,94){\scriptsize{$79.31$}}
        \put(8,89){\scriptsize{$79.28$}}
        \put(8,84){\scriptsize{$566$}}
    \end{overpic}
\end{minipage} 
\hfill
\begin{minipage}[t]{6in}
  \vspace{2mm}
    \begin{overpic}[width=2.8in]{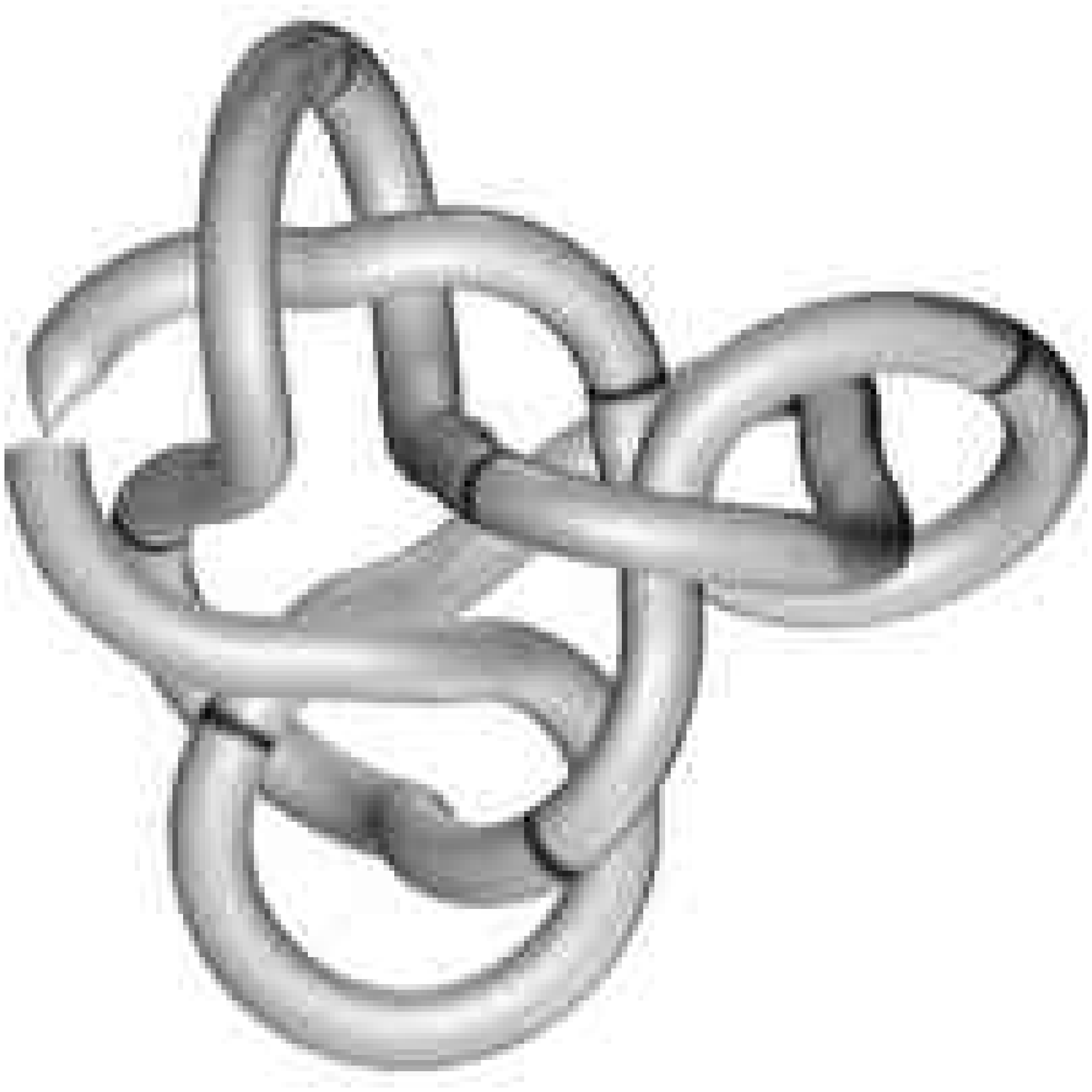}
        \put(-10,90){\large{$9_{10}$}}
    \end{overpic}
      \hspace{7mm}
    \begin{overpic}[width=2.8in]{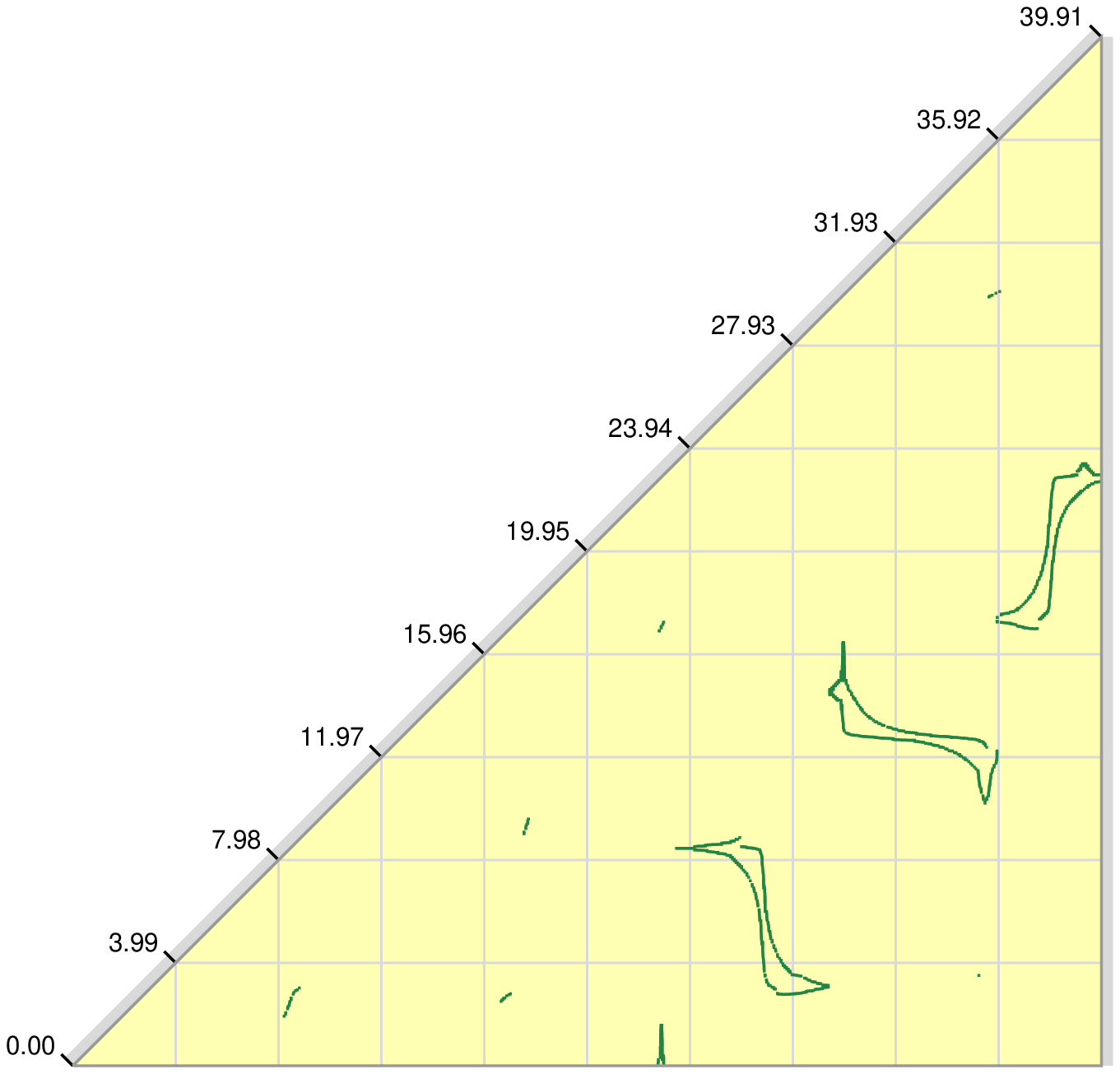}
        \put(8,94){\scriptsize{$79.82$}}
        \put(8,89){\scriptsize{$79.80$}}
        \put(8,84){\scriptsize{$570$}}
    \end{overpic}
\end{minipage} 
\hfill
\end{figure}
\clearpage
\pagebreak
\begin{figure}
\begin{minipage}[t]{6in}
  \vspace{2mm}
    \begin{overpic}[height=2.8in]{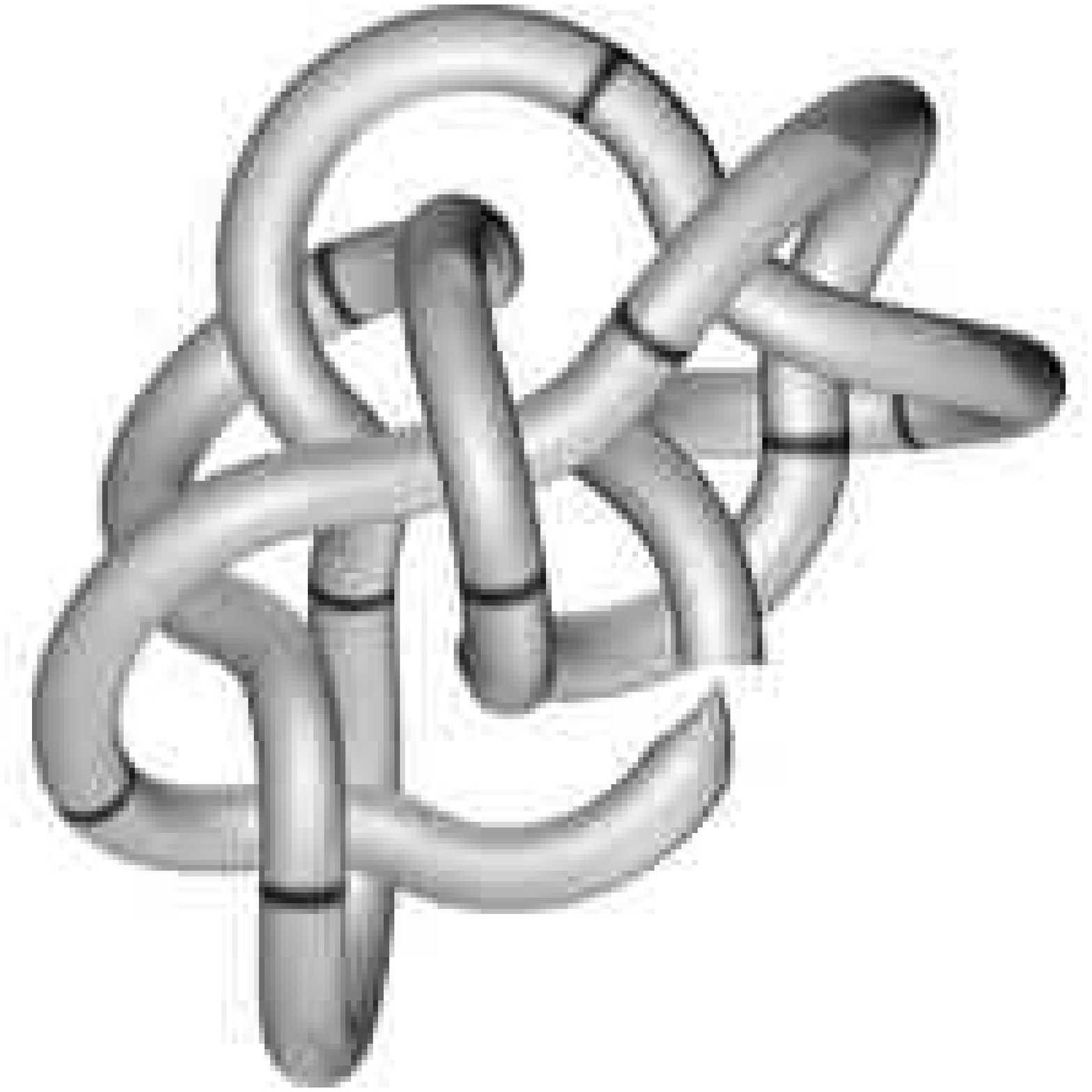}
        \put(-10,90){\large{$9_{11}$}}
    \end{overpic}
      \hspace{7mm}
    \begin{overpic}[width=2.8in]{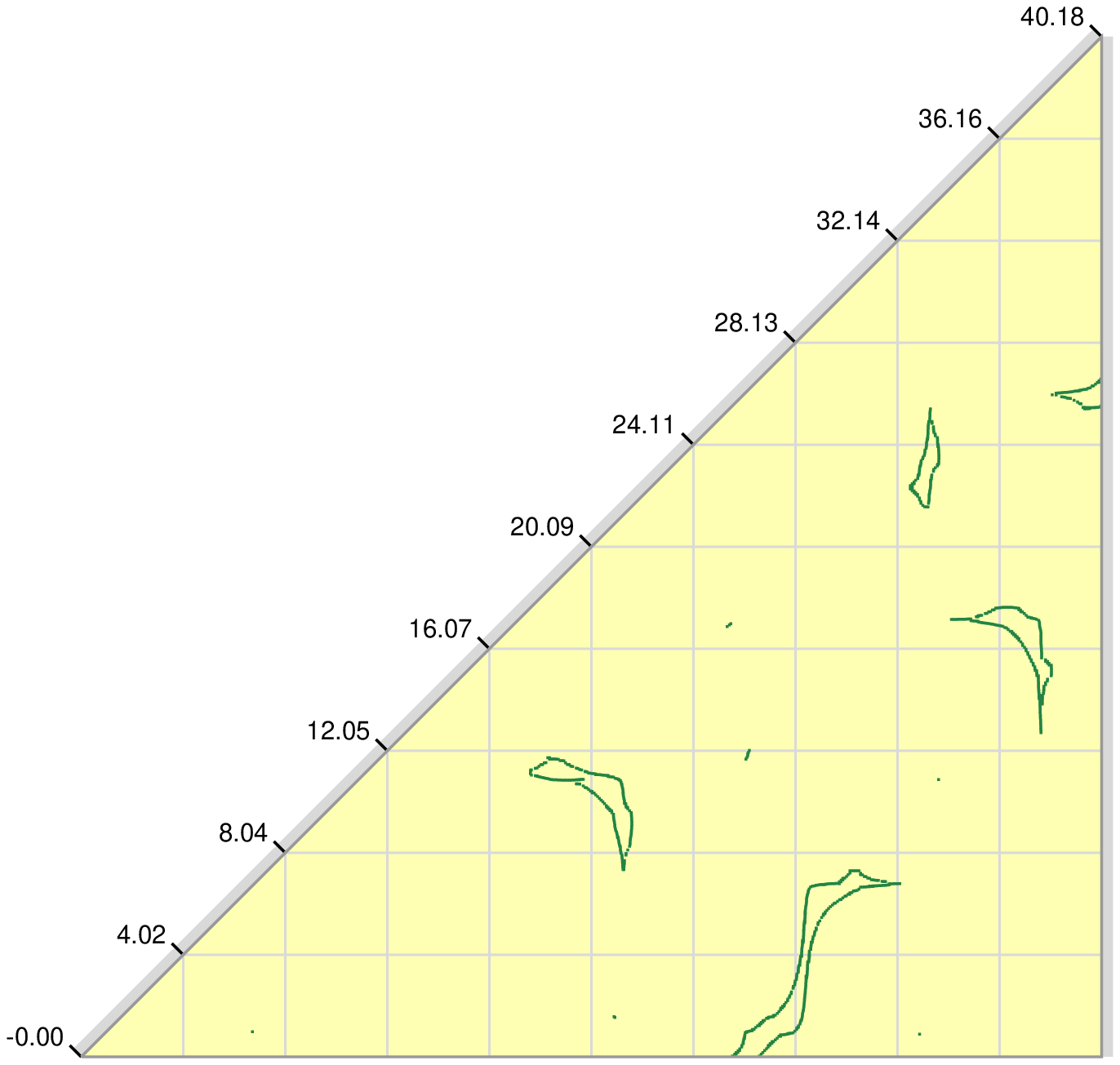}
        \put(8,94){\scriptsize{$80.37$}}
        \put(8,89){\scriptsize{$80.35$}}
        \put(8,84){\scriptsize{$574$}}
    \end{overpic}
\end{minipage} 
\hfill
\begin{minipage}[t]{6in}
  \vspace{2mm}
    \begin{overpic}[height=2.8in]{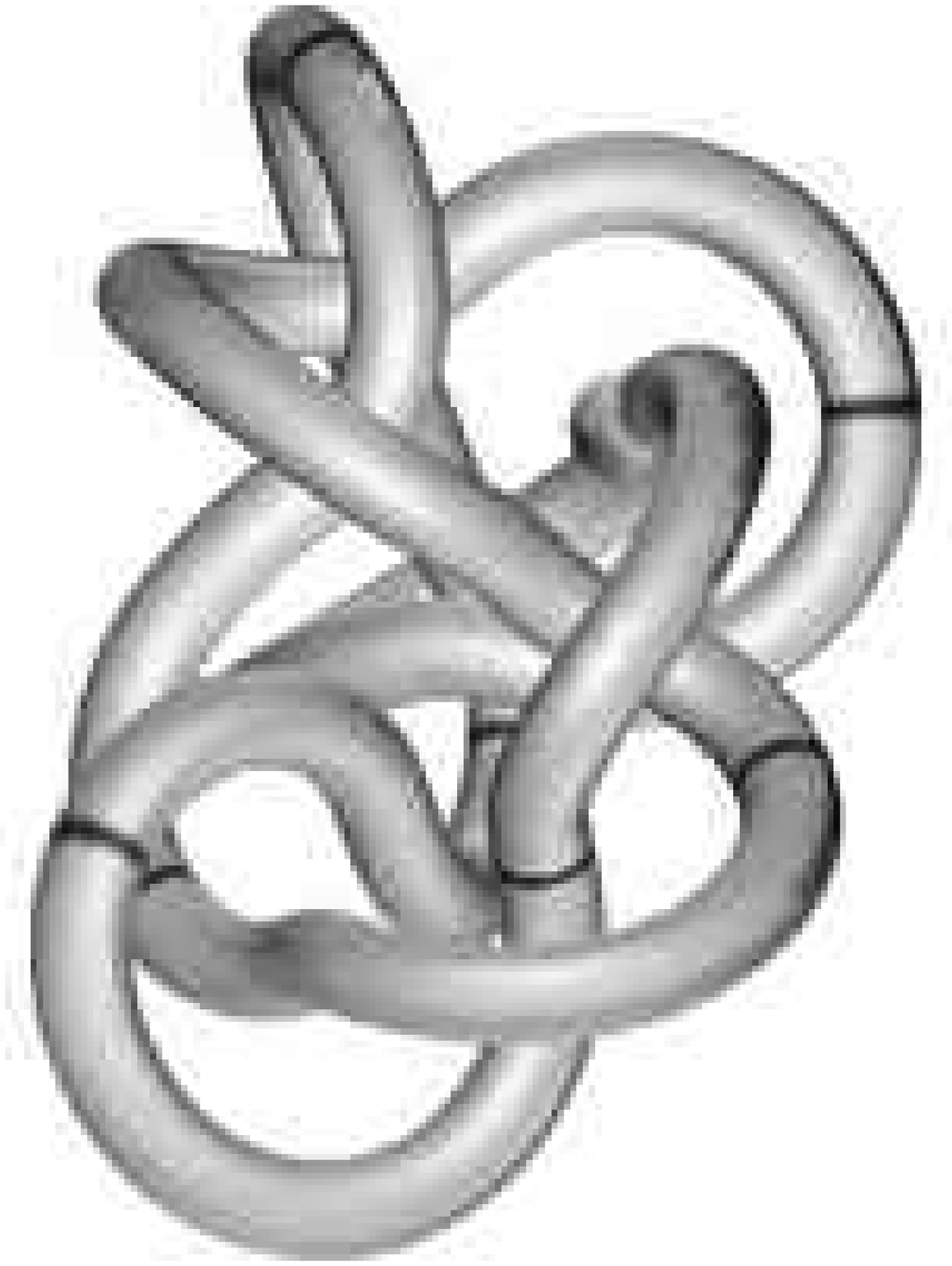}
        \put(-10,90){\large{$9_{12}$}}
    \end{overpic}
      \hspace{7mm}
    \begin{overpic}[width=2.8in]{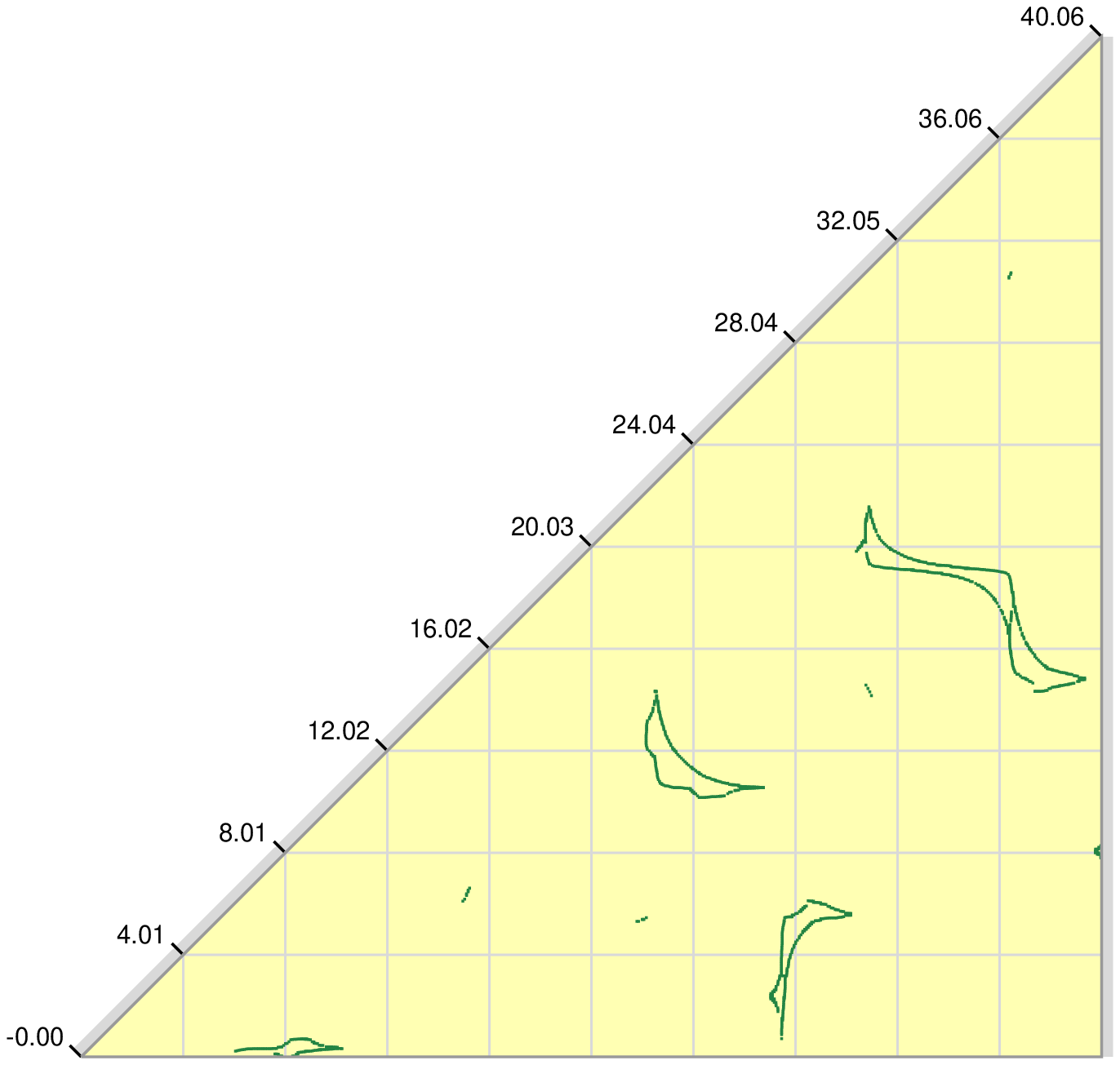}
        \put(8,94){\scriptsize{$80.13$}}
        \put(8,89){\scriptsize{$80.11$}}
        \put(8,84){\scriptsize{$572$}}
    \end{overpic}
\end{minipage} 
\hfill
\begin{minipage}[t]{6in}
  \vspace{2mm}
    \begin{overpic}[width=2.8in]{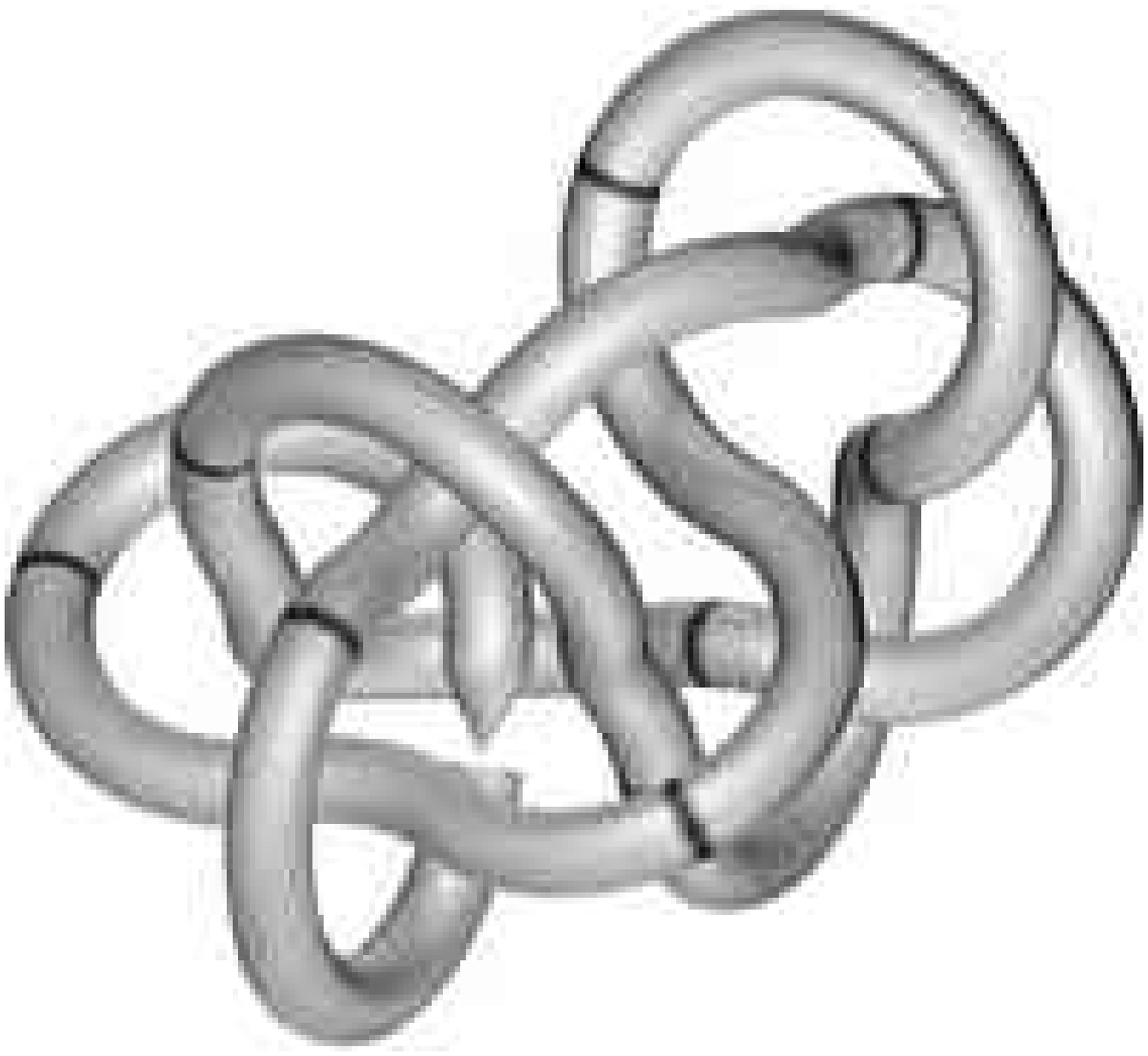}
        \put(-10,90){\large{$9_{14}$}}
    \end{overpic}
      \hspace{7mm}
    \begin{overpic}[width=2.8in]{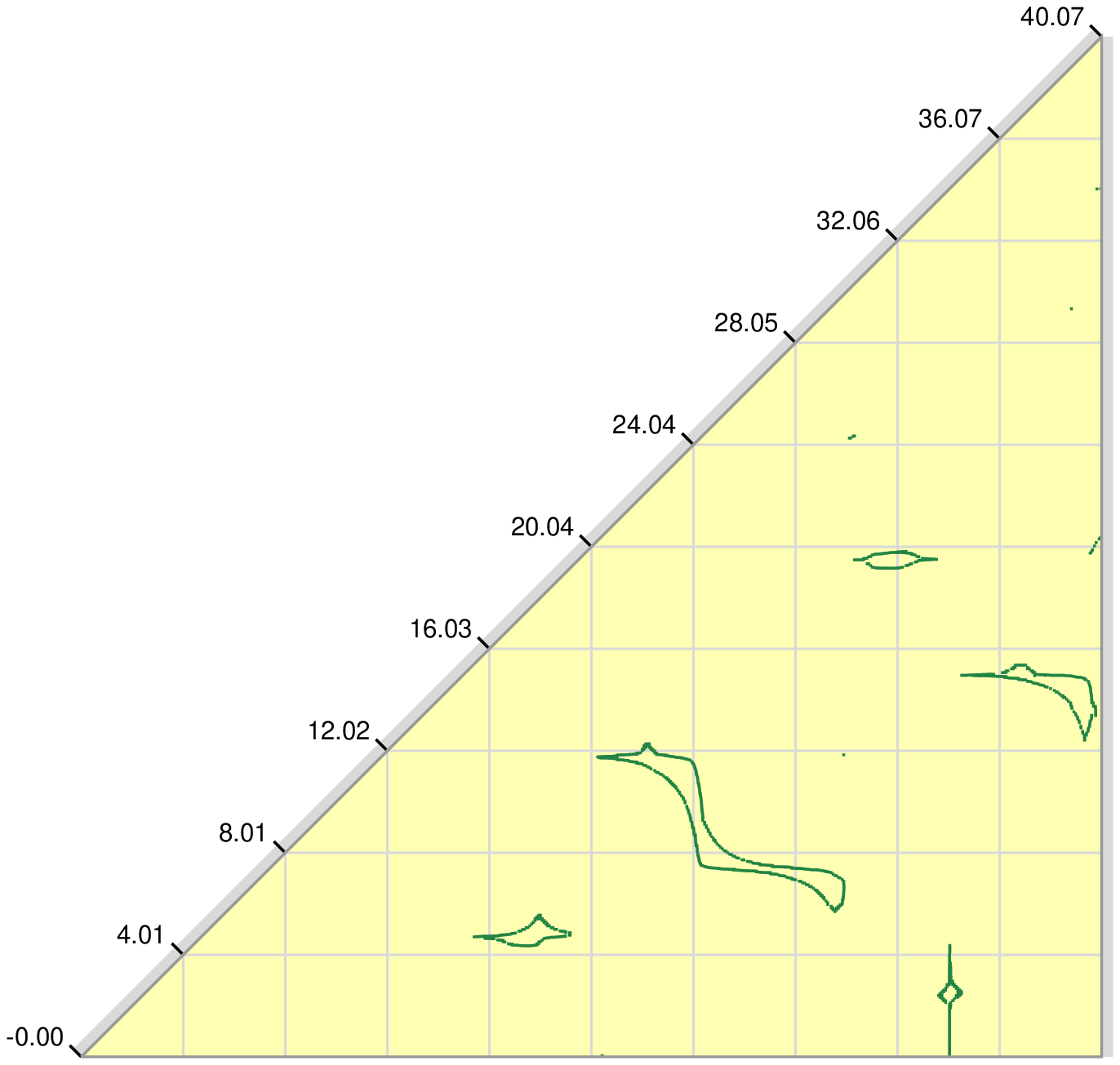}
        \put(8,94){\scriptsize{$80.16$}}
        \put(8,89){\scriptsize{$80.13$}}
        \put(8,84){\scriptsize{$572$}}
    \end{overpic}
\end{minipage} 
\hfill
\end{figure}
\clearpage
\pagebreak
\begin{figure}
\begin{minipage}[t]{6in}
  \vspace{2mm}
    \begin{overpic}[height=2.8in]{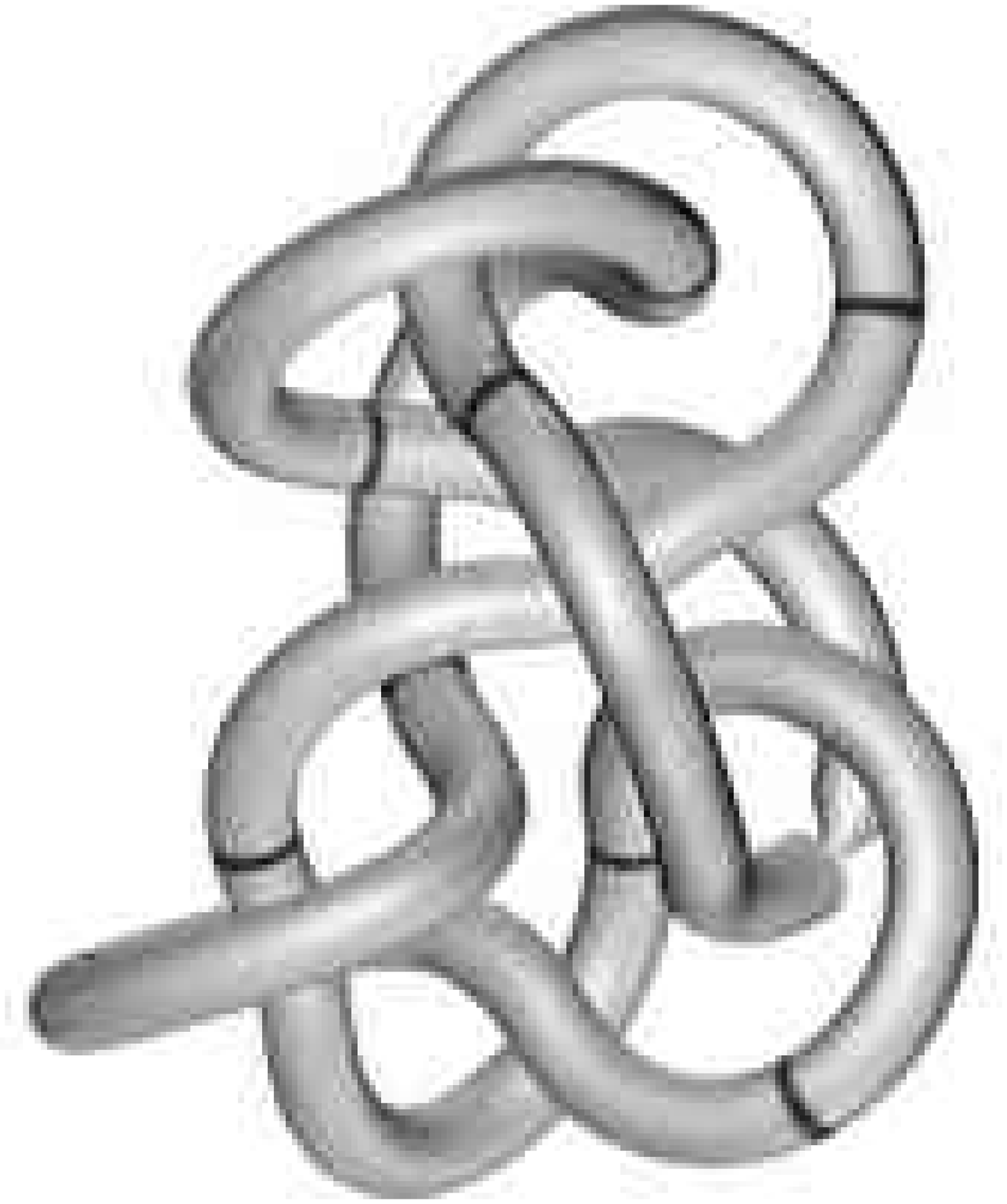}
        \put(-10,90){\large{$9_{15}$}}
    \end{overpic}
      \hspace{7mm}
    \begin{overpic}[width=2.8in]{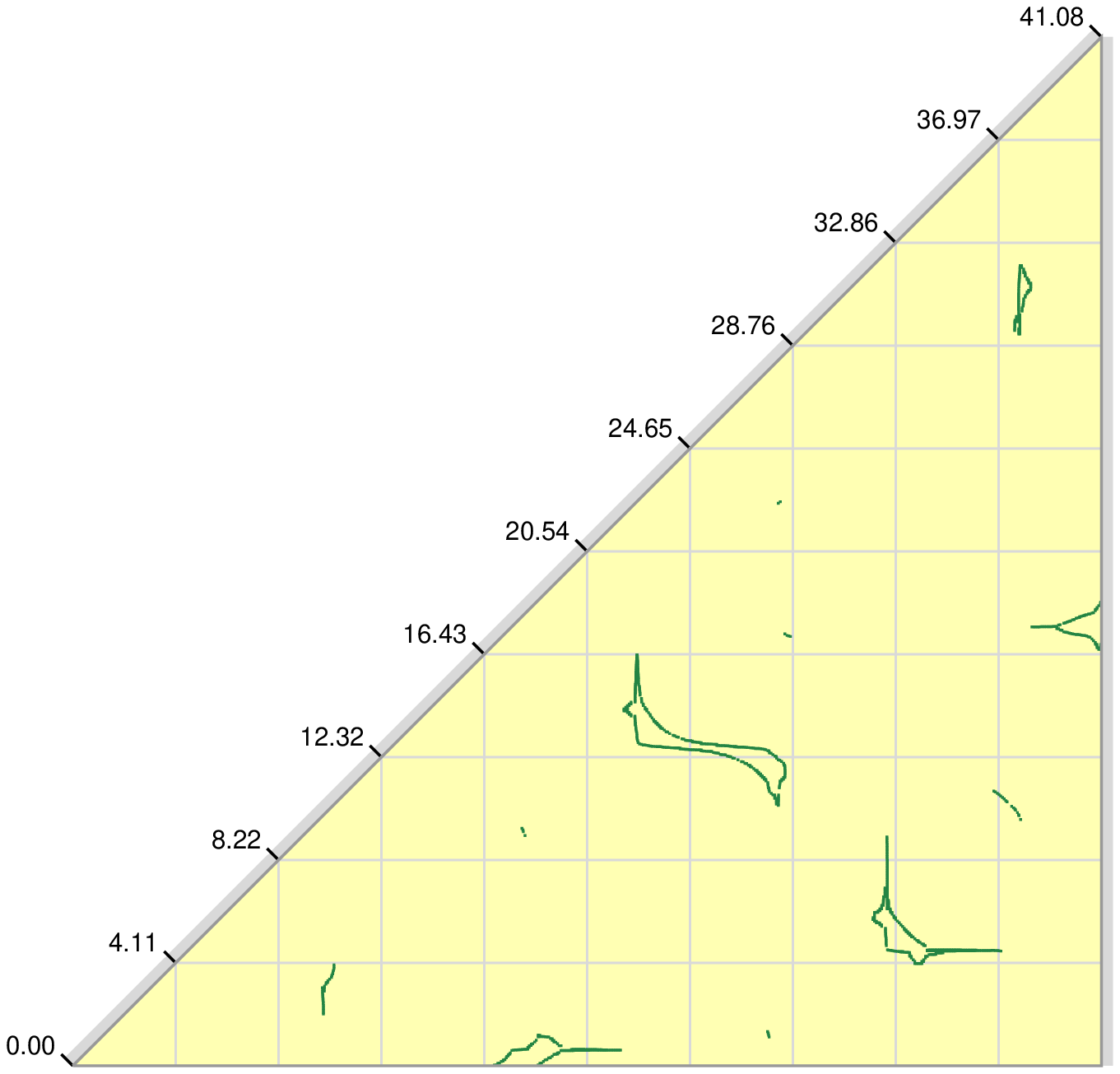}
        \put(8,94){\scriptsize{$82.17$}}
        \put(8,89){\scriptsize{$82.14$}}
        \put(8,84){\scriptsize{$587$}}
    \end{overpic}
\end{minipage} 
\hfill
\begin{minipage}[t]{6in}
  \vspace{2mm}
    \begin{overpic}[height=2.8in]{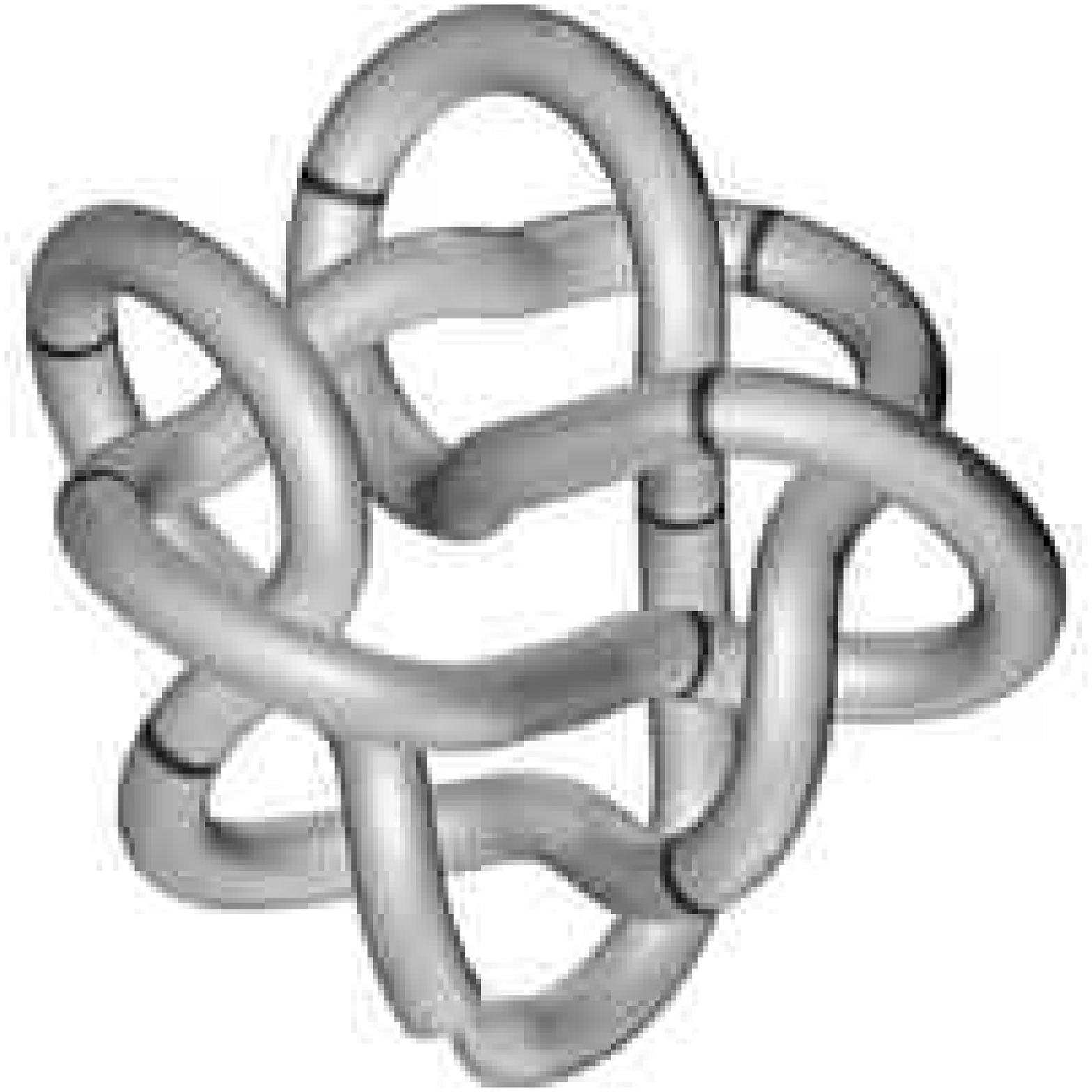}
        \put(-10,90){\large{$9_{18}$}}
    \end{overpic}
      \hspace{7mm}
    \begin{overpic}[width=2.8in]{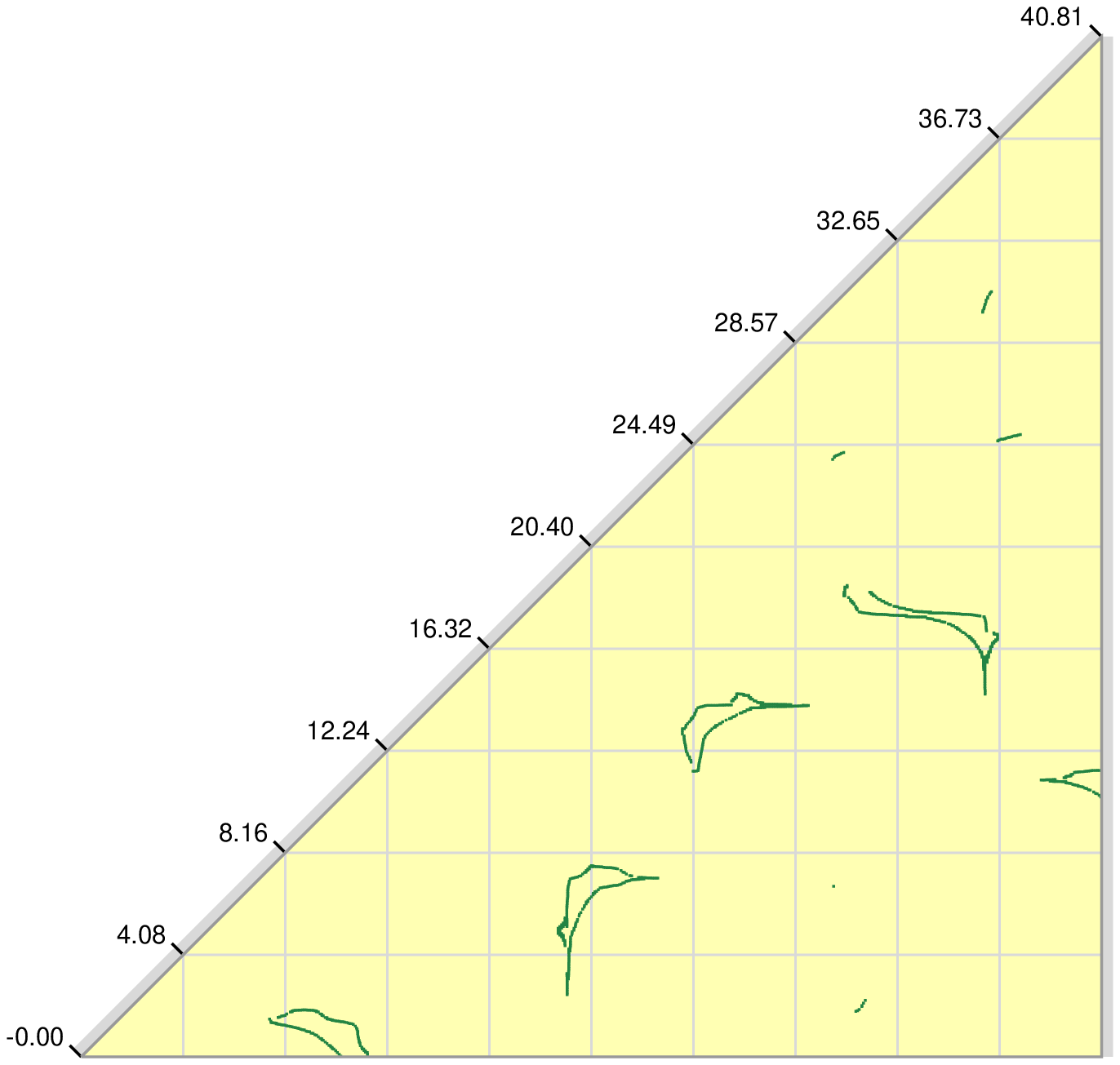}
        \put(8,94){\scriptsize{$81.63$}}
        \put(8,89){\scriptsize{$81.60$}}
        \put(8,84){\scriptsize{$583$}}
    \end{overpic}
\end{minipage} 
\hfill
\begin{minipage}[t]{6in}
  \vspace{2mm}
    \begin{overpic}[width=2.8in]{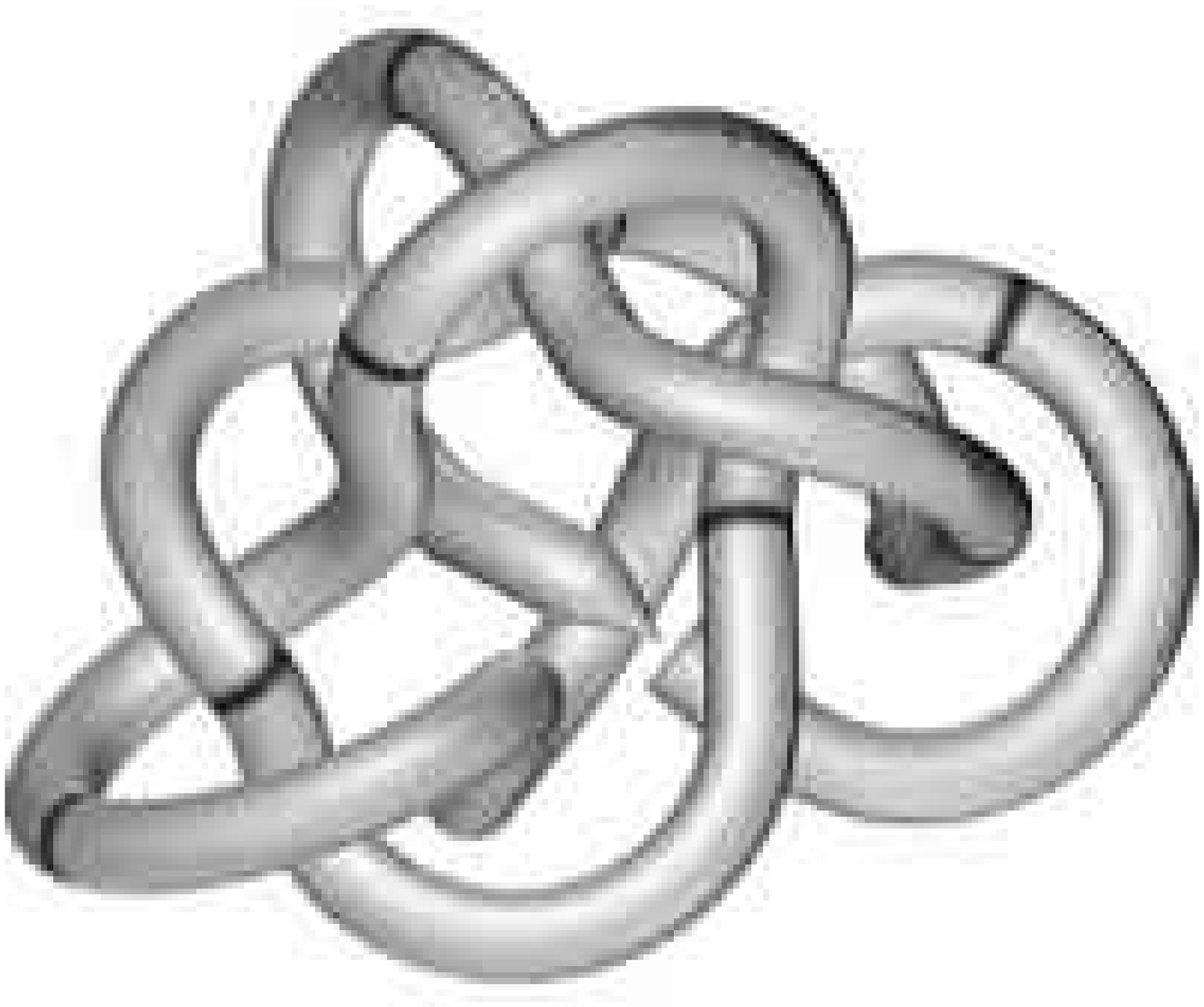}
        \put(-10,90){\large{$9_{22}$}}
    \end{overpic}
      \hspace{7mm}
    \begin{overpic}[width=2.8in]{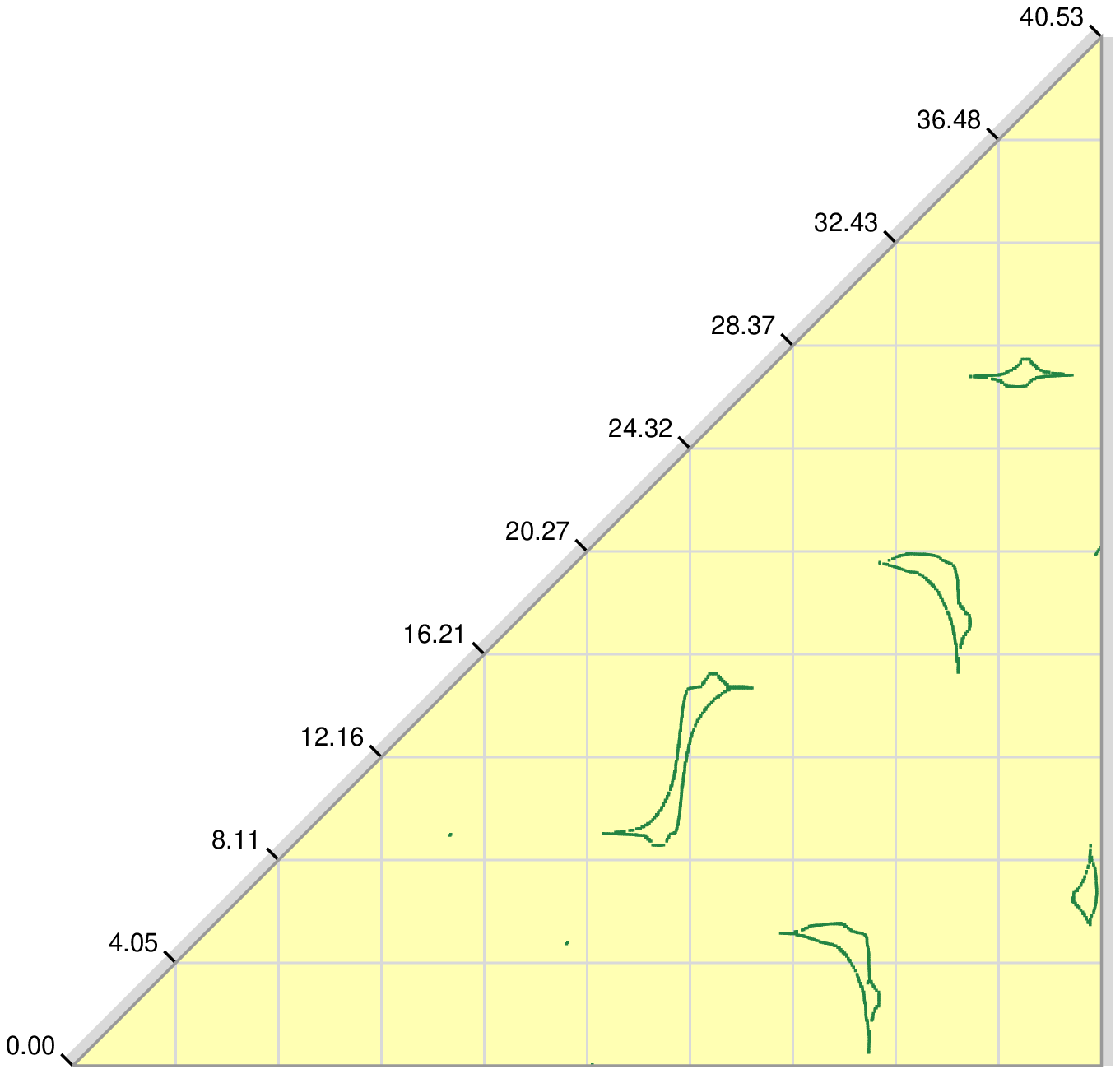}
        \put(8,94){\scriptsize{$81.08$}}
        \put(8,89){\scriptsize{$81.05$}}
        \put(8,84){\scriptsize{$579$}}
    \end{overpic}
\end{minipage} 
\hfill
\end{figure}
\clearpage
\pagebreak
\begin{figure}
\begin{minipage}[t]{6in}
  \vspace{2mm}
    \begin{overpic}[height=2.8in]{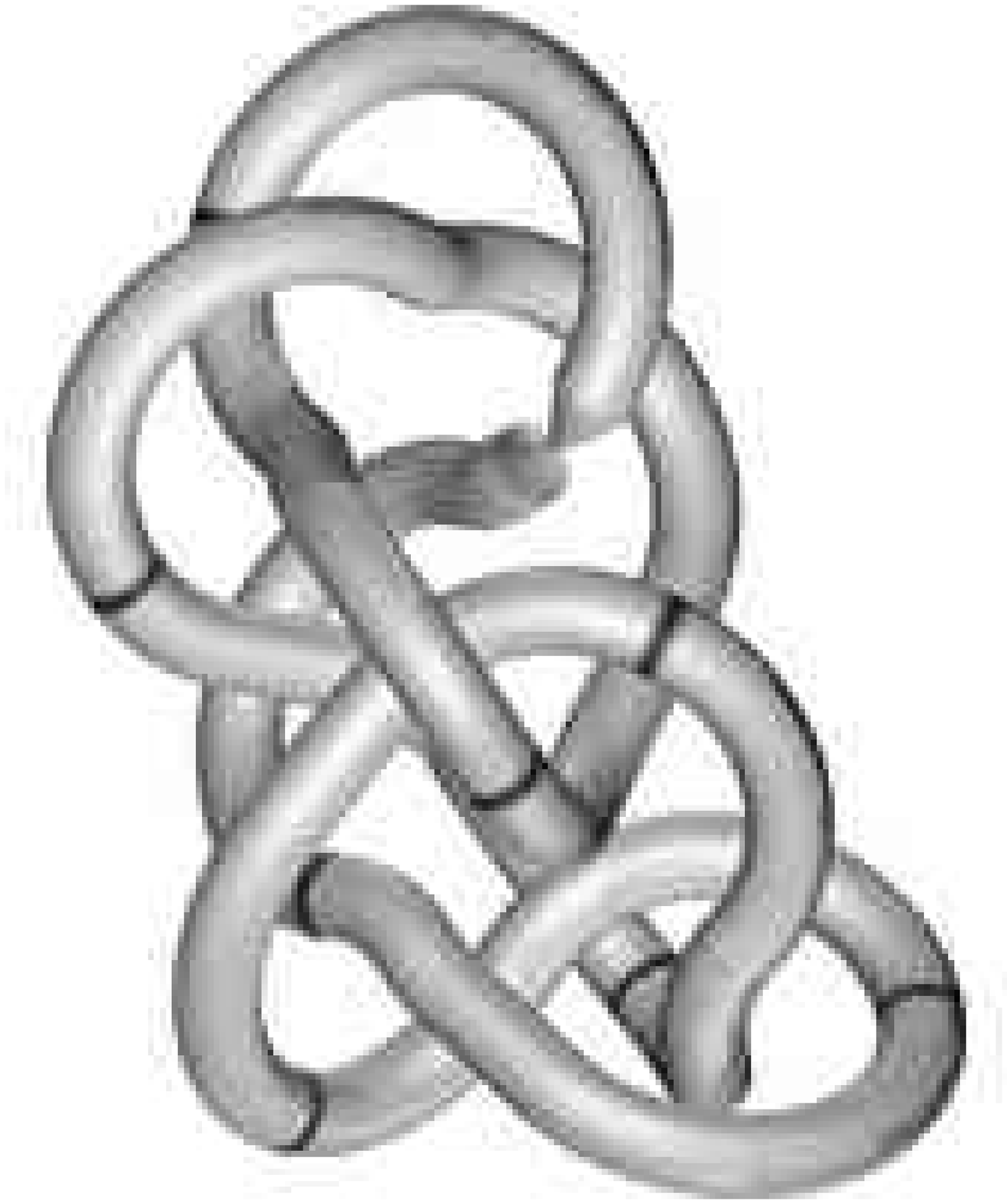}
        \put(-10,90){\large{$9_{24}$}}
    \end{overpic}
      \hspace{7mm}
    \begin{overpic}[width=2.8in]{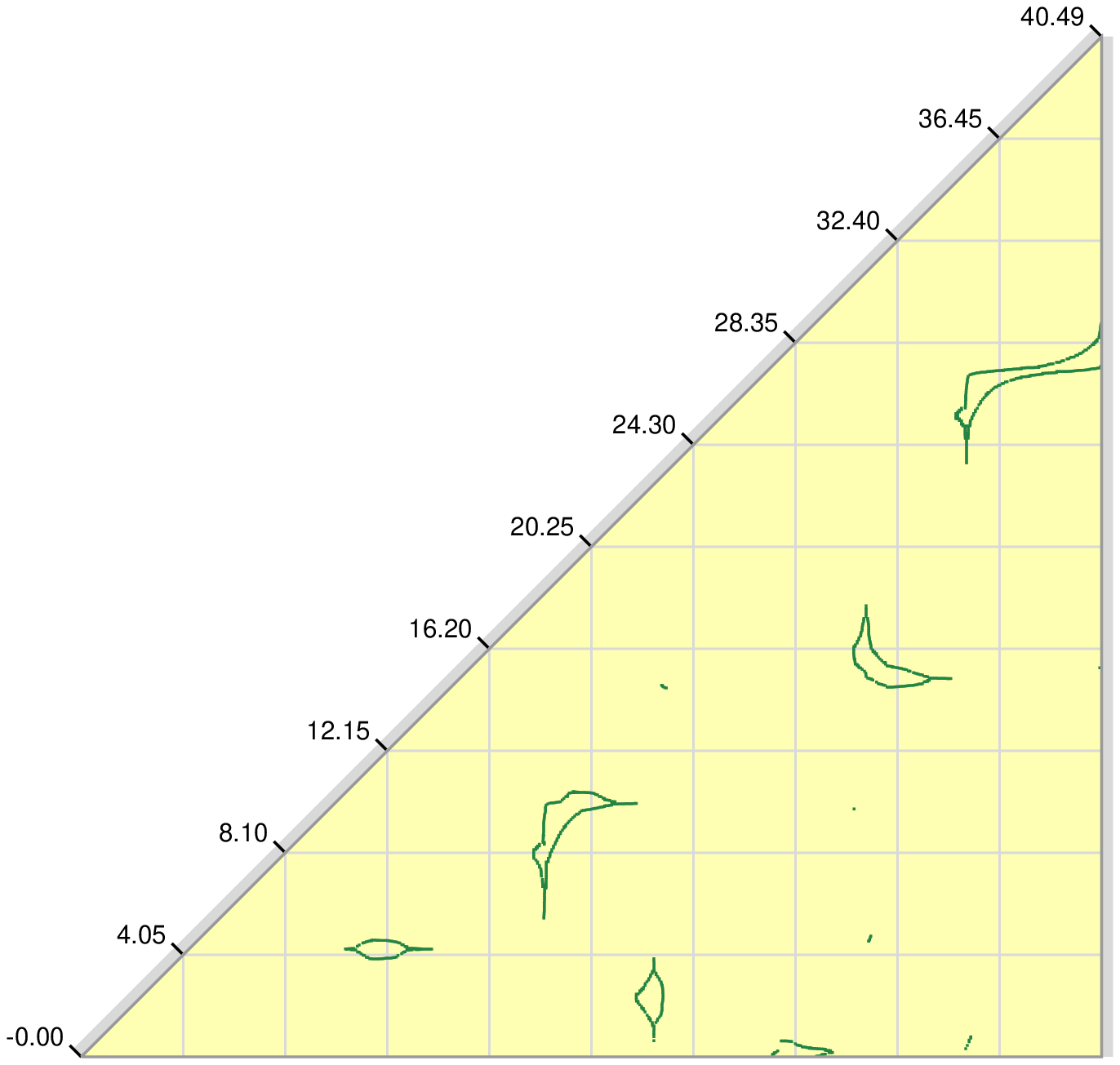}
        \put(8,94){\scriptsize{$81.00$}}
        \put(8,89){\scriptsize{$80.98$}}
        \put(8,84){\scriptsize{$578$}}
    \end{overpic}
\end{minipage} 
\hfill
\begin{minipage}[t]{6in}
  \vspace{2mm}
    \begin{overpic}[height=2.8in]{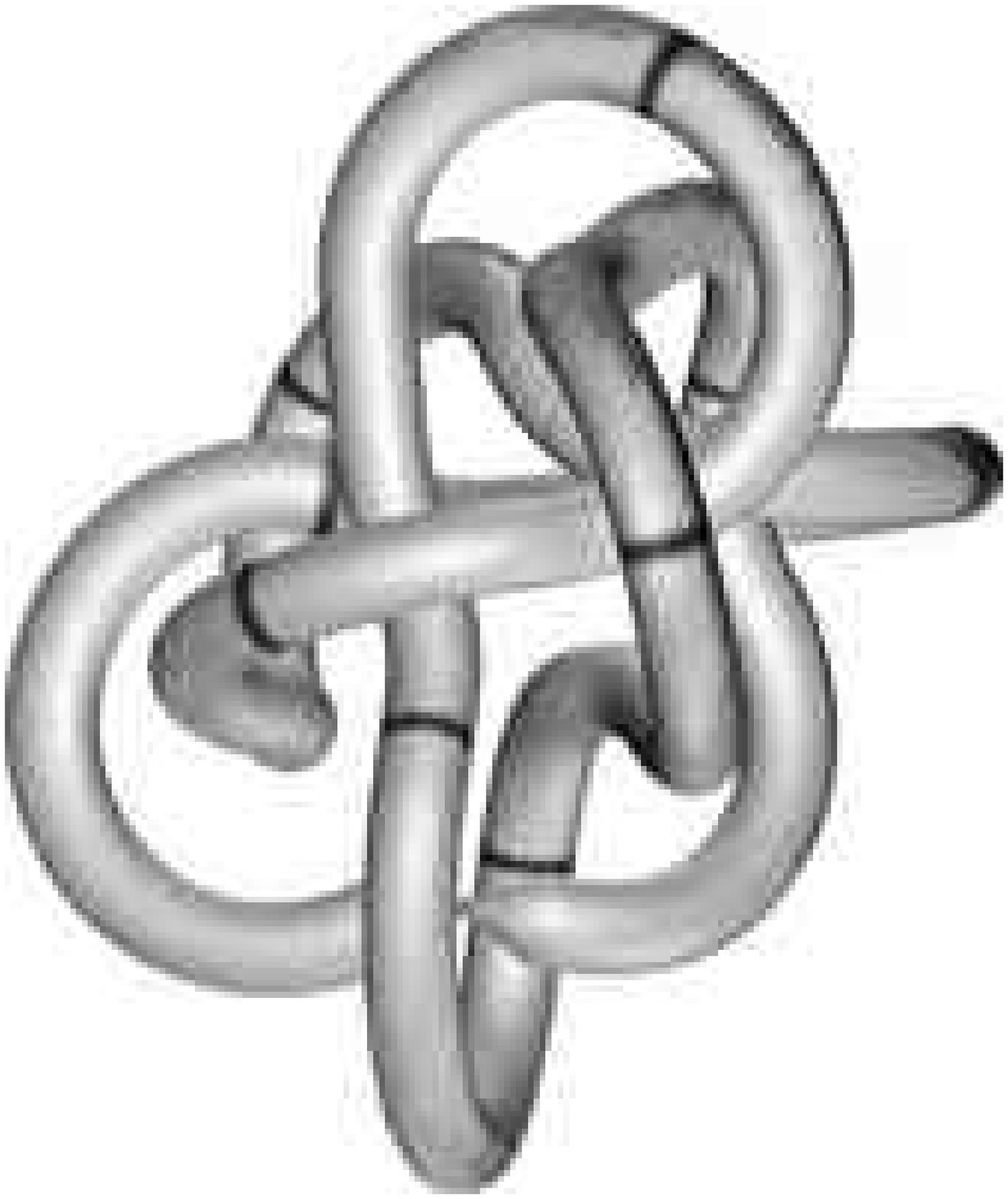}
        \put(-10,90){\large{$9_{25}$}}
    \end{overpic}
      \hspace{7mm}
    \begin{overpic}[width=2.8in]{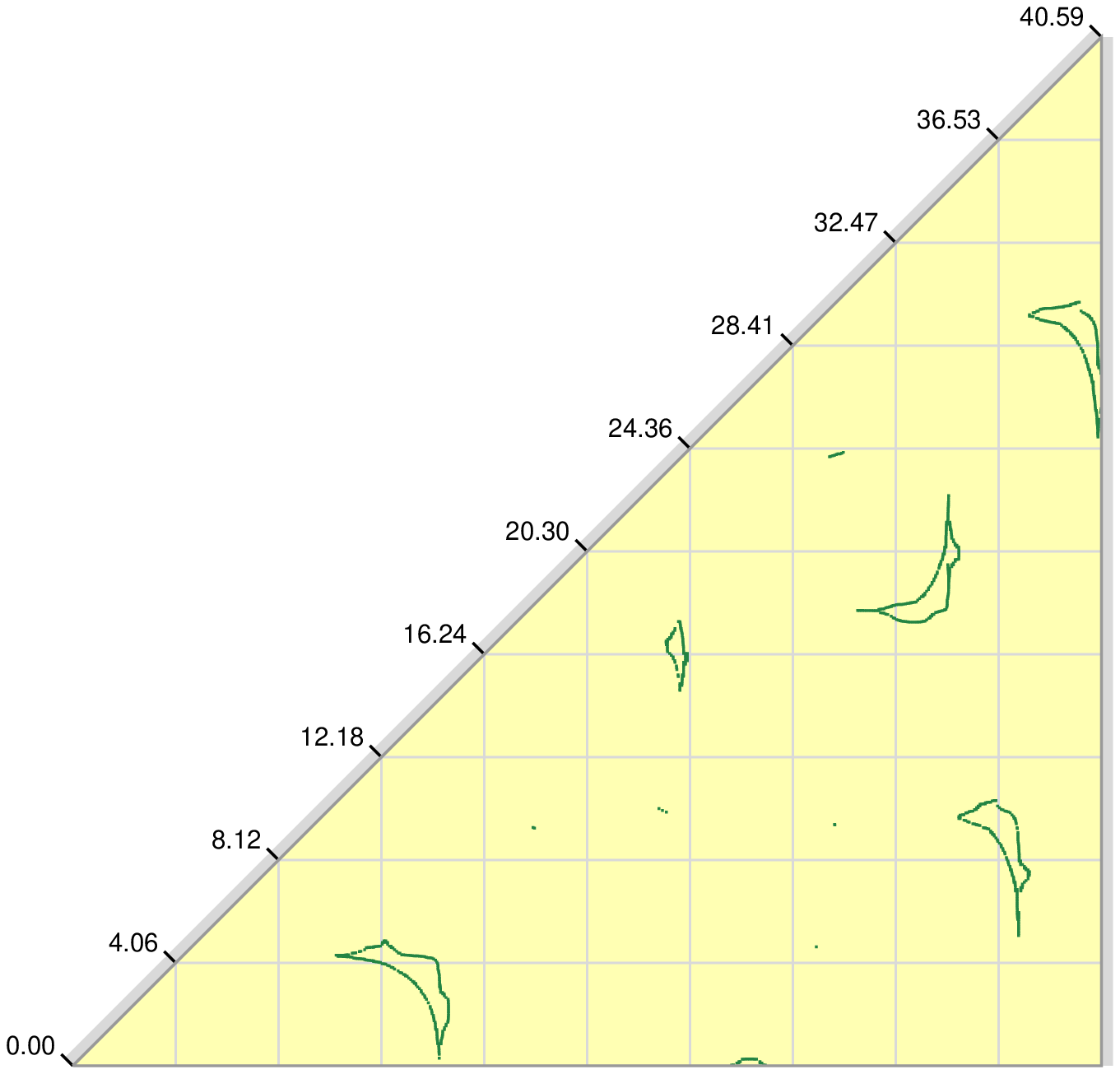}
        \put(8,94){\scriptsize{$81.19$}}
        \put(8,89){\scriptsize{$81.17$}}
        \put(8,84){\scriptsize{$580$}}
    \end{overpic}
\end{minipage} 
\hfill
\begin{minipage}[t]{6in}
  \vspace{2mm}
    \begin{overpic}[width=2.8in]{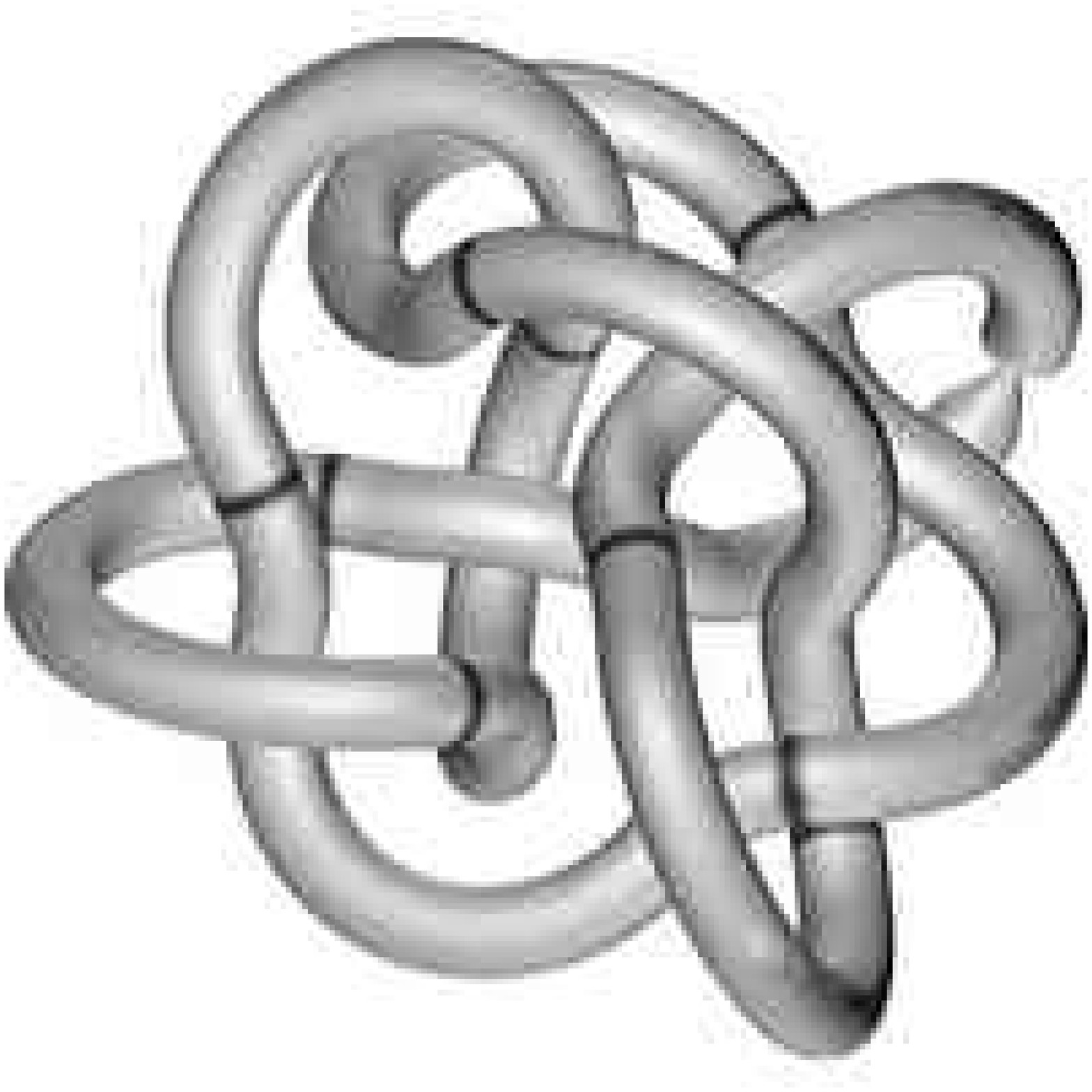}
        \put(-10,90){\large{$9_{30}$}}
    \end{overpic}
      \hspace{7mm}
    \begin{overpic}[width=2.8in]{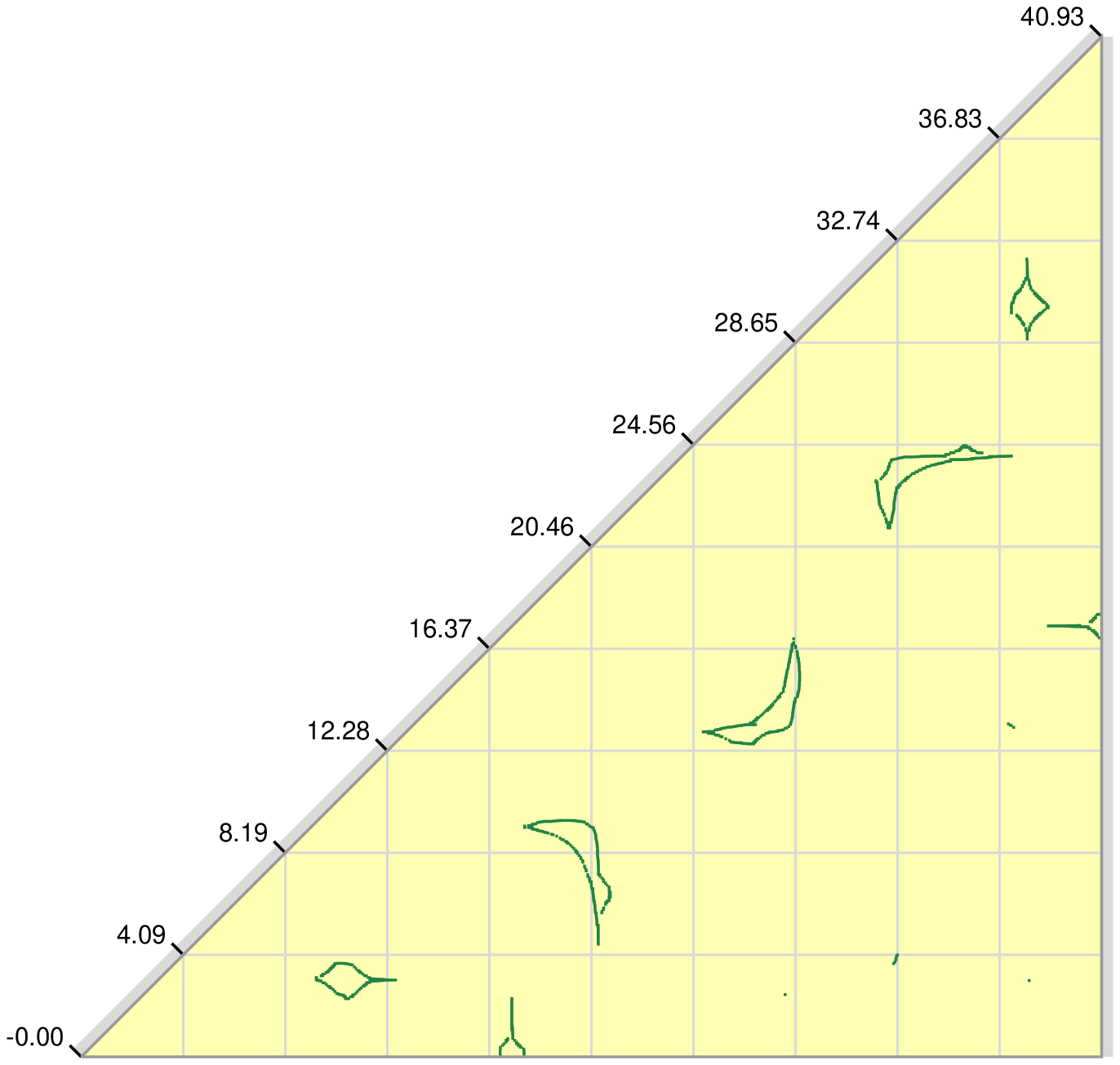}
        \put(8,94){\scriptsize{$81.86$}}
        \put(8,89){\scriptsize{$81.84$}}
        \put(8,84){\scriptsize{$585$}}
    \end{overpic}
\end{minipage} 
\hfill
\end{figure}
\clearpage
\pagebreak
\begin{figure}
\begin{minipage}[t]{6in}
  \vspace{2mm}
    \begin{overpic}[height=2.8in]{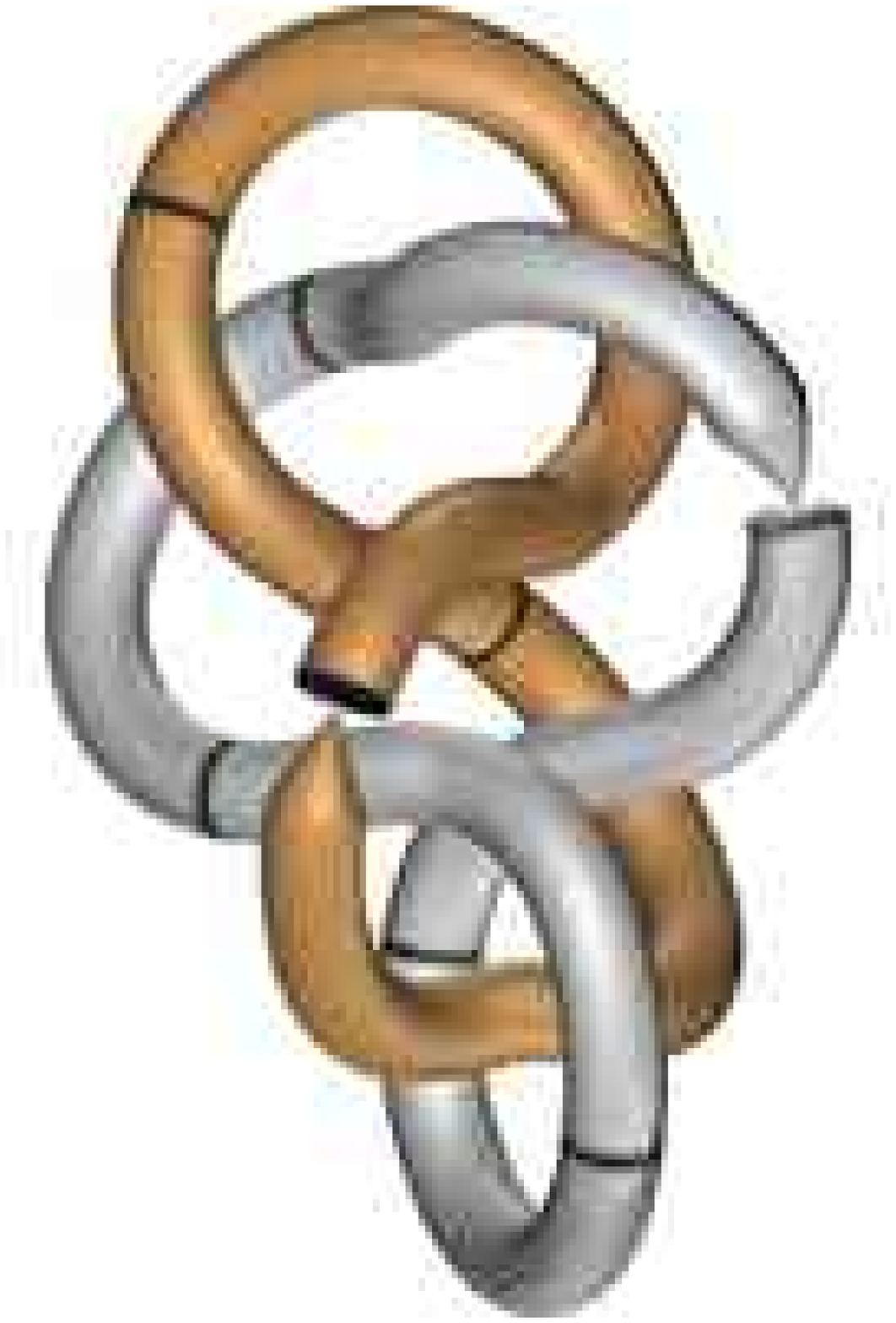}
        \put(-10,90){\large{$6^{2}_{2}$}}
    \end{overpic}
      \hspace{7mm}
    \begin{overpic}[width=2.8in]{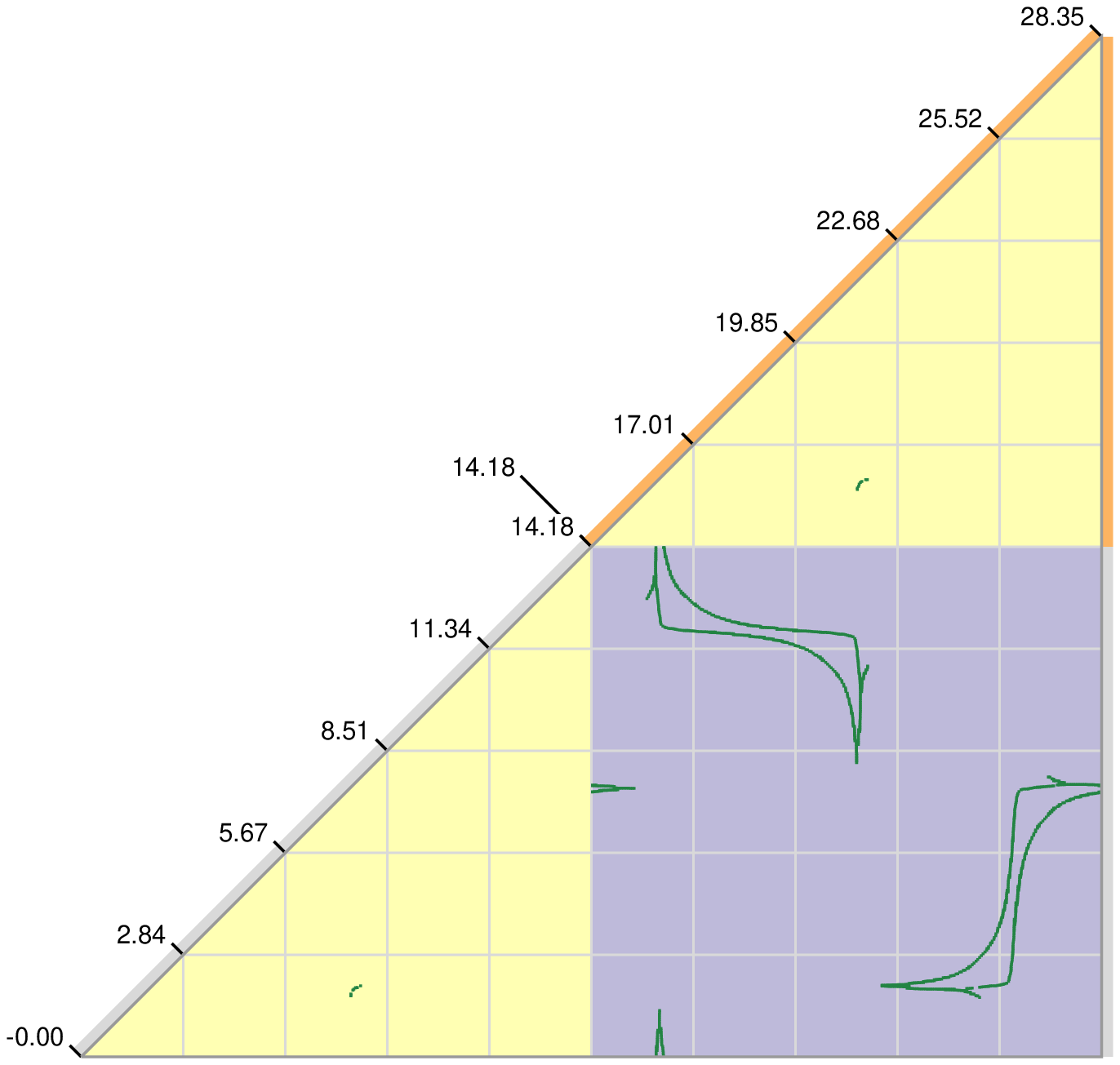}
        \put(8,94){\scriptsize{$56.72$}}
        \put(8,89){\scriptsize{$56.71$}}
        \put(8,84){\scriptsize{$562$}}
    \end{overpic}
\end{minipage} 
\hfill
\begin{minipage}[t]{6in}
  \vspace{2mm}
    \begin{overpic}[width=2.8in]{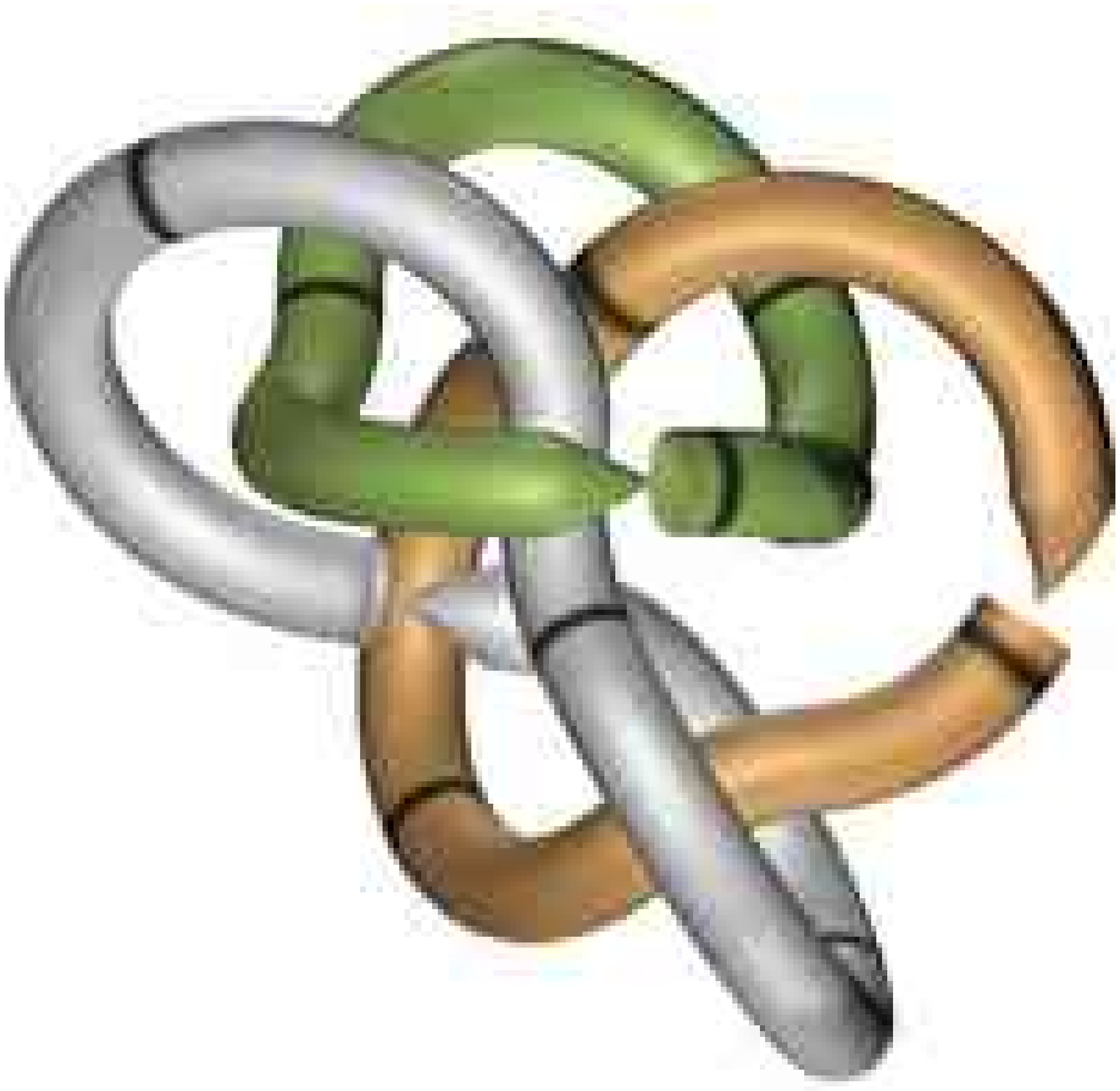}
        \put(-10,90){\large{$6^{3}_{1}$}}
    \end{overpic}
      \hspace{7mm}
    \begin{overpic}[width=2.8in]{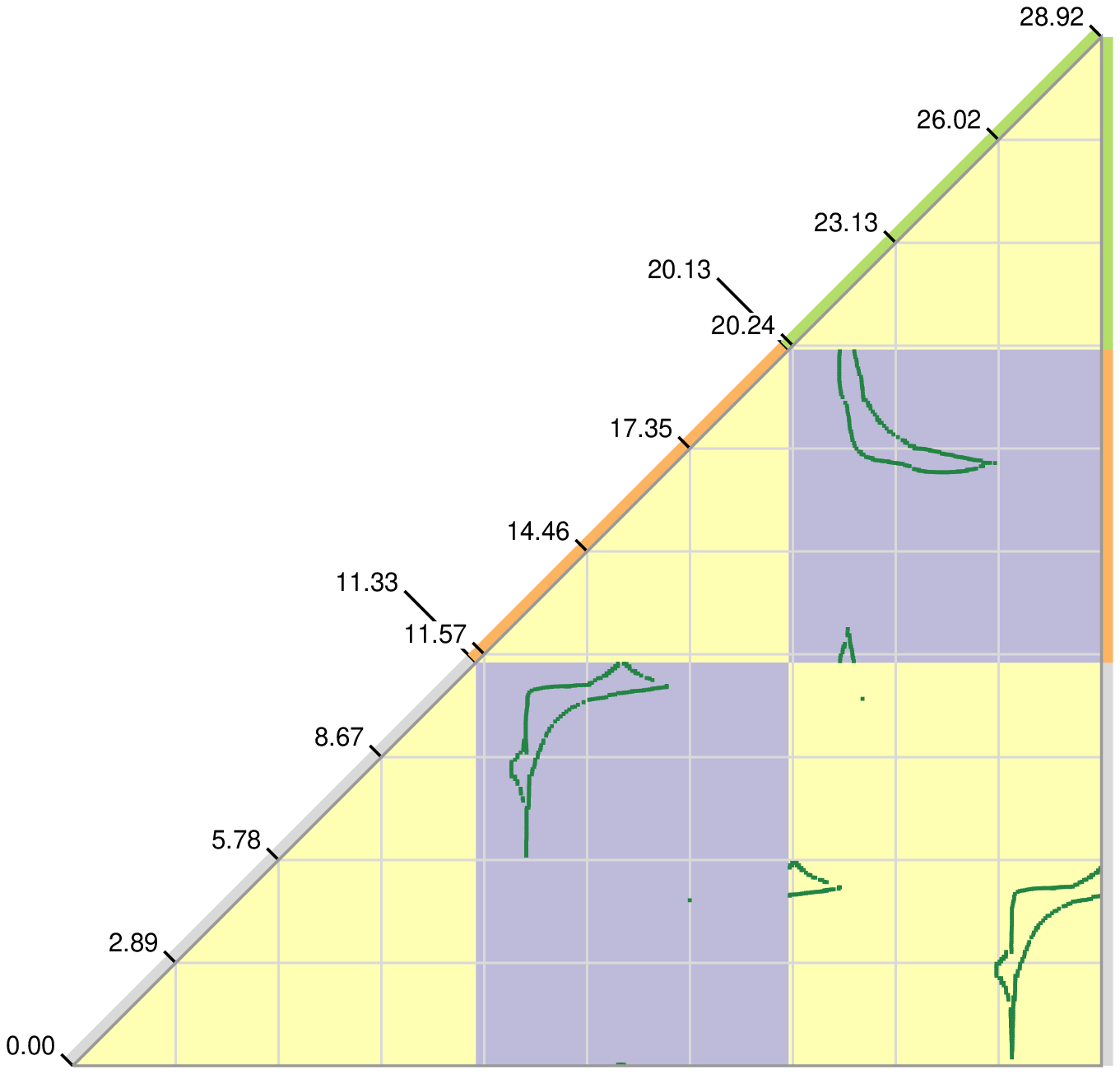}
        \put(8,94){\scriptsize{$57.84$}}
        \put(8,89){\scriptsize{$57.82$}}
        \put(8,84){\scriptsize{$413$}}
    \end{overpic}
\end{minipage} 
\hfill
\begin{minipage}[t]{6in}
  \vspace{2mm}
    \begin{overpic}[height=2.8in]{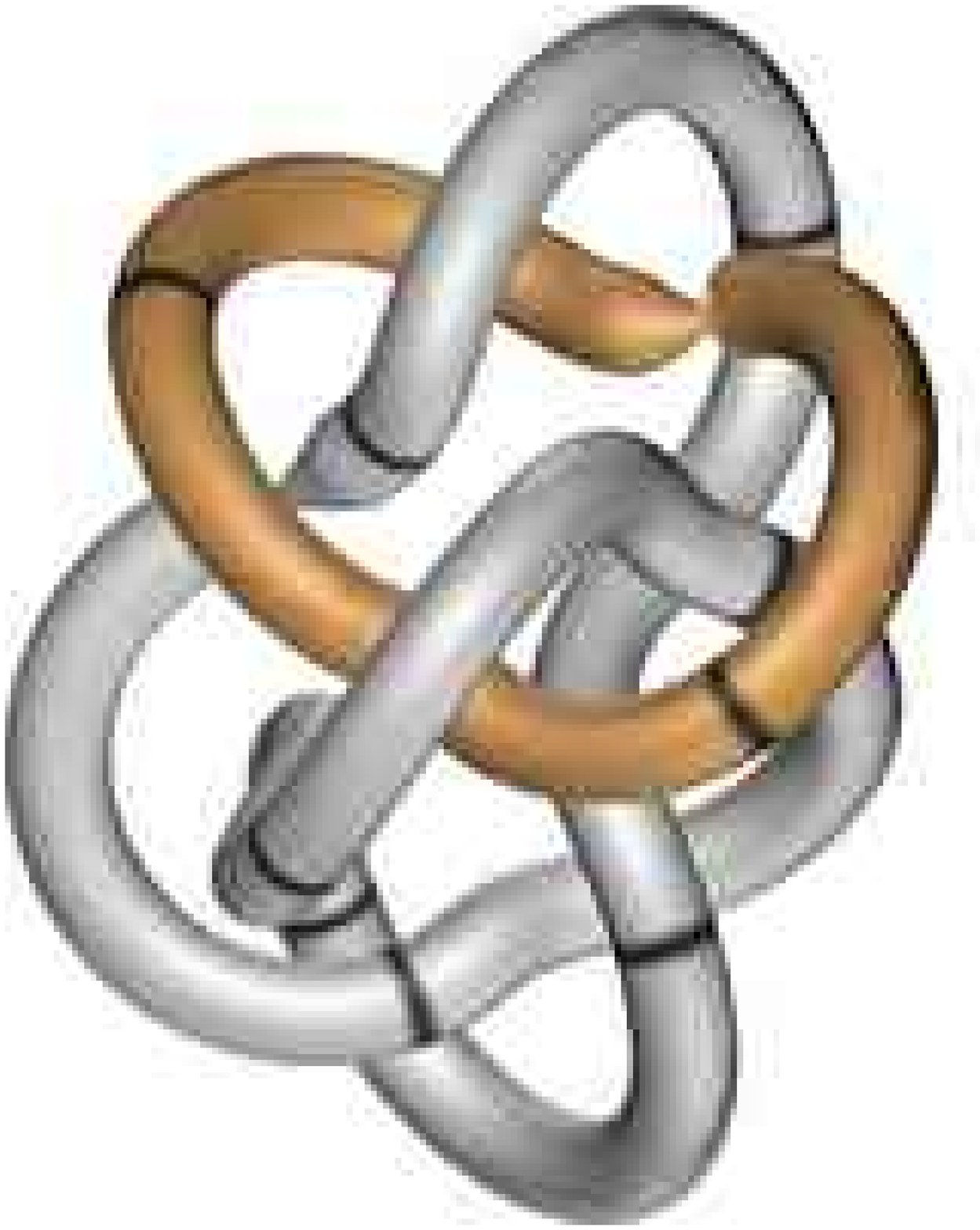}
        \put(-10,90){\large{$7^{2}_{4}$}}
    \end{overpic}
      \hspace{7mm}
    \begin{overpic}[width=2.8in]{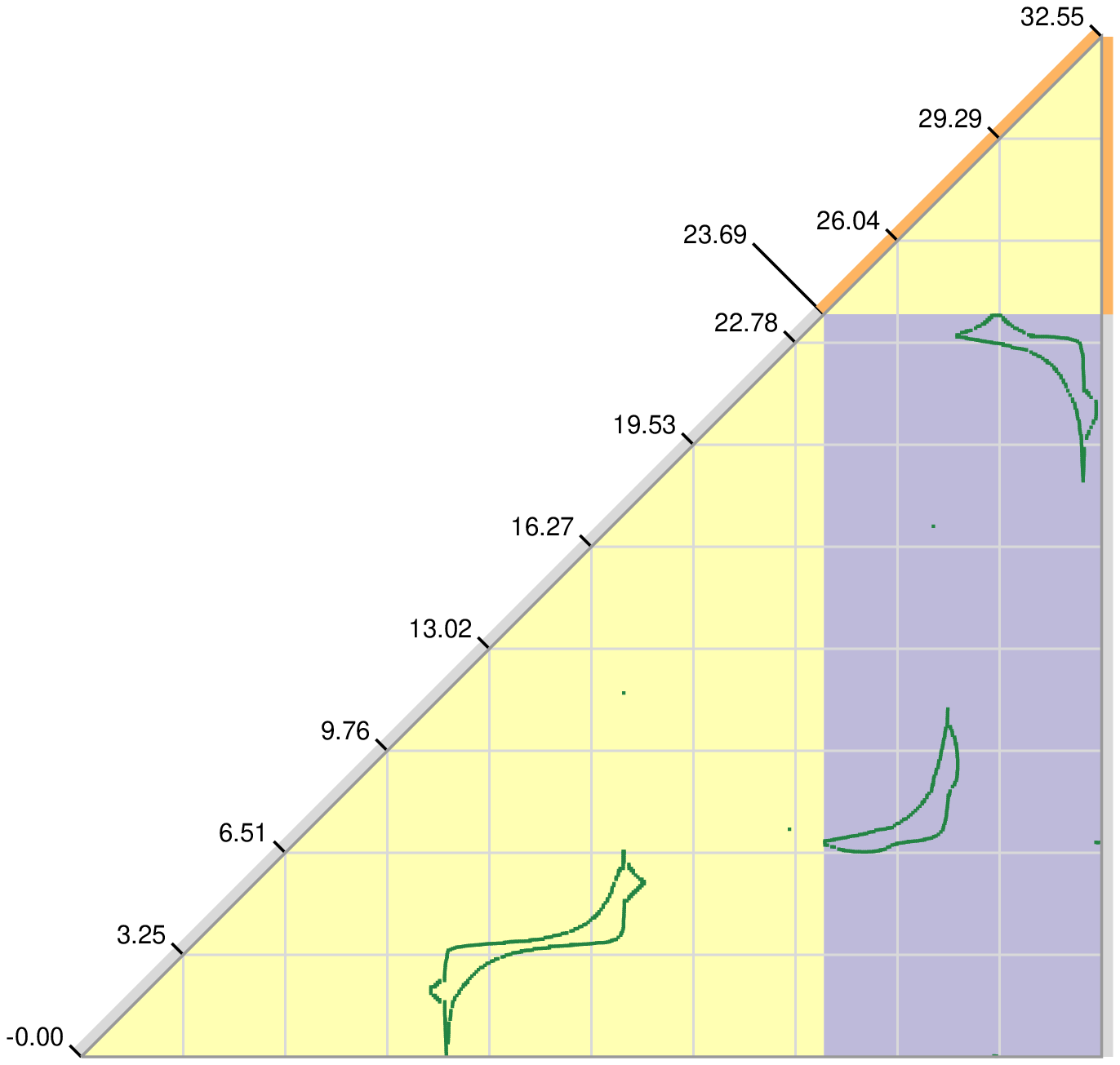}
        \put(8,94){\scriptsize{$65.10$}}
        \put(8,89){\scriptsize{$65.08$}}
        \put(8,84){\scriptsize{$465$}}
    \end{overpic}
\end{minipage} 
\hfill
\end{figure}
\clearpage
\pagebreak
\begin{figure}
\begin{minipage}[t]{6in}
  \vspace{2mm}
    \begin{overpic}[width=2.8in]{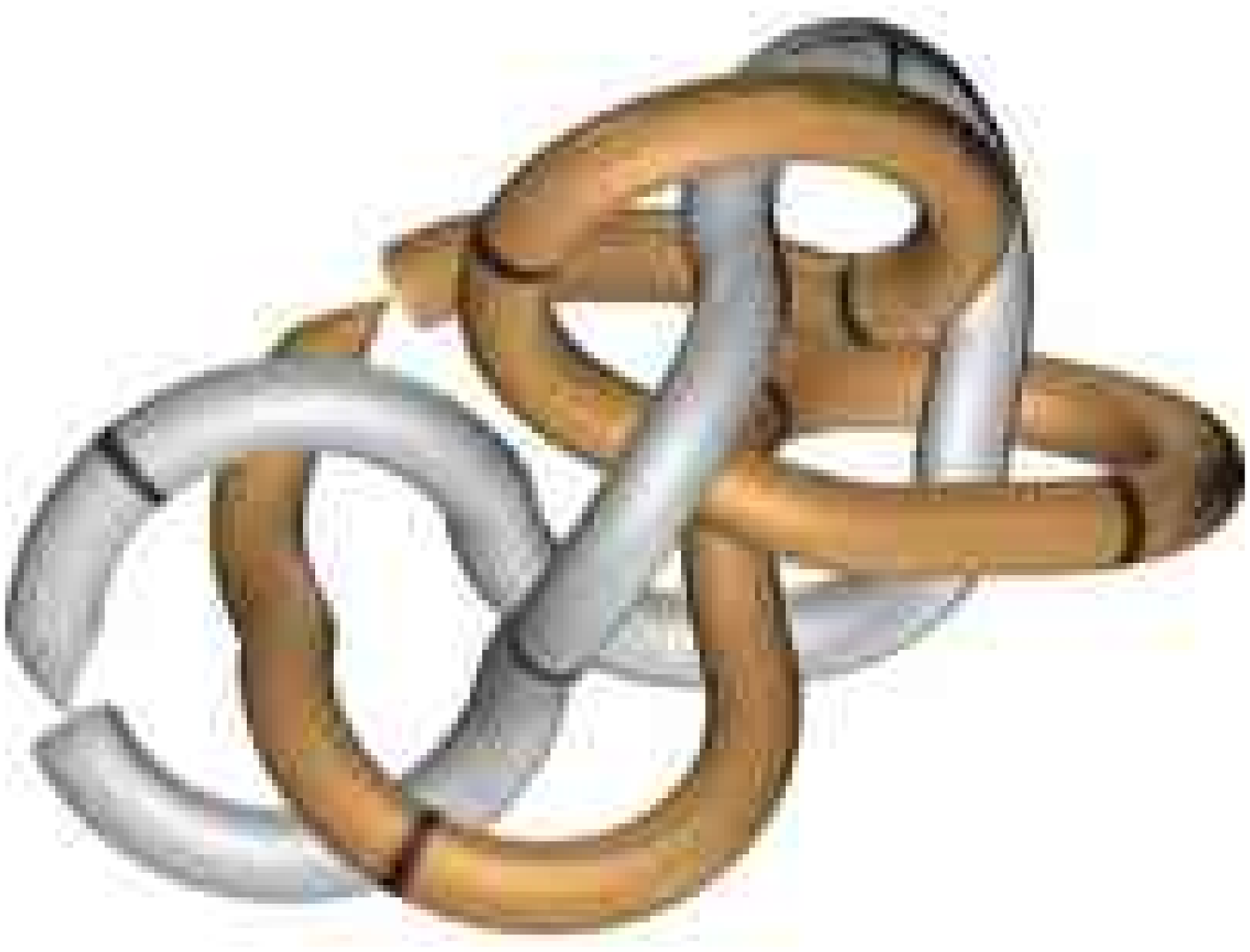}
        \put(-10,90){\large{$8^{2}_{4}$}}
    \end{overpic}
      \hspace{7mm}
    \begin{overpic}[width=2.8in]{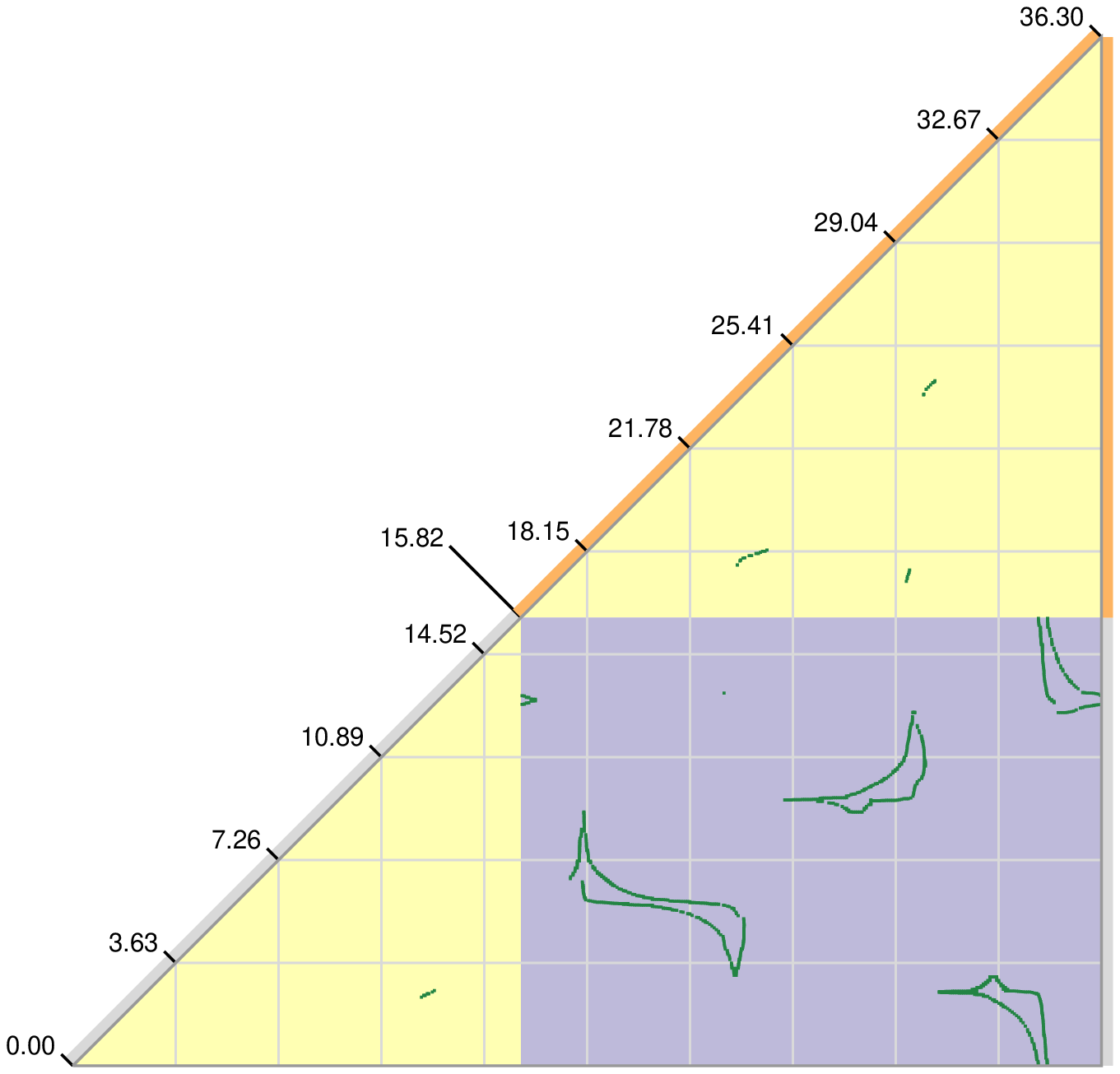}
        \put(8,94){\scriptsize{$72.62$}}
        \put(8,89){\scriptsize{$72.60$}}
        \put(8,84){\scriptsize{$518$}}
    \end{overpic}
\end{minipage} 
\hfill
\begin{minipage}[t]{6in}
  \vspace{2mm}
    \begin{overpic}[height=2.8in]{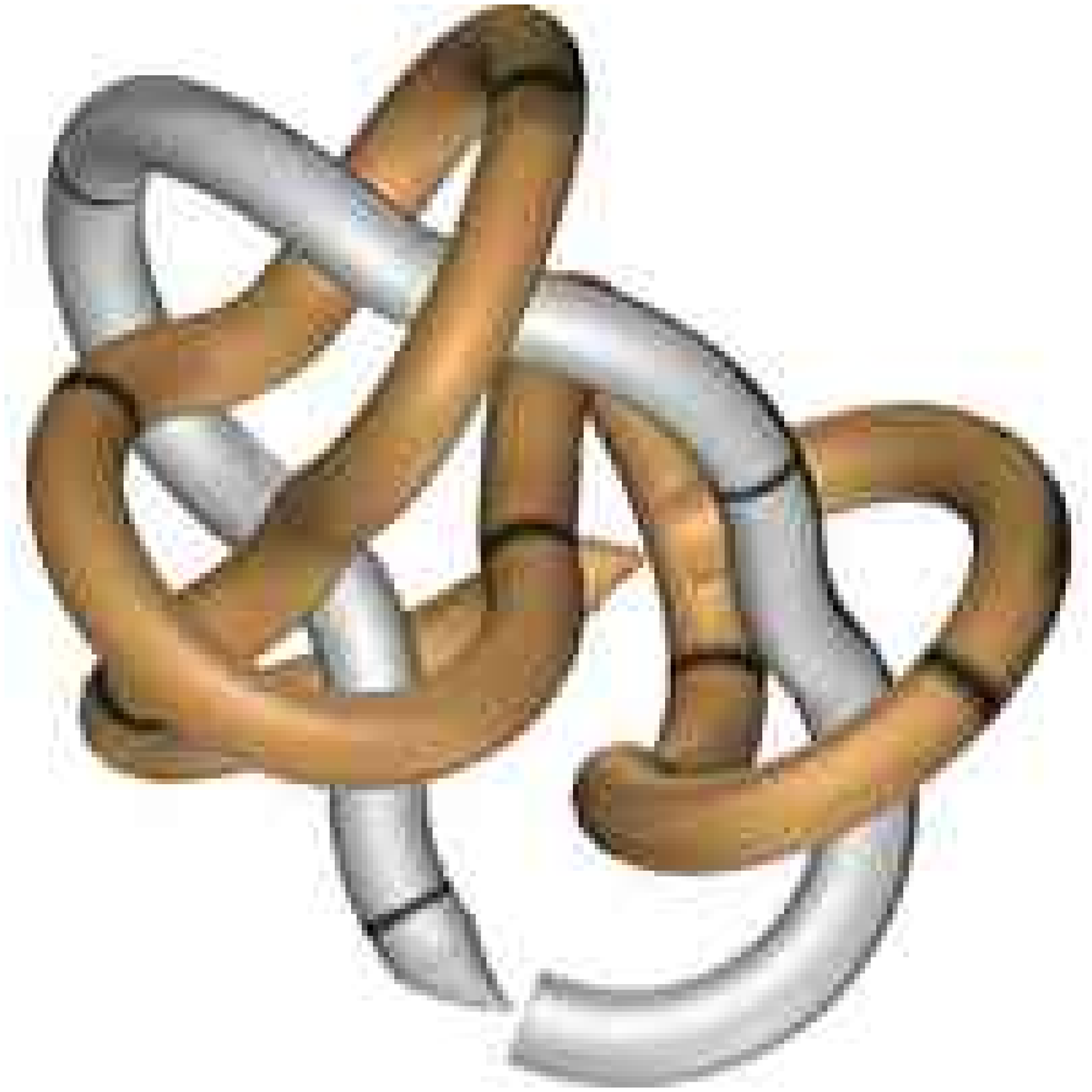}
        \put(-10,90){\large{$8^{2}_{7}$}}
    \end{overpic}
      \hspace{7mm}
    \begin{overpic}[width=2.8in]{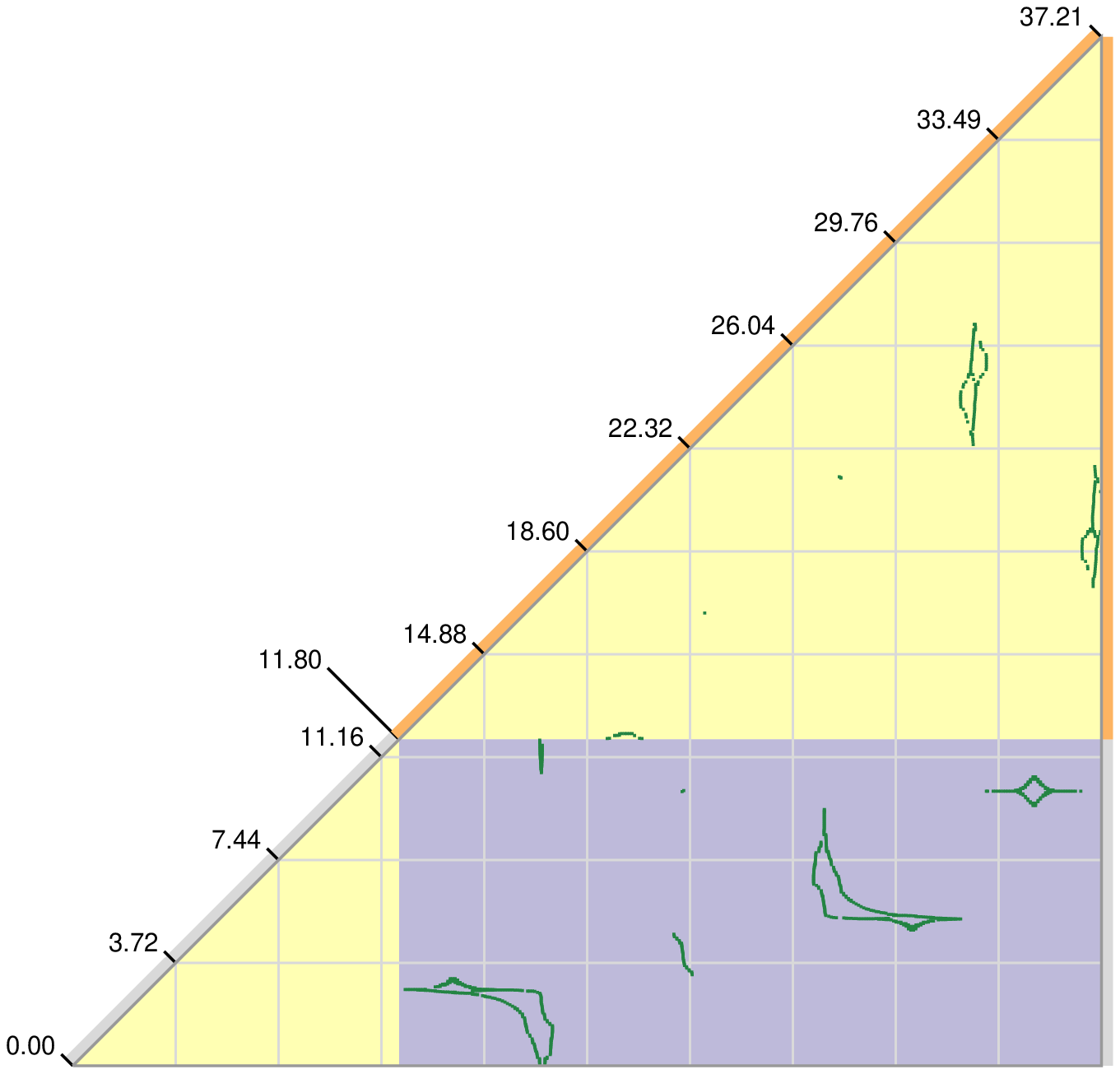}
        \put(8,94){\scriptsize{$74.42$}}
        \put(8,89){\scriptsize{$74.39$}}
        \put(8,84){\scriptsize{$531$}}
    \end{overpic}
\end{minipage} 
\hfill
\begin{minipage}[t]{6in}
  \vspace{2mm}
    \begin{overpic}[height=2.8in]{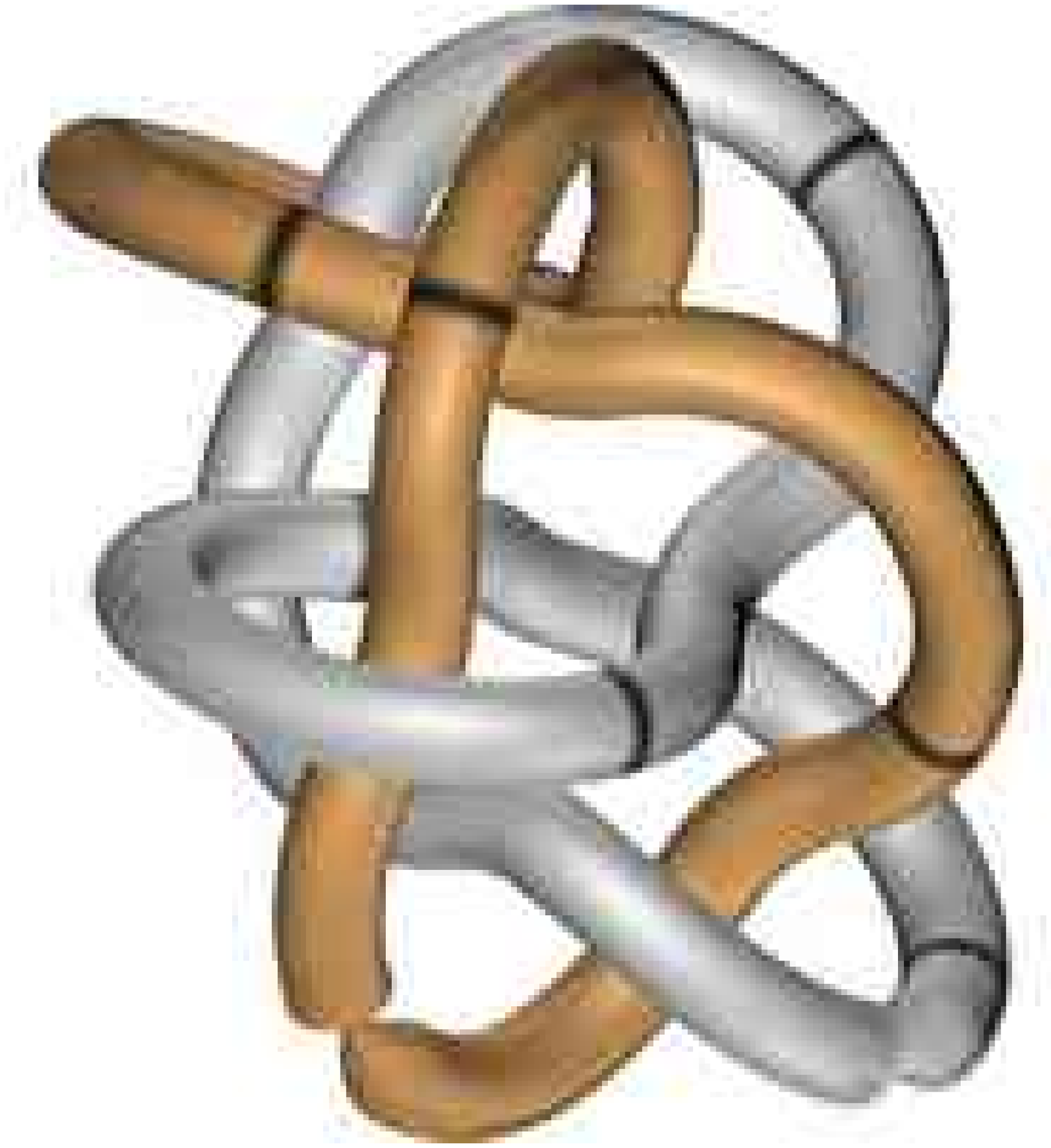}
        \put(-10,90){\large{$8^{2}_{8}$}}
    \end{overpic}
      \hspace{7mm}
    \begin{overpic}[width=2.8in]{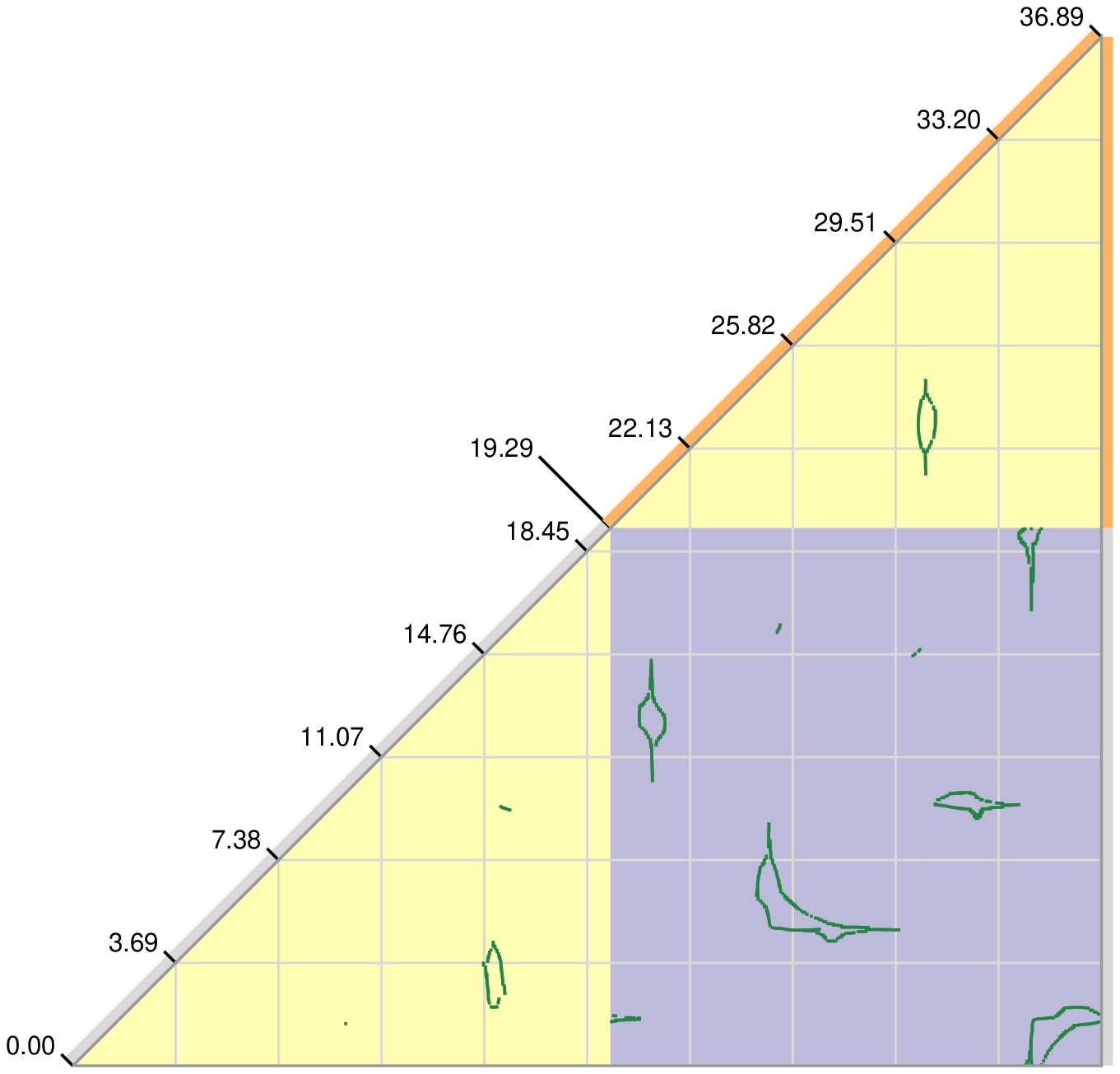}
        \put(8,94){\scriptsize{$73.79$}}
        \put(8,89){\scriptsize{$73.78$}}
        \put(8,84){\scriptsize{$527$}}
    \end{overpic}
\end{minipage} 
\hfill
\end{figure}
\clearpage
\pagebreak
\begin{figure}
\begin{minipage}[t]{6in}
  \vspace{2mm}
    \begin{overpic}[height=2.8in]{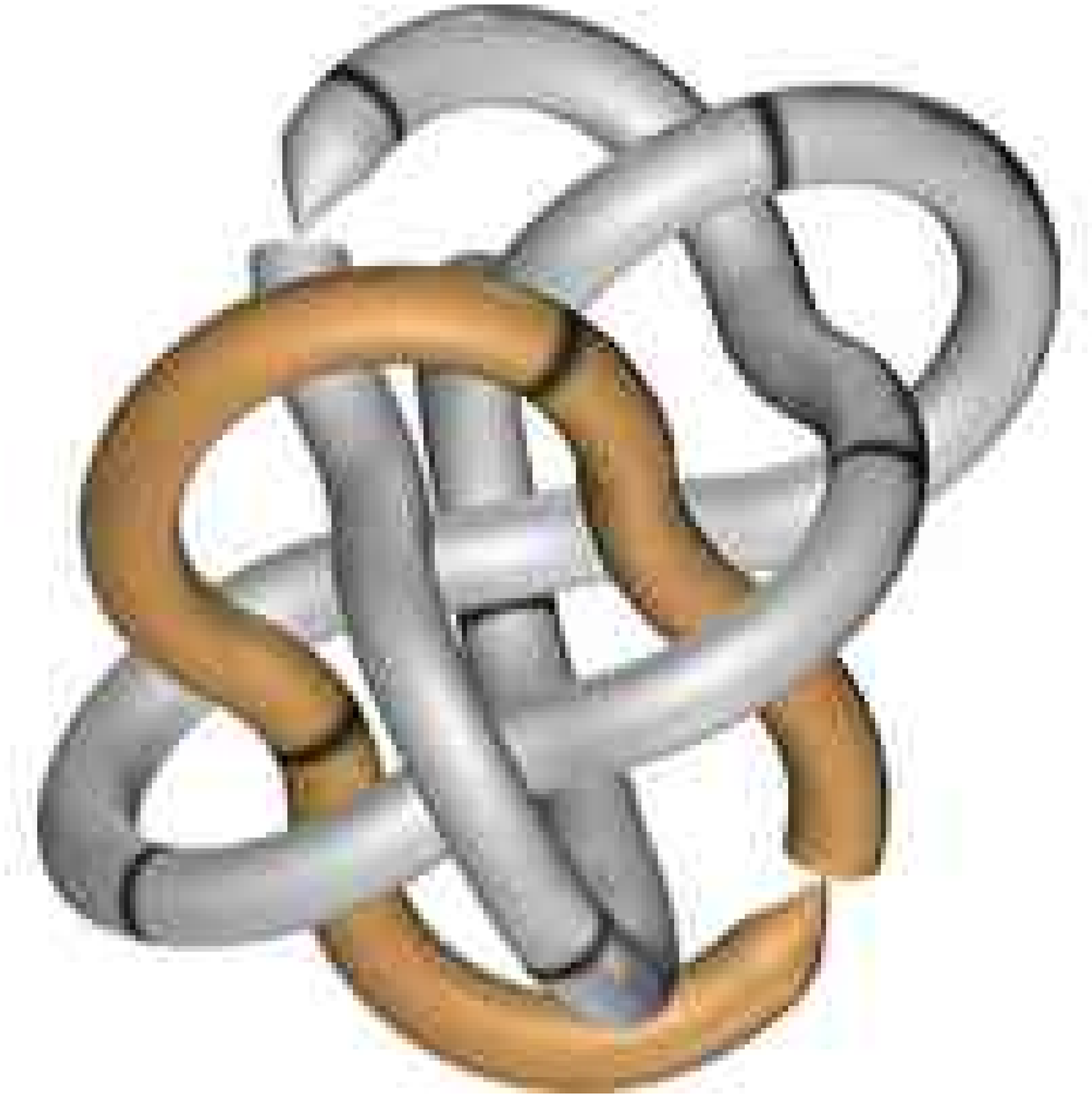}
        \put(-10,90){\large{$8^{2}_{13}$}}
    \end{overpic}
      \hspace{7mm}
    \begin{overpic}[width=2.8in]{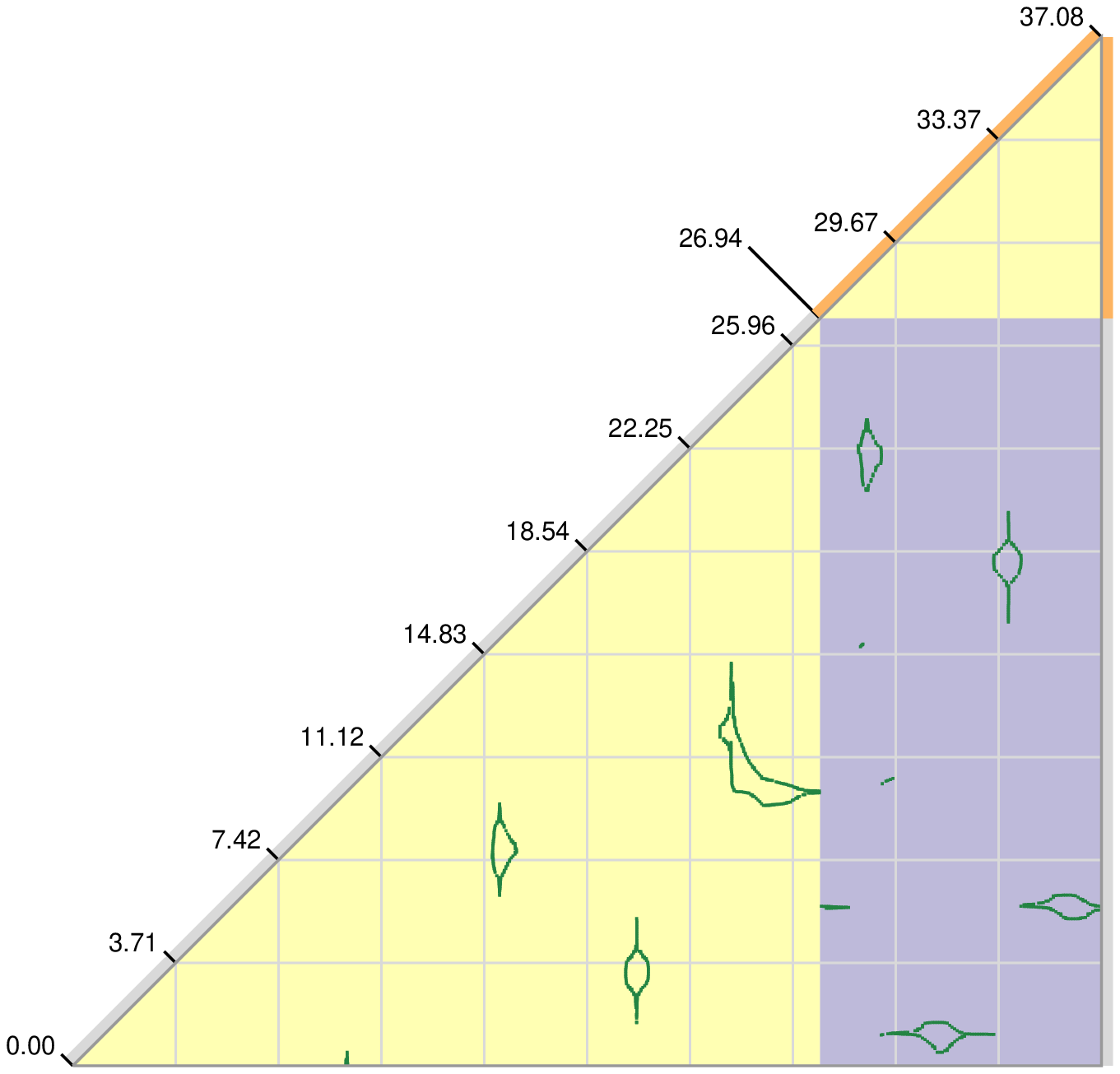}
        \put(8,94){\scriptsize{$74.17$}}
        \put(8,89){\scriptsize{$74.16$}}
        \put(8,84){\scriptsize{$529$}}
    \end{overpic}
\end{minipage} 
\hfill
\begin{minipage}[t]{6in}
  \vspace{2mm}
    \begin{overpic}[width=2.8in]{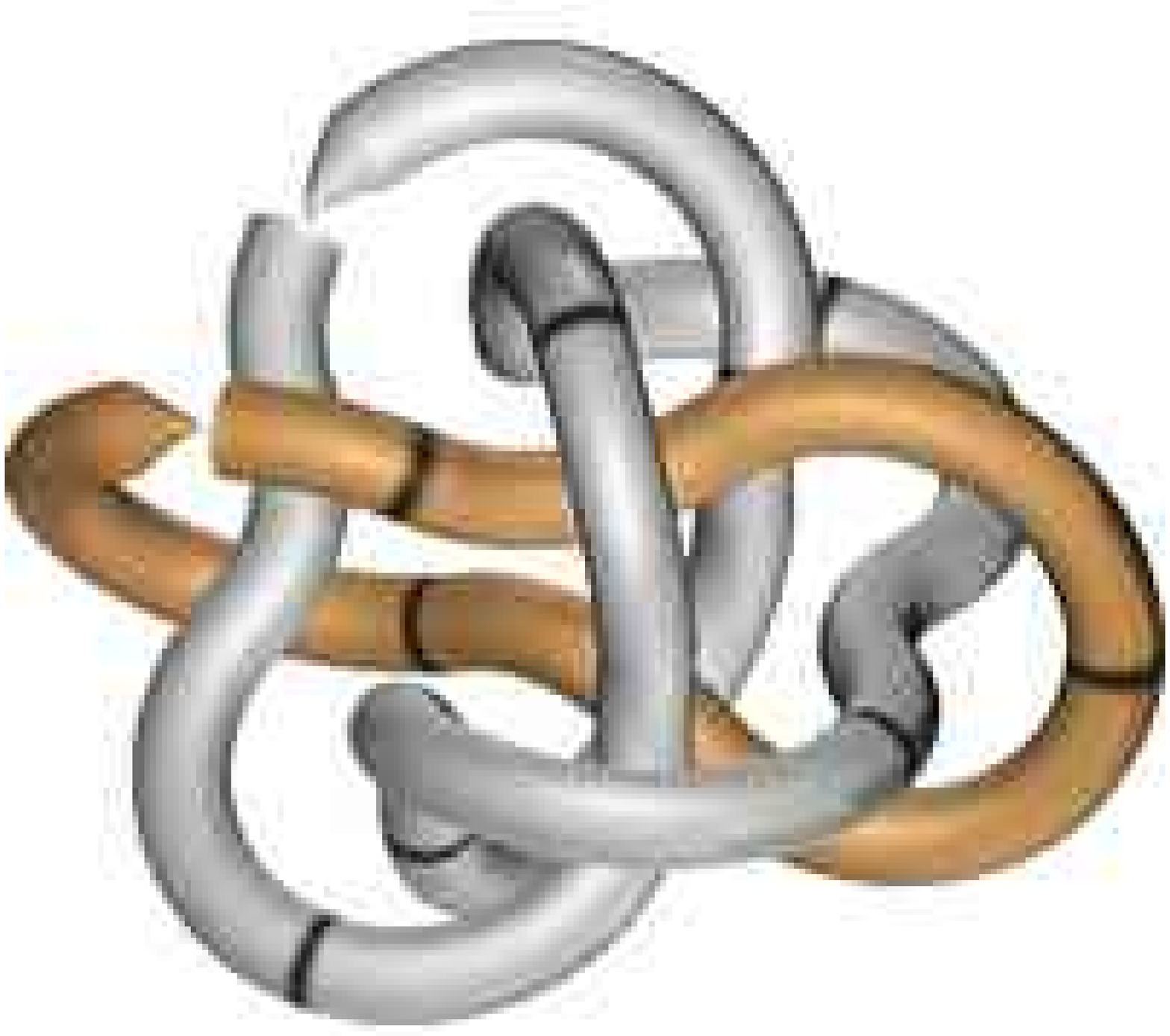}
        \put(-10,90){\large{$8^{2}_{14}$}}
    \end{overpic}
      \hspace{7mm}
    \begin{overpic}[width=2.8in]{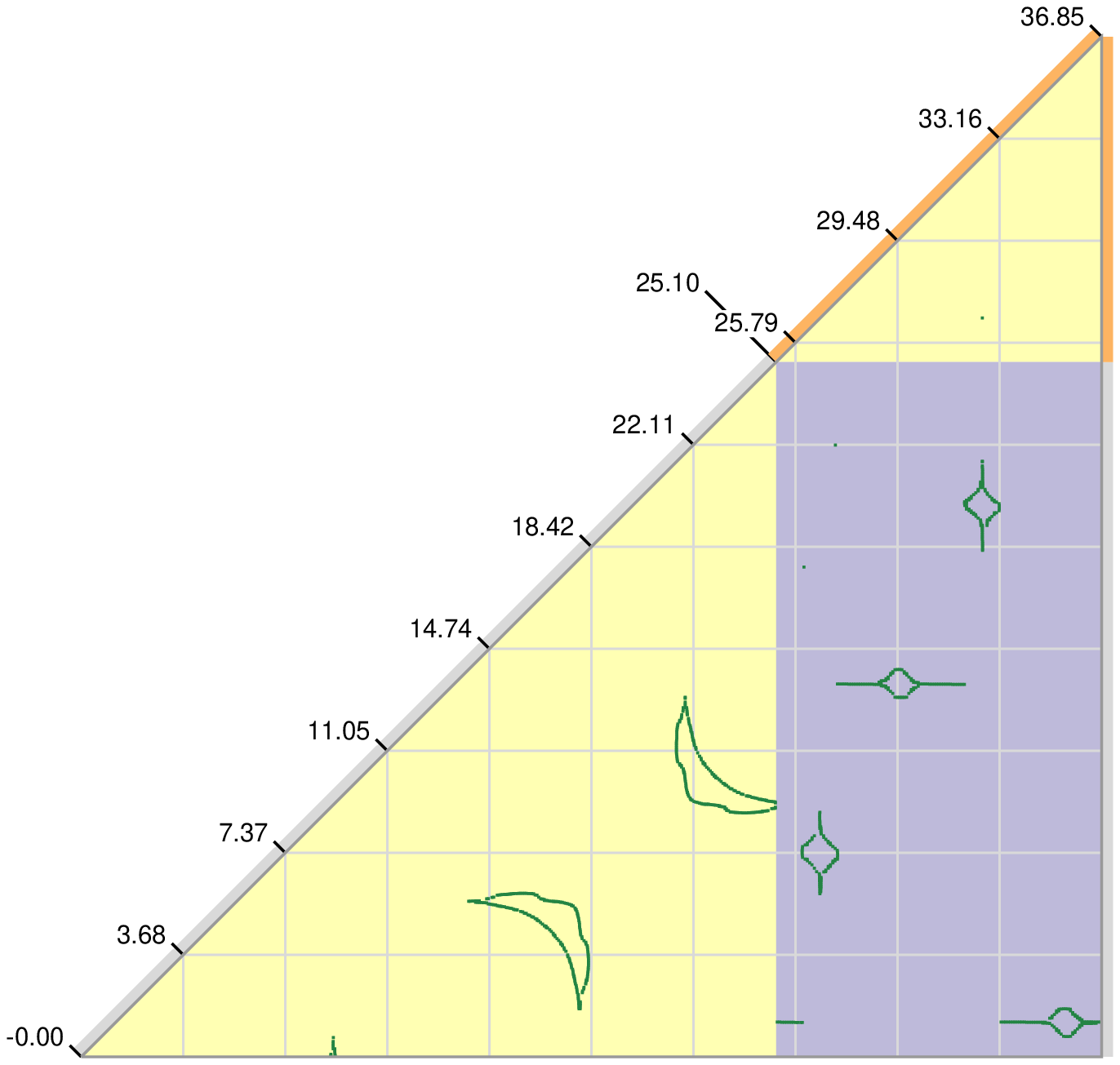}
        \put(8,94){\scriptsize{$73.70$}}
        \put(8,89){\scriptsize{$73.68$}}
        \put(8,84){\scriptsize{$526$}}
    \end{overpic}
\end{minipage} 
\hfill
\begin{minipage}[t]{6in}
  \vspace{2mm}
    \begin{overpic}[height=2.8in]{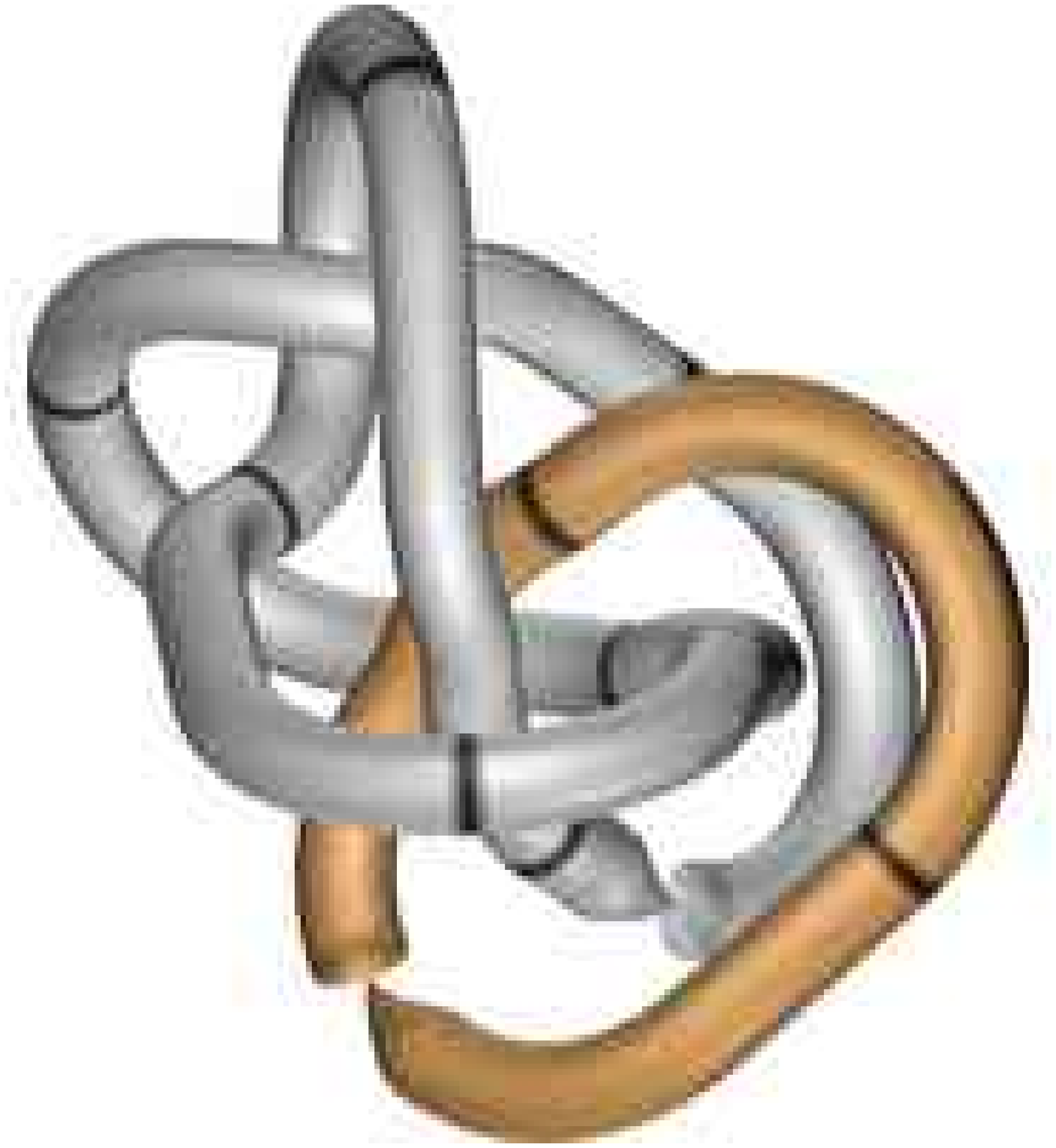}
        \put(-10,90){\large{$8^{2}_{15}$}}
    \end{overpic}
      \hspace{7mm}
    \begin{overpic}[width=2.8in]{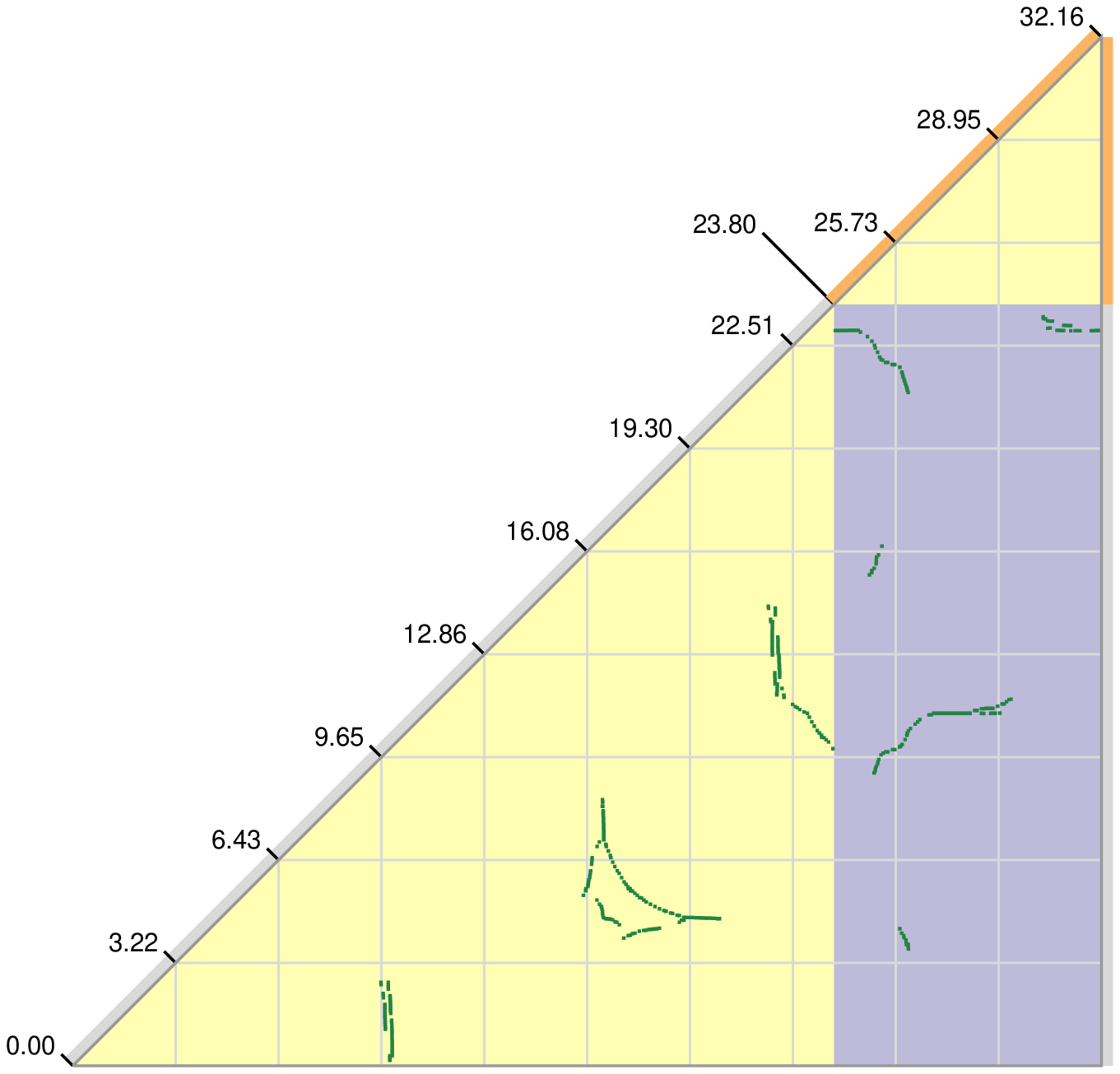}
        \put(8,94){\scriptsize{$64.34$}}
        \put(8,89){\scriptsize{$64.31$}}
        \put(8,84){\scriptsize{$459$}}
    \end{overpic}
\end{minipage} 
\hfill
\end{figure}
\clearpage
\pagebreak
\begin{figure}
\begin{minipage}[t]{6in}
  \vspace{2mm}
    \begin{overpic}[height=2.8in]{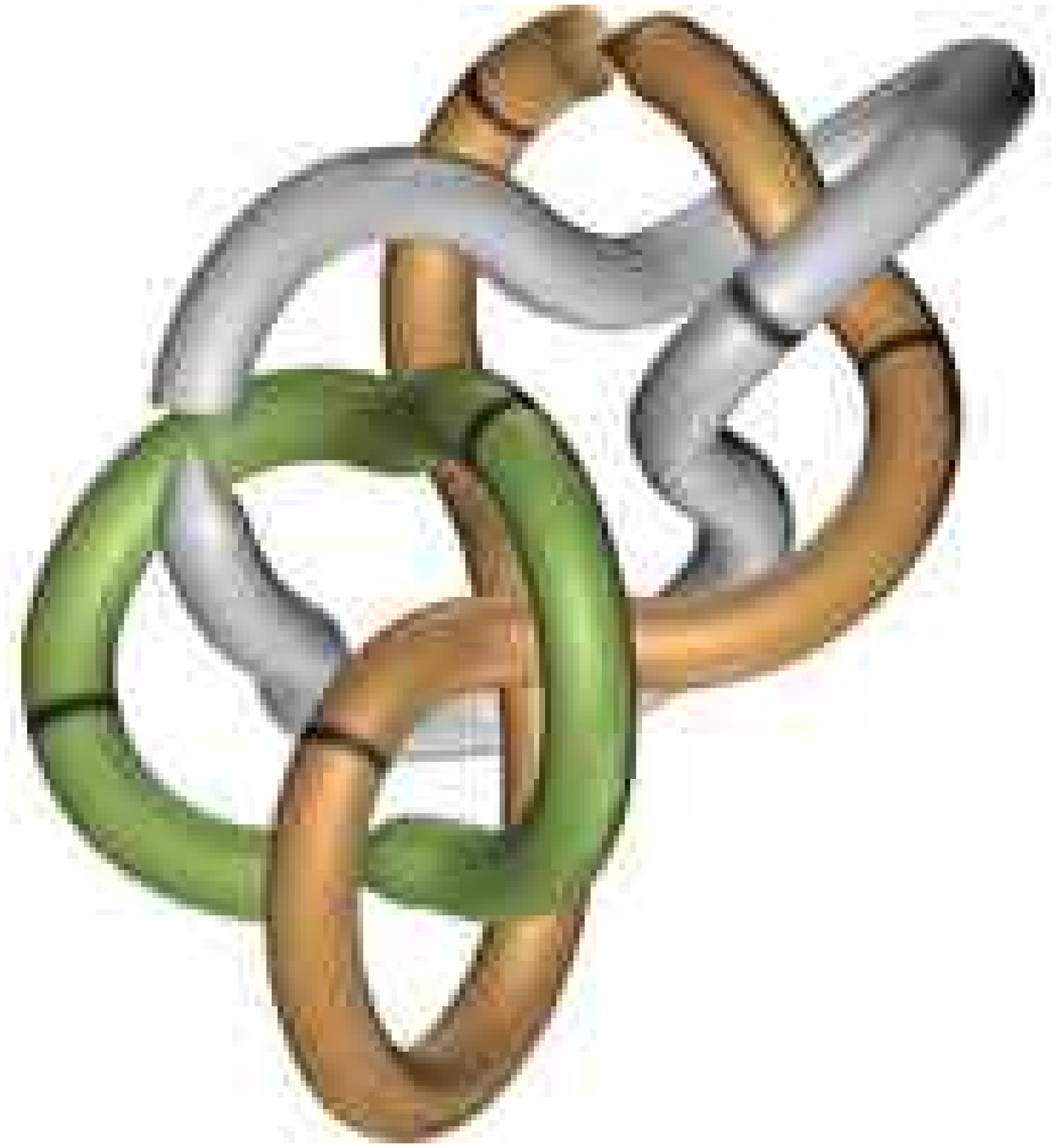}
        \put(-10,90){\large{$8^{3}_{2}$}}
    \end{overpic}
      \hspace{7mm}
    \begin{overpic}[width=2.8in]{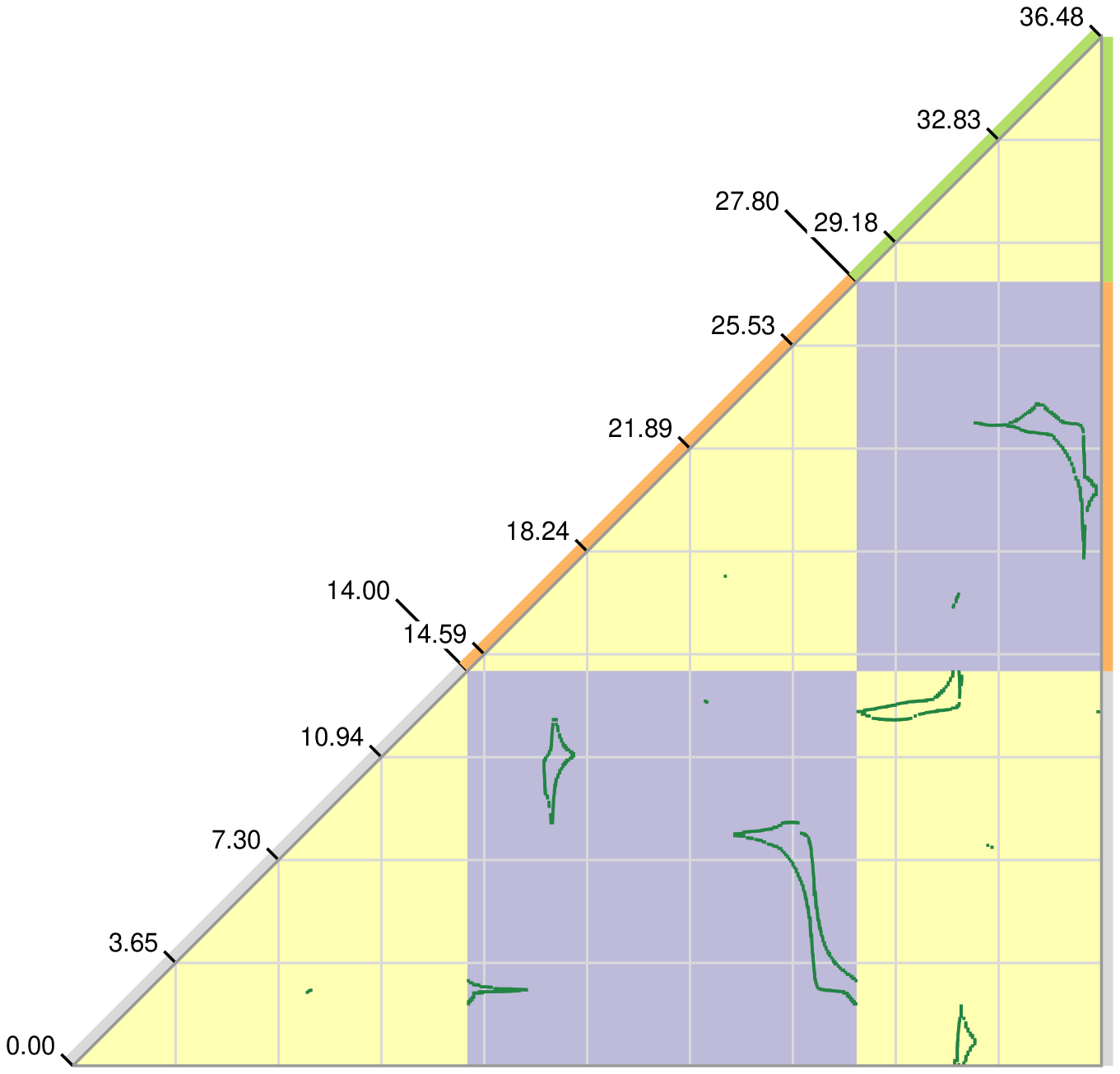}
        \put(8,94){\scriptsize{$72.96$}}
        \put(8,89){\scriptsize{$72.94$}}
        \put(8,84){\scriptsize{$521$}}
    \end{overpic}
\end{minipage} 
\hfill
\begin{minipage}[t]{6in}
  \vspace{2mm}
    \begin{overpic}[height=2.8in]{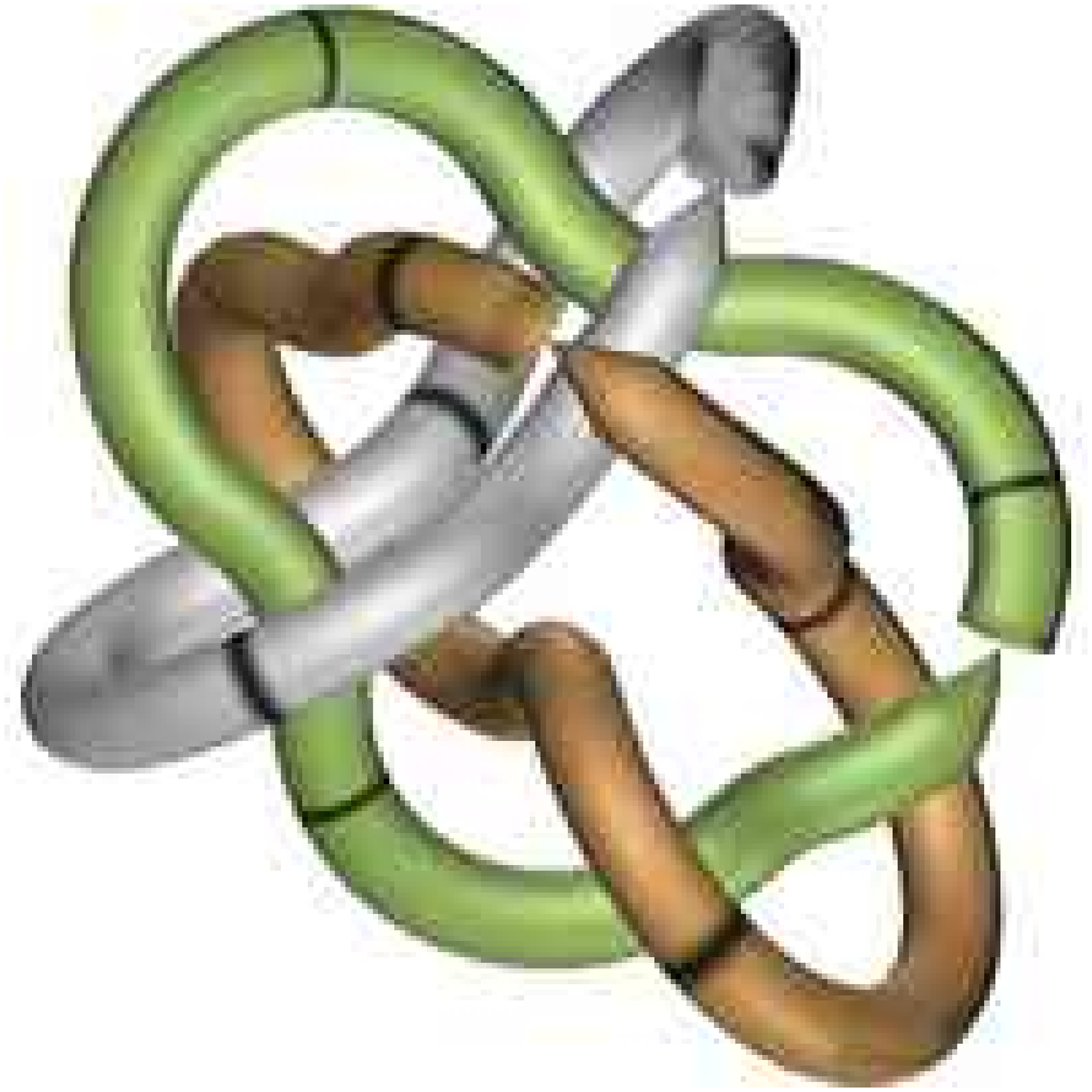}
        \put(-10,90){\large{$8^{3}_{5}$}}
    \end{overpic}
      \hspace{7mm}
    \begin{overpic}[width=2.8in]{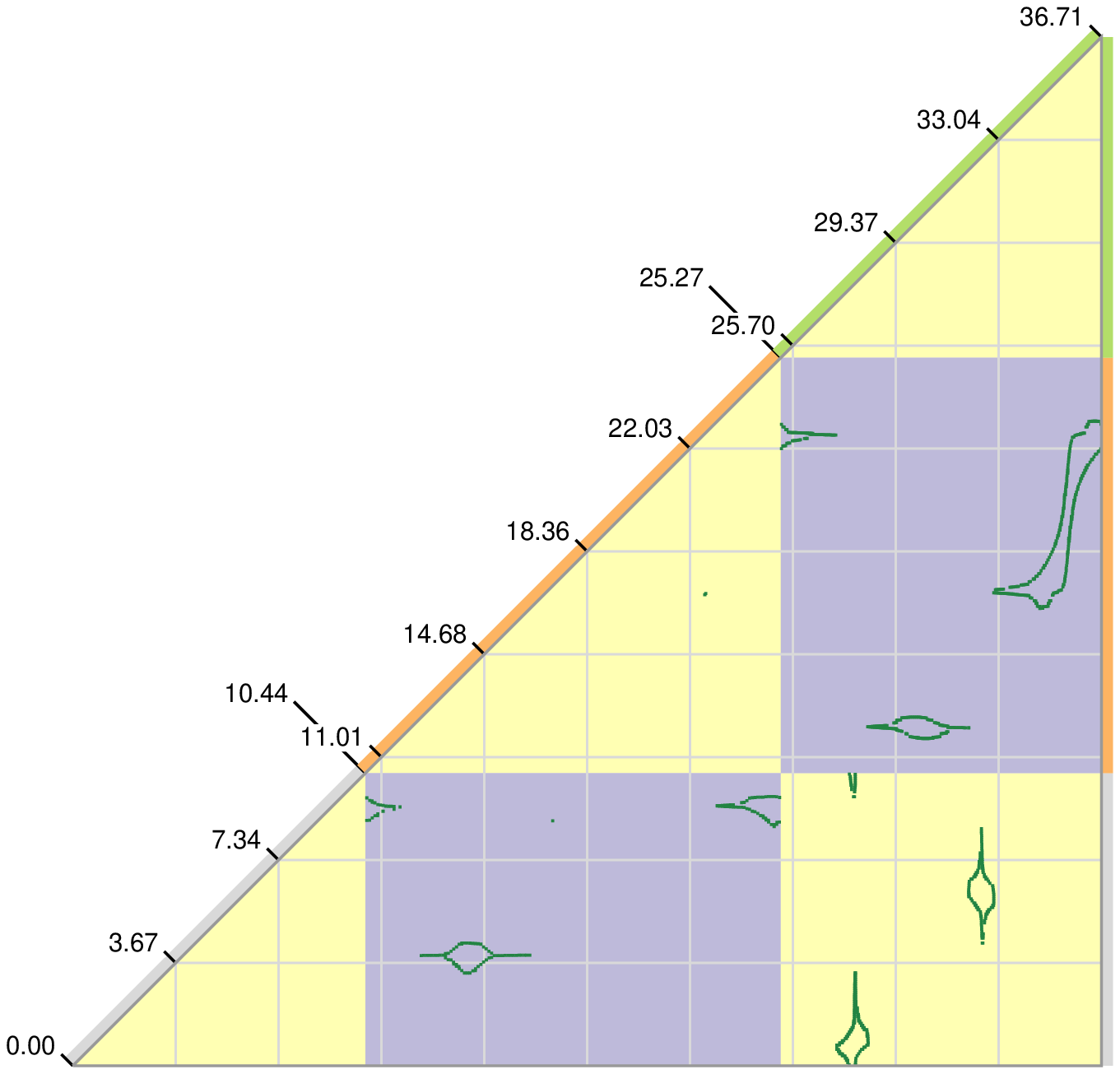}
        \put(8,94){\scriptsize{$73.43$}}
        \put(8,89){\scriptsize{$73.41$}}
        \put(8,84){\scriptsize{$524$}}
    \end{overpic}
\end{minipage} 
\hfill
\begin{minipage}[t]{6in}
  \vspace{2mm}
    \begin{overpic}[width=2.8in]{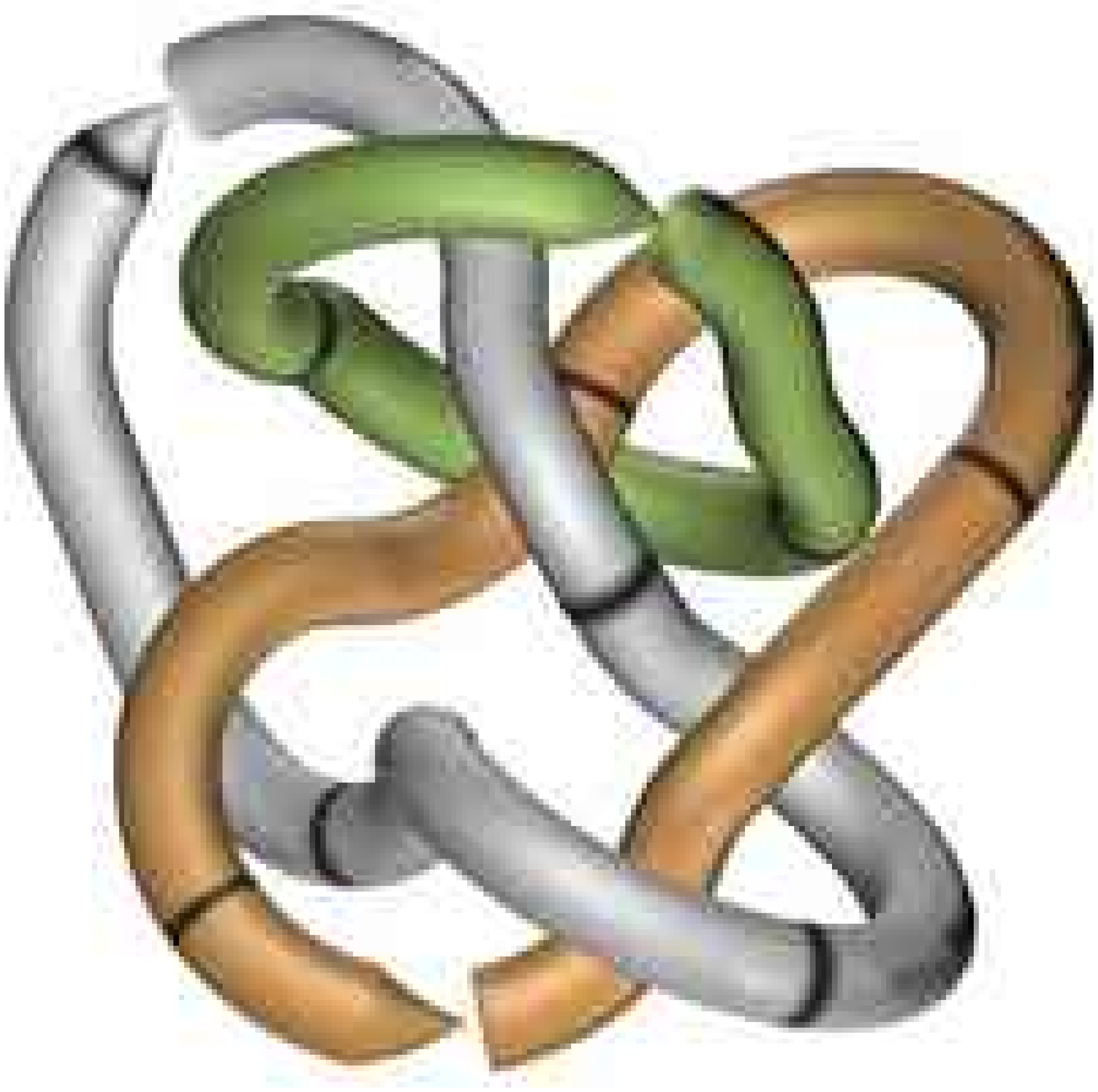}
        \put(-10,90){\large{$8^{3}_{8}$}}
    \end{overpic}
      \hspace{7mm}
    \begin{overpic}[width=2.8in]{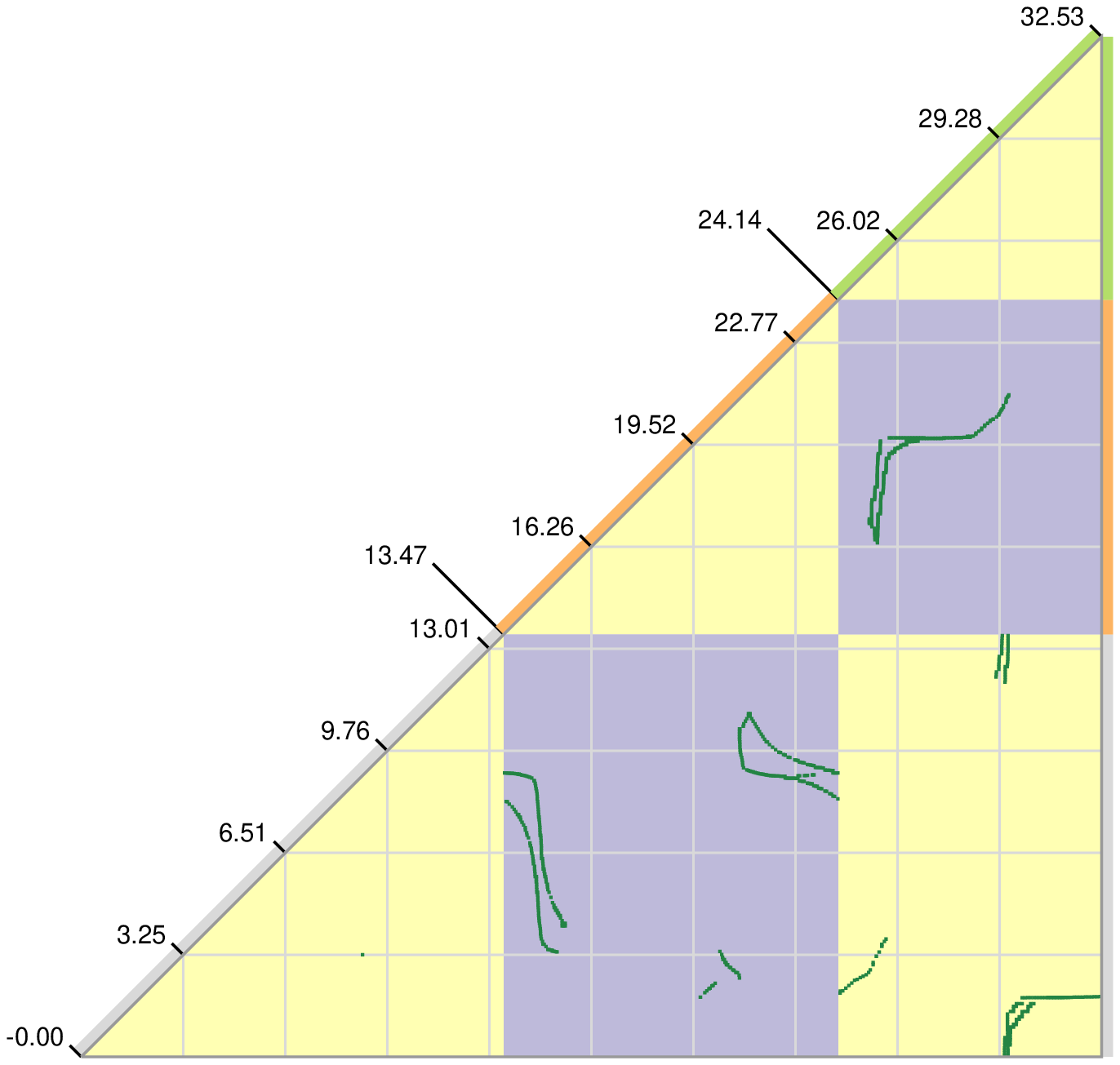}
        \put(8,94){\scriptsize{$65.07$}}
        \put(8,89){\scriptsize{$65.05$}}
        \put(8,84){\scriptsize{$464$}}
    \end{overpic}
\end{minipage} 
\hfill
\end{figure}
\clearpage
\pagebreak
\begin{figure}
\begin{minipage}[t]{6in}
  \vspace{2mm}
    \begin{overpic}[width=2.8in]{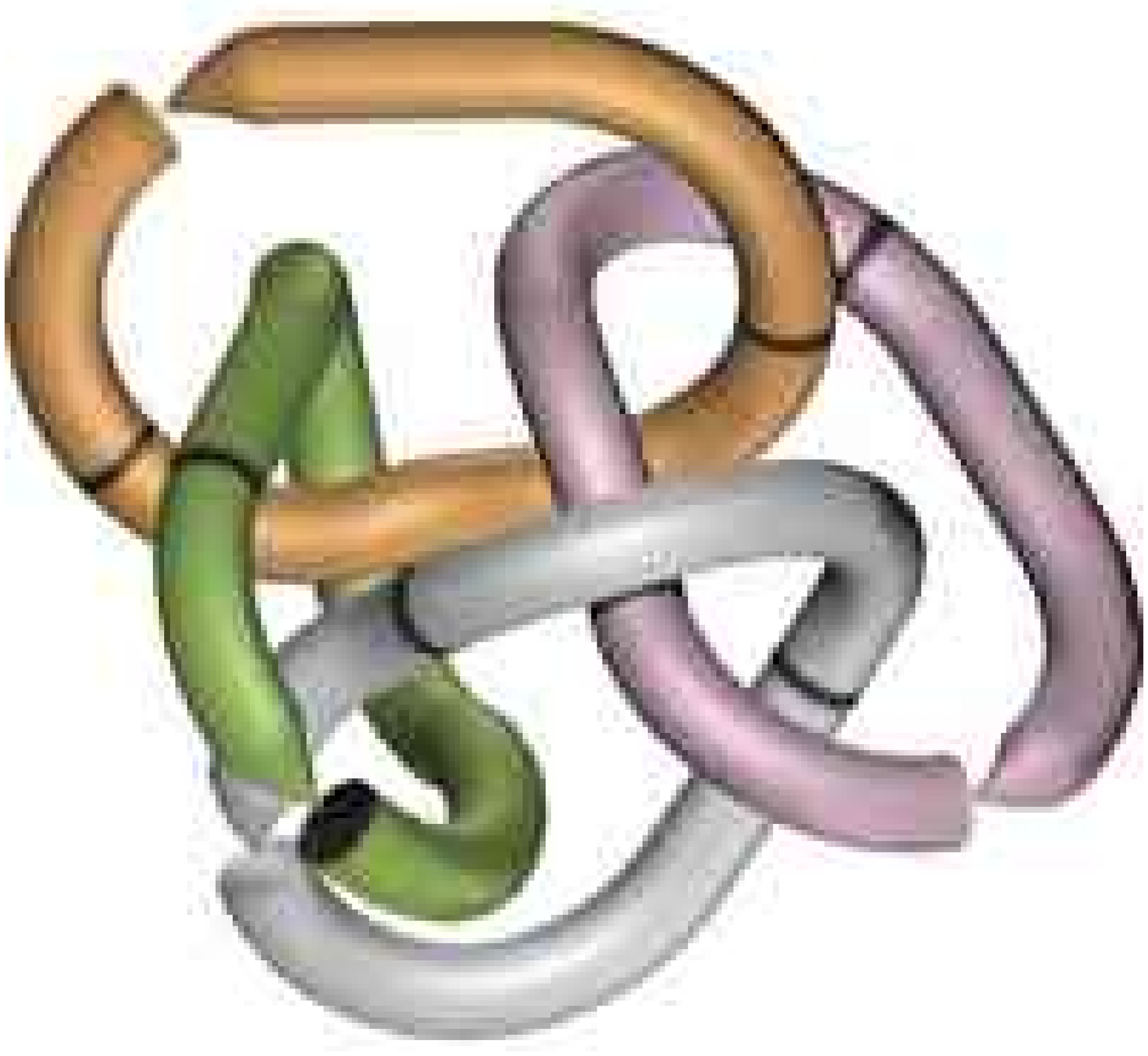}
        \put(-10,90){\large{$8^{4}_{2}$}}
    \end{overpic}
      \hspace{7mm}
    \begin{overpic}[width=2.8in]{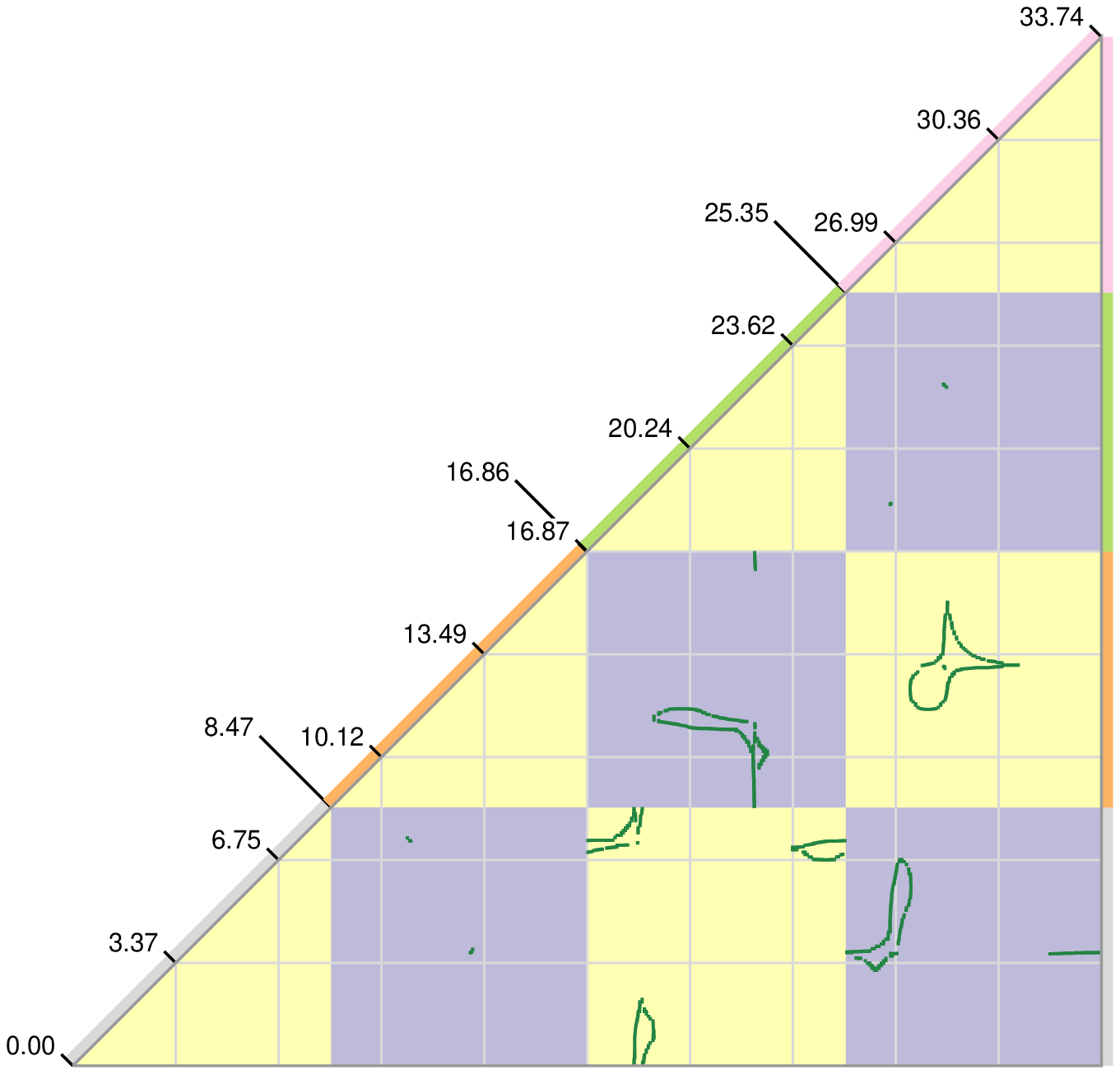}
        \put(8,94){\scriptsize{$67.48$}}
        \put(8,89){\scriptsize{$67.46$}}
        \put(8,84){\scriptsize{$482$}}
    \end{overpic}
\end{minipage} 
\hfill
\begin{minipage}[t]{6in}
  \vspace{2mm}
    \begin{overpic}[width=2.8in]{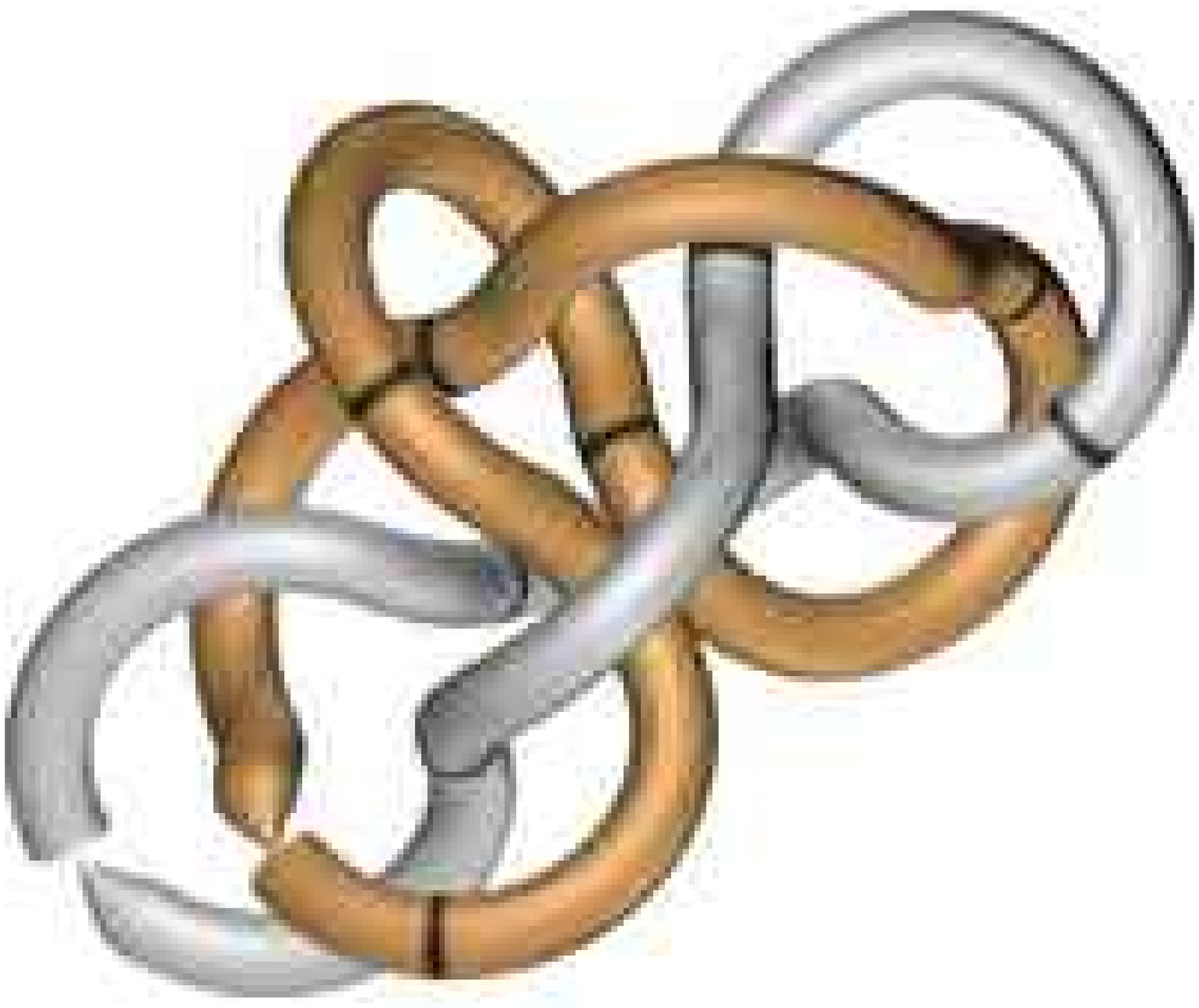}
        \put(-10,90){\large{$9^{2}_{4}$}}
    \end{overpic}
      \hspace{7mm}
    \begin{overpic}[width=2.8in]{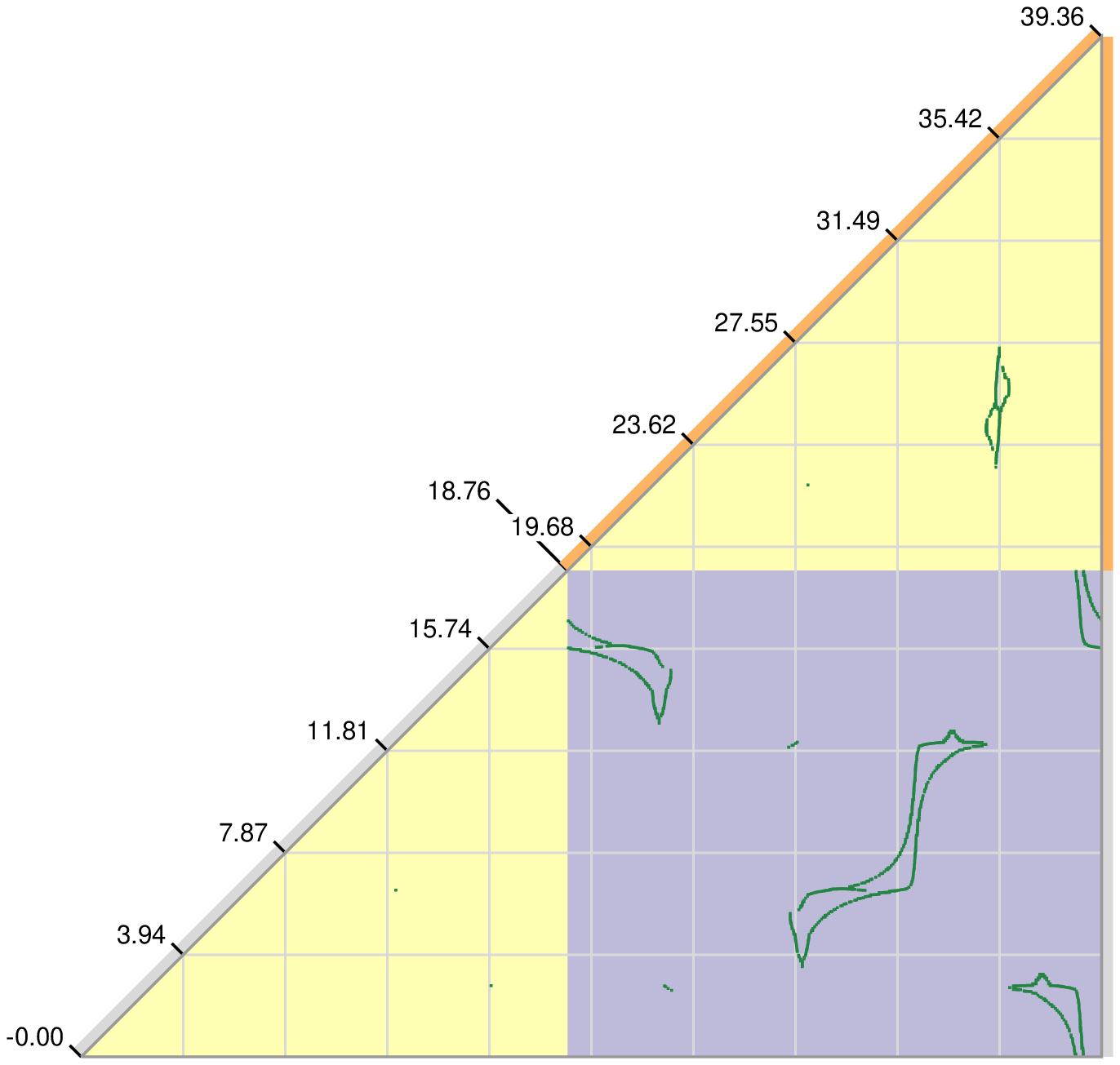}
        \put(8,94){\scriptsize{$78.73$}}
        \put(8,89){\scriptsize{$78.70$}}
        \put(8,84){\scriptsize{$562$}}
    \end{overpic}
\end{minipage} 
\hfill
\begin{minipage}[t]{6in}
  \vspace{2mm}
    \begin{overpic}[height=2.8in]{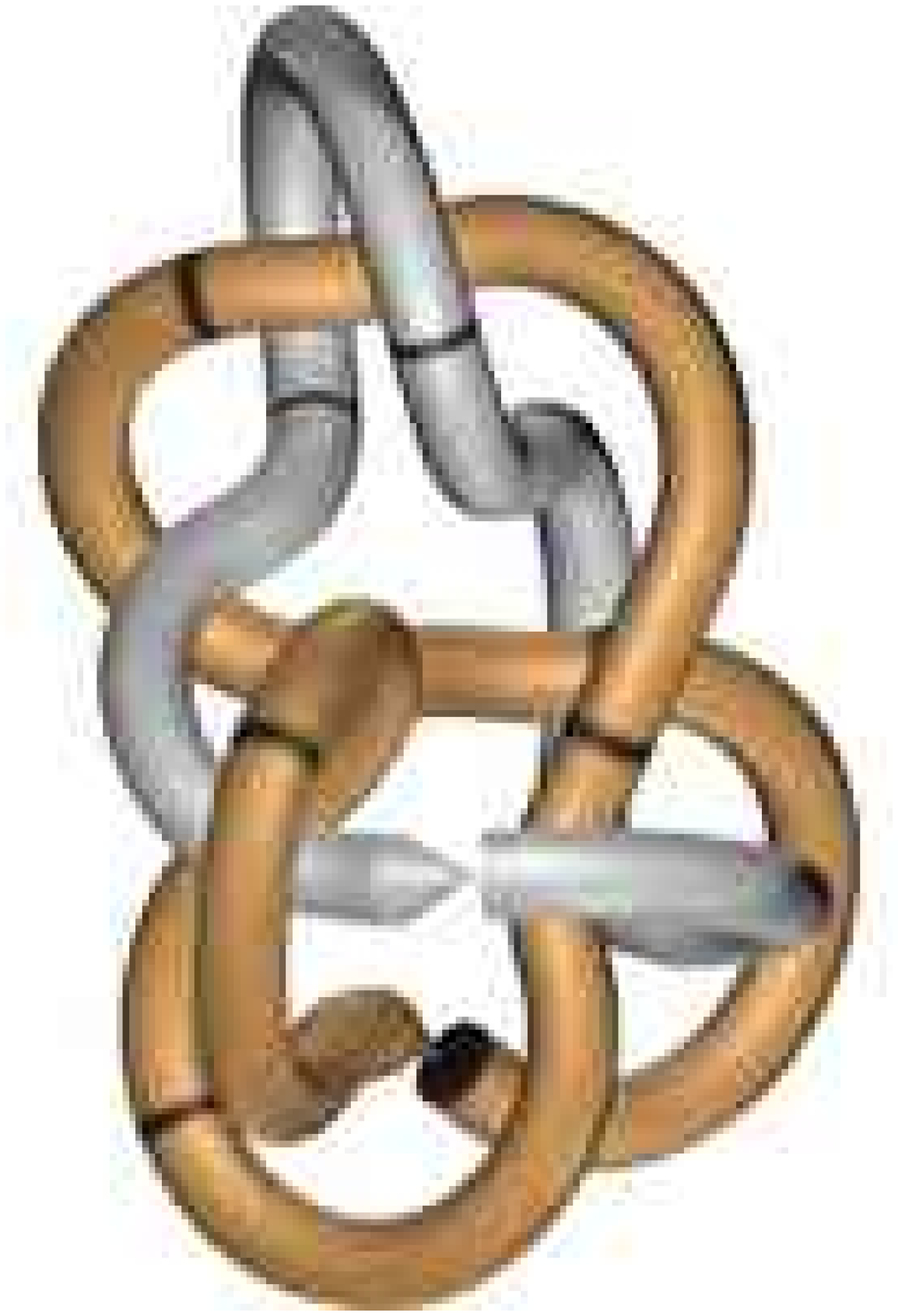}
        \put(-10,90){\large{$9^{2}_{7}$}}
    \end{overpic}
      \hspace{7mm}
    \begin{overpic}[width=2.8in]{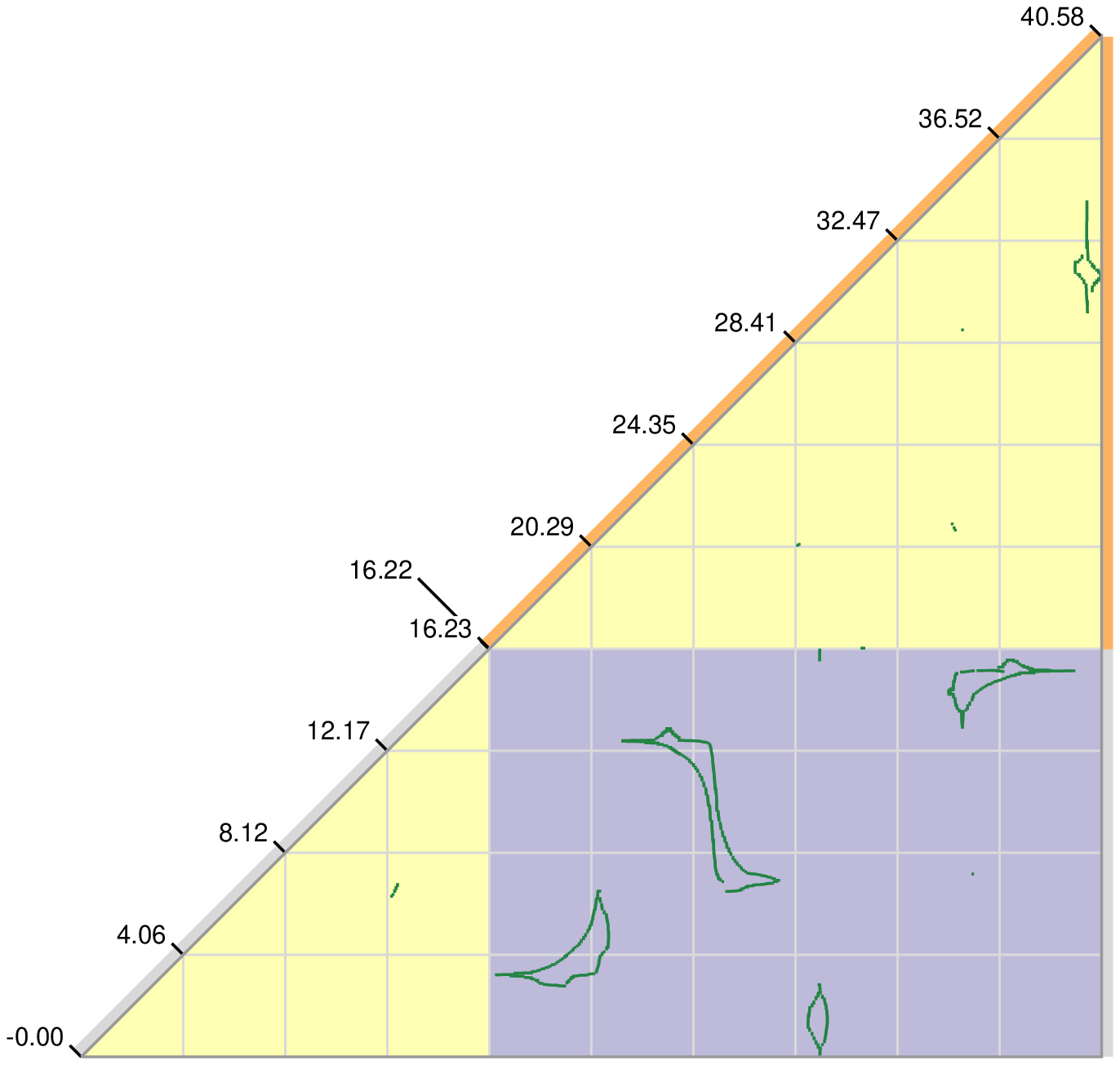}
        \put(8,94){\scriptsize{$81.18$}}
        \put(8,89){\scriptsize{$81.15$}}
        \put(8,84){\scriptsize{$580$}}
    \end{overpic}
\end{minipage} 
\hfill
\end{figure}
\clearpage
\pagebreak
\begin{figure}
\begin{minipage}[t]{6in}
  \vspace{2mm}
    \begin{overpic}[height=2.8in]{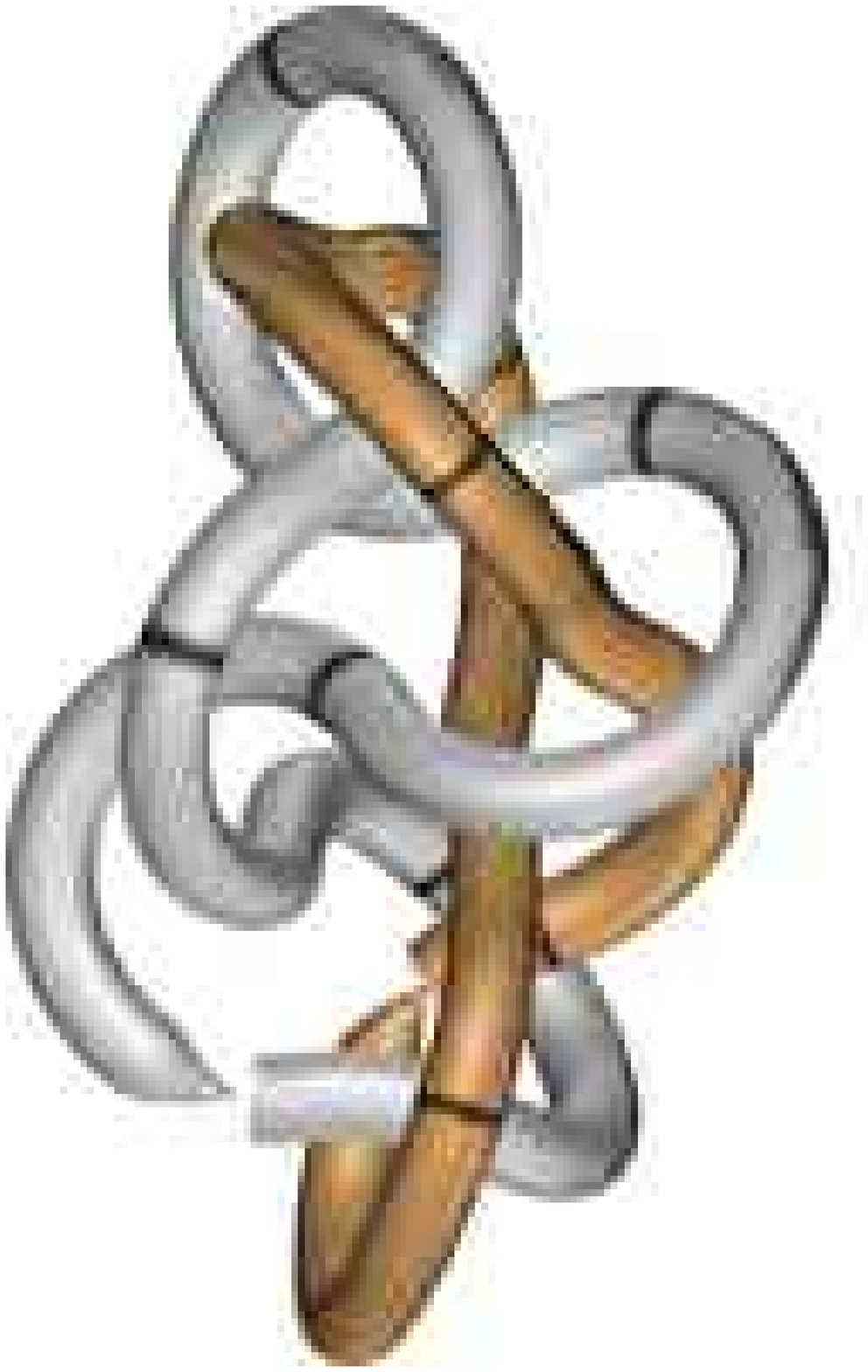}
        \put(-10,90){\large{$9^{2}_{8}$}}
    \end{overpic}
      \hspace{7mm}
    \begin{overpic}[width=2.8in]{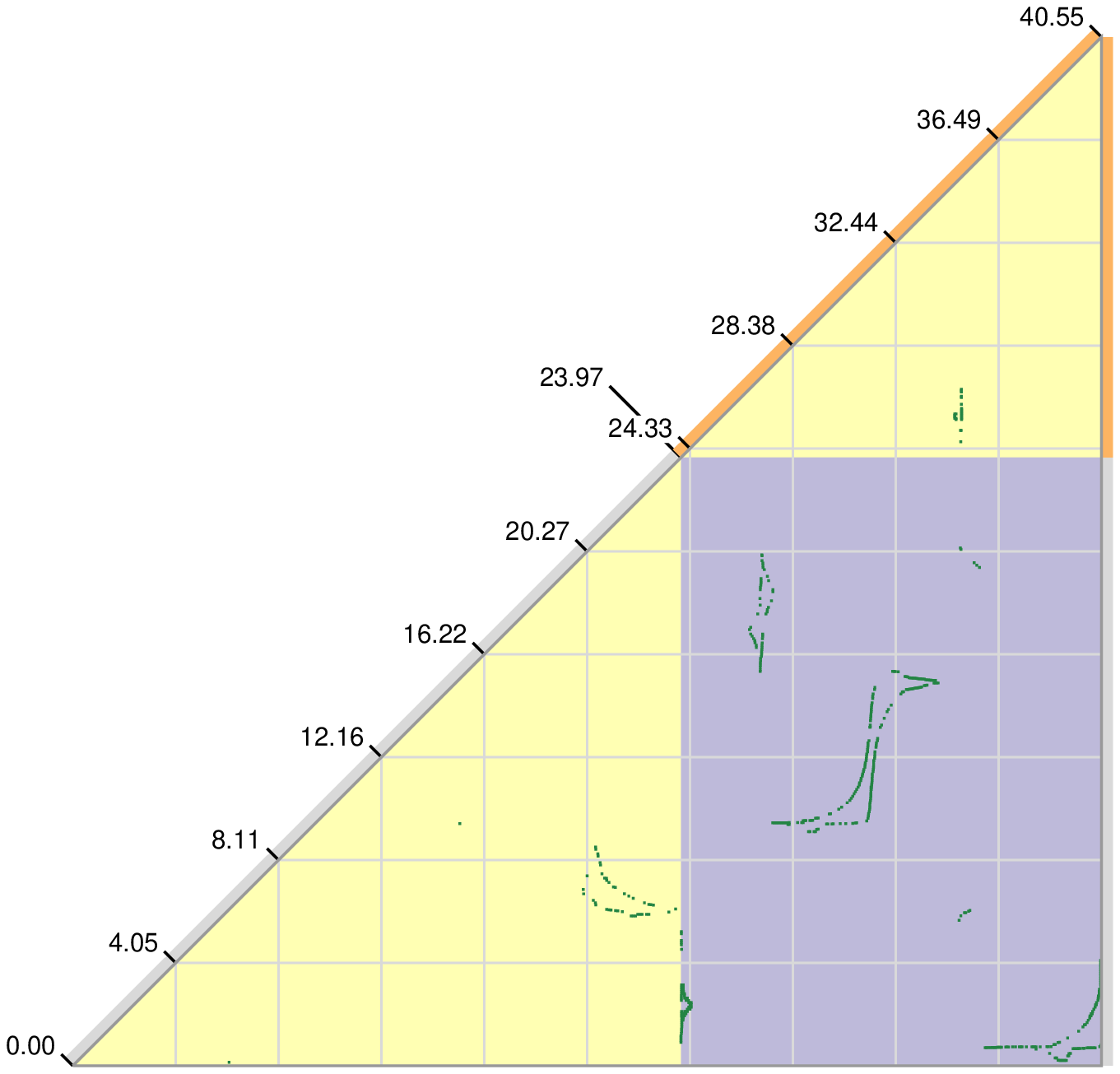}
        \put(8,94){\scriptsize{$81.10$}}
        \put(8,89){\scriptsize{$81.08$}}
        \put(8,84){\scriptsize{$579$}}
    \end{overpic}
\end{minipage} 
\hfill
\begin{minipage}[t]{6in}
  \vspace{2mm}
    \begin{overpic}[width=2.8in]{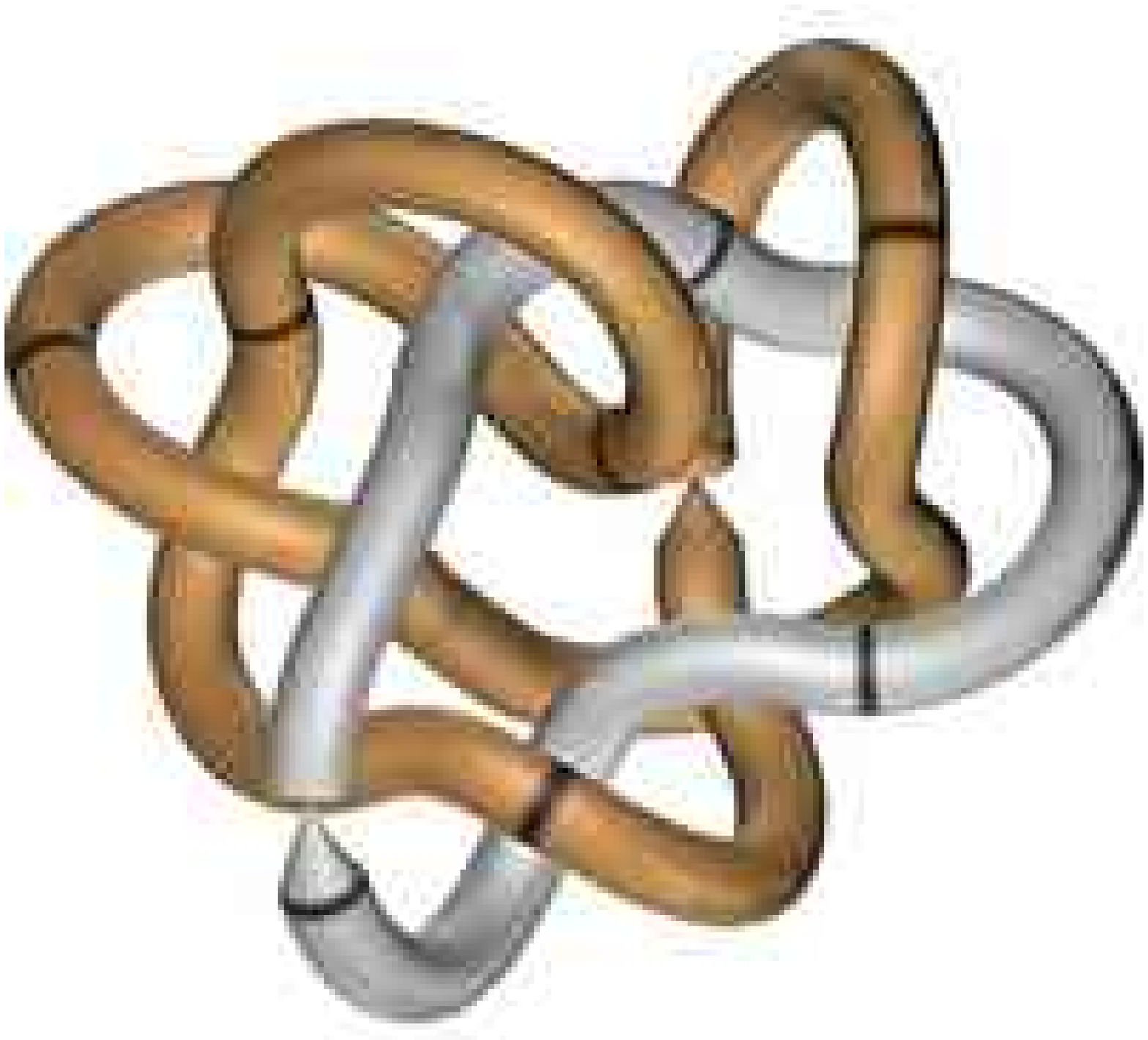}
        \put(-10,90){\large{$9^{2}_{11}$}}
    \end{overpic}
      \hspace{7mm}
    \begin{overpic}[width=2.8in]{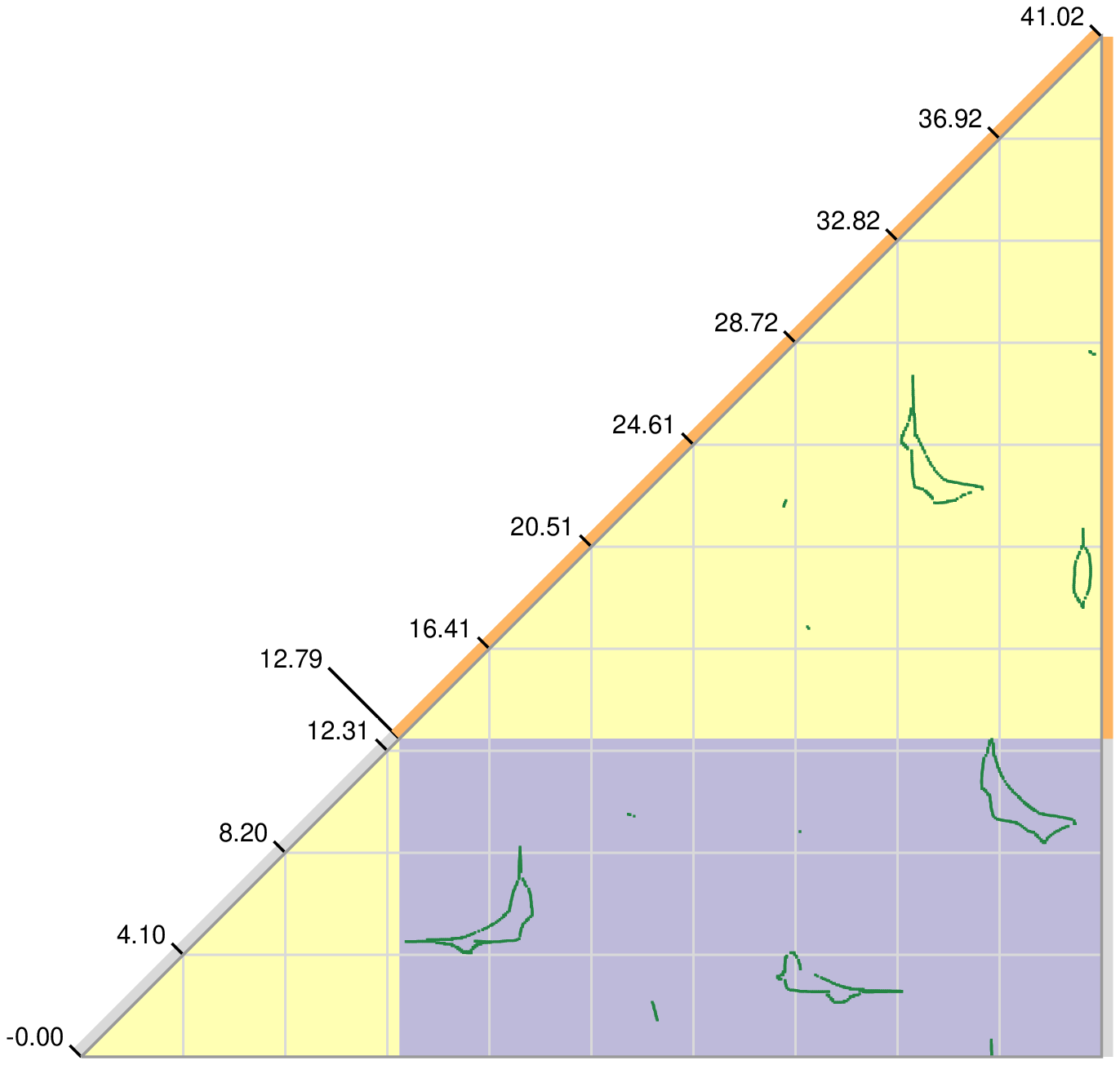}
        \put(8,94){\scriptsize{$82.06$}}
        \put(8,89){\scriptsize{$82.03$}}
        \put(8,84){\scriptsize{$586$}}
    \end{overpic}
\end{minipage} 
\hfill
\begin{minipage}[t]{6in}
  \vspace{2mm}
    \begin{overpic}[width=2.8in]{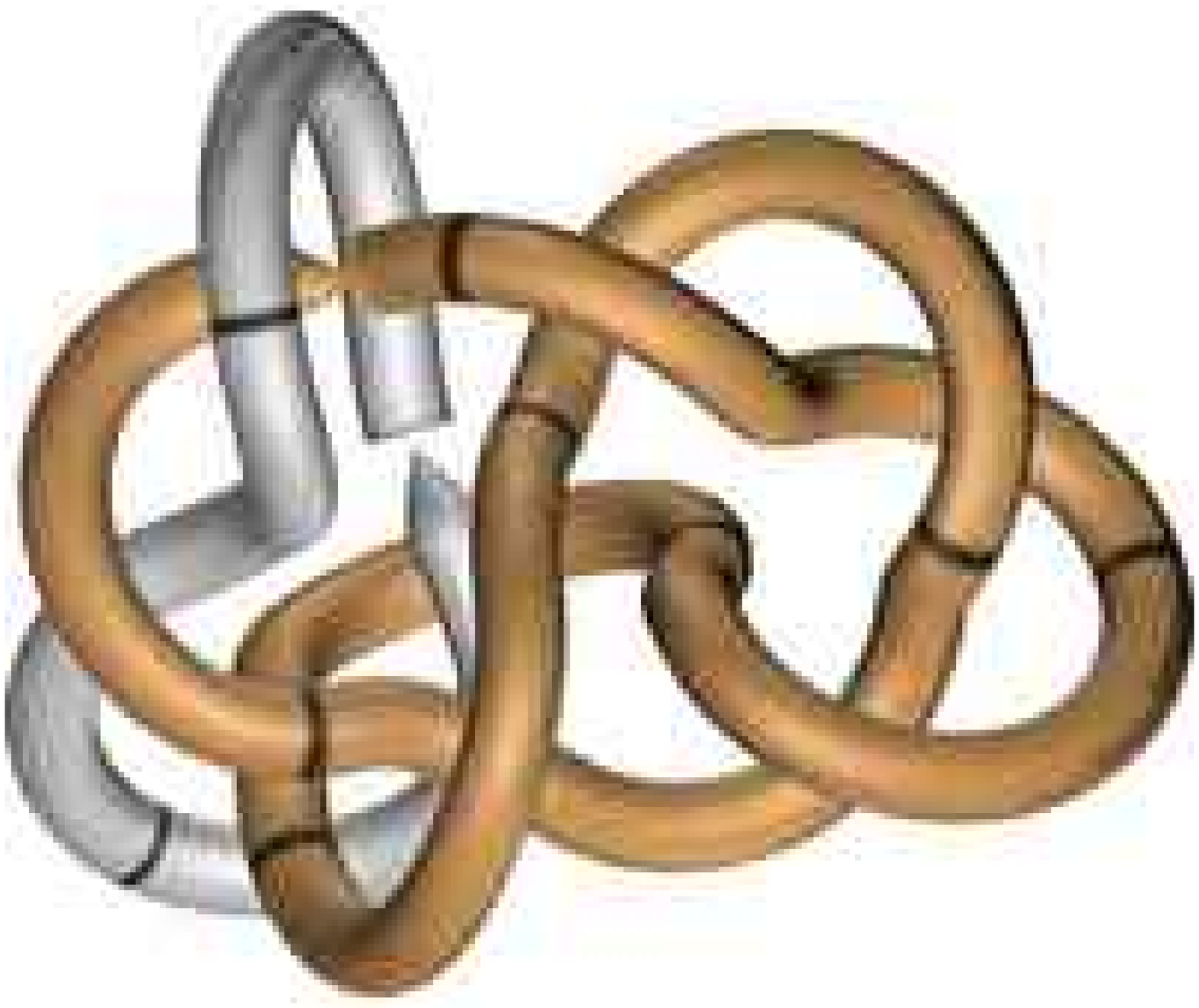}
        \put(-10,90){\large{$9^{2}_{14}$}}
    \end{overpic}
      \hspace{7mm}
    \begin{overpic}[width=2.8in]{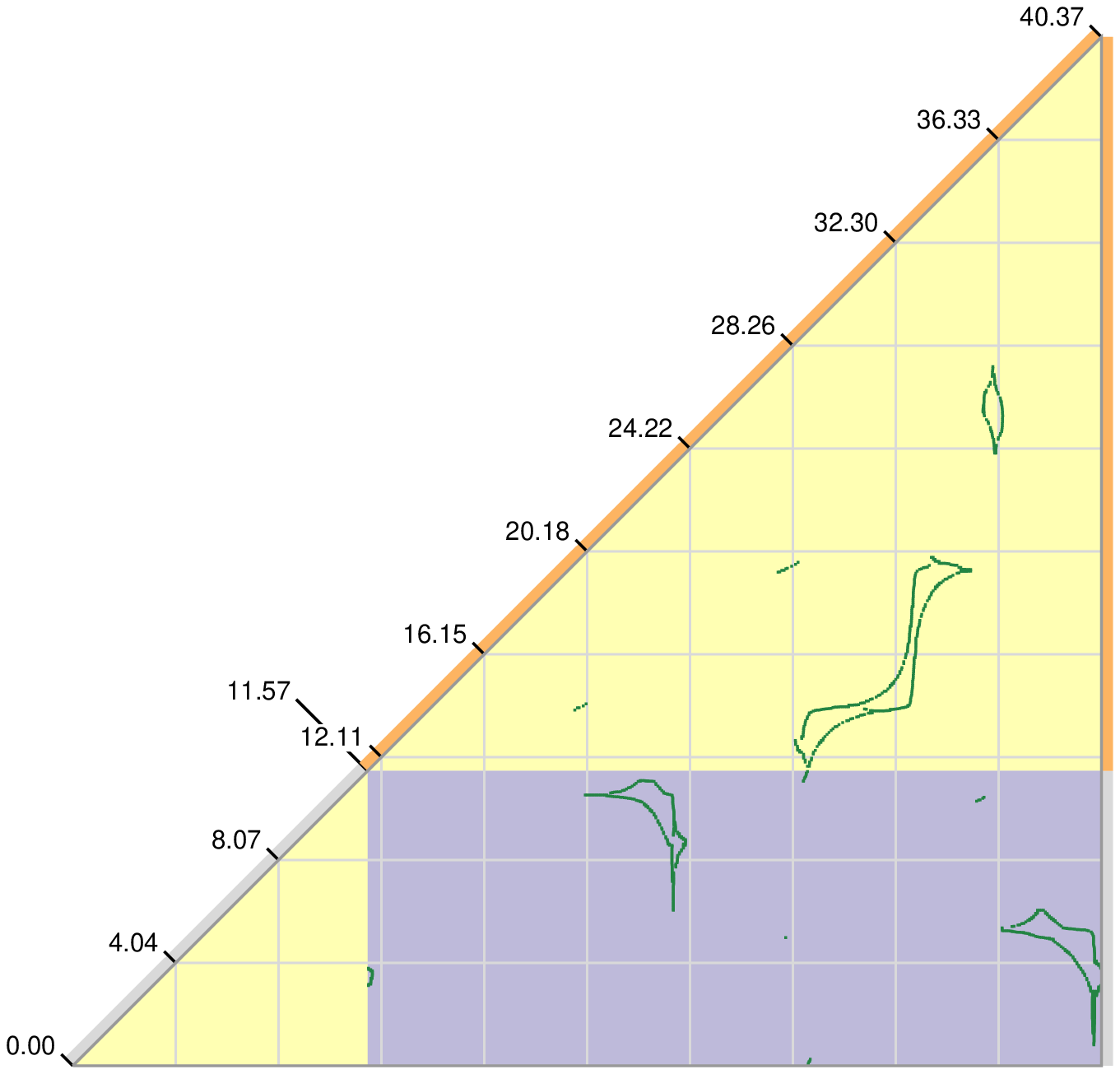}
        \put(8,94){\scriptsize{$80.75$}}
        \put(8,89){\scriptsize{$80.72$}}
        \put(8,84){\scriptsize{$577$}}
    \end{overpic}
\end{minipage} 
\hfill
\end{figure}
\clearpage
\pagebreak
\begin{figure}
\begin{minipage}[t]{6in}
  \vspace{2mm}
    \begin{overpic}[height=2.8in]{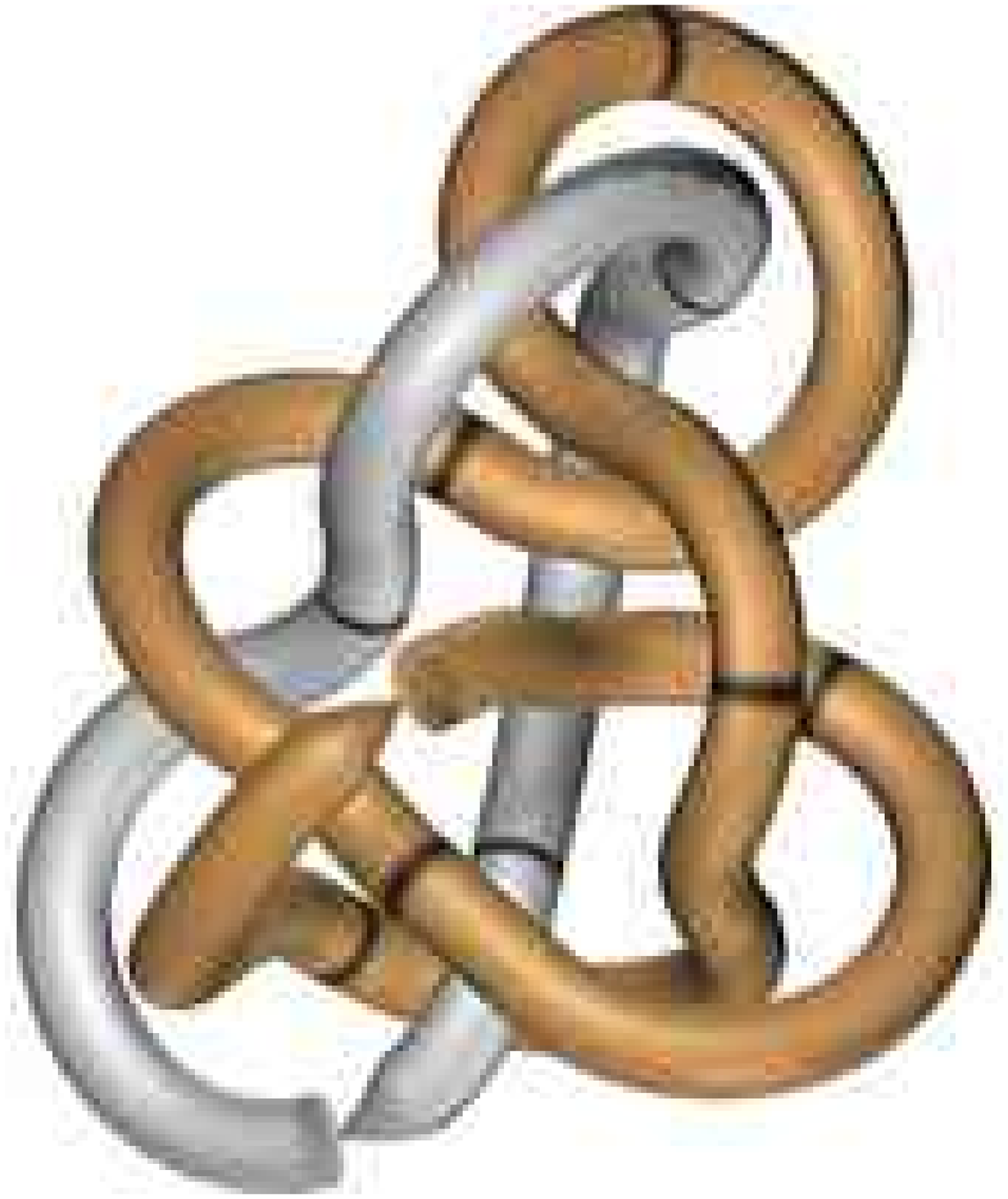}
        \put(-10,90){\large{$9^{2}_{22}$}}
    \end{overpic}
      \hspace{7mm}
    \begin{overpic}[width=2.8in]{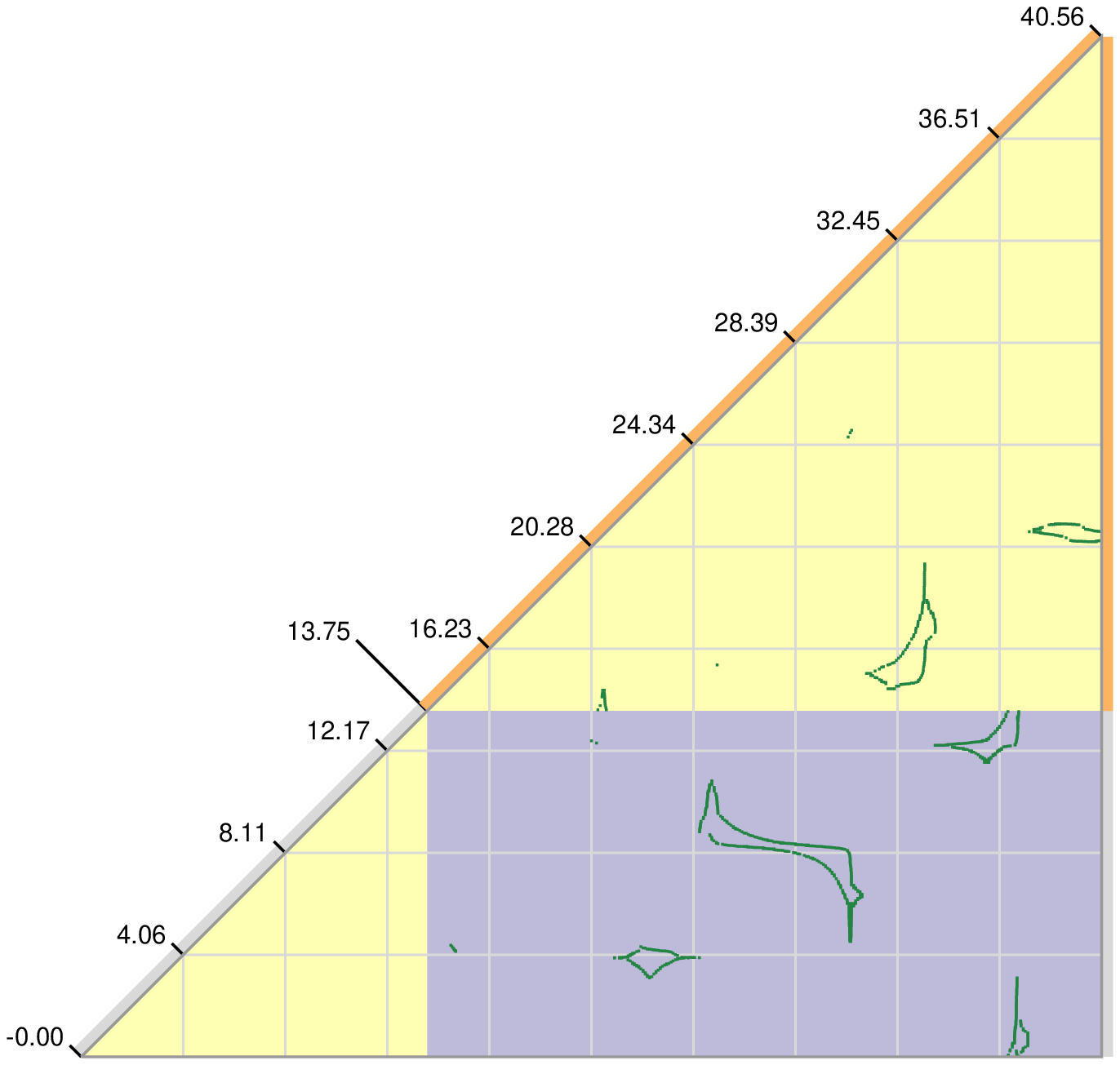}
        \put(8,94){\scriptsize{$81.14$}}
        \put(8,89){\scriptsize{$81.11$}}
        \put(8,84){\scriptsize{$579$}}
    \end{overpic}
\end{minipage} 
\hfill
\begin{minipage}[t]{6in}
  \vspace{2mm}
    \begin{overpic}[height=2.8in]{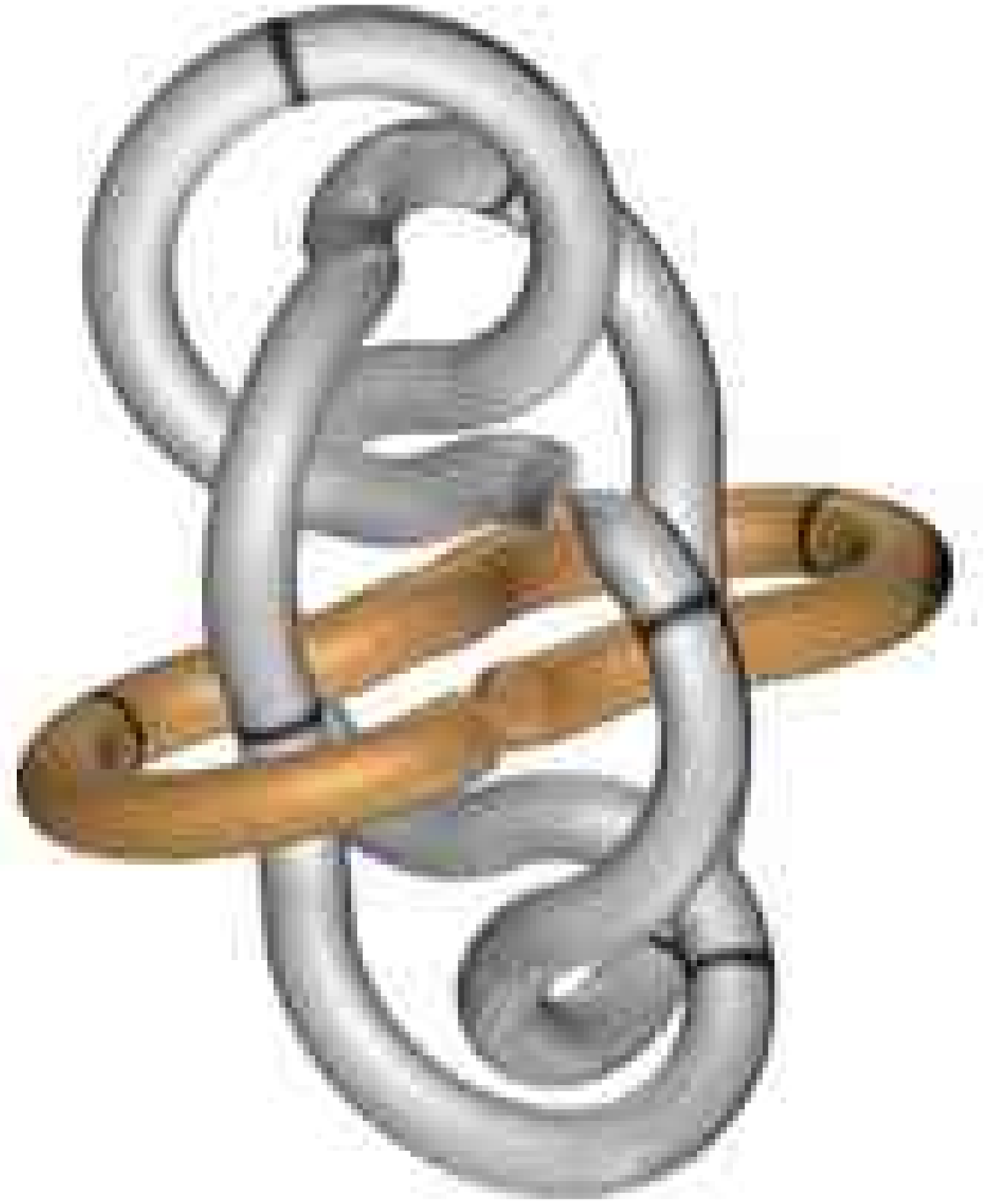}
        \put(-10,90){\large{$9^{2}_{37}$}}
    \end{overpic}
      \hspace{7mm}
    \begin{overpic}[width=2.8in]{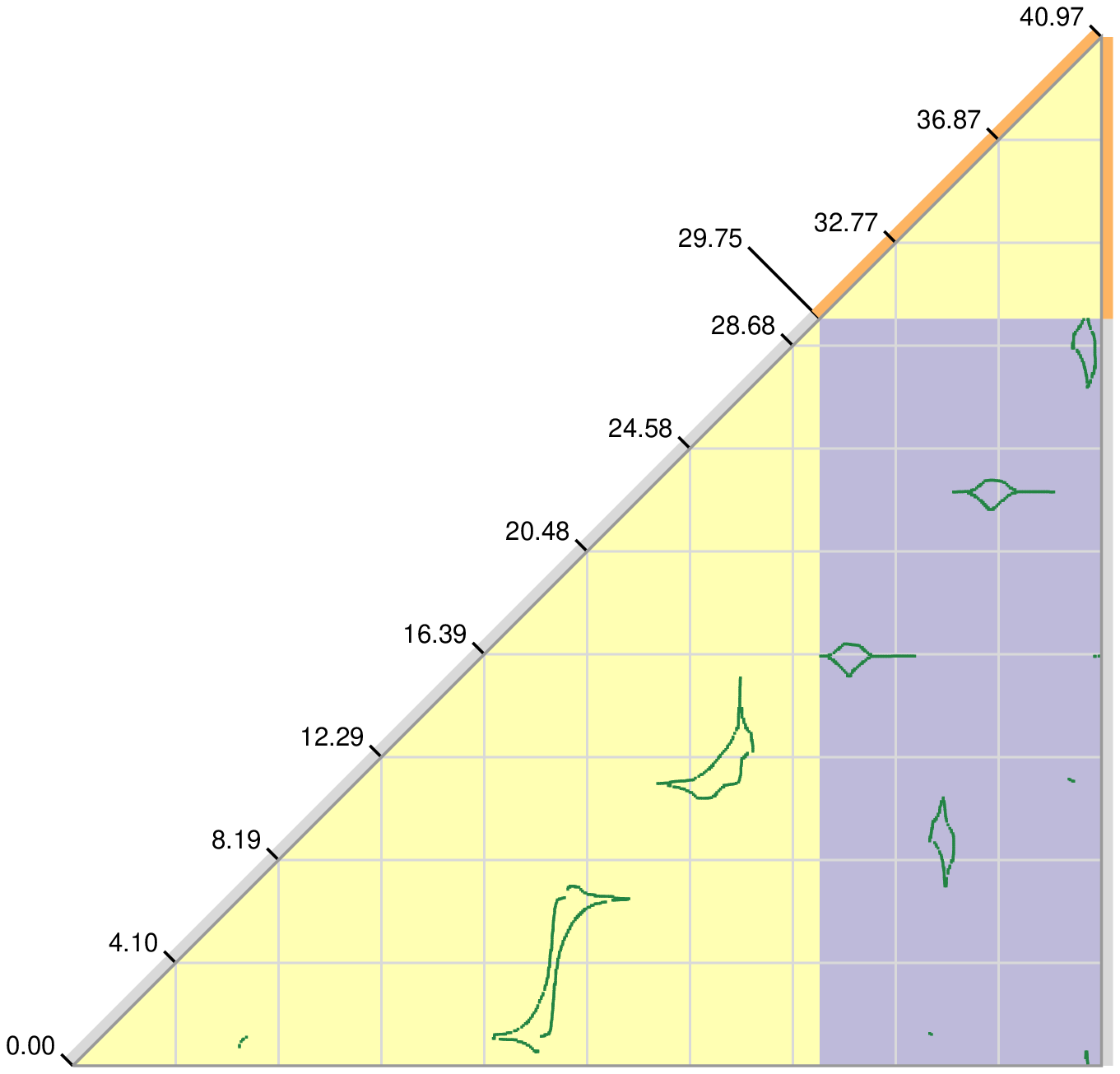}
        \put(8,94){\scriptsize{$81.94$}}
        \put(8,89){\scriptsize{$81.92$}}
        \put(8,84){\scriptsize{$585$}}
    \end{overpic}
\end{minipage} 
\hfill
\begin{minipage}[t]{6in}
  \vspace{2mm}
    \begin{overpic}[width=2.8in]{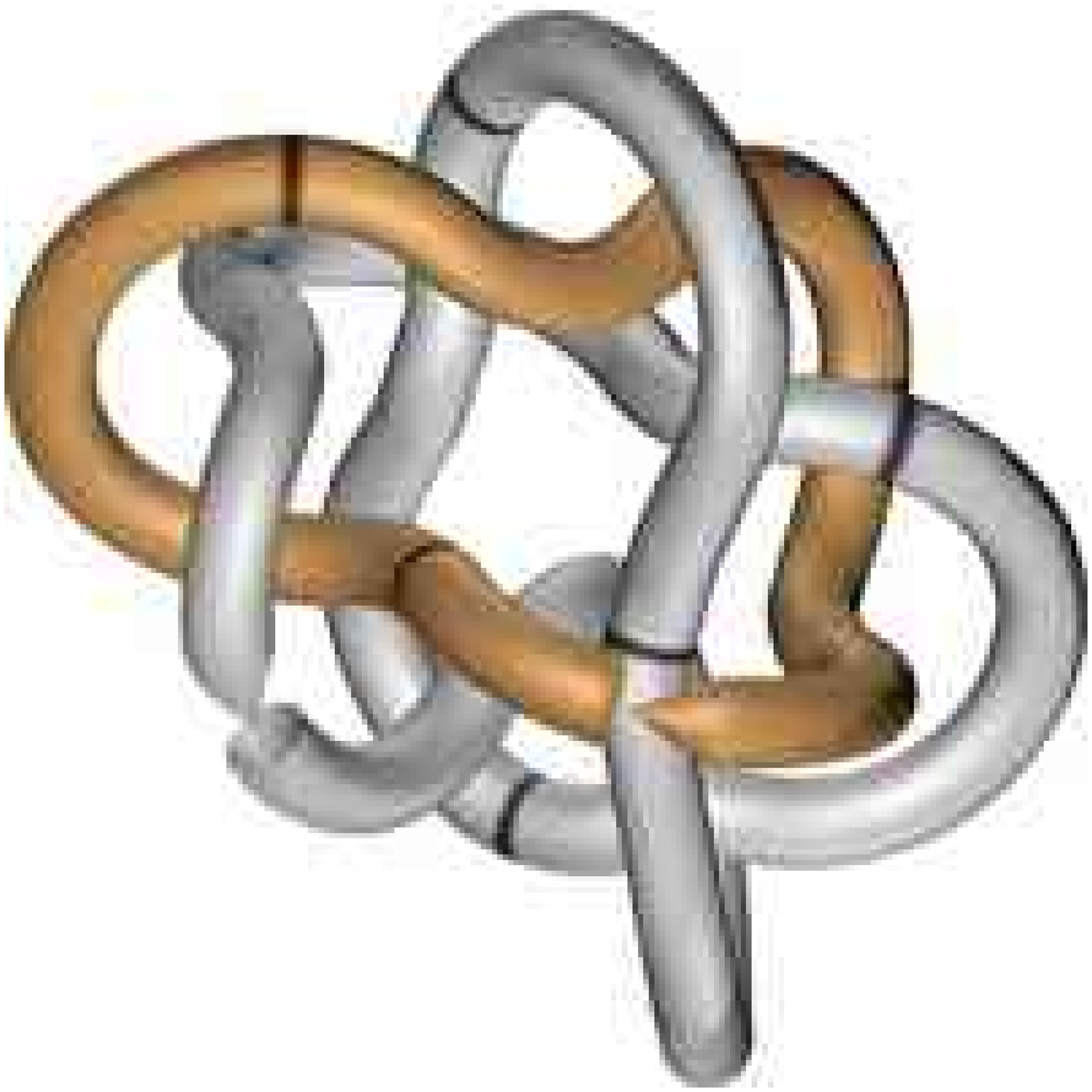}
        \put(-10,90){\large{$9^{2}_{39}$}}
    \end{overpic}
      \hspace{7mm}
    \begin{overpic}[width=2.8in]{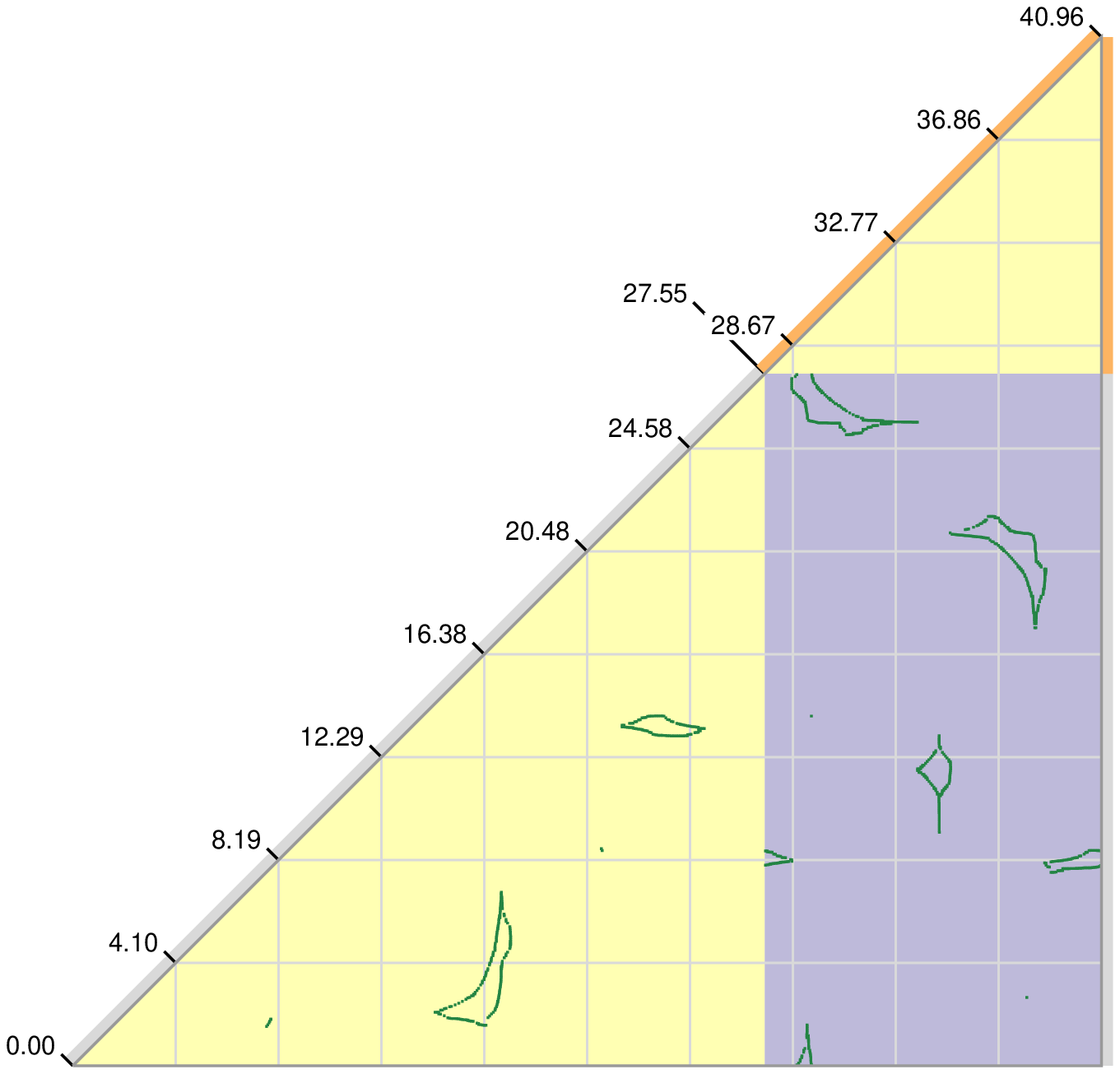}
        \put(8,94){\scriptsize{$81.93$}}
        \put(8,89){\scriptsize{$81.90$}}
        \put(8,84){\scriptsize{$585$}}
    \end{overpic}
\end{minipage} 
\hfill
\end{figure}
\clearpage
\pagebreak
\begin{figure}
\begin{minipage}[t]{6in}
  \vspace{2mm}
    \begin{overpic}[height=2.8in]{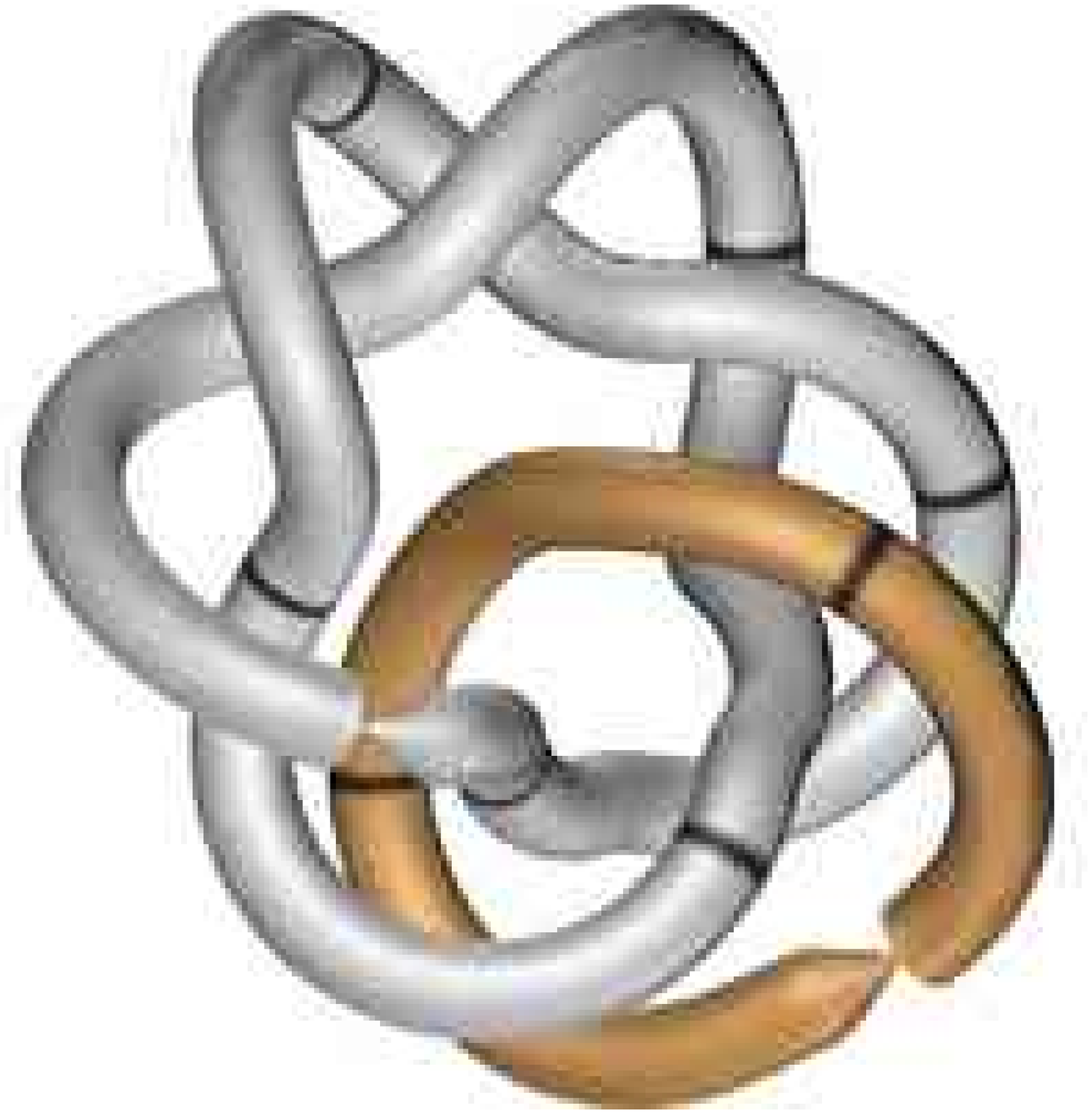}
        \put(-10,90){\large{$9^{2}_{43}$}}
    \end{overpic}
      \hspace{7mm}
    \begin{overpic}[width=2.8in]{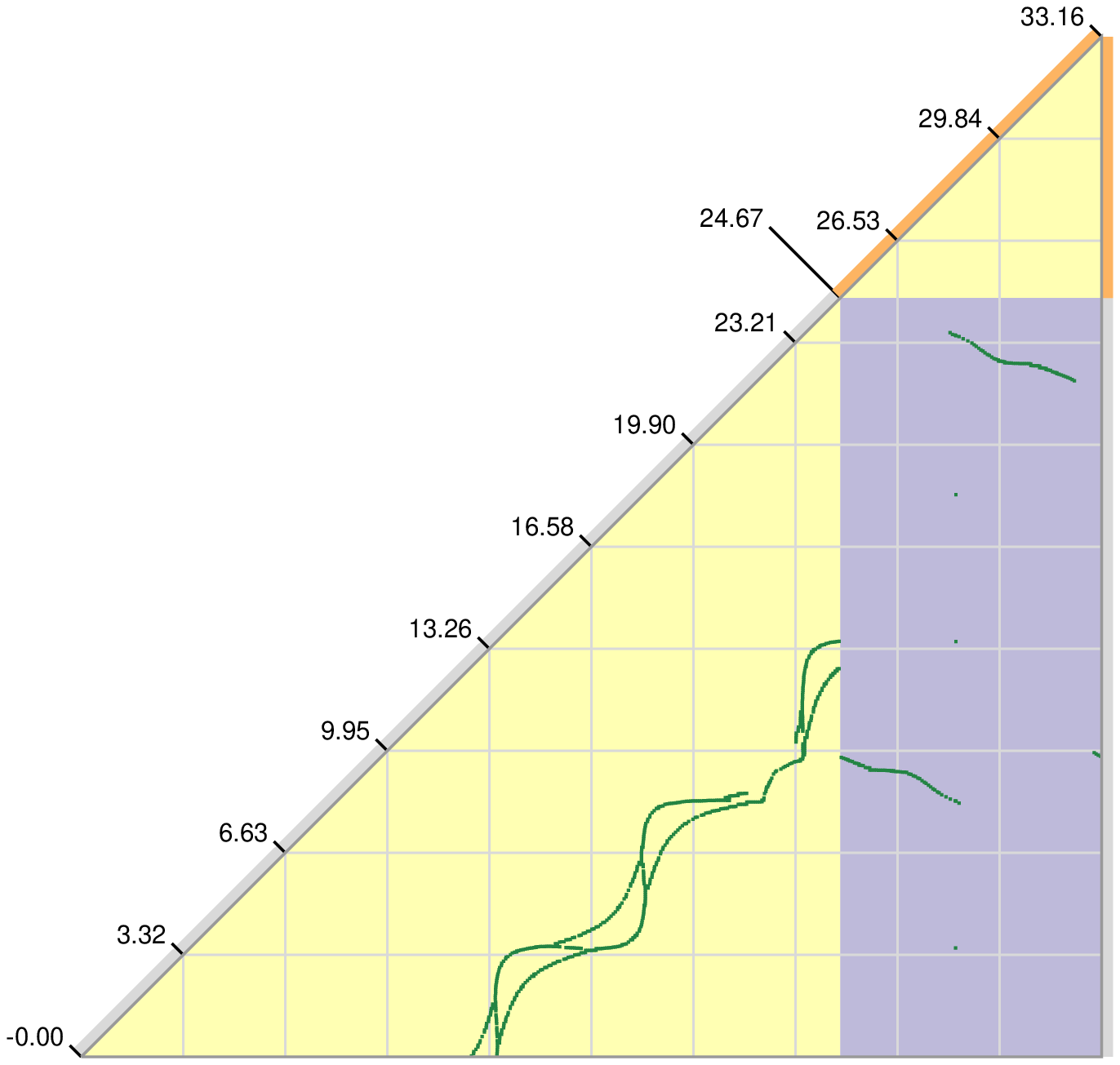}
        \put(8,94){\scriptsize{$66.33$}}
        \put(8,89){\scriptsize{$66.31$}}
        \put(8,84){\scriptsize{$473$}}
    \end{overpic}
\end{minipage} 
\hfill
\begin{minipage}[t]{6in}
  \vspace{2mm}
    \begin{overpic}[height=2.8in]{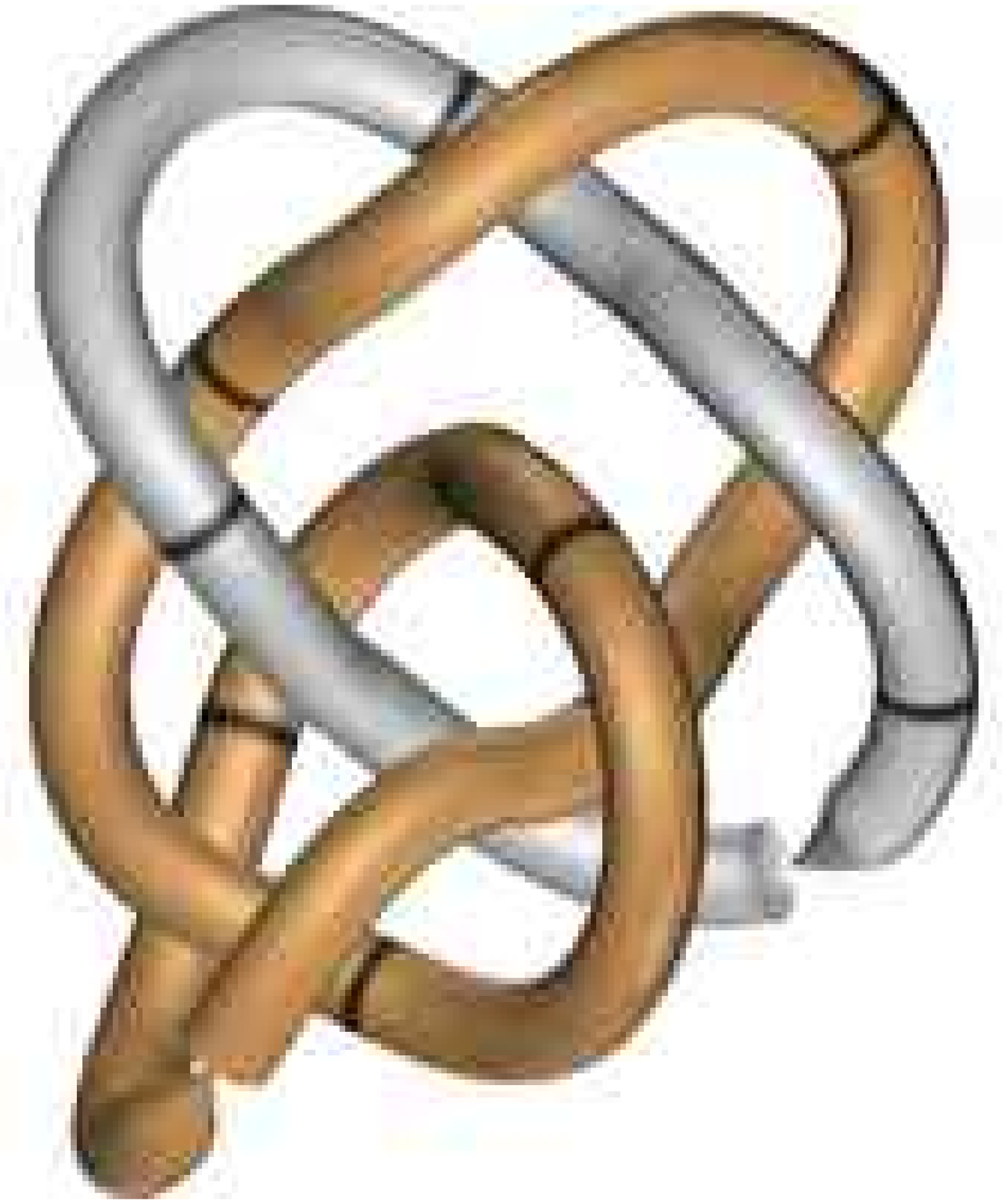}
        \put(-10,90){\large{$9^{2}_{49}$}}
    \end{overpic}
      \hspace{7mm}
    \begin{overpic}[width=2.8in]{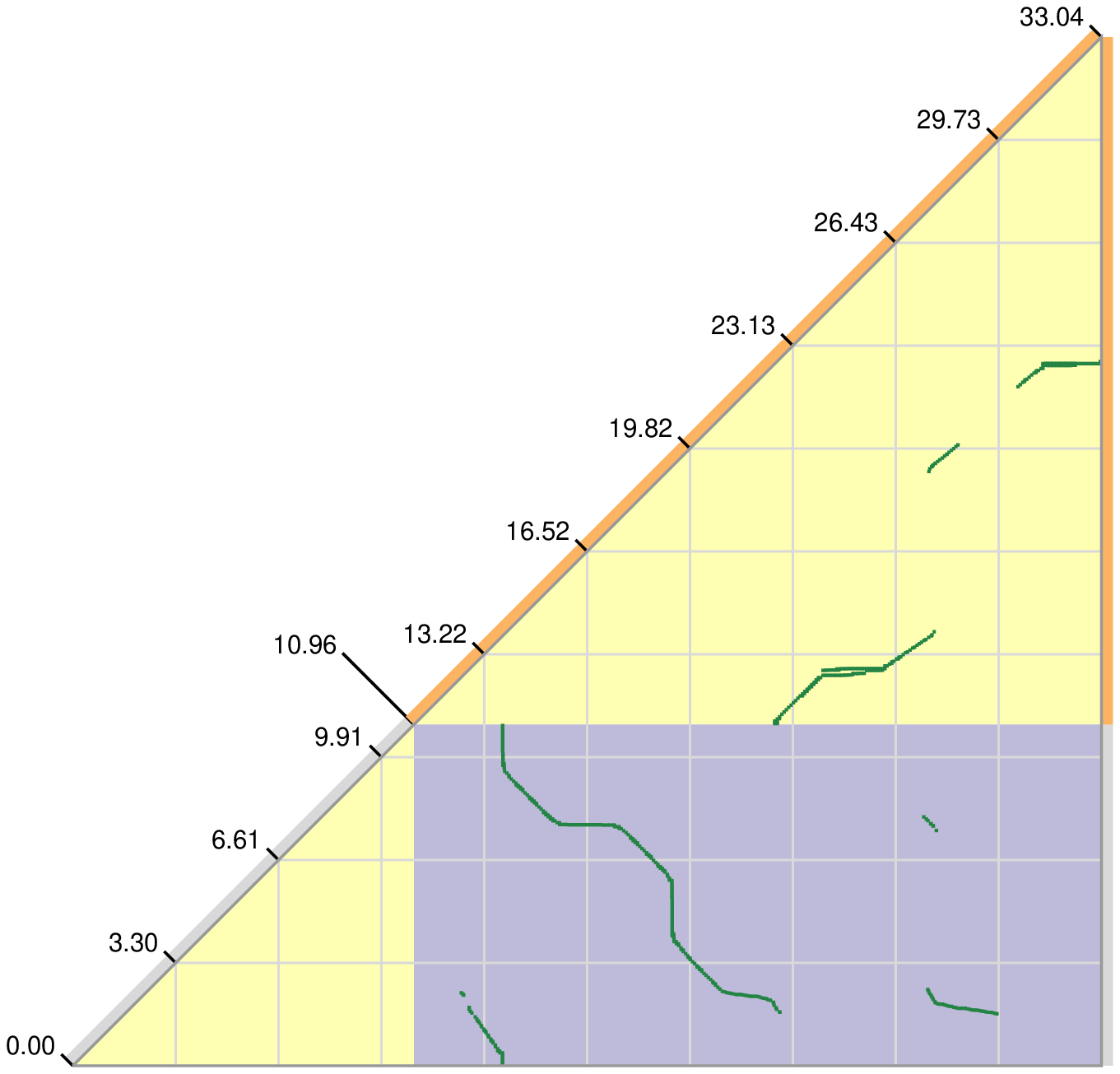}
        \put(8,94){\scriptsize{$66.09$}}
        \put(8,89){\scriptsize{$66.07$}}
        \put(8,84){\scriptsize{$472$}}
    \end{overpic}
\end{minipage} 
\hfill
\begin{minipage}[t]{6in}
  \vspace{2mm}
    \begin{overpic}[height=2.8in]{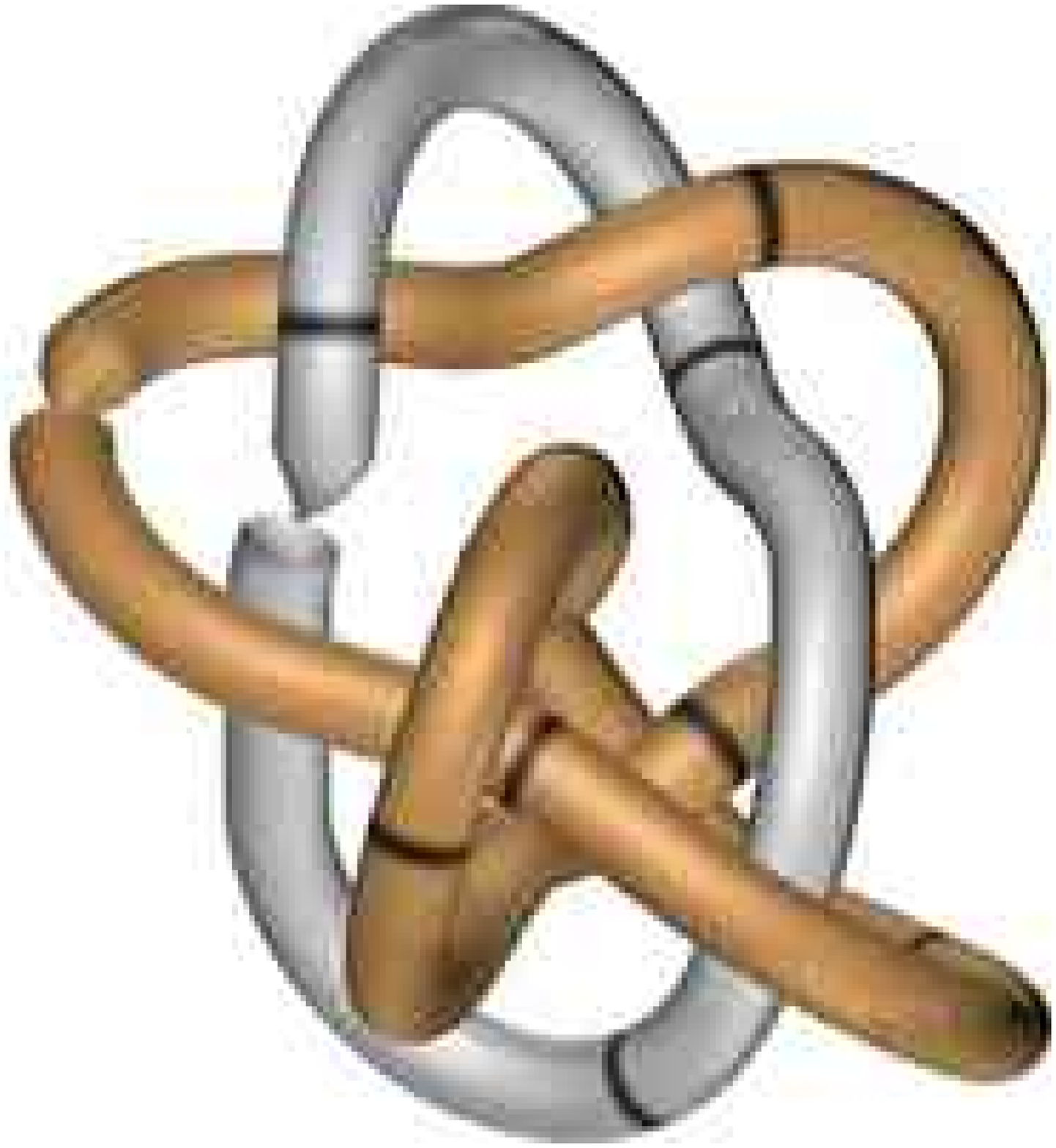}
        \put(-10,90){\large{$9^{2}_{50}$}}
    \end{overpic}
      \hspace{7mm}
    \begin{overpic}[width=2.8in]{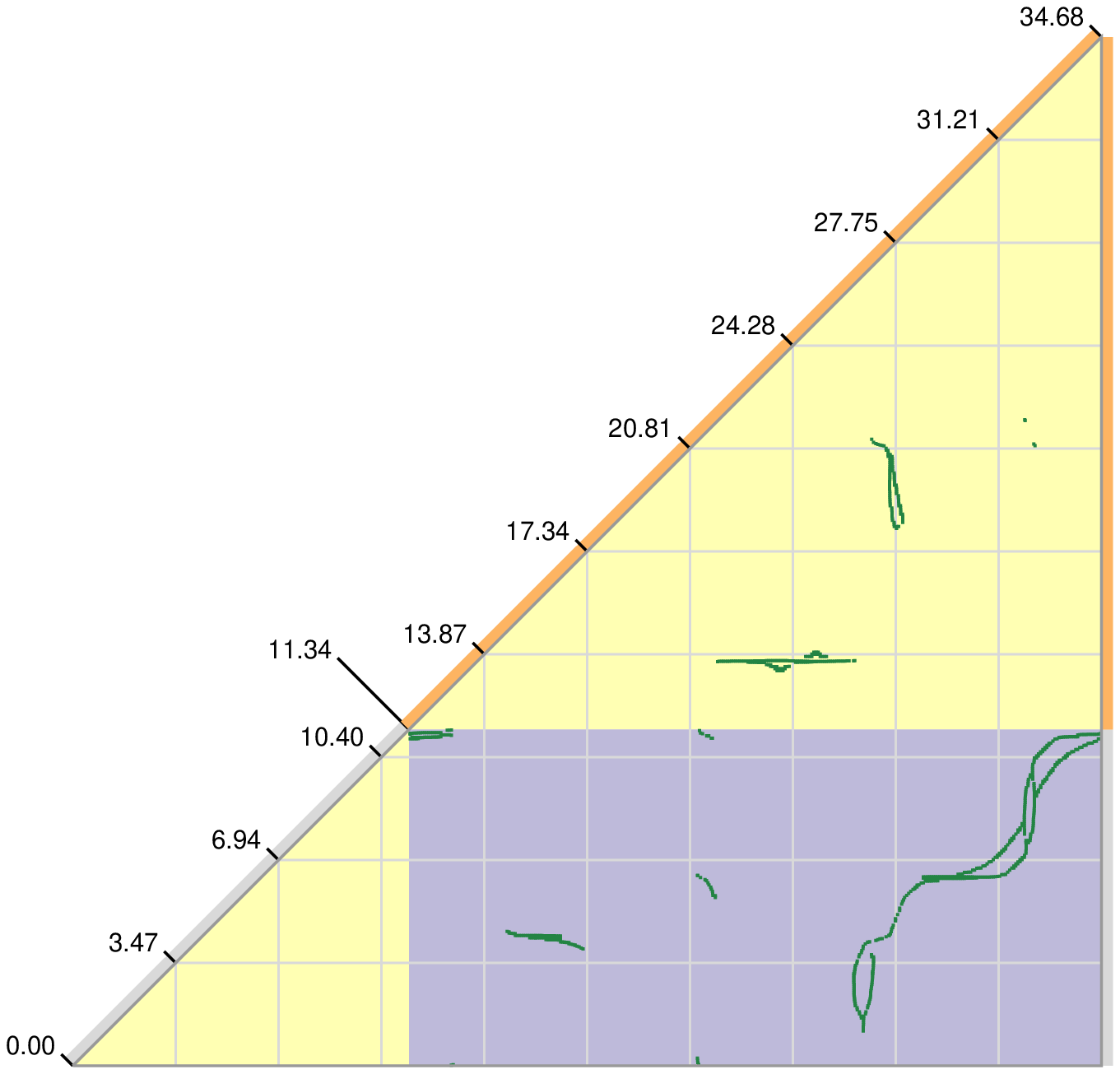}
        \put(8,94){\scriptsize{$69.37$}}
        \put(8,89){\scriptsize{$69.35$}}
        \put(8,84){\scriptsize{$495$}}
    \end{overpic}
\end{minipage} 
\hfill
\end{figure}
\clearpage
\pagebreak
\begin{figure}
\begin{minipage}[t]{6in}
  \vspace{2mm}
    \begin{overpic}[height=2.8in]{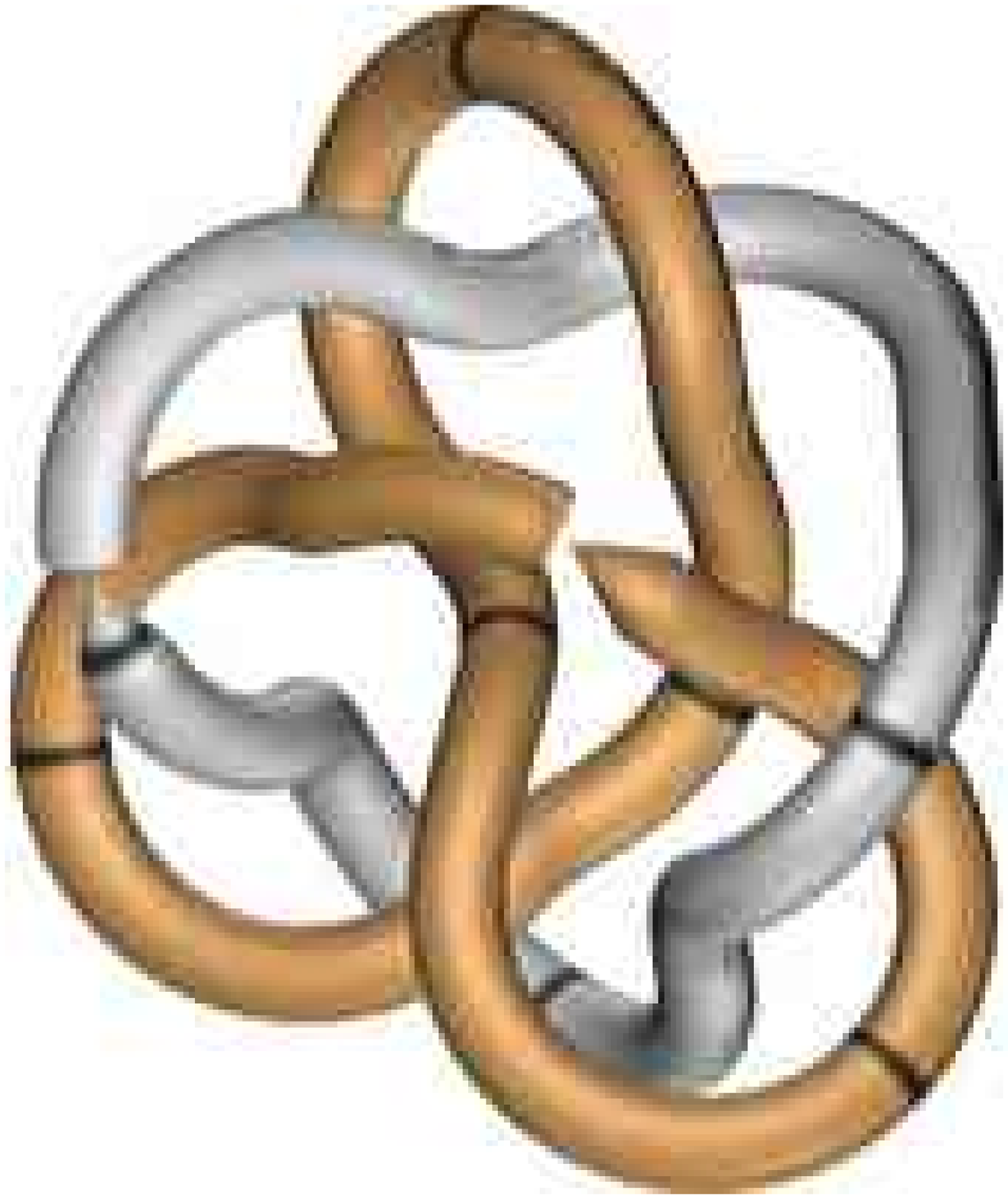}
        \put(-10,90){\large{$9^{2}_{54}$}}
    \end{overpic}
      \hspace{7mm}
    \begin{overpic}[width=2.8in]{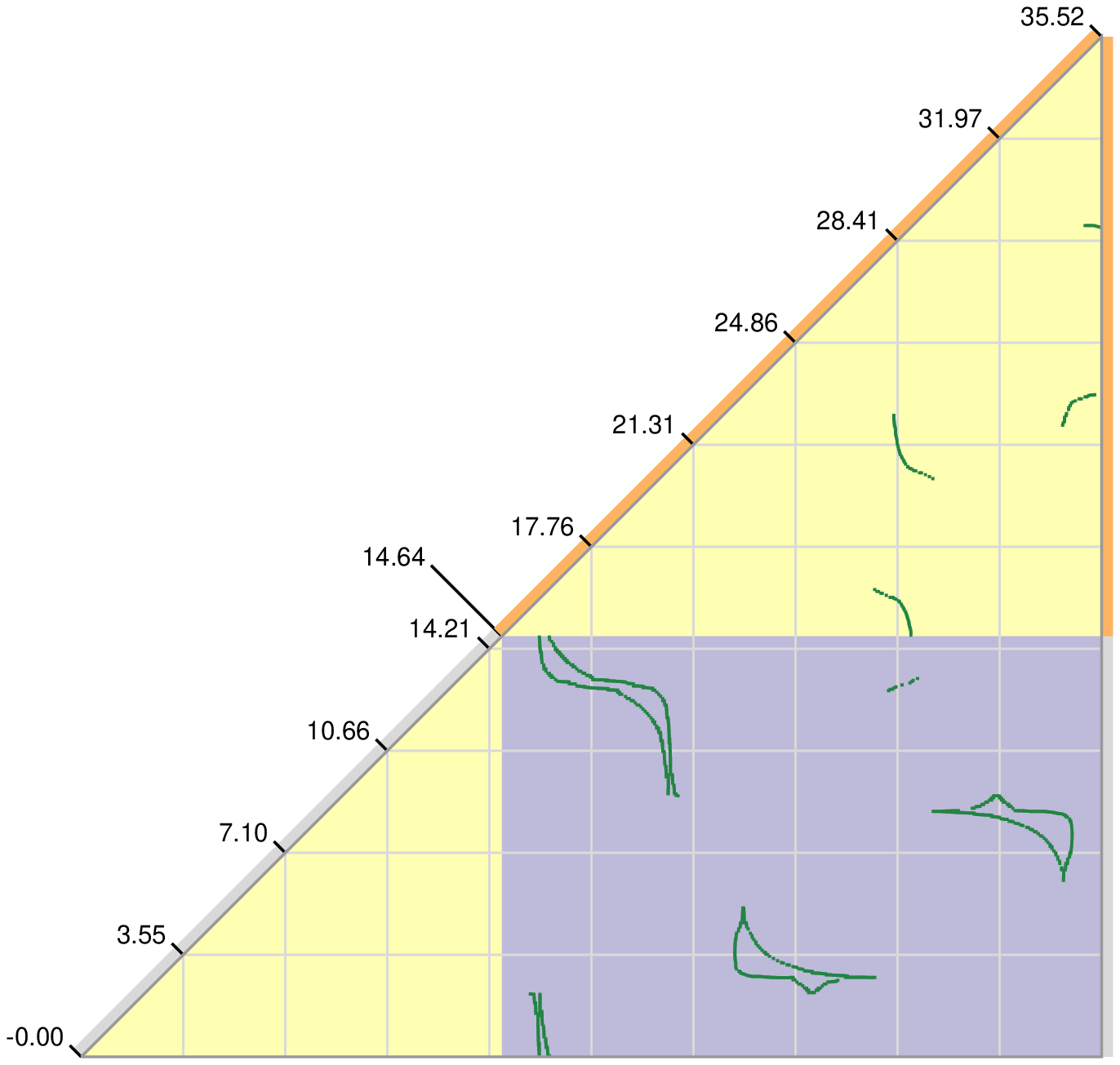}
        \put(8,94){\scriptsize{$71.04$}}
        \put(8,89){\scriptsize{$71.02$}}
        \put(8,84){\scriptsize{$507$}}
    \end{overpic}
\end{minipage} 
\hfill
\begin{minipage}[t]{6in}
  \vspace{2mm}
    \begin{overpic}[height=2.8in]{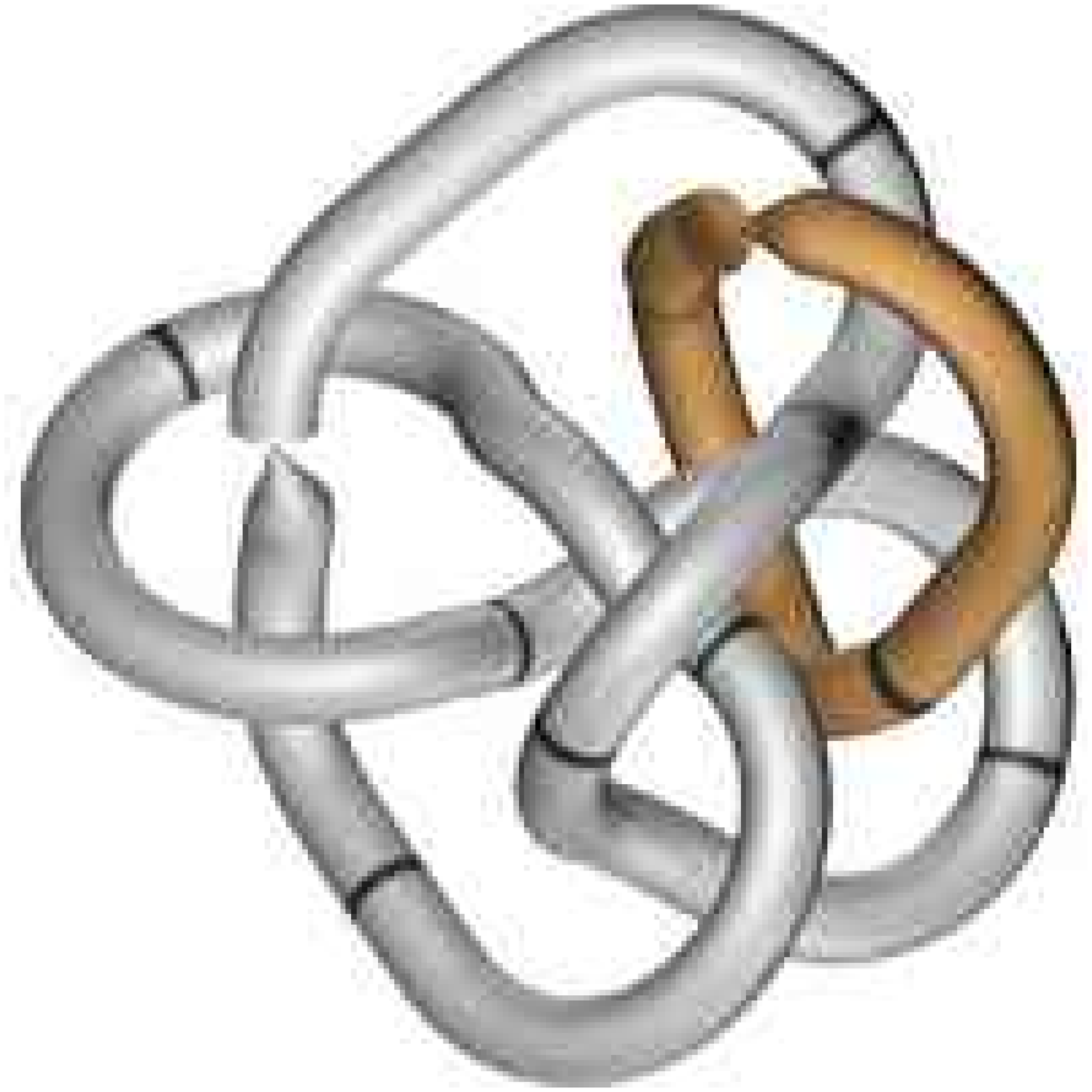}
        \put(-10,90){\large{$9^{2}_{58}$}}
    \end{overpic}
      \hspace{7mm}
    \begin{overpic}[width=2.8in]{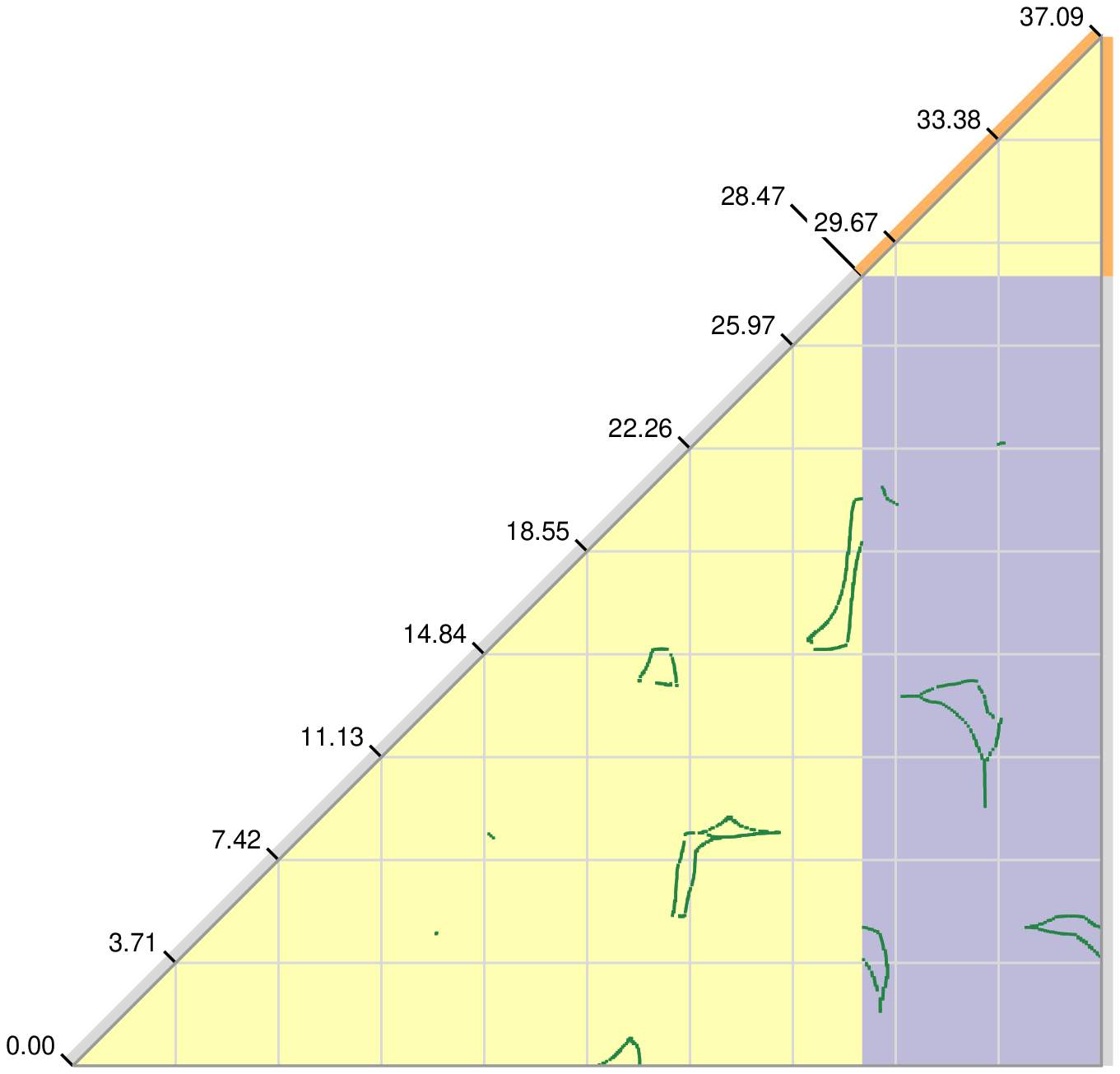}
        \put(8,94){\scriptsize{$74.20$}}
        \put(8,89){\scriptsize{$74.17$}}
        \put(8,84){\scriptsize{$530$}}
    \end{overpic}
\end{minipage} 
\hfill
\begin{minipage}[t]{6in}
  \vspace{2mm}
    \begin{overpic}[height=2.8in]{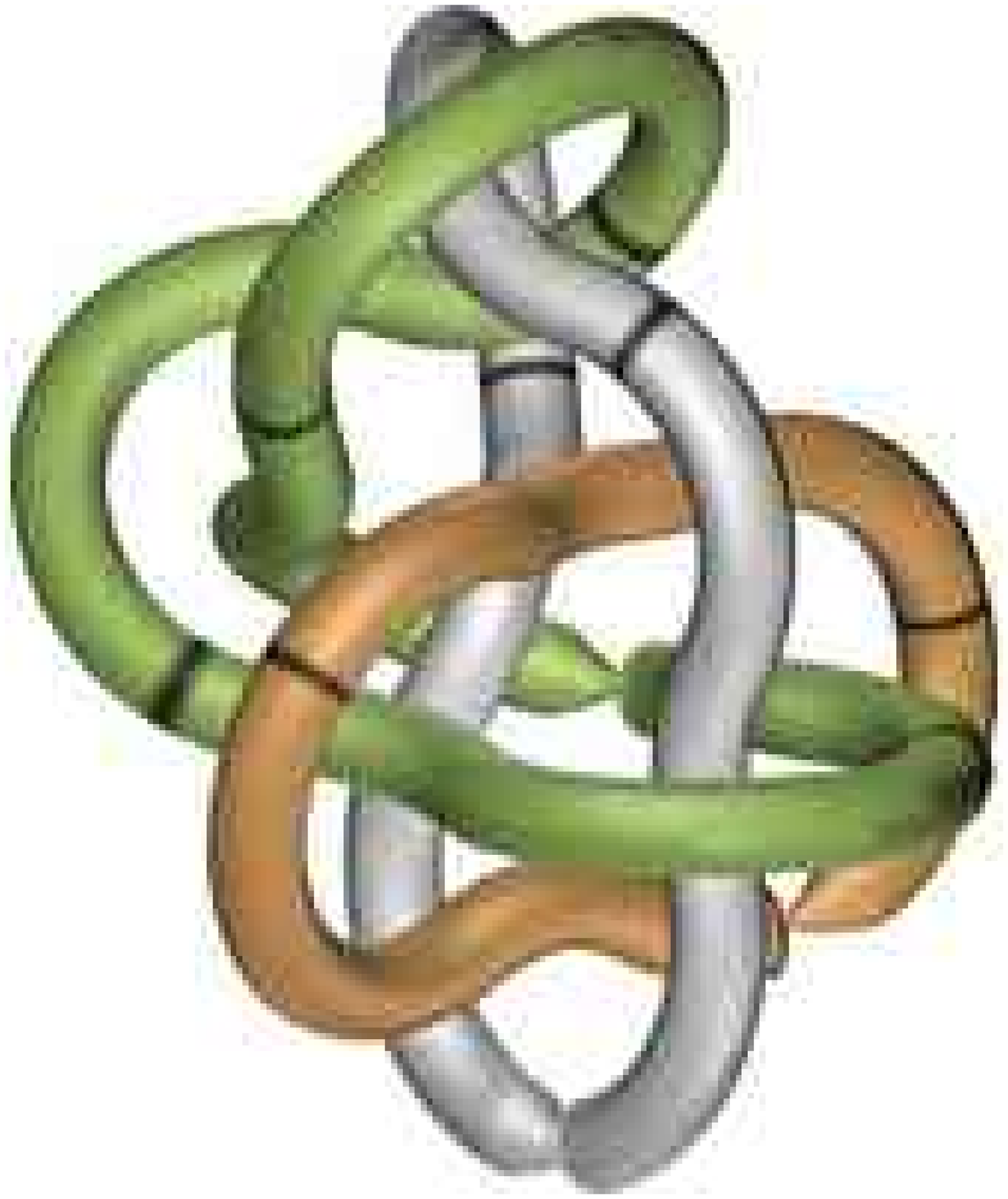}
        \put(-10,90){\large{$9^{3}_{10}$}}
    \end{overpic}
      \hspace{7mm}
    \begin{overpic}[width=2.8in]{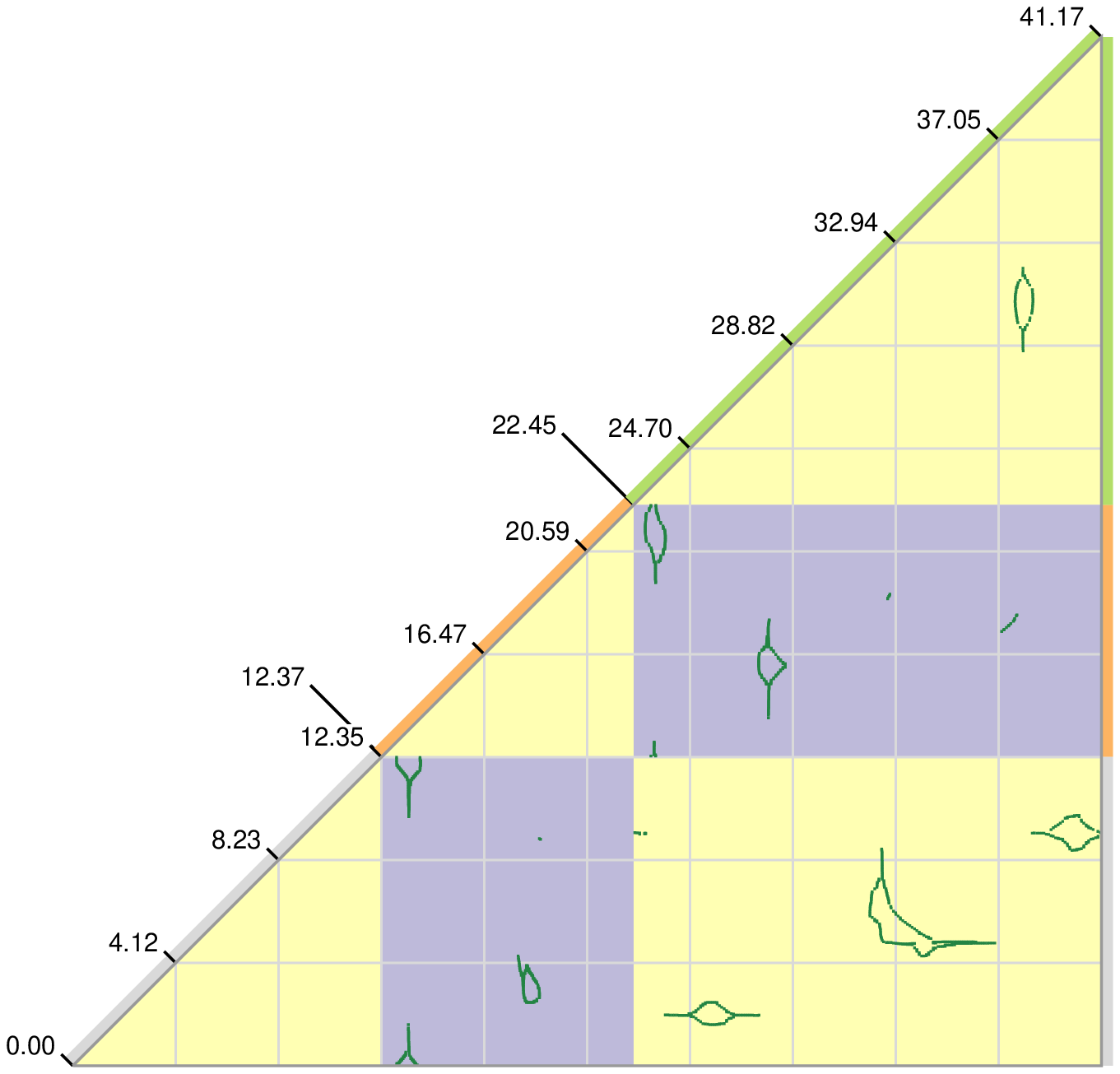}
        \put(8,94){\scriptsize{$82.35$}}
        \put(8,89){\scriptsize{$82.33$}}
        \put(8,84){\scriptsize{$588$}}
    \end{overpic}
\end{minipage} 
\hfill
\end{figure}
\clearpage
\pagebreak
\begin{figure}
\begin{minipage}[t]{6in}
  \vspace{2mm}
    \begin{overpic}[height=2.8in]{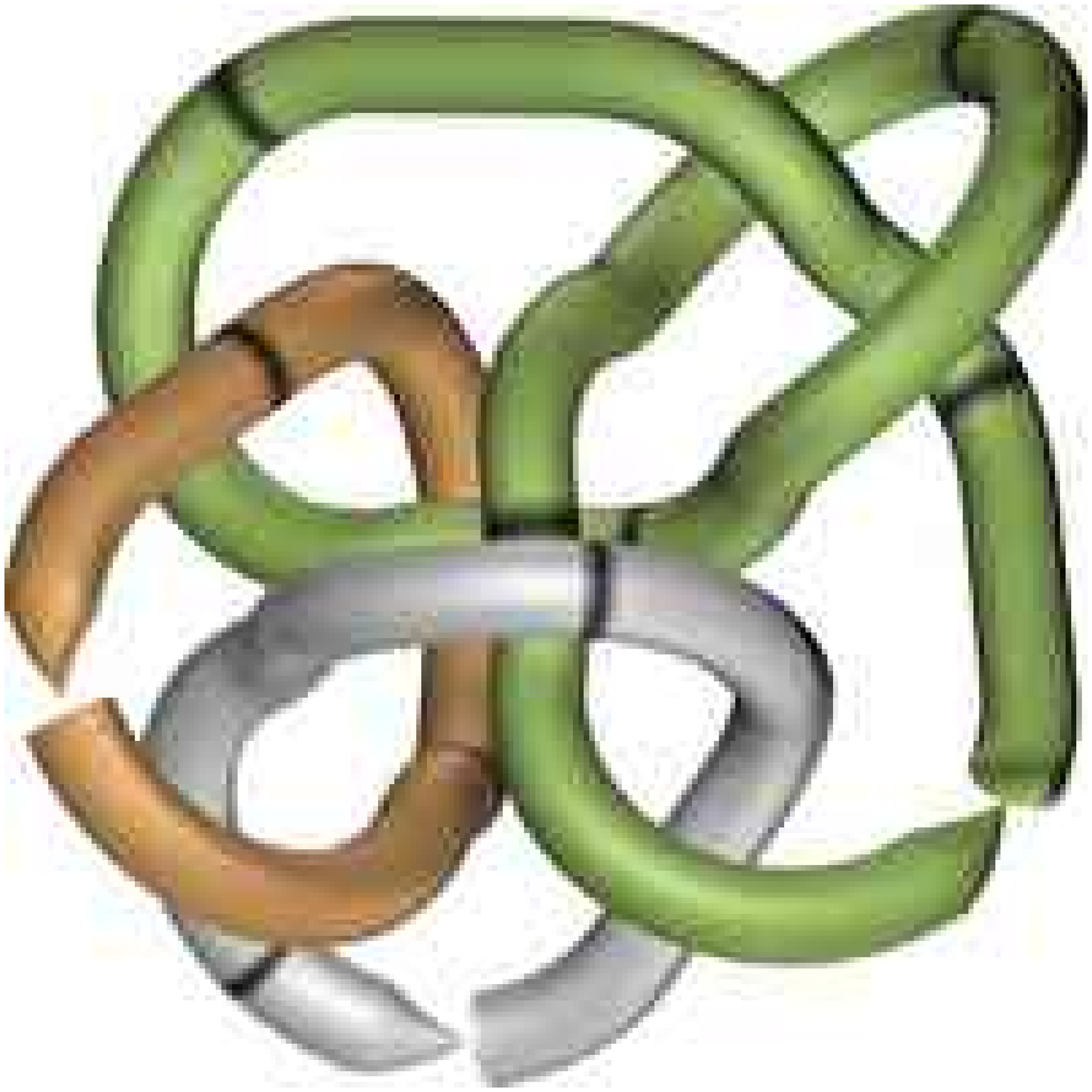}
        \put(-10,90){\large{$9^{3}_{16}$}}
    \end{overpic}
      \hspace{7mm}
    \begin{overpic}[width=2.8in]{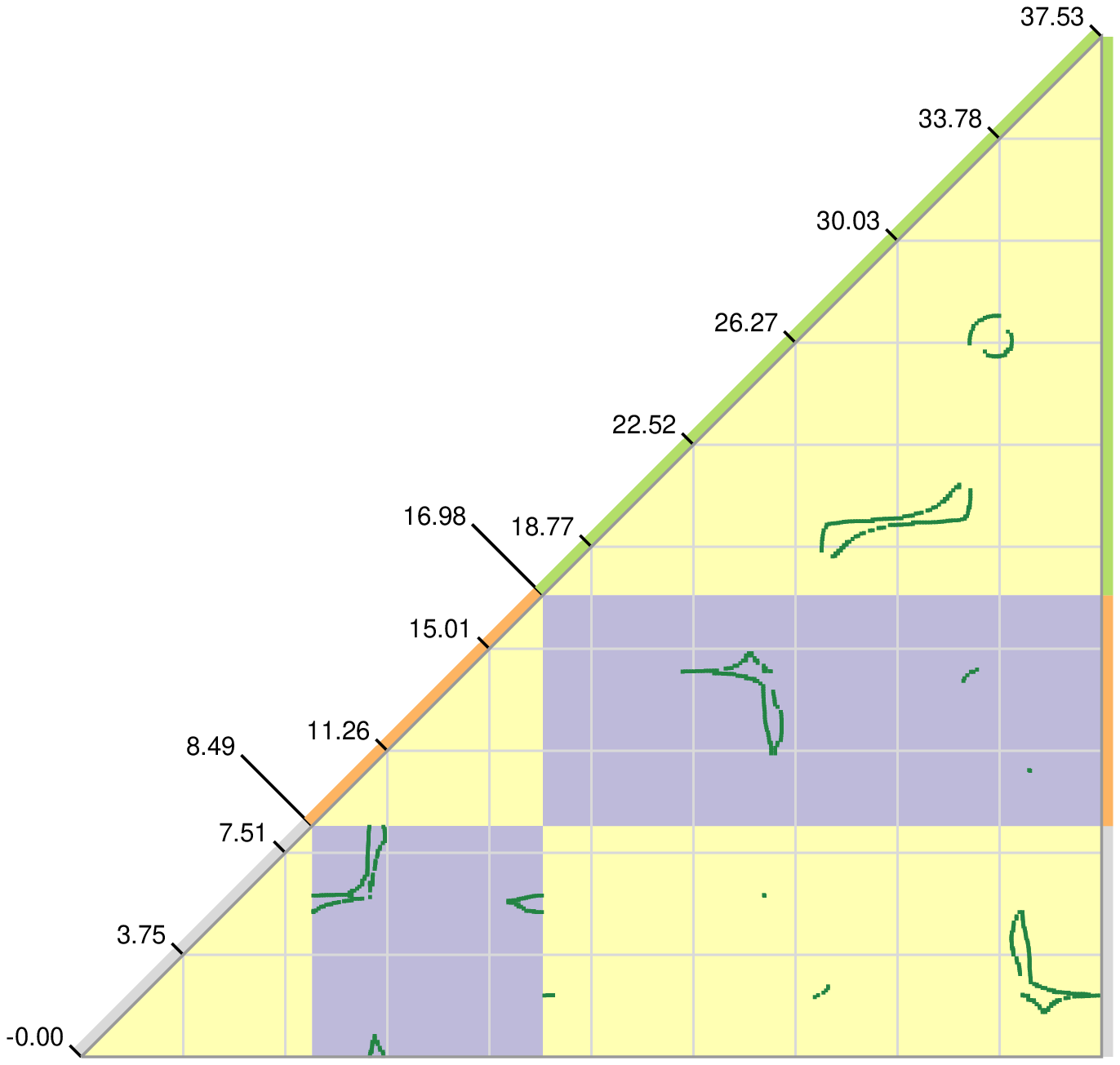}
        \put(8,94){\scriptsize{$75.07$}}
        \put(8,89){\scriptsize{$75.04$}}
        \put(8,84){\scriptsize{$400$}}
    \end{overpic}
\end{minipage} 
\hfill
\begin{minipage}[t]{6in}
  \vspace{2mm}
    \begin{overpic}[width=2.8in]{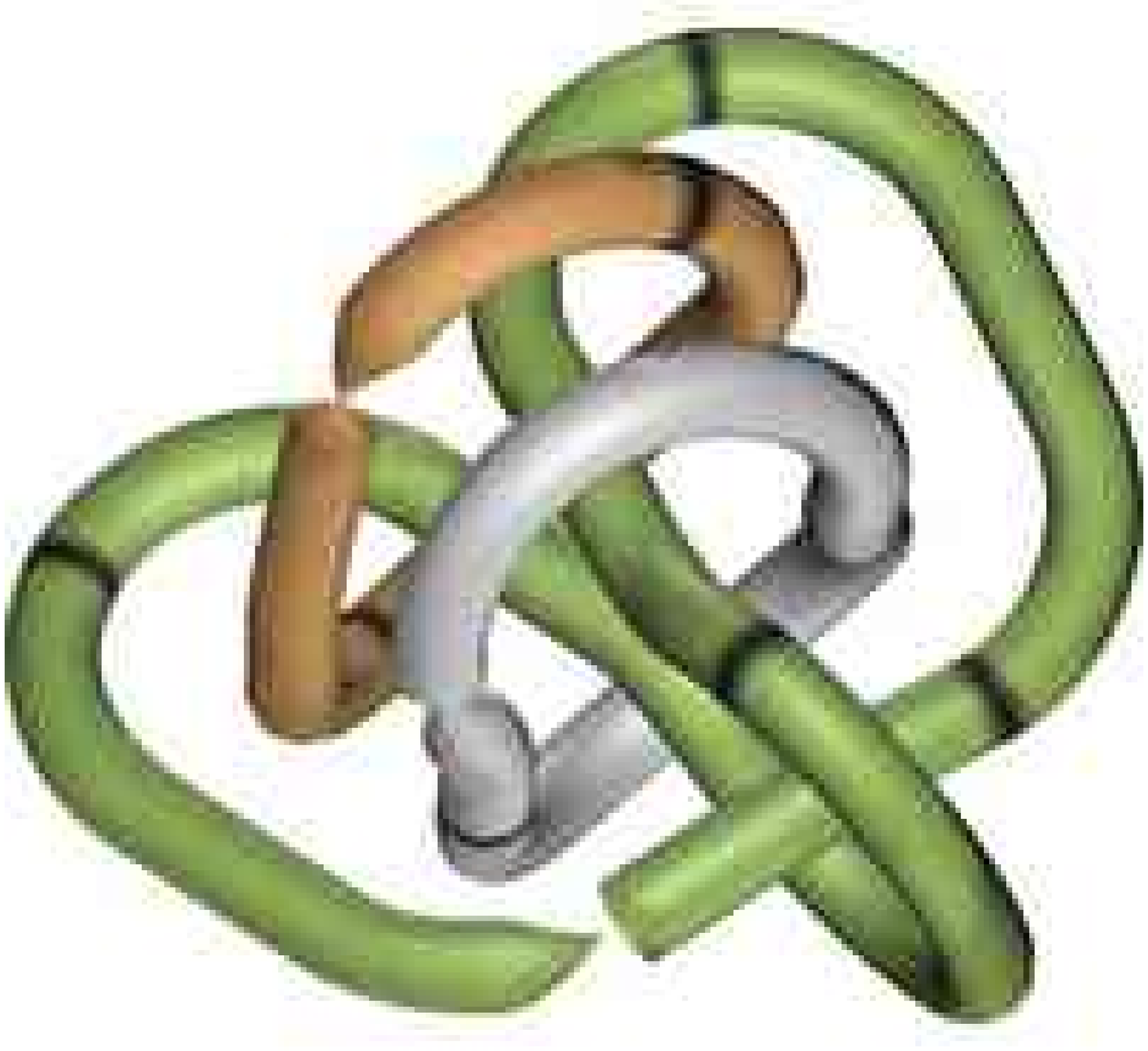}
        \put(-10,90){\large{$9^{3}_{21}$}}
    \end{overpic}
      \hspace{7mm}
    \begin{overpic}[width=2.8in]{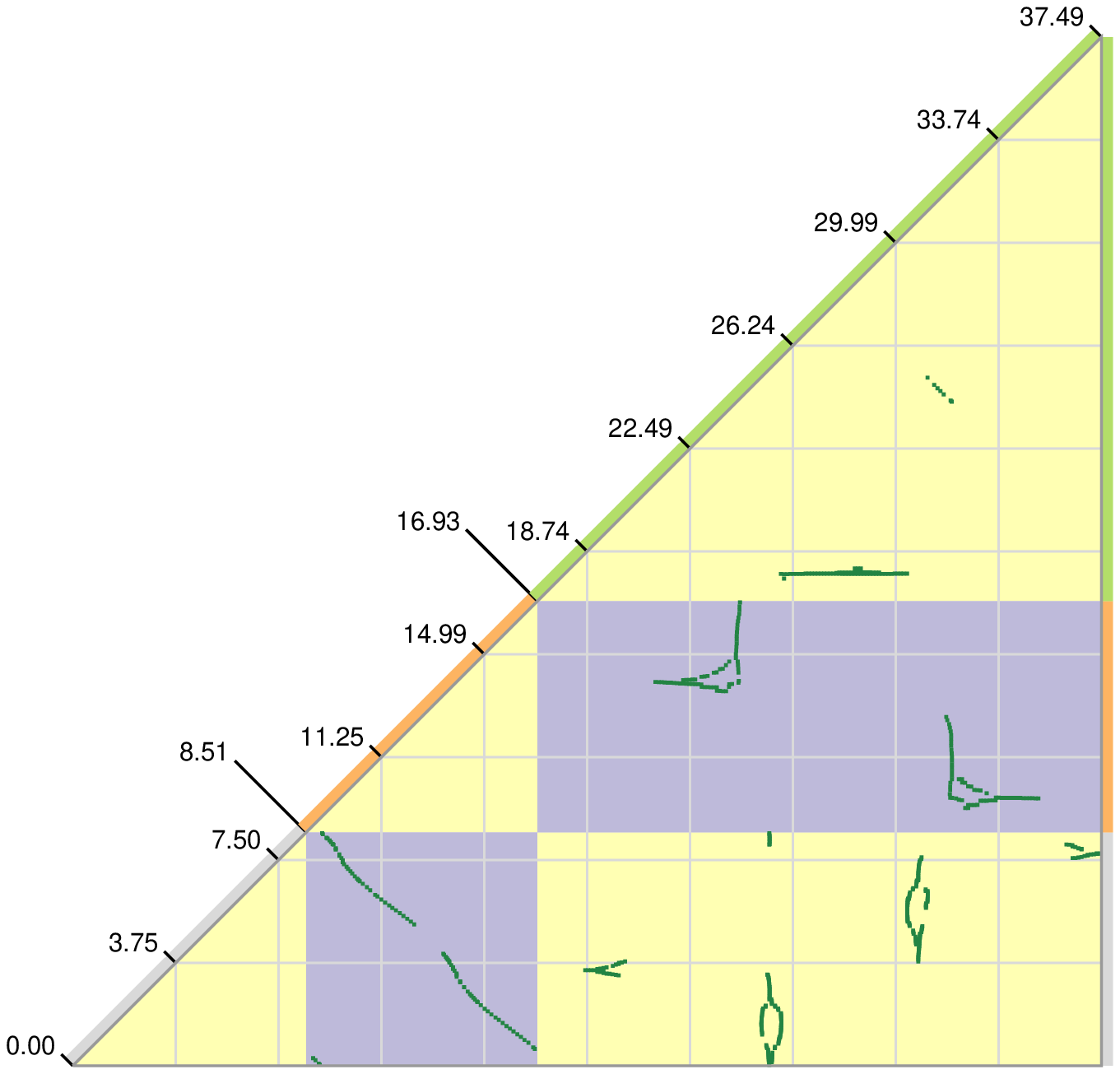}
        \put(8,94){\scriptsize{$74.98$}}
        \put(8,89){\scriptsize{$74.97$}}
        \put(8,84){\scriptsize{$400$}}
    \end{overpic}
\end{minipage} 
\hfill
\end{figure}

\end{center}

\end{document}